\documentclass[12pt,letterpaper]{report}

\usepackage{amsfonts}
\usepackage{amsmath,amsthm,amssymb,bbm,mathrsfs,mathtools,xfrac} 
\allowdisplaybreaks
\usepackage{breqn}
\usepackage{nicefrac}
\usepackage{comment}

\usepackage{apacite} 
\usepackage[authoryear]{natbib}   

\usepackage[toc,page]{appendix}
\usepackage{caption}
\usepackage{tcolorbox} 
\usepackage{color,colortbl,xcolor}
\usepackage{float} 
\usepackage{framed}
\usepackage[letterpaper, margin=1in]{geometry}		
\usepackage{graphicx}
\usepackage{listings}
\usepackage{placeins}
\usepackage{booktabs}
\usepackage{threeparttable}
\usepackage{tabulary}
\usepackage{acro}
\usepackage{subfig,multirow}
\usepackage{tabularx}
\usepackage[noend]{algpseudocode}
\usepackage{algorithm}
\usepackage{algorithmicx}

\usepackage{pbox}
\usepackage{makecell}
\usepackage{soul}

\usepackage{lscape}
  \usepackage{bm}
  \usepackage{amsfonts}
  \usepackage{amsmath}
  \usepackage{mathrsfs}
  \usepackage{tikz}
   \usetikzlibrary{arrows,arrows.meta}
    \usepackage{prodint}

    \usepackage{graphicx}
    \usepackage{lipsum}
    \usepackage[export]{adjustbox}

\definecolor{whitesmoke}{RGB}{245, 245, 245}

\def\@tab@fn#1{\ensuremath{^{\mbox{{\scriptsize #1}}}}}
\def\tabfnm#1{\rlap{\@tab@fn{#1}}}
\def\tabfnt#1#2{\raggedright\@tab@fn{#1}#2}

\def\@tab@fn#1{\ensuremath{^{\mbox{{\scriptsize #1}}}}}
\def\tabfnm#1{\rlap{\@tab@fn{#1}}}
\def\tabfnt#1#2{\raggedright\@tab@fn{#1}#2}

\usepackage[utf8]{inputenc}
\usepackage[T1]{fontenc}
\usepackage{pifont}

\usepackage{siunitx}
\usepackage{array}
\usepackage[a-1b]{pdfx}
\hypersetup{
	unicode=false,
	pdftoolbar=true,
	pdfmenubar=true,
	pdffitwindow=false,     
	pdfstartview={FitH},    
	pdftitle={Contributions to Robust and Efficient Methods for Analysis of High-Dimensional Data},    
	pdfauthor={Kai Yang},     
	pdfsubject={Contributions to Robust and Efficient Methods for Analysis of High-Dimensional Data},   
	pdfcreator={Kai Yang},   
	pdfproducer={Kai Yang}, 
	pdfnewwindow=true,      
	colorlinks=true,       
	linkcolor=red,          
	citecolor=blue,        
	filecolor=black,      
	urlcolor=blue           
}

\usepackage{cleveref}

\usepackage[parfill]{parskip} 
\usepackage{setspace}
\usepackage{subfig}
\usepackage{pgfplots}

\makeatletter
\def\maxwidth{ %
	\ifdim\Gin@nat@width>\linewidth
	\linewidth
	\else
	\Gin@nat@width
	\fi
}
\makeatother

\definecolor{fgcolor}{rgb}{0.345, 0.345, 0.345}

\usepackage{framed}
\makeatletter
{\par\unskip\endMakeFramed%
	\at@end@of@kframe}
\makeatother

\definecolor{shadecolor}{rgb}{.97, .97, .97}
\definecolor{messagecolor}{rgb}{0, 0, 0}
\definecolor{warningcolor}{rgb}{1, 0, 1}
\definecolor{errorcolor}{rgb}{1, 0, 0}

\usepackage{alltt}

\newcommand{\figref}[1]{\Cref{#1}}
\DeclarePairedDelimiter\abs{\lvert}{\rvert}%
\DeclarePairedDelimiter\norm{\lVert}{\rVert}%

\makeatletter
\let\oldabs\abs
\def\abs{\@ifstar{\oldabs}{\oldabs*}}
\let\oldnorm\norm
\def\norm{\@ifstar{\oldnorm}{\oldnorm*}}
\makeatother

\IfFileExists{upquote.sty}{\usepackage{upquote}}{}

\newtheoremstyle{def}
{10pt}
{20pt}
{}
{}
{\rmfamily\mdseries\scshape}
{.}
{.5em}
{}
\theoremstyle{def}

\newtheorem{thm}{Theorem}
\newtheorem{cor}[thm]{Corollary}

\DeclareAcronym{ABIDE}{
  short=ABIDE,
  long=Autism Brain Imaging Data Exchange,
  tag = abbrev
}
\DeclareAcronym{KL divergence}{
  short=KL divergence,
  long=Kullback--Leibler divergence,
  tag = abbrev
}
\DeclareAcronym{kNN}{
  short=$k$NN,
  long=$k$-Nearest Neighbors,
  tag = abbrev
}
\DeclareAcronym{DFT}{
  short=DFT,
  long=Discrete Fourier Transform,
  tag = abbrev
}
\DeclareAcronym{FFT}{
  short=FFT,
  long=Fast Fourier Transform,
  tag = abbrev
}
\DeclareAcronym{KDE}{
  short=KDE,
  long=Kernel Density Estimation,
  tag = abbrev
}
\DeclareAcronym{FFTKDE}{
  short=FFTKDE,
  long=Fast Fourier Transform-based Kernel Density Estimation,
  tag = abbrev
}
\DeclareAcronym{AG}{
  short=AG,
  long=Accelerated Gradient,
  tag = abbrev
}
\DeclareAcronym{CG}{
  short=CG,
  long=Conjugate Gradient,
  tag = abbrev
}
\DeclareAcronym{GMRES}{
  short=GMRES,
  long=Generalized Minimal Residuals,
  tag = abbrev
}

\DeclareAcronym{SCAD}{
  short=SCAD,
  long=Smoothly Clipped Absolute Deviation,
  tag = abbrev
}
\DeclareAcronym{MCP}{
  short=MCP,
  long=Minimax Concave Penalty,
  tag = abbrev
}

\DeclareAcronym{PCA}{
  short=PCA,
  long=Principal Components Analysis,
  tag = abbrev
}

\DeclareAcronym{ICA}{
  short=ICA,
  long=Independent Component Analysis,
  tag = abbrev
}

\DeclareAcronym{GWAS}{
  short=GWAS,
  long=Genome-Wide Association Studies,
  tag = abbrev
}

\DeclareAcronym{GMM}{
  short=GMM,
  long=Generalized Method of Moments,
  tag = abbrev
}

\DeclareAcronym{MRI}{
  short=MRI,
  long=Magnetic Resonance Images,
  tag = abbrev
}
\DeclareAcronym{CRLB}{
  short=CRLB,
  long=Cramer--Rao Lower Bound,
  tag = abbrev
}
\DeclareAcronym{NYSE}{
  short=NYSE,
  long=New York Stock Exchange,
  tag = abbrev
}
\DeclareAcronym{NASDAQ}{
  short=NSADAQ,
  long=National Association of Securities Dealers Automated Quotations,
  tag = abbrev
}
\DeclareAcronym{BSM}{
  short=BSM,
  long=Black--Scholes--Merton,
  tag = abbrev
}
\DeclareAcronym{LASSO}{
  short=LASSO,
  long=Least Absolute Shrinkage and Selection Operator,
  tag = abbrev
}
\DeclareAcronym{SNR}{
  short=SNR,
  long=Signal-to-Noise Ratio,
  tag = abbrev
}

\usepackage{bold-extra}
\renewcommand{\bf}{\bfseries}

\providecommand{\lemmaname}{Lemma}
\providecommand{\remarkname}{Remark}
\theoremstyle{plain}
\newtheorem{lem}[thm]{\protect\lemmaname}
\theoremstyle{remark}
\newtheorem{rem}[thm]{\protect\remarkname}

\newenvironment{QZ@Cent}{\centering}{\par}%

\newboolean{SetDSpace}%
\setboolean{SetDSpace}{true}%

\newcommand{\BaseDiff}{0}

\newcommand{\GoSingle}{\renewcommand{\baselinestretch}{1}
	\normalfont\tiny\normalsize
}

\newcommand{\GoDouble}{\renewcommand{\baselinestretch}{1.655}
	\renewcommand{\BaseDiff}{0.655}\normalfont\tiny\normalsize
}

\newcommand*{\SetTitle}[1]{\renewcommand*{\Title}{#1}}
\newcommand*{\Title}{No Title Given}

\newcommand*{\SetAuthor}[1]{\renewcommand*{\FullName}{#1}}
\newcommand*{\FullName}{Please Define Your Name}

\newcommand*{\SetDegreeType}[1]{\renewcommand*{\DegreeType}{#1}}
\newcommand*{\DegreeType}{UNDEFINED DEGREE}

\newcommand*{\SetDepartment}[1]{\renewcommand*{\ETDDepartment}{#1}}
\newcommand*{\ETDDepartment}{UNDEFINED DEPARTMENT}

\newcommand*{\SetUniversity}[1]{\renewcommand*{\ETDUniversity}{#1}}
\newcommand*{\ETDUniversity}{McGill University}

\newcommand*{\SetUniversityAddr}[1]{\renewcommand*{\ETDUniversityAddr}{#1}}
\newcommand*{\ETDUniversityAddr}{Montreal, Quebec}

\newcommand*{\SetThesisDate}[1]{\renewcommand*{\ETDThesisDate}{#1}}
\newcommand*{\ETDThesisDate}{Date111}

\newcommand*{\SetRequirements}[1]{\renewcommand*{\ETDRequirements}{#1}}
\newcommand*{\ETDRequirements}{Date222}

\newcommand*{\SetCopyright}[1]{\renewcommand*{\ETDCopyright}{#1}}
\newcommand*{\ETDCopyright}{Date223333}

\newenvironment{cent}{\centering}{\par}

\renewcommand{\maketitle}{
	\clearpage
	\thispagestyle{empty}
	\begin{cent}
		\Title
		\vfill
		\GoSingle
		\normalsize\normalfont
		\FullName\normalsize\normalfont \\*[\BaseDiff\baselineskip]
		\vfill
		\DegreeType\normalsize\normalfont \\*[\BaseDiff\baselineskip]
		\vfill
		\ETDDepartment\normalsize\normalfont \\*[\BaseDiff\baselineskip]
		\vfill
		\ETDUniversity\normalsize\normalfont  \\*[\BaseDiff\baselineskip]
		\ETDUniversityAddr\normalsize\normalfont \\*[\BaseDiff\baselineskip]
		\ETDThesisDate\normalsize\normalfont \\*[\BaseDiff\baselineskip]
		\vfill
		\ETDRequirements\normalsize\normalfont \\*[\BaseDiff\baselineskip]
		\textcircled{c} \ETDCopyright\normalsize\normalfont \\*[\BaseDiff\baselineskip]
	\end{cent}
	\vspace*{0.5in}
	\clearpage
}

\newcommand*{\ETDDedicationName}{Dedication}

\newcommand*{\SetDedicationText}[1]{\renewcommand*{\ETDDedicationText}{#1}}
\newcommand*{\ETDDedicationText}{Dedication text goes here!}

\newenvironment{simpleenv}[4]{\clearpage}{\clearpage}

\newcommand{\Dedication}{
	\begin{simpleenv}{}{}{}{}
		\pagestyle{plain}
		\GoSingle
		\begin{cent}
			\bfseries{\ETDDedicationName}
		\end{cent}
		\vspace*{0.5in}
		\par
		\GoDouble
		\ETDDedicationText
	\end{simpleenv}
}

\newcommand*{\ETDAcknowledgeName}{Acknowledgements}

\newcommand*{\SetAcknowledgeText}[1]{\renewcommand*{\ETDAcknowledgeText}{#1}}
\newcommand*{\ETDAcknowledgeText}{Acknowledgements text goes here!}

\newcommand{\Acknowledge}{
	\begin{simpleenv}{}{}{}{}
		\pagestyle{plain}
		\GoSingle
		\begin{cent}
			\bfseries{\ETDAcknowledgeName}
		\end{cent}
		\vspace*{0.5in}
		\par
		\GoDouble
		\ETDAcknowledgeText
	\end{simpleenv}
}

\newcommand*{\ETDPrefaceName}{Preface }

\newcommand*{\SetPrefaceText}[1]{\renewcommand*{\ETDPrefaceText}{#1}}
\newcommand*{\ETDPrefaceText}{Preface text goes here!}

\newcommand{\Preface}{
	\begin{simpleenv}{}{}{}{}
		\GoSingle
		\begin{cent}
			\bfseries{\ETDPrefaceName}
		\end{cent}
		\vspace*{0.5in}
		\par
		\GoDouble
		\ETDPrefaceText
	\end{simpleenv}
}

\newenvironment{romanPagenumber}[1]
{\setcounter{page}{#1}\renewcommand{\thepage}{\roman{page}}}
{\pagenumbering{arabic}}

\newcommand*{\ETDAbstractEnName}{Abstract}

\newcommand*{\SetAbstractEnText}[1]{\renewcommand*{\ETDAbstractEnText}{#1}}
\newcommand*{\ETDAbstractEnText}{Abstract text goes here!}

\newcommand*{\AbstractEn}{
    \clearpage
    \pagestyle{plain}
    \GoSingle
	\begin{cent}
		\bfseries{\ETDAbstractEnName}
	\end{cent}
    \par
    \GoDouble
    \ETDAbstractEnText
}

\newcommand*{\SetAbstractFrName}[1]{\renewcommand*{\ETDAbstractFrName}{#1}}
\newcommand*{\ETDAbstractFrName}{Abr\'{e}g\'{e}}

\newcommand*{\SetAbstractFrText}[1]{\renewcommand*{\ETDAbstractFrText}{#1}}
\newcommand*{\ETDAbstractFrText}{Abstract text goes here!}

\newcommand*{\AbstractFr}{
    \clearpage
    \pagestyle{plain}
    \GoSingle
	\begin{cent}
		\bfseries{\ETDAbstractFrName}
	\end{cent}
    \par
    \GoDouble
    \ETDAbstractFrText
}

\newcommand*{\bibHeading}[1]{
	\renewcommand{\bibname}{#1}
}

\newcommand*{\LOFHeading}[1]{
    \renewcommand{\listfigurename}{#1}
}

\definecolor{dkgreen}{rgb}{0,0.6,0}
\definecolor{gray}{rgb}{0.5,0.5,0.5}
\definecolor{mauve}{rgb}{0.58,0,0.82}

\newcommand*{\ETDAppendix}[2]{%
	\begin{simpleenv}{}{}{}{}%
		\setboolean{SetDSpace}{false}%
		\pagestyle{plain}%
		\GoSingle%
		\begin{QZ@Cent}%
			\bfseries{#1} 
		\end{QZ@Cent}%
		\vspace*{0.5in}%
		\par%
		\GoDouble%
		\addcontentsline{toc}{extrachapter}{#1}%
		#2%
	\end{simpleenv}%
}

\usepackage{tablefootnote}
\pgfplotsset{compat=1.15}

\SetTitle{\huge{Contributions to Robust and Efficient Methods for Analysis of High--Dimensional Data}}%
\SetAuthor{Kai Yang}
\SetDegreeType{Doctor of Philosophy}%
\SetDepartment{\begin{figure}[H]
\centering{}\includegraphics[width=0.6\textwidth]{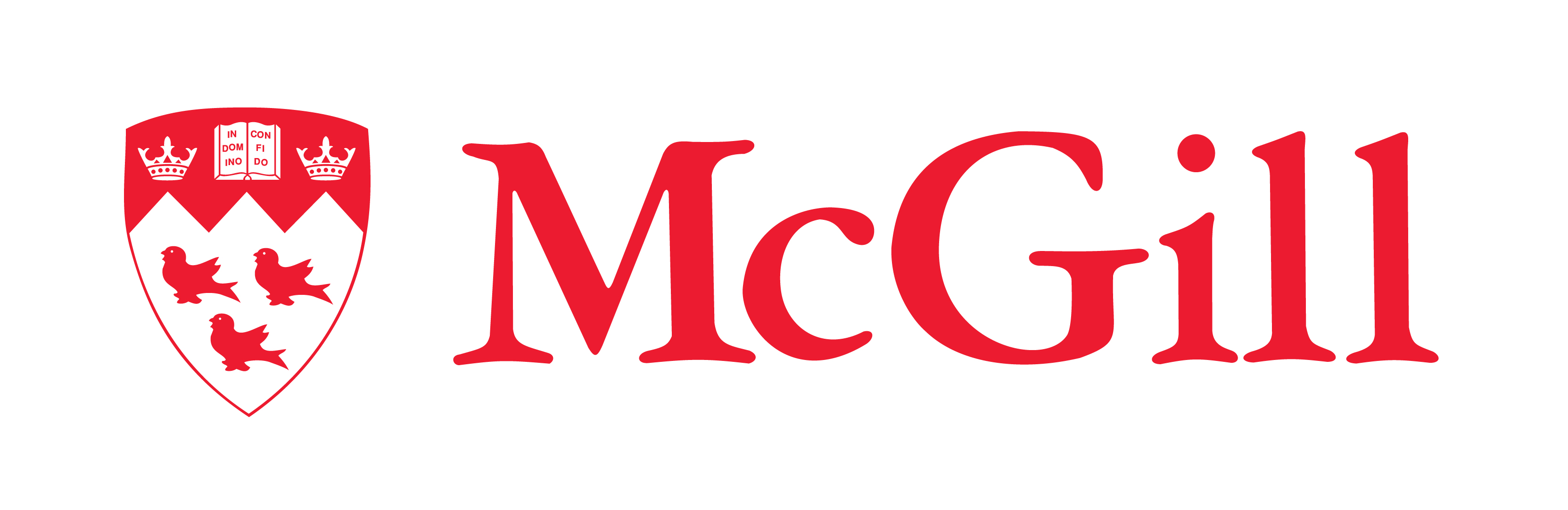}\end{figure}
\vspace{.8in}
Department of Epidemiology, Biostatistics and Occupational Health}
\SetUniversity{McGill University}%
\SetUniversityAddr{Montr\'eal, Qu\'ebec}%
\SetThesisDate{July 2024}
\SetRequirements{A thesis submitted to McGill University in partial fulfillment of the requirements of the degree of Doctor of Philosophy}%
\SetCopyright{Copyright Kai Yang, 2024}

\begin{document}

\maketitle%

\begin{romanPagenumber}{2}%

\SetDedicationText{
I tread paths laid by giants, whose towering achievements guide me \emph{philosophically and academically.}\\
\begin{quote}
\textit{Wir m\"ussen wissen, wir werden wissen.}\\
\hfill \textit{--- David Hilbert, 8 September 1930}\\
\textit{G{\"o}del's incompleteness theorems \citep{Goedel1931}}\\
\hfill \textit{--- Kurt G{\"o}del}
\end{quote}
I dedicate this thesis to you, the reader, who will navigate through my approximately 200 pages of writing with statements that can be deeply traced back to ZFC; I also dedicate this to us, the people, living in the time post--$\text{sub}\left(n,n,17\right)$.\\
} 
\Dedication%

\SetAcknowledgeText{\textcolor{black}
{This dissertation could not have been completed without the invaluable contributions and support from a host of dedicated individuals. I am immensely grateful for their assistance throughout this journey.
\newline
I extend my deepest appreciation to my supervisors, Celia M.T. Greenwood, Masoud Asgharian, Sahir Bhatnagar, for their invaluable guidance, persistent support, and expert advice throughout the course of this research. Special appreciation is due to Celia Greenwood for translating the English abstract into French, consistently providing timely and detailed feedback, and fulfilling the duties of a supervisor with exceptional dedication, even while facing significant health challenges. Additionally, Masoud Asgharian made substantial contributions to the discussion chapter by suggesting numerous avenues for future research, which have greatly enriched the scope and depth of this thesis, and offered a profound course in advanced nonparametric statistics that bridged theory with application together with a concise and invaluable reference summary of concepts in probability theory, contained within just a few pages. Their expertise and encouragement have been crucial to the success of this academic endeavor, and they have consistently offered prompt and profound insights throughout the research process.
\newline
I also like to thank Adam Oberman, Tim Hoheisel, Courtney Paquette, Gantumur Tsogtgerel, Jean--Christophe Nave, Jean--Philippe Lessard, and Russell Davidson for their exceptional teaching in mathematical machine learning, convex analysis, functional analysis, numerical analysis, dynamical systems, and stochastic differential equations. Special thanks to Jean--Philippe Lessard for sharing his pre--published book \emph{Ordinary Differential Equations: A Constructive Approach} \citep{Berg2023}, which greatly enhanced my learning.
\newline
Thanks also to Shayda Asgharian, Masoud Asgharian's daughter, for her assistance in translating the English abstract to French.
}}%
\Acknowledge%

%
%
\SetPrefaceText{
The work presented, including the introduction, literature review, bridging texts, discussion, and conclusion, was authored by myself, Kai Yang, and significantly enhanced under the diligent guidance and thorough revisions provided by my supervisors, Celia Greenwood and Masoud Asgharian.\\
The author contributions to each of the three manuscripts included in this thesis are as follows:\\
{\bf Manuscript 1:}
\begin{itemize}
\item Kai Yang (student): Conceptualization (introduced mutual information estimation using FFT Kernel Density Estimation for variable screening); Formal analysis, Methodology, Investigation, and Software (carried out all mathematical proofs and developments, developed the screening method and Python package, designed and executed simulations and case studies); Visualization; Writing - Original Draft (wrote the initial draft); Writing - Review \& Editing (revisions).
\item Masoud Asgharian (supervisor): Conceptualization (proposed the Linfoot measure concept for variable screening); Investigation (collaborated on design and execution of simulations and case studies); Writing - Review \& Editing (manuscript revisions); Supervision.
\item Nikhil Baghwat: Data Curation (handled the preprocessing of the ABIDE data).
\item Jean-Baptiste Poline: Resources (provided the ABIDE data).
\item Celia Greenwood (supervisor): Investigation (assisted in designing simulations and case studies); Writing - Review \& Editing (manuscript revisions); Resources (data provision); Supervision; Funding acquisition.
\end{itemize}
{\bf Manuscript 2:}
\begin{itemize}
\item Kai Yang (student): Formal analysis, Methodology, and Investigation (carried out all mathematical proofs and developments; developed the theoretical part, designed and carried out simulation studies); Visualization; Conceptualization (proposed the accelerated gradient approach); Writing - Original Draft (composed the initial draft); Writing - Review \& Editing (draft revisions).
\item Masoud Asgharian (supervisor): Formal analysis, Methodology, and Investigation (contributed to the proof of Theorem 2 on $O(1/k)$ convergence by proposing to use HM-GM inequality approach, designed the simulation studies); Writing - Review \& Editing (edited the manuscript); Supervision.
\item Sahir Bhatnagar (supervisor): Conceptualization (proposed SCAD/MCP in the original problems to be solved); Investigation (designed the simulation studies); Writing - Review \& Editing (edited the manuscript); Supervision.
\end{itemize}
{\bf Manuscript 3:}
\begin{itemize}
    \item Kai Yang (student): Formal analysis and Methodology (carried out all mathematical proofs and developments, developed the theoretical part, developed the optimization framework based on variational and nonsmooth analysis and the conjugate gradient method); Conceptualization (proposed to use the conjugate gradient method); Writing - Original Draft (wrote the initial draft); Writing - Review \& Editing (subsequent revisions).
    \item Masoud Asgharian (supervisor): Conceptualization (proposed the use of Tsallis entropy and $q$Gaussian distribution for modeling); Writing - Review \& Editing (reviewed and revised the manuscript); Supervision.
    \item Celia Greenwood (supervisor): Writing - Review \& Editing (reviewed and revised the manuscript); Supervision; Funding acquisition.
\end{itemize}
This doctoral thesis presents original scholarship and distinct contributions to knowledge, specifically in the area of statistical computing and robust statistical methods for high-dimensional data analysis. The core contributions of this work, which advance knowledge within the field, are the development of new theories and methodologies, comprehensive simulation results, and data analyses. These contributions are thoroughly detailed in the chapters within.
}
\Preface




\SetAbstractEnText{ 


A ubiquitous feature of biological data of our era, such as brain
functional magnetic resonance imaging or genetic data, is their extra-large
sizes and dimensions. However, analyzing such high--dimensional biological
data poses significant challenges, since the feature dimension is
often much larger than the sample size. This thesis introduces robust
and computationally efficient methods to address several common challenges associated
with high--dimensional data.

In my first manuscript, I propose a coherent approach to variable
screening that can accommodate nonlinear associations. I develop a
novel variable screening method that transcends traditional linear
assumptions by leveraging mutual information, with an intended application
in neuroimaging data. This approach allows for a more accurate identification
of important variables by capturing nonlinear as well as linear relationships
between the outcome and the covariates. This strategy proves to be
transformative in the analysis of neuroimaging data, as demonstrated
through a detailed examination of the prepossessed Autism Brain Imaging
Data Exchange dataset \citep{Cameron2013, Barry2020}.

Then, building on this foundation, I develop new computing techniques
for sparse estimation using nonconvex penalties in my second manuscript.
These methods address notable challenges in current statistical computing
practices, facilitating computationally efficient and robust analyses of complex datasets.
While my study in the second manuscript is mainly motivated by computational
challenges in sparse estimation using nonconvex penalties, the proposed
method can be applied to a considerably general class of optimization
problems.

In my third manuscript, I contribute to the development of robust
modeling of high--dimensional correlated observations by relaxing some
of the underlying assumptions for the analysis of such data. I develop
a $q$Gaussian linear mixed-effects model, designed to surpass the constraints
of conventional Gaussian linear mixed-effects models by accommodating
a broader class of distributions that are more robust toward outliers.
For correlated observations, this $q$Gaussian model enhances the robustness
and flexibility of statistical analyses, providing a more comprehensive
tool for modeling the widely-correlated observations frequently encountered
in biological and medical studies.

Collectively, these contributions aim at addressing the multifaceted
challenges of high--dimensional biological data analysis and paving
the way for deeper insights into complex biological systems by seamlessly
integrating solutions to nonlinearity, nonconvex nonsmooth optimization,
and the need for more robust and adaptable models.

}
\AbstractEn%

\SetAbstractFrName{\MakeUppercase{ABR\'{E}G\'{E}}}%
\SetAbstractFrText{ 


Une caract\'{e}ristique omnipr\'{e}sente des donn\'{e}es biologiques de notre \'{e}poque, telles que l’imagerie par r\'{e}sonance magn\'{e}tique fonctionnelle du cerveau ou les donn\'{e}es g\'{e}n\'{e}tiques, est leur taille et leur dimension extra-larges. Cependant, l’analyse de ces donn\'{e}es biologiques \`{a} haute dimension pose des d\'{e}fis importants, car la dimension des caract\'{e}ristiques est souvent beaucoup plus grande que la taille de l’\'{e}chantillon. Cette th\`{e}se introduit des m\'{e}thodes robustes et efficaces pour r\'{e}pondre \`{a} plusieurs d\'{e}fis communs associ\'{e}s aux donn\'{e}es de haute dimension.

Dans mon premier manuscrit, je propose une approche coh\'{e}rente de la s\'{e}lection des variables qui peut prendre en compte les associations non lin\'{e}aires. Je d\'{e}veloppe une nouvelle m\'{e}thode de s\'{e}lection des variables qui transcende les hypoth\`{e}ses lin\'{e}aires traditionnelles en utilisant le concept de l’information mutuelle. Cette approche permet une identification plus pr\'{e}cise des variables importantes en capturant les relations non seulement lin\'{e}aires mais aussi non lin\'{e}aires entre le r\'{e}sultat et les covariables. Cette strat\'{e}gie s’av\`{e}re transformatrice dans l’analyse des donn\'{e}es de neuro-imagerie, comme le montre mon analyse de donn\'{e}es d’imagerie c\'{e}r\'{e}brale sur l’autisme \citep{Cameron2013, Barry2020}.

Ensuite, en m’appuyant sur cette base, je d\'{e}veloppe de nouvelles techniques de calcul pour la s\'{e}lection de variables parcimonieuse en utilisant des p\'{e}nalit\'{e}s non convexes dans mon deuxi\`{e}me manuscrit. Ces m\'{e}thodes abordent des d\'{e}fis notables dans les pratiques actuelles de calcul statistique, facilitant des analyses efficaces et robustes d’ensembles de donn\'{e}es complexes. Alors que mon \'{e}tude dans le deuxi\`{e}me manuscrit est principalement motiv\'{e}e par les d\'{e}fis de calcul dans l’estimation parcimonieuse utilisant des p\'{e}nalit\'{e}s non convexes, la m\'{e}thode propos\'{e}e peut être appliqu\'{e}e \`{a} une classe tr\`{e}s g\'{e}n\'{e}rale de probl\`{e}mes d’optimisation.

Dans mon troisi\`{e}me manuscrit, je contribue au d\'{e}veloppement d’une mod\'{e}lisation robuste d’observations corr\'{e}l\'{e}es en haute dimension, en assouplissant certaines des hypoth\`{e}ses sous-jacentes pour l’analyse de telles donn\'{e}es. Je d\'{e}veloppe un mod\`{e}le lin\'{e}aire mixte $q$Gaussien, conçu pour d\'{e}passer les contraintes des mod\`{e}les lin\'{e}aires mixtes Gaussiens conventionnels. Ceci permet une adaptation \`{a} une classe plus large de distributions qui sont plus robustes vis-\`{a}-vis des valeurs aberrantes. Pour les observations corr\'{e}l\'{e}es, ce mod\`{e}le $q$Gaussien est robuste et flexible , fournissant un outil plus complet pour mod\'{e}liser les observations largement corr\'{e}l\'{e}es, qui sont fr\'{e}quemment rencontr\'{e}es dans les \'{e}tudes biologiques et m\'{e}dicales.

Collectivement, ces contributions visent \`{a} relever les d\'{e}fis \`{a} multiples aspects de l'analyse des donn\'{e}es biologiques \`{a} haute dimension, et \`{a} ouvrir la voie \`{a} une meilleure compr\'{e}hension des syst\`{e}mes biologiques complexes en int\'{e}grant de mani\`{e}re transparente des solutions \`{a} la non-lin\'{e}arit\'{e}, \`{a} l'optimisation non convexe et non lisse, et \`{a} la n\'{e}cessit\'{e} de mod\`{e}les plus robustes et adaptables.

}%
\AbstractFr%


\tableofcontents 

\LOFHeading{List of Figures}
\listoftables 
\listoffigures

\newpage

\printacronyms[include=abbrev, name=Abbreviations]
\newpage

\end{romanPagenumber}


\doublespacing

\chapter{Introduction}\label{ch:intro}

\acresetall






In the domain of biostatistics, the prevalence of high--dimensional biological data stands as a testament to the field's intricate relationship with complex datasets, notably within genetic research and brain neuroimaging. The breadth and complexity of these data landscapes emphasize the vital role of biostatistics in deciphering meaningful scientific insights from extra large high--dimensional datasets. 

High--dimensional genetic data reflect a wealth of information about individual susceptibilities to diseases, physiological traits, and other critical biological attributes. This intricate dataset has been the foundation for numerous groundbreaking studies aimed at deciphering the molecular underpinnings of diseases, subsequently leading to innovative therapeutic approaches. \ac{GWAS} have been instrumental in identifying genetic factors that contribute to the biology of diseases, thus paving the way for new therapeutic developments \citep{Visscher2017}. Moreover, the interpretative analysis of statistical genetic models has enriched our understanding of heritability \citep{Yang2017}. The application of genetic information to identify individuals at an elevated risk of specific diseases enhances disease screening strategies \citep{Chatterjee2016}, while genome analysis initiatives have refined diagnostic and screening processes for complex disorders, illustrating the capacity of genetic data to revolutionize healthcare practices \citep{Khera2017, Pashayan2015}. Technology advances in the last decade have dramatically increased the volume and complexity of genetic datasets. For example, the UK Biobank project, with its extensive collection of genetic variants from approximately half a million individuals, embodies this evolution, presenting more than $800,000$ attributes of unique genetic markers \citep{Bycroft2018}. 

In parallel, neuroimaging data emerge as another example of high--dimensional biomedical datasets. The complexity and high dimension of neuroimaging data have catalyzed advances in variable selection techniques, as evidenced by a notable increase in related research publications: \citep{Adeli2017,Fan2016,Febles2022,GomezVerdejo2019,Hao2020,He2018,Hunt2014,Ivanoska2021,Mohr2006,Martino2008,Pereda2018,Roy2021,Schloegl2002,Sofer2014,Suresh2022}. Acquisition of magnetic resonance images (MRI) produces data on an unprecedented scale, capturing measurements in millions of voxels \citep{Bell2018, Liang2022, Linn2016, Fan2016}. The advent of multiple imaging modalities has introduced multiple sets of high--dimensional features, each providing different insights into brain function and exhibiting complex correlation patterns. This multiplicity of data accentuates the critical need for sophisticated analytical techniques capable of managing and interpreting the intricate details captured within and across these modalities.

The analysis of large high--dimensional biological datasets, common in fields such as genomics and neuroimaging, presents ultimate challenges in statistical computing. Often, these datasets are so voluminous that they exceed available memory capacity, necessitating strategies for dimension reduction to perform statistical analysis on the data. In this context, feature selection emerges as a crucial technique. Unlike other dimension reduction methods such as \ac{PCA} and \ac{ICA}, univariate variable screening stands out for its computational efficiency. Additionally, univariate variable screening adapts to limited memory resources, as it processes only the outcome and a single covariate at each iteration, making it especially suitable for analyzing extensive datasets. Moreover, it offers the advantage of straightforward interpretability; the variables selected through this process directly correspond to features of interest, providing clear insights without the obfuscation that can accompany other dimensionality reduction techniques.

Another critical benefit of univariate variable screening is its compatibility with parallel computing frameworks. This adaptability allows the simultaneous processing of data segments, significantly speeding up the variable screening step of high--dimensional large datasets that are typical in genetic research and neuroimaging studies. Such computational efficiency is crucial in these fields, where rapid and effective interpretation of data can lead to significant scientific advancements.

Furthermore, when comparing univariate screening with multivariable selection methods, univariate approaches maintain consistency in variable selection. This consistency stems from the fact that the calculated measure of the association between each covariate and the outcome is independent of the influence of other covariates. This feature ensures that the introduction of additional covariates into the analysis does not necessitate a re-evaluation of existing associations, a requirement that multivariable approaches cannot circumvent. In scenarios where new covariates are added to the dataset, univariate screening only requires the calculation of associations with the outcome for these new covariates, whereas multivariable dimension reduction or variable selection methods would need to reassess the entire dataset, including both established and newly incorporated covariates. This distinction underscores the practicality and computational efficiency of univariate variable screening in the dynamic environment of high--dimensional data analysis, making it an invaluable tool for researchers navigating the complexities of genetic studies and neuroimaging data.

Hence, in my first manuscript, I introduce a coherent approach to univariate variable screening that is robust to nonlinear associations. The variable screening methods in my first manuscript are incorporated in a Python pacakge {\tt fastHDMI}, which stands for \emph{Fast Mutual Information Estimation for high--dimensional Data}. This innovative tool consists of three mutual information estimation techniques for variable selection within neuroimaging analyses. Using extensive simulation studies based on the preprocessed \emph{\ac{ABIDE} dataset} \citep{Cameron2013, Barry2020}, my screening methods are evaluated under various conditions, highlighting the superiority of mutual information estimation through \emph{\ac{FFTKDE}} for variable screening when the continuous outcome is nonlinearly associated with the covariates, as well as the advantage of variable screening using mutual information estimation by binning continuous variables when the binary outcome is nonlinearly associated with the covariates. Furthermore, based on case studies to predict the continuous outcome age and the binary outcome autism diagnosis, my research showcases the package's capability in variable screening by comparing the performance of various predictive models built using the selected variables from screening, demonstrating {\tt fastHDMI}'s significant contribution to enhancing neuroimaging data analysis and expanding the repertoire of variable screening tools for researchers when it comes to high--dimensional data prevalent in biomedical studies.

Building upon the foundational work presented in my first manuscript, my second manuscript ventures into the realm of developing new statistical computing techniques for sparse estimation, specifically addressing the challenges posed by nonconvex penalties. These innovative methods tackle significant obstacles encountered in current statistical computing paradigms, enhancing the computational efficiency of analyzing high-dimensional large datasets. Central to this exploration is the adaptation of Nesterov's \ac{AG} method \citep{Nesterov1983, Nesterov2004a} to nonconvex nonsmooth settings --- a notable departure from its conventional application to convex nonsmooth penalties such as $\ell_{1}$ penalty \citep{Tibshirani1996} or the elastic net penalty \citep{Zou2005}. This adaptation is particularly crucial given the convergence challenges associated with nonconvex penalties such as \ac{SCAD} \citep{Fan2001} and \ac{MCP} \citep{Zhang2010a}. This adaption is established upon the methodologies outlined in \citep{Ghadimi2015}, setting a foundation for the algorithmic analysis and development presented in this manuscript. 

My second manuscript details a sophisticated algorithm focused on the selection of critical optimization hyperparameters, pivotal for its practical implementation. It delves into the intricacies of selecting these hyperparameters, proposing a strategy based on complexity upper bounds to accelerate convergence, thereby making a significant contribution to sparse learning in a high--dimensional context. Furthermore, by establishing the rate of convergence and presenting a novel bound to describe the optimal damping sequence, this work not only underscores the algorithm's theoretical underpinnings but also demonstrates its superior performance over existing methods through comprehensive simulation studies by nonconvex penalized linear and logistic models. This manuscript, while primarily motivated by computational challenges in sparse estimation with nonconvex penalties, ultimately presents a methodology with broad applicability across a diverse spectrum of optimization problems, marking a significant step forward in the field of statistical computing. This manuscript has now been recognized and disseminated through its publication in the journal \emph{Statistics and Computing}, an achievement that highlights its contribution to the field \citep{Yang2024}. 

Biostatistical datasets often feature correlated observations, a notable example being genetic data, which inherently embodies structured correlation between observations \citep{Bycroft2018}. Neglecting population structure often leads to a considerable lack of fit: previous research demonstrates that the predictions obtained by the expectations of linear models do not predict as accurately as the maximum a posteriori (MAP) predictions obtained by linear mixed models (LMM), with the latter incorporating population structure \citep{Bhatnagar2019}. The population structure can also be a confounder for the phenotype and the genetic data; hence, it might cause spurious correlations discovered if not accounted for. Specific to variable selection, not accounting for population structure might cause some population-related variables falsely selected when they are not, in fact, related to the phenotype --- in this view, it might even cause true variables not selected. The motivation behind my third manuscript is driven by the need to address this issue, proposing a linear mixed-effects model based on the idea of Tsallis entropy maximization. This method effectively handles the correlation among observations, while also incorporating variable selection for fixed-effects covariates, utilizing sparse penalties that function as regularizers when the dimensionality of the design matrix surpasses the number of observations.

The developed $q$Gaussian linear mixed effects model marks a significant advance in statistical sparse learning, providing an approach to analyze high--dimensional and correlated observations robust to outliers and the underlying distributional assumption. This innovation addresses the limitations inherent in traditional Gaussian distribution assumptions that have historically constrained statistical analysis. Based on the principle of maximizing Tsallis entropy, the $q$Gaussian model excels in navigating the complexities of biostatistical data, characterized by correlated observations and heterogeneity of variances, a scenario frequently encountered in genetic and longitudinal data.

In my third manuscript, I re-derive the multivariate probability density function from Tsallis entropy maximization. This allows for statistical modeling using the likelihood-ist approach, overcoming the constraints imposed by conventional Gaussian assumptions, which often fall short in robustness towards outliers and the accurate representation of underlying distributional shapes. Furthermore, I introduce a novel framework that leverages numerous numerical methods originally designed to find equilibria in flows, thus addressing the composite optimization problems characteristic of statistical sparse learning. The framework is further applied to the state-of-the-art Hager-Zhang conjugate gradient algorithm \citep{Hager2005}, which yields a numerically stable and computationally efficient algorithm for sparse statistical learning. 



In essence, through the development of robust and computationally efficient methods, this thesis enhances the ability to model and predict using large high-dimensional datasets frequently encountered in biostatistics, such as in neuroimaging and genetics. The groundwork laid by this research promises to propel forward in statistical computing and robust modeling, setting the stage for future investigations that delve deeper into rich, uncharted territories of biomedical data.

\newpage

\chapter{Literature review}\label{ch:litt}

\noindent 
   
In this section, a summary of pertinent literature related to the thesis is provided. For an in-depth exploration of the literature, please consult the literature review sections within each of the three manuscripts included in this thesis. A motivating factor for the research presented in this dissertation stems from the challenge posed by high-dimensional datasets, where the number of features often surpasses the number of observations. This results in a row rank deficiency in the design matrix $\mathbf{X}$, leading to the null space $\text{null}\left(\mathbf{X}\right)\neq \emptyset$. The foundation of many statistical learning methods is the linear predictor $\mathbf{X}\boldsymbol{\beta}$, with the estimation of $\boldsymbol{\beta}$ parameters typically achieved through the minimization of an objective function. Such functions include least-square loss, robust objective functions such as Huber loss function, (negative) $\log-$likelihood, (negative) partial $\log-$likelihood, and the \ac{GMM}, among others. The existence of a nonempty null space indicates that the solutions to these minimization problems with respect to $\boldsymbol{\beta}$ are not uniquely defined, rendering the problem ill-posed. To address this, regularization via a strongly convex function, dimension reduction, or variable selection can be employed. For dimension reduction methods such as \ac{PCA}, \ac{ICA}, and autoencoders \citep{Hinton2006}, as discussed in the Introduction chapter, these methods exhibit certain limitations compared to variable selection. As elucidated previously in the text, variable selection through the application of penalties can function as localized regularization effects under certain conditions.

\section{Mutual Information}

Mutual information is defined as the Kullback--Leibler divergence between the joint distribution of two variables and their outer product distribution, as described in \eqref{eq:MI-defn-thesis-textbook-review}.

\begin{equation}
    I\left(X,Y\right)\coloneqq\mathbb{E}_{\left[X,Y\right]^{T}}\left[-\log\frac{p\left(X\otimes Y\right)}{p\left(\left[X,Y\right]^{T}\right)}\right]
    \label{eq:MI-defn-thesis-textbook-review}
\end{equation}

This concept can quantify dependencies without assuming a linear relationship between variables, in contrast to conventional methods like Pearson's correlation, as shown in \eqref{eq:pearson-corr-defn-thesis-textbook-review} for vector $x,y$ being the realizations drawn from the random variables $X,Y$, which assumes linearity and does not perform well when there are nonlinear associations.

\begin{equation}
    r\left(x,y\right)\coloneqq\left\langle \frac{x}{\left\Vert x\right\Vert _{2}},\frac{y}{\left\Vert y\right\Vert _{2}}\right\rangle 
    \label{eq:pearson-corr-defn-thesis-textbook-review}
\end{equation}

In comparison, mutual information-based variable screening approaches are robust towards nonlinear associations. Previous research has introduced some methods to evaluate the association between outcome and covariates \citep{Renyi1959, Reshef2011, Speed2011}, where the measures of association detailed in these studies are all monotonically increasing functions of mutual information. Hence, variable screening using any of these association measures yield identical results to variable screening using mutual information.

In the domain of neuroimaging data analysis, mutual information has been extensively employed, demonstrating its versatility and efficacy in deciphering the intricacies of neural datasets. Noteworthy applications include the use of Gaussian copula for mutual information estimation in continuous datasets \citep{Ince2016, Magri2009}, the application of mutual information in fMRI data analysis \citep{Nemirovsky2023, Tsai1999}, and the exploration of EEG-based brain-computer interfaces through mutual information \citep{Schloegl2002}. Despite its widespread application, there remains a gap in neuroimaging research with respect to the use of mutual information for feature screening when the variables are continuous, primarily due to the challenges associated with estimating mutual information for continuous variables.

Estimating mutual information for discrete variables is straightforward; however, estimating mutual information for continuous variables involves a variety of methodologies. Previous studies have considered techniques such as the binning of continuous variables to transform them into discrete variables \citep{Ross2014}, kernel density estimation (KDE) \citep{Steuer2002, Moon1995, Khan2007, Gao2014}, and $k$-nearest-neighbor estimation ($k$NN) \citep{Faivishevsky2008,Kraskov2004,Victor2002,Pal2010,Lord2018,Gao2014} strategies. KDE-based methods, in particular, have shown superior performance in mutual information estimation, especially in settings with small sample sizes and high noise levels \citep{Steuer2002, Moon1995, Khan2007, Gao2014}. For binning estimation of mutual information, choosing the number of bins is critical. Previous literature \citep{Birge2006} suggests an approach to find the optimal number of bins based on Castellan's bounds on risk of penalized maximum likelihood estimators \citep{Castellan2000}, which therefore enables a data-driven number of bins for the estimation of mutual information. Furthermore, for a detailed description of mutual information estimation using the \emph{Fast Fourier Transform based Kernel Density Estimation (FFTKDE)} or $k$-nearest-neighbor ($k$NN), please refer to Appexdix \ref{apdx:method-contribution} of my first manuscript. 

\section{$\ell_1$-induced Sparse Learning}

Sparse learning, also called variable selection, has always been a major approach in multivariable statistical analysis of high--dimensional data. This approach assumes that only a small number of predictors are relevant to the outcome. The resulting statistical models usually perform better in terms of predictive accuracy and possibly also interpretability. For these reasons, sparse learning has received much attention in the statistical literature over the past two decades (for example, \citep{Tibshirani1996, Zou2005, Buehlmann2014}). Sparse learning is commonly accomplished by adding sparse penalties to the objective function, either in the main problem or the subproblems, to produce sparsity in the estimation of coefficients. Let $\ell_{p}$ denote the sequential space endowed by the sequential norm $\left\Vert x\right\Vert _{p}\coloneqq\left(\sum_{j\in\mathbb{N}_{>0}}\left|x_{j}\right|^{p}\right)^{\frac{1}{p}}$. In a context of sparse learning, the $\ell_p$ penalty refers to the penalty term being $\left\Vert \cdot \right\Vert _{p}$ multiplied by a positive penalty hyperparameter $\lambda$ to control the level of penalization. Specifically, $\ell_1$ is usually used to achieve sparsity, as Lagrangian duality reveals a geometric interpretation that the solution of $\ell_1$ penalized problems will be on the boundary of the $\ell_1$ ball of some radius. for any smooth\footnote{In this thesis, smooth denotes {\em first-order} continuous differentiable unless otherwise specified} function $f\left(\boldsymbol{\beta}^\prime \right)$ with $L$-Lipschitz continuous gradients, when penalized by the convex nonsmooth $\ell_1$ penalty, the resulting estimator, denoted by $\hat{\boldsymbol{\beta}^\prime}$, will need to satisfy the first-order (necessary) optimality condition for the objective function $f\left(\boldsymbol{\beta}^\prime\right)+\lambda \left\Vert \boldsymbol{\beta}^\prime\right\Vert _{1}$:
\begin{equation}\label{subgrad-fermat}
    -\nabla_{\boldsymbol{\beta}^{\prime}}f\left(\hat{\boldsymbol{\beta}^{\prime}}\right)\in\frac{\partial}{\partial\boldsymbol{\beta}^{\prime}}\left(\lambda\left\Vert \boldsymbol{\beta}^{\prime}\right\Vert _{1}\right)\left(\hat{\boldsymbol{\beta}^{\prime}}\right)
\end{equation}
where $\lambda>0$ controls the amount of penalization and $\frac{\partial}{\partial\boldsymbol{\beta}^{\prime}}\left(\lambda\left\Vert \boldsymbol{\beta}^{\prime}\right\Vert _{1}\right)$ denotes the subgradient set operator for the nonsmooth convex function $\lambda\left\Vert \boldsymbol{\beta}^{\prime}\right\Vert _{1}$. Let $\Tilde{\boldsymbol{\beta}^\prime}$ denote the estimator obtained by minimizing $f\left(\boldsymbol{\beta}^\prime\right)$ itself, we have that $\Tilde{\boldsymbol{\beta}^\prime}$ satisfies 
\begin{equation}\label{smooth-fermat}
    \nabla_{\boldsymbol{\beta}^{\prime}}f\left(\tilde{\boldsymbol{\beta}^{\prime}}\right)=0=-\nabla_{\boldsymbol{\beta}^{\prime}}f\left(\tilde{\boldsymbol{\beta}^{\prime}}\right)
\end{equation}
For each coordinate $\beta^{\prime}_i$, if the partial derivative satisfies $\left|\frac{\partial}{\partial\beta_{i}^{\prime}}f\left(0\right)\right|\leq\lambda$, we'll have $\hat{\beta^{\prime}_i}=0$ \citep{Tibshirani2011}. By Lipschitz continuity of the gradients, a sufficient condition for $\hat{\beta^{\prime}_i}=0$ is: $$\left|\tilde{\beta_{i}^{\prime}}\right|\leq\frac{\lambda}{L}$$ In view of \eqref{smooth-fermat}, as an interpretation, should a coefficient be close to $0$, the $\ell_1$ sparse penalty will set the coefficient estimator zero. The mechanism above to induce the sparsity of $\ell_1$ penalization is, in fact, used in almost all sparse penalties for statistical learning, which makes the objective function nonsmooth. Furthermore, \cite{Nikolova2000} suggested that to achieve sparsity by optimizing a penalized problem, the first derivative of the objective function must be discontinuous. In a statistical context, this implies the nonsmoothness of the sparse penalty.

\section{Penalties with Oracle Property}

{\em Oracle property}, originally proposed by \citeauthor{Fan2001}, is a useful property of statistical estimators for sparse learning \citep{Fan2001}: as an interpretation, the oracle property demonstrates that the asymptotic distribution of the estimator yielded by penalized loss function is the same as the asymptotic distribution of the unpenalized estimator based on the loss function fitted only on the true support. In view of \eqref{subgrad-fermat}, the estimator yielded by the $\ell_1$ penalty, $\hat{\boldsymbol{\beta}^\prime}$, will be biased even when $\beta_i$ is in the true support --- this implies that when applied to statistical learning problems, the penalty $\ell_1$ cannot yield oracle estimators. 

The vast majority of penalties used for sparse learning consist of $\ell_1$ as one of its components to induce sparsity. On the other hand, to satisfy the oracle property, a necessary condition is that the estimator yielded by penalized MLE is unbiased for large $\beta_i$ \citep{Nikolova2000}. That is, for a sparse penalty $p\left(\boldsymbol{\beta}_{-0}\right)\coloneqq\lambda\left\Vert \boldsymbol{\beta}_{-0}\right\Vert _{1}+q\left(\boldsymbol{\beta}_{-0}\right)$, we need to have \citep{Nikolova2000} $$\frac{\partial}{\partial\beta_{i}}p\left(\boldsymbol{\beta}_{-0}\right)=\lambda\cdot\text{sgn}\left(\beta_{i}\right)+\frac{\partial}{\partial\beta_{i}}q\left(\boldsymbol{\beta}_{-0}\right)=0\text{ as }\left|\beta_{i}\right|\rightarrow\infty$$ This suggests that data-independent sparse penalties with oracle property must be nonconvex. Two famous penalties, {\em smoothly clipped absolute deviation (SCAD)} \citep{Fan2001} and {\em minimax concave penalty (MCP)} \citep{Zhang2010}, have been shown to possess oracle properties under certain conditions. 

Based on the above discussion, for a data-independent sparse penalty to possess oracle property, the penalty must be nonsmooth and nonconvex. As a result, computation for the penalized MLE must accommodate for this deduced nonconvexity and nonsmoothness. 

\section{Past Approaches to Solve Nonconvex Nonsmooth Penalties}
In this paragraph, we summarize the past approaches of computation methods proposed for SCAD/MCP and other nonconvex nonsmooth penalties; that is, problem \eqref{eq:penalized-nonconvex-nonsmooth}, where $p_{\lambda,\gamma}$ is the penalty function depending on penalty hyper-parameters $\lambda,\gamma$, $f$ is the unpenalized loss function -- the intercept coefficient, $\beta_{0}$, is not penalized.

\begin{equation}
\min_{\boldsymbol{\beta}\in\mathbb{R}^{q+1}}f\left(\boldsymbol{\beta}\right)+\sum_{j=1}^{q}p_{\lambda,\gamma}\left(\beta_{j}\right)
\label{eq:penalized-nonconvex-nonsmooth}
\end{equation}

\citeauthor{Zou2008} proposed to perform a local linear approximation, which yields a descent majorization-minimization (MM) algorithm \citep{Zou2008}. \citeauthor{Breheny2011} proposed to use the coordinate descent method to carry out the estimation for linear models with least-square loss or logistic regression, penalized by SCAD and MCP \citep{Breheny2011}. \citeauthor{Mazumder2011} implemented in an {\tt R} package, {\tt sparsenet}, which carries out root-finding process in a coordinate manner \citep{Mazumder2011}. \citeauthor{Kim2008} discussed difference of convex programming (DCP) method for OLS estimators penalized by SCAD penalty \citep{Kim2008}, which was later generalized by \citeauthor{Wang2013} to a class of nonconvex penalties \citep{Wang2013}. \citeauthor{Lee2016} developed a modified second-order method originally designed for the loss function of \ac{LASSO} with extension to SCAD and MCP \citep{Lee2016}, this attempt was later extended to generalized linear models such as logistic or Poisson, and Cox's proportional hazard model \citep{Kim2018}. There are also a few other attempts to apply the quasi-Newton method or a mixture of first and second order descant method on the objective function with nonconvex penalties \citep{Ibrahim2012, Ghosh2016}. Rigorous proof of global convergence and the rate of convergence have rarely been established for these approaches; rather, most of them illustrated their convergence properties using simulated studies. Furthermore, for high--dimensional problems, second-order methods suffer from computational inefficiency when accounting for the computational cost in the evaluation of the secant condition. The first-order methods proposed above are prone to the behavior of ``zigzagging'' when the problem is ill--conditioned \citep{Watt2020}. For a smooth ill--conditioned problem, in view of a local quadratic approximation, the search direction experiences oscillations along the direction of eigenvectors corresponding to a greater absolute eigenvalue while moving very slowly towards the direction of eigenvectors corresponding to a less absolute eigenvalue, resulting in steps that ``zigzag,'' thus converges in a much slower speed numerically. To address this issue of ``zigzagging,'' accelerated gradient methods have been initially developed for smooth objective functions \citep{Polyak1964, Nesterov1983}, subsequently extended to nonsmooth convex problems \citep{Nesterov2004}, and more recently adapted to nonconvex and nonsmooth problems \citep{Ghadimi2015}. 



\section{$q$Gaussian Distribution}
The vast majority of distributions used in biostatistics can be considered as derived by maximizing Shannon's entropy under certain constraints \citep{Cover2006}. Specifically, the Gaussian distribution can be derived by maximizing Shannon's entropy with first-moment and second central-moment constraints. The Gaussian distribution is widely used in statistical machine learning, but suffers from several disadvantages, notably its exponential tail decay and lack of a shape parameter, which compromise robustness towards outliers and limits distribution shape representation. The $q$Gaussian distribution is derived from the maximization of the Tsallis entropy, and is a generalization of bell-curve distributions. Therefore, it emerges as a robust alternative capable of accurately modeling the diverse shapes of bell curve distributions and accounting for heavy-tailed characteristics, with wide usage, such as in modeling financial return data \citep{Borland2002,Borland2002a,Domingo2017}. Despite its proven advantages in finance, the $q$Gaussian distribution's application within biostatistics and statistical sparse learning is limited. 

\section{Existing Algorithms for Optimizing Sparse Learning Problems for Linear Mixed-effects Models}
When it comes to variable selection of the fixed-effects covariates in a context of linear mixed-effects models, there have been multiple approaches related to penalized LMMs computation. Most of these previous approaches were based on convex penalties, including \ac{LASSO} \citep{Tibshirani1996}, Adaptive \ac{LASSO}, or elastic net \citep{Zou2005} (for which the penalty is a linear combination of $\ell_1$ and $\ell_2$ norms of fixed effects coefficients). \citeauthor{Xiong2019} performed coordinate descent for adaptive \ac{LASSO} \citep{Xiong2019}, \citeauthor{Schelldorfer2011} built the {\tt R} package {\tt lmm\ac{LASSO}} for \ac{LASSO}/adaptive \ac{LASSO} based on proximal quasi-Newton method with Armijo rule \citep{Schelldorfer2011}, \citeauthor{Wang2018a} chose to use proximal gradient descent method for sparse penalties \citep{Wang2018a}. \citeauthor{Pan2018} used proximal Newton-Ralphson method for adaptive \ac{LASSO} \citep{Pan2018}. However, convex methods might not retain global convergence when applied to nonconvex problems. 

Proximal methods emerge as pivotal strategies for their effectiveness in handling sparse-induced nonsmooth optimization problems. These methods achieve superior numerical performance over alternative nonsmooth optimization methods (see, for example, \citep{Yu2017,Li2016}). 

One of the most basic proximal methods, the proximal gradient methods, is formulated as
\[
\min_x g(x) + h(x),
\]
where $g$ is a globally smooth function and $h$ is a convex, possibly nonsmooth function. At each iteration $k$, the method updates the variable $x$ according to the equation:
\[
x^{(k+1)} = \text{prox}_{\alpha^{\left(k\right)} h}(x^{(k)} - \alpha^{\left(k\right)} \nabla g(x^{(k)})),
\]

where \( \alpha^{\left(k\right)} \) is the step size and \( \text{prox}_{\alpha^{\left(k\right)} g} \) is the proximal operator of \( g \) parameterized by \( \alpha^{\left(k\right)} \). For a more detailed review of variational and nonsmooth analysis and proximal operators, please refer to Section \ref{subsec:variational-analysis} of the third manuscript. Simultaneously, various numerical algorithms are examined and utilized in the context of dynamic systems, particularly to find the equilibria of flows. Please refer to Section \ref{sec:dynamical-system} for a brief introduction of dynamical systems. Such capabilities are extensively documented in various scholarly resources on numerical analysis \citep{Quarteroni2007,Atkinson1989,Lubich2006,Hubbard1995,Helmke1994,Ross2019,Riahi2018}. In my third manuscript, I develop a method to transform the vast majority of numerical methods to find equilibria of flows to a optimization algorithm for nonconvex composite problems with the smooth component being globally Lispchitz-smooth. Moreover, the global convergence of such numerical method is preserved under this transformation. 

At the same time, Krylov subspace methods are recognized as the foundational pillars in numerical analysis, providing computationally efficient solutions for large-scale optimization problems \citep{Saad2003}, with excellent convergence acceleration and enhanced numerical stability. The Krylov subspace methods aim to solve the linear system $Ax=b$, which involve constructing a sequence of subspaces, known as Krylov subspaces, which are defined as:
\[
\mathcal{K}_r(A, b) = \text{span}\{b, Ab, A^2b, \dots, A^{r}b\},
\]
where \(r\) is the order of the Krylov subspace; $K_r$ is known as an order-$r$ Krylov subspace. For a detailed explanation of Krylov subspace methods, please refer to Section \ref{sec:Krylov-subspace}. These methods solve linear system \(Ax = b\) by iteratively improving an estimate of the solution or eigenvalue/eigenvector, leveraging the properties of the Krylov subspace to minimize computational effort while increasing the accuracy of the solution with each iteration.

The nonlinear conjugate gradient methods, which were developed based on the conjugate gradient method, are highlighted for their excellence in smooth optimization. This is attributed to its computational and memory efficiency, scalability, and numerical stability, making it a method of choice in the optimization landscape. Numerous previous literature has proposed various conjugate gradient methods; among these, the Hager-Zhang conjugate gradient method \citep{Hager2005}, when applied to a globally smooth problem, not necessarily convex, achieves global convergence and has shown good numerical results \citep{Hager2006}. Building on this algorithm, I applied the framework mentioned above to Hager-Zhang conjugate gradient method, which can achieve fast and numerically stable results for composite optimization problems very often encountered in a sparse learning/variable selection context. 

\section{Krylov Subspace Methods \label{sec:Krylov-subspace}}

This section provides an overview of Krylov subspace methods, which aids in comprehending the nonlinear conjugate gradient method discussed in Manuscript 3. A substantial portion of the content presented in this section is based on the notes taken during my numerical analysis course, wherein \citep{Trefethen2022} served as one of the primary references.

\paragraph{Arnoldi Iteration }

The Krylov subspace methods are best understood based on the idea
of projecting onto Krylov subspaces. Given matrix $A\in\mathbb{R}^{m\times m}$
and vector $b$, \emph{Krylov sequence }refers to the set of vectors
$b,Ab,A^{2}b,\dots$, and the corresponding \emph{Krylov subspaces}
of order $r$ is then defined as the space spanned by the first $r$
terms of Krylov sequence. \emph{Arnoldi iteration} then can be interpreted
as performing (modified) Gram-Schmidt on the Krylov matrix 
\[
\mathcal{K}_{n}\coloneqq\left[\begin{array}{ccccc}
b & Ab & A^{2}b & \cdots & A^{n-1}b\end{array}\right]\in\mathbb{R}^{m\times n}
\]
to construct its orthonormal basis. A matrix is in Hessenberg form
if it is \textquotedbl almost\textquotedbl{} triangular: all elements
below the first sub-diagonal are zero. In view of $A$ itself, Arnoldi
iteration can be considered as an Hessenberg-ized method analogous
to Gram-Schmidt, see Table \ref{tab:Householder's-reflection-vs}
-- one similarity is, they both can stop at any iteration with a
sequence of triangular/Hessenberg factors and a partial orthogonalized
factor $Q^{\left(k\right)}$, therefore serves as a better iterative
method. \emph{Householder's reflector} is a method to numerically
compute QR-decomposition. Different from how givens rotations method
\emph{rotates} the vector to zeroing an entry, Householder's method
will \emph{reflect} the vector by a hyperplane $H$ such that the
reflection can point to the desired direction -- reflecting one column
vector of $A$ at a time. For example, for the first column vector
of $A$, $a_{1}\in\mathbb{R}^{n}$, we try to reflect $a_{1}$ to
the direction of $e_{1}$ by left-multiplying an orthogonal matrix
$Q_{1}$ such that $Q_{1}a_{1}=\left\Vert a_{1}\right\Vert e_{1}$,
where $e_{1}$ denotes the vector with the first entry being $1$
and the rest of entries being $0$ per usual; the hyperplane $H$
is set orthogonal to $v\coloneqq\left\Vert a_{1}\right\Vert e_{1}-a_{1}$,
therefore the orthogonal matrix can be constructed by 
\[
Q_{1}\coloneqq I-2\frac{vv^{T}}{v^{T}v}
\]
where $Q_{1}$ is orthogonal.

\begin{table}[H]
\caption{\label{tab:Householder's-reflection-vs}Householder's reflection vs Gram-Schmidt/Arnoldi process}
\centering
\begin{tabularx}{\textwidth}{|>{\centering\arraybackslash}X|>{\centering\arraybackslash}X|>{\centering\arraybackslash}X|}
\hline 
 & \textbf{QR factorization} $A=QR$ & \textbf{Hessenberg formation} $A=QHQ^{*}$\\
\hline 
\textbf{Householder's reflection} & Orthogonal triangularization & Orthogonal Hessenberg formation\\
\hline
\textbf{Gram-Schmidt/Arnoldi process} & Triangular orthogonalization & Hessenbergized orthogonalization\\
\hline 
\end{tabularx}
\end{table}

For iterative methods we consider $m$ to be large or infinite, so
we only consider the first $n$ columns of $AQ=QH$. Let $Q_{n}\in\mathbb{R}^{m\times n}$
denote the first $n$ columns of $Q$; and let $\tilde{H}_{n}\in\mathbb{R}^{\left(n+1\right)\times n}$
be the submatrix located at the upper-left corner of $H$, which will
also be a Hessenberg matrix itself. Then we'll have $AQ_{n}=Q_{n+1}\tilde{H}_{n}$
as the first $n$ columns of $AQ=QH$. \emph{And equating the $n$th
column of both sides gives us $Aq_{n}=h_{1n}q_{1}+\cdots+h_{nn}q_{n}+h_{n+1,n}q_{n+1}$,
which is a recurrence relation for $q_{n+1}$ -- Arnoldi iteration
follows directly on this recurrence relation: let $q_{1}=\frac{b}{\left\Vert b\right\Vert }$
be the initializer, and choose $h_{kn}$ such that $h_{kn}q_{k}$
is a projection of $q_{k}$ on $Aq_{n}$ for $k=1,2,\dots,n$; as
an interpretation, the updating step first subtracts the projections
of the built orthogonal bases from $Aq_{n}$, then normalizing the
reminder with $h_{n+1,n}$ to ensure $\left\Vert q_{n+1}\right\Vert =1$.}
Because the recurrence formula states that each $q_{n}$ is formed
by a linear combination of $Aq_{n-1}$ and $q_{1},q_{2},\dots,q_{n-2}$,
each $q_{n}$ is therefore a degree-$\left(n-1\right)$ polynomial
of $A$ times $b$; hence $q_{1},q_{2},\dots,q_{n}$ form an orthonormal
basis for the Krylov subspace 
\[
\mathcal{K}_{n}\coloneqq\left\langle b,Ab,\dots,A^{n-1}b\right\rangle 
\]
\emph{(i). In this view, Arnoldi process can be considered as systematic
construction of orthonormal bases for successive Krylov subspaces
$\mathcal{K}_{1},\mathcal{K}_{2},\mathcal{K}_{3},\dots$.} Because
Arnoldi iteration constructs orthonormal basis in a Gram-Schmidt manner,
the $Q_{n}$ here will be exactly the same as the $Q_{n}$ present
in the Gram-Schmidt QR factorization of $K_{n}$, while here $K_{n}$
and $R$ \emph{per se} are never explicitly constructed. And it is
called \emph{modified} Gram-Schmidt because at iteration $k$, we
subtract projections of constructed bases $q_{1},q_{2},\dots,q_{k}$
from the vector $Aq_{k}$ instead of the ``original'' vector $A^{k}b$. 

\emph{(ii). Another view of Arnoldi process is a computation of projections
onto successive Krylov subspaces.} Note that $Q_{n}^{*}Q_{n+1}$ is
a $n\times\left(n+1\right)$ matrix with $1$ on the diagonal and
$0$ elsewhere; then from $AQ_{n}=Q_{n+1}\tilde{H}_{n}$ we have 
\[
\underbrace{Q_{n}^{*}Q_{n+1}\tilde{H}_{n}}_{\eqqcolon H_{n}}=Q_{n}^{*}AQ_{n}.
\]
Apparently, $H_{n}$ here will be the $n\times n$ submatrix located
at the upper-left corner of $H$. This is an analogue to change of
basis, with $Q_{n}$ not orthogonal but of shape $m\times n$ --
and the resulting interpretation is: given some $v\in\mathcal{K}_{n}$,
applying $A$ to it, then orthogonally project $Av$ back to $\mathcal{K}_{n}$. 

Note that here $H_{n}$ and $A$ are \emph{pseudo-}similar. Intuitively,
one might then consider the eigenvalues of $H_{n}$ as estimates for
the eigenvalues of $A$ -- for this reason, they are called \emph{Arnoldi
eigenvalue estimates} (at step $n$) or \emph{Ritz values }(wrt. $\mathcal{K}_{n}$). 

Consider a vector $x\in\mathcal{K}_{n}$, such a vector can then be
written as a linear combination of Krylov's vectors $b,Ab,\dots,A^{n-1}b$,
put in polynomial form, it will be 
\[
x=q\left(A\right)b
\]
Now consider $P^{n}\coloneqq\left\{ \text{monic polynomials of degree }n\right\} $,
the famous \emph{Arnoldi-Lanczos approximation problem} is proposed
as 
\[
\min_{p^{n}\in P^{n}}\left\Vert p^{n}\left(A\right)b\right\Vert 
\]
and the Arnoldi iteration solves this problem exactly (if it doesn't
break down ofc...) -- \emph{the minimizer $\bar{p}^{n}$ is uniquely
given by the characteristic polynomial of $H_{n}$.} As a proof, let
$y\coloneqq A^{n}b-p^{n}\left(A\right)b\in\mathcal{K}_{n}$, then
the problem can be considered as minimizing $\left\Vert A^{n}b-y\right\Vert $
wrt. $y$; i.e., minimizing the distance from $A^{n}b$ to $\mathcal{K}_{n}$
-- thus the minimization can be characterized by $p^{n}\left(A\right)b\perp\mathcal{K}_{n}\Leftrightarrow Q_{n}^{*}p^{n}\left(A\right)b=0$
as $q_{1},q_{2},\dots,q_{n}$ are a basis of $\mathcal{K}_{n}$. Now
consider $A=QHQ^{*}$; where $Q\coloneqq\left[\begin{array}{cc}
Q_{n} & U\end{array}\right]$ such that $Q$ is a orthogonal matrix extended from $Q_{n}$, and
$H\coloneqq\left[\begin{array}{cc}
H_{n} & X_{2}\\
X_{1} & X_{3}
\end{array}\right]$, where the entries of $X_{1}$ is all $0$ besides its upper-right
entry and $X_{3}$ is Hessenberg -- due to the Hessenberg structure
of $H$. Then we have 
\begin{align}
 & Q_{n}^{*}p^{n}\left(A\right)b=0\nonumber \\
\Leftrightarrow & Q_{n}^{*}Qp^{n}\left(H\right)Q^{*}b=0\nonumber \\
\Leftrightarrow & \left[\begin{array}{cc}
I_{n} & 0\end{array}\right]p^{n}\left(H\right)e_{1}\left\Vert b\right\Vert =0\label{eq:arnoldi-lanczo-prob}
\end{align}
and \eqref{eq:arnoldi-lanczo-prob} follows from $q_{1}=\frac{b}{\left\Vert b\right\Vert }$.
The interpretation of last equation is, the minimization characterization
now becomes that the first $n$ entries in the first column of $p^{n}\left(H\right)$
are $0$. Due to the Hessenberg structure of $H$, the first $n$
entries in the first column of $p^{n}\left(H\right)$ are exactly
the first column of $p^{n}\left(H_{n}\right)$ -- in view of this,
it is \emph{sufficient} to make $p^{n}\left(H_{n}\right)=0$: by Cayley-Hamilton
theorem, if $p^{n}$ is the characteristic polynomial of $H_{n}$,
$p^{n}\left(H_{n}\right)=0$. Proof of uniqueness uses contradiction:
if uniqueness is voided, taking difference of two distinct degree-$n$
monic polynomials that both minimize $\left\Vert p^{n}\left(A\right)b\right\Vert $
will then result in a non-zero polynomial $q\left(A\right)$ of degree
$\leq n-1$ such that $q\left(A\right)b=0$ -- this contradicts the
assumption that $K_{n}$ is of full-rank. 

Based on this finding, \emph{(iii). the Ritz values generated by Arnoldi
iteration are the roots of the optimal polynomial to the Arnoldi-Lanczos
approximation problem.} And this gives the Ritz values some invariant
properties: 
\begin{itemize}
\item \emph{(translation invariance)} If $A$ is changed to $A+\sigma I$
for some $\sigma\in\mathbb{R}$, and $b$ is left unchanged, the Ritz
values $\left\{ \theta_{j}\right\} $ at each step will be changed
to $\left\{ \theta_{j}+\sigma\right\} $ 
\item \emph{(scale invariance)} If $A$ is changed to $\sigma A$ for some
$\sigma\in\mathbb{R}$, and $b$ is left unchanged, the Ritz values
$\left\{ \theta_{j}\right\} $ at each step will be changed to $\left\{ \sigma\theta_{j}\right\} $ 
\item \emph{(unitary similarity transformation invariance)} If $A$ is changed
to $UAU^{*}$ for some unitary $U$, and $b$ is changed to $Ub$,
the Ritz values do not change 
\end{itemize}

\paragraph{\ac{GMRES} }

\emph{\ac{GMRES}} is a method \emph{using
Arnoldi iteration to solve a linear system $Ax=b$, the resulting
mechanic is to use $x_{n}\in\mathcal{K}_{n}$ at step $n$ to approximate
the root by formulating the problem:} 
\begin{align}
 & \min_{x_{n}\in\mathcal{K}_{n}}\left\Vert Ax_{n}-b\right\Vert \nonumber \\
\Leftrightarrow & \min_{c\in\mathbb{R}^{n}}\left\Vert AK_{n}c-b\right\Vert \nonumber \\
\Leftrightarrow & \min_{y\in\mathbb{R}^{n}}\left\Vert AQ_{n}y-b\right\Vert \nonumber \\
\Leftrightarrow & \min_{y\in\mathbb{R}^{n}}\left\Vert Q_{n+1}\tilde{H}_{n}y-b\right\Vert \nonumber \\
\Leftrightarrow & \min_{y\in\mathbb{R}^{n}}\left\Vert \tilde{H}_{n}y-Q_{n+1}^{*}b\right\Vert \label{eq:GMRES-1}\\
\Leftrightarrow & \min_{y\in\mathbb{R}^{n}}\left\Vert \tilde{H}_{n}y-\left\Vert b\right\Vert e_{1}\right\Vert \label{eq:GMRES-2}
\end{align}
where \eqref{eq:GMRES-1} is because that $b$ is in the column space
of $Q_{n+1}$ (because $q_{1}\coloneqq\frac{b}{\left\Vert b\right\Vert }$),
therefore left multiplication of $Q_{n+1}^{*}$ does not change the
norm. Furthermore, note that $Q_{n+1}^{*}b=\left\Vert b\right\Vert e_{1}$,
which gives us \eqref{eq:GMRES-2}. 

On another note, the initial assumption for \ac{GMRES} of $x_{n}\in\mathcal{K}_{n}$
is equivalent to $x_{n}=q_{n}\left(A\right)b$ for some degree-$\left(n-1\right)$
polynomial $q_{n}$, with coefficients being $c$ mentioned in above
equations. Then the residual satisfies $b-Ax_{n}=\left(I-Aq_{n}\left(A\right)\right)b$;
let $p_{n}\left(z\right)\coloneqq1-zq\left(z\right)$, then \ac{GMRES}
in fact solves problem $\min_{p_{n}\in P_{n}}\left\Vert p_{n}\left(A\right)b\right\Vert $,
but with $$P_{n}\coloneqq\left\{ \text{degree}\leq n\text{ polynomials }p\text{ with }p\left(0\right)=1\right\}.$$

\paragraph{Lanczos Iteration and Conjugate Gradient}

If $A$ is symmetric, or Hermitian over the complex space, the Arnoldi iteration will be redundant to find
eigenvalues of $A$ -- a method called Lanczos iteration was introduced
as a simplification of Arnoldi iteration (\emph{mainly simplified
by noticing that $H_{n}$ becomes tri-diagonal now}). With a similar
simplification idea, if $A$ is symmetric positive definite, solving
$\min_{x}\left\Vert Ax-b\right\Vert $ using \ac{GMRES} will in fact not
be efficient -- \emph{\ac{CG}} was then introduced
based on \emph{minimizing the A-norm of the error}; where the \emph{$A$-norm
of $e_{n}\coloneqq x^{*}-x_{n}$ }is defined as $e_{n}^{T}Ae_{n}$.
Specifically, the famous CG is proposed as Algorithm \ref{alg:Conjugate-Gradient-(CG)}. 

\begin{algorithm}[h]
\caption{\label{alg:Conjugate-Gradient-(CG)}Conjugate Gradient (CG) }

\begin{algorithmic}[1] 	
\Require $A\in \mathbb{R}^{m\times m}\succ 0$, $b\in \mathbb{R}^m$ 	
\Ensure $x_n$ -- the solution of linear system $Ax=b$ 	
\State Set $x_0\leftarrow0$, $r_0\leftarrow b$, $p_0\leftarrow r_0$ \Comment{Initialization} 	
\While{not converged} 	
\State $\alpha_k \leftarrow \frac{r_{k-1}^Tr_{k-1}}{p_{k-1}^TAp_{k-1}}$ \Comment{calculate step length} 	
\State $x_k \leftarrow x_{k-1}+\alpha_kp_{k-1}$ \Comment{approximate solution} 	
\State $r_k \leftarrow r_{k-1} - \alpha_kAp_{k-1}$ \Comment{calculate residual} 	
\State $\beta_k \leftarrow \frac{r_k^Tr_k}{r_{k-1}^Tr_{k-1}}$ \Comment{calculate improvement from this step} 
\State $p_k \leftarrow r_k + \beta_kp_{k-1}$ \Comment{calculate next step's search direction} 	
\EndWhile 
\end{algorithmic}
\end{algorithm}

And induction on $n$ can show that: 
\begin{enumerate}
\item \emph{(identity of subspaces) }
\begin{align*}
\mathcal{K}_{n} & =\left\langle x_{1},x_{2},\dots,x_{n}\right\rangle =\left\langle p_{0},p_{1},\dots,p_{n-1}\right\rangle \\
 & =\left\langle r_{0},r_{1},\dots,r_{n-1}\right\rangle =\left\langle b,Ab,\dots,A^{n-1}b\right\rangle 
\end{align*}
\item \emph{(orthogonal residuals) }
\[
r_{i}^{T}r_{j}=0,\ \forall i\neq j
\]
\item \emph{($A$-conjugate search directions) }
\[
p_{i}^{T}Ap_{j}=0,\ \forall i\neq j
\]
\end{enumerate}
Following results above, for iteration $n$, we can show that $x_{n}$
is the unique point in $\mathcal{K}_{n}$ that minimizes $\left\Vert e_{n}\right\Vert _{A}$;
and the convergence is monotonic (descent property), i.e., 
\[
\left\Vert e_{n}\right\Vert _{A}\leq\left\Vert e_{n-1}\right\Vert _{A}
\]
and $e_{n}=0$ is achieved for some $n\leq m$. The first statement
follows some simple calculations, the monotonicity follows $\mathcal{K}_{n}\subset\mathcal{K}_{n+1}$. 

As CG minimizes $A$-norm of the error in an iterative manner, this
enables us to view CG as an optimization algorithm -- simple calculations
allow us to formulate the following problem for CG: 
\[
\min_{x\in\mathbb{R}^{m}}\frac{1}{2}x^{T}Ax-x^{T}b
\]
Lastly, similar to how we build the connection between Arnoldi iteration
and \ac{GMRES} in a polynomial minimization manner at the end of last section,
it is similar for CG: CG approximation problem can be formulated as
$\min_{p_{n}\in P_{n}}\left\Vert p_{n}\left(A\right)e_{0}\right\Vert _{A}$;
where $e_{0}\coloneqq x^{*}-x_{0}$ denotes the initial error, and
$P_{n}\coloneqq\left\{ \text{degree}\leq n\text{ polynomials }p\text{ with }p\left(0\right)=1\right\} $,
same as before. To conclude this section, it is worth noting that a plethora of nonlinear conjugate gradient methods have been derived from this original linear conjugate gradient \citep{Hager2006}.

\section{Brief Introduction on Dynamical Systems \label{sec:dynamical-system}}

This section presents a concise overview of certain dynamical system concepts, which enhances the understanding of the topics discussed in Manuscript 3. A substantial portion of the content presented in this section is based on the notes taken during the dynamical system course, wherein \citep{Berg2023} served as the primary reference.

For a topological space $X$, a \emph{flow} is a continuous map $\phi:\mathbb{R}\times X\mapsto X$
such that $\forall x\in X$, $t,s\in\mathbb{R}$, 
\begin{enumerate}
\item $\phi\left(0,x\right)=x$
\item $\phi\left(t,\varphi\left(s,x\right)\right)=\phi\left(t+s,x\right)$
\end{enumerate}
As mentioned in the third manuscript, the gradient flow is the flow
generated by the ordinary differential equation 
\[
\dot{x}=-\nabla f\left(x\right)
\]
for some smooth objective function $f$. The equilibrium point of
a dynamical system, $\bar{x}$, describes the steady state by setting
$\phi\left(\bar{x}\right)=0$. Consider two flows, $\phi_{1}:\mathbb{R}\times X_{1}\mapsto X_{1}$
and $\phi_{2}:\mathbb{R}\times X_{2}\mapsto X_{2}$, the homeomorphism
is a function $\varphi:X_{1}\mapsto X_{2}$ such that $\varphi$ is
bijective, continuous with continuous inverse, and $\forall t\in\mathbb{R},x\in X_{1}$,
$\varphi\left(\phi_{1}\left(t,x\right)\right)=\phi_{2}\left(t,\varphi\left(x\right)\right)$,
known as \emph{flow intertwining property.} Homeomorphism demonstrates
that two flows, or dynamical systems in general, are topologically
equivalent, which makes them important in analyzing the behavior of
dynamical systems. 

For the numerous numerical methods to find the equilibria of flow,
the foundation is Cauchy-Lipschitz theorem, as known as Picard--Lindel\"{o}f
theorem, the theorem states let $D\subseteq\mathbb{R}\times\mathbb{R}^{n}$
be closed and let $\left(t_{0},y_{0}\right)\in\text{int }D$. Let
$f:D\mapsto\mathbb{R}^{n}$ be a function that is continuous in $t$
and Lipschitz continuous in $y$; then $\exists\varepsilon>0$ such
that the initial value problem (IVP) 
\begin{align*}
\dot{y} & =f\left(t,y\left(t\right)\right)\\
y\left(t_{0}\right) & =y_{0}
\end{align*}
has a unique solution $y\left(t\right)$ on $\overline{B\left(t_{0},\varepsilon\right)}$, the closed ball centered at $t_{0}$ with radius $\varepsilon$.
This theorem ensures the existence and uniqueness of the underlying
dynamical system governed by an IVP. The vast majority of differential
equations can not be solved analytically; in this view, this theorem
gives the condition that the underlying dynamical system is unique,
further allowing various numerical methods to solve the system or to find
the equilibria of the dynamical system.

\section{Conclusion of Literature Review }

To conclude this chapter, it is worth reiterating that a more detailed literature review is available in each of the three manuscripts included in this thesis. Specifically, the \emph{Fast Fourier Transform-based Kernel Density Estimation (FFTKDE)} method as well as the $k$-nearest-neighbour method for estimating mutual information is elaborated in the Appendix \ref{apdx:method-contribution} of the first manuscript. Furthermore, a foundational overview of variational and nonsmooth analysis is presented in Section \ref{subsec:variational-analysis} of the third manuscript, providing the necessary theoretical underpinning.

\chapter{{\tt fastHDMI}: Fast Mutual Information Estimation for High-Dimensional Data}\label{ch:paper1}

\indent \textbf{Preamble to Manuscript 1.}
\subsection*{Introduction to the Study and Its Place in the Workflow:}
Manuscript 1 introduces {\tt fastHDMI}, a Python package that carries out variable screening, representing a significant advancement in the initial stage of high-dimensional data analysis. This study is crucial because it directly addresses the challenge of efficiently analyzing complex neuroimaging datasets, which are characterized by their high dimensionality. The manuscript's approach, which is based on three different mutual information estimation methodologies, ensures that only the most relevant variables are selected for further analysis while maintaining robustness to nonlinear association, thereby allowing the subsequent stages of modeling and interpretation.

\subsection*{Interconnection with Subsequent Research Phases:}
The methodologies developed in this manuscript provide a foundational tool for the entire analysis workflow outlined in the thesis. By successfully identifying key variables through {\tt fastHDMI}, researchers can ensure that the modeling phase, discussed in subsequent manuscripts, is based on the most pertinent data, thereby enhancing the efficacy and computational efficiency of the models. This initial screening is particularly vital given the nonlinear associations often present in biological data, which traditional linear screening methods might miss.

\subsection*{Enhancement of Neuroimaging Data Analysis:}
The utilization of the preprocessed Autism Brain Imaging Data Exchange (ABIDE) dataset to confirm the efficiency and computational efficiency of these screening techniques underlines the practical relevance of {\tt fastHDMI}. Through a thorough evaluation of the various mutual information estimation techniques included in the package, the manuscript demonstrates its suitability for the analysis of real--world neuroimaging data.

\subsection*{Contribution to the Broader Research Goals:}
The findings from Manuscript 1 significantly contribute to the overall objective of the thesis by enhancing our understanding of how to efficiently perform a screening of variables that are robust to nonlinear associations. This step is essential for the efficient processing and analysis of large datasets prevalent in biostatistics and underpins the subsequent methodological advances explored in Manuscripts 2 and 3. By establishing a robust and computationally efficient approach to variable screening, this manuscript ensures that the data fed into more complicated statistical models at the later stage is of the highest relevance and quality, thereby facilitating more robust and insightful analyses.

\subsection*{Transition to Manuscript 2:}
Building on the computational efficiency achieved in Manuscript 1, Manuscript 2 expands these concepts into the realm of sparse estimation using nonconvex penalties. The ability to screen variables effectively sets the stage for these advanced computational techniques, which are designed to handle the challenges in statistical computing posed by the high--dimensional data structures. This natural progression underscores the interconnectedness of the manuscripts, as each builds upon the previous findings to enhance the overall efficacy of biostatistical data analyses.


\newpage
\acresetall

\vspace*{2cm}

\begin{center}
	\Large{{\tt fastHDMI}: Fast Mutual Information Estimation for High-Dimensional Data}
\end{center}

\vspace*{3cm}

\begin{center} 
	Kai~Yang$^{1}$, Masoud~Asgharian$^{2}$, Nikhil~Baghwat$^{4}$, Jean-Baptiste~Poline$^{4}$, Celia~M.~T.~Greenwood$^{1,3}$
\end{center}

\vspace*{1cm}

\begin{center} 
	$^{1}$\textit{Department of Epidemiology, Biostatistics, and Occupational Health, McGill University}\\
        $^{2}$\textit{Department of Mathematics and Statistics, McGill University}\\        
        $^{3}$\textit{Lady Davis Institute for Medical Research, Montr\'{e}al}\\
        $^{4}$\textit{Department of Neurology and Neurosurgery, McGill University}\\
\end{center}

\vspace*{1cm}


\newpage

\section*{Abstract}
In this paper, we introduce {\tt fastHDMI}, a Python package for the efficient execution of variable screening for high--dimensional datasets, including neuroimaging datasets. This study marks the inaugural application of three distinct mutual information estimation methodologies for variable selection in the context of neuroimaging analysis, a novel contribution implemented through {\tt fastHDMI}. Such advancements are critical for dissecting the complex architectures inherent in neuroimaging datasets, offering refined mechanisms for variable selection against the backdrop of high dimensionality. Employing the preprocessed Autism Brain Imaging Data Exchange (ABIDE) dataset \citep{Cameron2013, Barry2020} as a foundation, we assess the efficacy of these variable screening methodologies through extensive simulation studies. These evaluations encompass a diverse set of conditions, including linear and nonlinear associations, alongside continuous and binary outcomes. The results delineate the {\em Fast Fourier Transform Kernel Density Estimation (FFTKDE)}-based mutual information estimation approach as preeminent for feature screening with continuous nonlinear outcomes, while the binning-based methodology is identified as superior for binary outcomes contingent on nonlinear underlying probability preimage. For linear simulations, a parity in performance is observed for continuous outcomes between the absolute Pearson correlation and FFTKDE-based mutual information estimation, with the former also exhibiting dominance in binary outcomes predicated on linear underlying probability preimage. A comprehensive case analysis utilizing the preprocessed Autism Brain Imaging Data Exchange (ABIDE) dataset further illuminates the applicative potential of {\tt fastHDMI}, demonstrating the predictive capabilities of models constructed from variables selected through our implemented screening methods. This research not only substantiates the computational prowess and methodological robustness of {\tt fastHDMI}, but also contributes significantly to the arsenal of analytical tools available for neuroimaging research. 

\newpage

\section{Introduction }\label{sec:introduction}

The question of how to best select a subset of variables from a large set is a commonly investigated topic in high--dimensional model fitting \citep{Chandrashekar2014}. This topic is often called ``variable selection'' in statistics, or ``feature selection'' in the machine learning world. Feature selection may be necessary either to fit a particular statistical model or, in some situations, because the data are too large for memory. Neuroimaging data provide a good example of such challenges. For example, \ac{MRI} result in measurements at millions of voxels \citep{Bell2018, Liang2022, Linn2016, Fan2016}, and the development of multiple imaging modalities is leading to multiple high--dimensional sets of features, each capturing a different aspect of brain function, that can show widespread correlation patterns within and between each modality.
These high dimensions in neuroimaging data have stimulated the development of variable selection methods; indeed, there has been a recent surge in publications: see, for example,
\citep{Adeli2017,Fan2016,Febles2022,GomezVerdejo2019,Hao2020,He2018,Hunt2014,Ivanoska2021,Mohr2006,Martino2008,Pereda2018,Roy2021,Schloegl2002,Sofer2014,Suresh2022}. These papers take a wide variety of strategies ranging from univariate to multivariate selection. Among these, \cite{Fan2016, Schloegl2002} considered absolute correlation or mutual information with respect to the outcome as a conventional univariate approach; selection based on sparse-inducing penalties on multi-variable models were proposed on the data \citep{Fan2016, Hao2020, Hunt2014, Roy2021} or transformed data \citep{Adeli2017}. Multivariate selection based on random forest variable importance \citep{Febles2022, Hao2020} or sign consistency from the support vector machine \citep{GomezVerdejo2019} has also been applied previously.  A ``potential support vector machine" was applied by \cite{Mohr2006}, an idea that rests on exchanging   the roles of data points and features. Another approach can be seen in \citep{Martino2008}, where they selected features recursively based on multivariate model fitting. Evidently, these papers take a wide variety of strategies ranging from simple methods like analyzing the direct absolute correlation between outcomes and features, to more complex approaches involving the use of univariate regression coefficients, univariate copulas, and techniques that leverage variable importance measures or sparse penalties in multivariate model fitting. Variable selection under a multivariable model generally requires certain assumptions, often including the assumption of linearity, which is not robust to misspecification. Furthermore, variable selection based on marginal associations demands less computational power and memory and can easily adapt to data inflow. Additionally, variable selection within a joint model framework allows for variable screening conditioned on other covariates, such as confounders. 

Although Pearson correlation is frequently used to measure the association between covariates and the outcome, in situations where nonlinearity may be present, a variety of strategies have been introduced to examine the relationship between the outcome and the covariates \citep{Renyi1959, Reshef2011, Speed2011}. These methods, when utilized for feature screening, effectively equate to screening via mutual information, as they are all deterministic monotonically increasing functions of mutual information. Among the strategies for feature selection, an entropy-based method, \emph{mutual information} has two appealing characteristics. As defined in \eqref{eq:mi-defn}, mutual information is defined as the \ac{KL divergence} between the joint distribution of two variables and their outer product distribution, effectively quantifying their dependency. This method can carry out model-independent feature selection, and is robust to non-linearity between the outcome and the features. For these reasons, mutual information has already been a popular choice for neuroimaging data. \cite{Ince2016} proposed to estimate mutual information based on the Gaussian Copula for continuous data, which works well for approximately Gaussian data, such as local field potentials and M/EEG data \citep{Magri2009}. \cite{Nemirovsky2023} advanced the analysis of functional MRI data by implementing integrated information theory, which is calculated based on the mutual information between the state of the conscious system over time and across the conscious system's partitions. 
\cite{Tsai1999} used mutual information to analyze functional MRI data to compute an activation map. \cite{Schloegl2002} used mutual information to study the EEG-based brain-computer interface. \cite{Chai2009} and \cite{Li2022} employed multivariate mutual information to study functional connectivity between brain regions in functional MRI data. \cite{Combrisson2022} proposed a nonparametric permutation-based framework for neurophysiological data to analyze cognitive brain networks.  

While mutual information estimation for discrete random variables is trivial, the estimation of mutual information for continuous random variables can be done using a few different approaches. 
One fundamental method is to estimate mutual information based on the binning of continuous variables to treat them as discrete variables. \cite{Steuer2002} reported improved performance using \ac{KDE} based methods. \ac{KDE}-based methods numerically calculate the mutual information estimation based on the estimated kernel density functions \citep{Moon1995}.  The \ac{kNN} approach was previously adapted to estimate mutual information \citep{Faivishevsky2008,Kraskov2004,Victor2002,Pal2010,Lord2018,Gao2014}. \cite{Khan2007} compared the performance of mutual information estimators based on \ac{kNN} and \ac{KDE} and concluded that \ac{KDE}-based mutual information estimators outperform \ac{kNN}--based estimators for small samples with a high noise level. \cite{Gao2014} argued that accurate estimation of mutual information of two strongly dependent variables using \ac{kNN}-based methods requires a prohibitively large sample size. As shown later in our simulation studies in Section \ref{subsec:ABIDE-simulation}, our \ac{KDE}-based mutual information screening method also outperforms the \ac{kNN}-based counterpart. Since kernel density estimation on large volume of data is a computationally challenging approach and that neuroimaging data is usually of large volume, variable screening based on mutual information has never been implemented for neuroimaging data to the best of our knowledge. In this paper, we implement variable screening methods using a few different approaches and carried out comprehensive simulation and real case studies using the preprocessed ABIDE data \citep{Cameron2013, Barry2020}. The variable screening functionality is encapsulated within our Python package, {\tt fastHDMI}, an acronym for \emph{Fast high--dimensional Mutual Information estimation}. This package is specifically designed to facilitate the effective processing and analysis of substantial volumes of neuroimaging data using a few different computationally efficient estimation methods. 


In Section \ref{sec:mi-estimation}, we will explore the concept of mutual information and provide an overview of the estimation methods. Subsequently, Section \ref{sec:sim-and-case-studies} assesses the efficacy of variable selection and the computational speed of the variable selection methods implemented in our {\tt fastHDMI} package. These methods encompass \emph{Fast Fourier Transform-based Kernel Density Estimation (FFTKDE)} mutual information estimation, mutual information estimation based on binning of continuous variables with the number of bins determined using the results of a previous study \citep{Birge2006} utilizing bounds on the risk of penalized maximum likelihood estimators due to Castellan \citep{Castellan2000}, \ac{kNN}-based mutual information estimation, and Pearson correlation. The \ac{kNN}-based mutual information estimation utilized in our work is adapted from the {\tt scikit-learn} library. We will begin by examining these variable screening methods within our {\tt fastHDMI} package through simulations in Section \ref{subsec:ABIDE-simulation}, then proceed to compare their computing speeds. Finally, in Section \ref{subsec:ABIDE-case-studies}, the performance of the predictive models created with the variables selected using our four implemented methods will be demonstrated. 

\section{Estimation of Mutual Information }\label{sec:mi-estimation}

The entropy-based screening methods are based on Shannon's
entropy \citep{Shannon1948}. Let $\mathbf{X}\in\mathbb{R}^{n}$ denote
a random variable residing in a probability space with probability
mass or density function $p\left(\mathbf{X}\right)$. Shannon's entropy
is defined as 
\begin{equation}
    \label{eq:shannon-entropy}
    H\left(\mathbf{X}\right)\coloneqq\mathbb{E}\left[-\log p\left(\mathbf{X}\right)\right].
\end{equation}
Furthermore, Lebesgue's decomposition theorem expands the above definition
for all other random variables. Relative entropy, also known as the {\em\ac{KL divergence}}, is a specific case of Bregman divergence applied to $-H$, the negative of Shannon’s entropy, which is a strictly convex functional: 
\begin{equation}
D_{KL}\left(\mathbf{X}_{1}\parallel\mathbf{X}_{2}\right)\coloneqq\mathbb{E}_{\mathbf{X}_{1}}\left[-\log\frac{p\left(\mathbf{X}_{2}\right)}{p\left(\mathbf{X}_{1}\right)}\right].
\end{equation}

Moreover, mutual information is defined as the \ac{KL divergence} from
the joint distribution $\left(\mathbf{X},\mathbf{Y}\right)$ to the
outer product distribution $\mathbf{X}\otimes\mathbf{Y}$, hence symmetric.
For random variables $\mathbf{X},\mathbf{Y}$, the mutual information
\begin{equation}
I\left(\mathbf{X},\mathbf{Y}\right)\coloneqq D_{KL}\left(\left(\mathbf{X},\mathbf{Y}\right)\parallel\mathbf{X}\otimes\mathbf{Y}\right).\label{eq:mi-defn}
\end{equation}
$\mathbf{X}$ and $\mathbf{Y}$ in \eqref{eq:mi-defn} are typically
univariate for variable screenings. The implementation of \ac{KDE}--based mutual information
estimation uses \ac{FFT} based \ac{KDE} methods from
the Python package {\tt KDEpy} \citep{Odland2018}. FFT-based \ac{KDE} was initially
proposed by \cite{Silverman1982} on Gaussian kernels with much faster
computing speed and much lower numerical errors.
As shown in the paper, such an approach significantly solves the computational
speed challenges that \ac{KDE} usually faces \citep{Silverman1982}. The performance of \ac{KDE} usually
depends on the bandwidth and kernel selection. While we leave it for
users to choose kernel and bandwidth, the default arguments are
set to be the state-of-the-art \emph{Improved Sheather-Jones }bandwidth
\citep{Botev2010} with Epanechnikov kernel \citep{Epanechnikov1969}. For a detailed explanation of the FFTKDE method for mutual information estimation, see Appendix \ref{apdx:method-contribution}.

At the same time, mutual information estimation using the \ac{kNN} method leverages the \ac{kNN} algorithm for entropy estimation, a technique introduced by \cite{L.F.Kozachenko1987}. This method estimates Shannon entropy, as detailed in equation \eqref{eq:shannon-entropy}, with the sample mean, alongside a trinomial distribution to estimate $\widehat{p\left(x_j\right)}$. The binning approach for mutual information estimation converts continuous variables into discrete variables through binning, with the optimal number of bins guided by findings from a previous study by \cite{Birge2006}, which derived the optimal number of bins based on the bounds on the risk of penalized maximum likelihood estimators due to \cite{Castellan2000}. Pearson correlation is calculated through the standardized inner product of outcomes and variables. Additionally, to drastically improve the processing speed for large-scale datasets, our package incorporates multiprocessing capabilities, enabling parallel processing across all employed methods. This adaptation to parallel computing significantly enhances the utility of our package, especially for extensive neuroimaging data analyses.

Previous studies demonstrated that the three density estimation methods discussed in this paper, \ac{KDE}, \ac{kNN}, and histogram--based methods, are consistent estimators under suitable conditions. The Lebesgue integral, as a linear operator, has its boundedness equivalent to continuity in a normed linear space. Since expectation is a linear operator, it is continuous under appropriate norms when it is bounded. By the continuous mapping Theorem, the mutual information estimated using these three density estimators is consistent, as the mutual information functional is continuous with respect to the joint likelihood, and continuity is preserved under finite composition. 

Furthermore, since mutual information is continuous with respect to the joint density, sufficiently small numerical errors will not significantly perturb the mutual information estimation. The numerical error associated with the \ac{FFT} procedure arises from multiple sources beyond numerical precision, including errors from using a finite number of \ac{DFT} terms --- such as discretization, truncation of frequencies, and aliasing; and errors from applying \ac{FFT} to a non--periodic function, including boundary effects, zero--padding, and interpolation. Notably, Fourier’s theorem implies that the error from \ac{FFT} for {\em periodic} functions vanishes asymptotically with respect to the number of \ac{DFT} terms. With a computational complexity of $O\left( n\log n \right)$, utilizing a sufficiently fine grid can mitigate these errors while maintaining high computational efficiency. Moreover, \ac{KDE} is inherently non-periodic. Consequently, errors due to boundary effects, zero–padding, and interpolation are influenced by the chosen interval for \ac{KDE} and will not asymptotically vanish with respect to the number of \ac{DFT} terms. The error due to the chosen bounded interval in which the data points reside presents a general challenge when evaluating mutual information numerically, not limited to the \ac{FFT} approach. Additionally, it is important to note that numerical errors, though generally insignificant when using a large number of \ac{DFT} terms, will not vanish asymptotically with respect to the number of data points in the dataset. In summary, \ac{FFT} is an efficient tool to perform \ac{KDE} while maintaining high computational efficiency, as evidenced by previous studies \citep{Silverman1982}.

\section{Simulation and Case Studies }\label{sec:sim-and-case-studies}
\emph{Autism Brain Imaging Data Exchange
(ABIDE) preprocessed Data} consists of preprocessed functional MRI
brain imaging data from \emph{$539$} individuals suffering from ASD
and $573$ typical controls \citep{Cameron2013}. In this paper, we used the preprocessed ABIDE data consisting of $149955$ brain imaging variables, together with age, biological sex, and diagnosis of autism for $508$ cases and $542$ controls \citep{Cameron2013, Barry2020}. The preprocessing was carried out exactly the same manner as the preprocessing performed earlier by \cite{Barry2020} (see also \citep{Fischl2012, Dale1999}): the T1-weighted Magnetic Resonance scans were processed through the FreeSurfer 6.0 pipeline \citep{Fischl2012} on the CBrain computing facility \citep{Sherif2014}. This pipeline delineates the cortical surface from magnetic resonance scans, allowing the quantification of the cortical thickness across the brain hemispheres \citep{Fischl2012, Dale1999}. The process involves several steps: affine registration to MNI305 space \citep{Collins1994}, bias field correction, removal of non-cortical regions, and the estimation of white matter and pial surfaces from intensity gradients, which are used to estimate cortical thickness. These cortical surfaces are projected into a common space (fsaverage) for comparison across individuals.

Brain MRI data has been used to predict age to study the brain aging process linked to diseases such as Alzheimer’s disease and Parkinson’s disease \citep{Jonsson2019,Jiang2020,Cole2017,Franke2010,Liem2017}. For the case studies based on the preprocessed ABIDE data \citep{Cameron2013, Barry2020} in Section \ref{subsec:ABIDE-case-studies}, we choose
age at the MRI scan as the continuous outcome and autism diagnosis
as the binary outcome. When using age at the MRI scan as the outcome,
we adjust for sex and autism diagnosis; we using autism diagnosis
as the outcome, we adjust for age and sex. 
We compare the few screening methods in our Python package {\tt fastHDMI}, including mutual information estimation using the FFTKDE and \ac{kNN} originally implemented in the {\tt scikit-learn} library, as well as Pearson correlation.


\subsection{Simulation based on the preprocessed ABIDE data \citep{Cameron2013, Barry2020} \label{subsec:ABIDE-simulation} }

We decided to simulate outcomes based on the preprocessed ABIDE MRI features in order to preserve the distribution patterns and the correlation structure in this high--dimensional dataset. Therefore, we simulated both nonlinear and linear outcomes from the preprocessed
ABIDE data \citep{Cameron2013, Barry2020}. Let $\mathbf{X}\in\mathbb{R}^{N\times p}$
denote the design matrix; i.e., all the MRI brain imaging variables from
the entire preprocessed ABIDE dataset. The simulation of the \emph{nonlinear}
outcomes proceeds in this manner -- the nonlinearity for continuous
outcomes comes from the quadratic manipulation, i.e., step 4: 

\vspace{.1in}

\noindent\fbox{%
    \parbox{\textwidth}{%
\vspace{.1in}
\begin{enumerate}
\item Pick the number of ``true'' covariates $p_{\text{true}}$, choose $p_{\text{true}}$ uniformly randomly from the full feature set; let
$\mathbf{X}_{\text{true}}\in\mathbb{R}^{N\times p_{\text{true}}}$ 
denote the corresponding design sub-matrix. 
\item Simulate the corresponding ``true'' coefficients $\boldsymbol{\beta}_{\text{true}}\in\mathbb{R}^{p_{\text{true}}}$ with $\boldsymbol{\beta}_{\text{true}}\sim N_{p_{\text{true}}}\left(1,\Sigma_{\boldsymbol{\beta}_{\text{true}}}\right)$ and $\Sigma_{\boldsymbol{\beta}_{\text{true}}}$ being a $0.6$ Toeplitz matrix. The correlation design aims to replicate the phenomenon of correlated brain signals. 
\item Standardize the design sub-matrix for the true features $\mathbf{X}_{\text{true}}$, to obtain
$\mathbf{X}_{\text{true},1}$. 
\item For nonlinear simulations only: take the element-wise square of $\mathbf{X}_{\text{true},1}$
and then standardize the matrix again to obtain $\mathbf{X}_{\text{true},2}$; the standardization here is to ensure that each feature impacts the simulated outcome proportionally. 
\item The continuous and binary outcomes are then simulated in this manner: 
\begin{enumerate}
\item To simulate continuous outcomes: 
\begin{enumerate}
\item Pick $\text{SNR}=3$; calculate $\sigma_{\text{true}}=\sqrt{\frac{\boldsymbol{\beta}_{\text{true}}^{T}\mathbf{X}_{\text{true},2}^{T}\mathbf{X}_{\text{true},2}\boldsymbol{\beta}_{\text{true}}}{\text{SNR}}}$; 
\item Simulate the error $\varepsilon_{j}\overset{i.i.d.}{\sim}N\left(0,\sigma_{\text{true}}^{2}\right)$;
\item The outcome is simulated as $\mathbf{y}=\mathbf{X}_{\text{true},2}\boldsymbol{\beta}_{\text{true}}+\boldsymbol{\varepsilon}$. 
\end{enumerate}
\item To simulate binary outcomes: 
\begin{enumerate}
\item Calculate $\boldsymbol{\tau}=\mathbf{X}_{\text{true},2}\boldsymbol{\beta}_{\text{true}}$; 
\item Standardize $\boldsymbol{\tau}$, obtain $\boldsymbol{\tau}^{\prime}$
-- this is to avoid the data being too centered, which will cause
all simulated binary outcomes in the same class; 
\item Take $\boldsymbol{\tau}^{\prime\prime}=\boldsymbol{\tau}^{\prime}+\text{arctanh}\sqrt{\frac{1}{3}}$
for \emph{translated} binary outcome simulations, or $\boldsymbol{\tau}^{\prime\prime}=\boldsymbol{\tau}^{\prime}$
for \emph{original} binary outcome simulations. The translated binary outcome
simulation is to make the logistic transformation of centered data in the next step
as nonlinear as possible, as $\text{\ensuremath{\pm}arctanh\ensuremath{\sqrt{\frac{1}{3}}}}$
is the location for the logistic transformation to achieve the greatest
absolute curvature value; 
\item The binary outcome is then simulated as $y_{j}\overset{\text{indep.}}{\sim}\text{Bern}\left(\text{logistic}\left(\tau_{j}^{\prime\prime}\right)\right)$. 
\end{enumerate}
\end{enumerate}
\end{enumerate}
For linear simulations, we omit step 4 and take $\mathbf{X}_{\text{true},2}\coloneqq\mathbf{X}_{\text{true},1}$
thereafter. 
\vspace{.05in}
    }%
}

\vspace{.05in}

The screening of features with respect to the simulated continuous and binary outcomes $\mathbf{y}$ are then carried out using the original entire design matrix $\mathbf{X}$. 
Variable selection performance is measured by \emph{Variable Selection
Area under Receiver Operating Curve (AUROC)}, which is the AUROC calculated with the true labels taking value $1$
for the simulated true coefficients and $0$ for other coefficients,
and the ranking of the coefficients follows the absolute value of
the three association measures, respectively; i.e., $\widehat{MI}$
based on FFTKDE and \ac{kNN}, as well as Pearson correlation. The top $p_{\text{true}}$ of the most associated covariates are then taken as selected covariates, which will take value $1$, and the others will take value $0$. Variable Selection
AUROC therefore measures the matching between the selected covariates and the simulated ``true'' covariates. Such
measures can differentiate distinct methods when the traditional measures
such as classification rate or adjusted Rand Index can not -- a scenario frequently occurs to variable selection for ultra-high--dimensional
data. 

\begin{figure}[h]
\centering{}\includegraphics[width=1.0\textwidth]{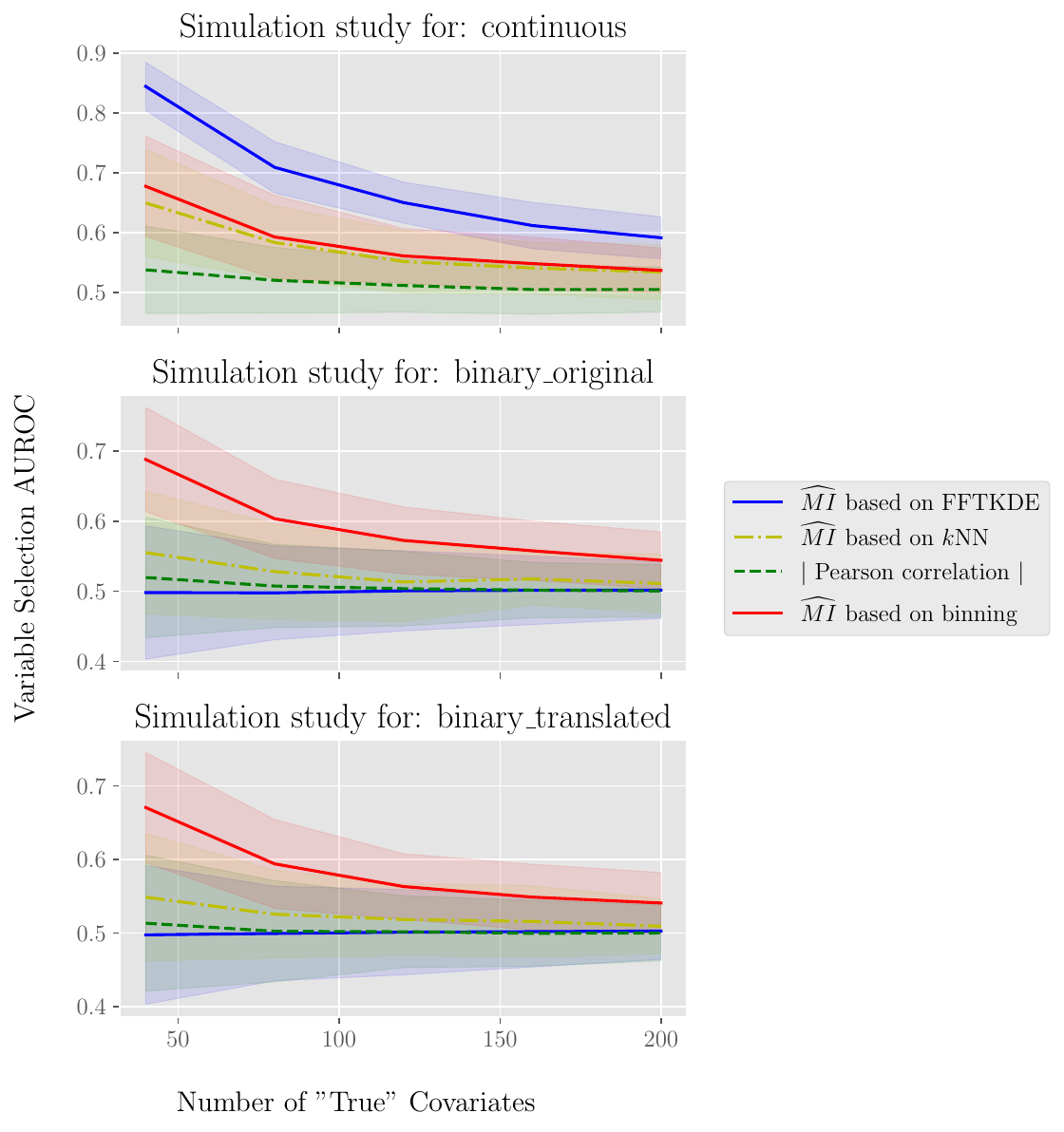}\caption{Variable selection AUROC on the simulated \emph{nonlinear} continuous
and original/translated binary outcomes; the horizontal axis is the
number of \textquotedblleft true\textquotedblright{} covariates used
in the outcome simulation. Means with their $95\%$ confidence intervals were plotted
for $100$ simulation replications. \label{fig:Variable-selection-AUROC-nonlinear}}
\end{figure}

\begin{figure}[h]
\centering{}\includegraphics[width=1.0\textwidth]{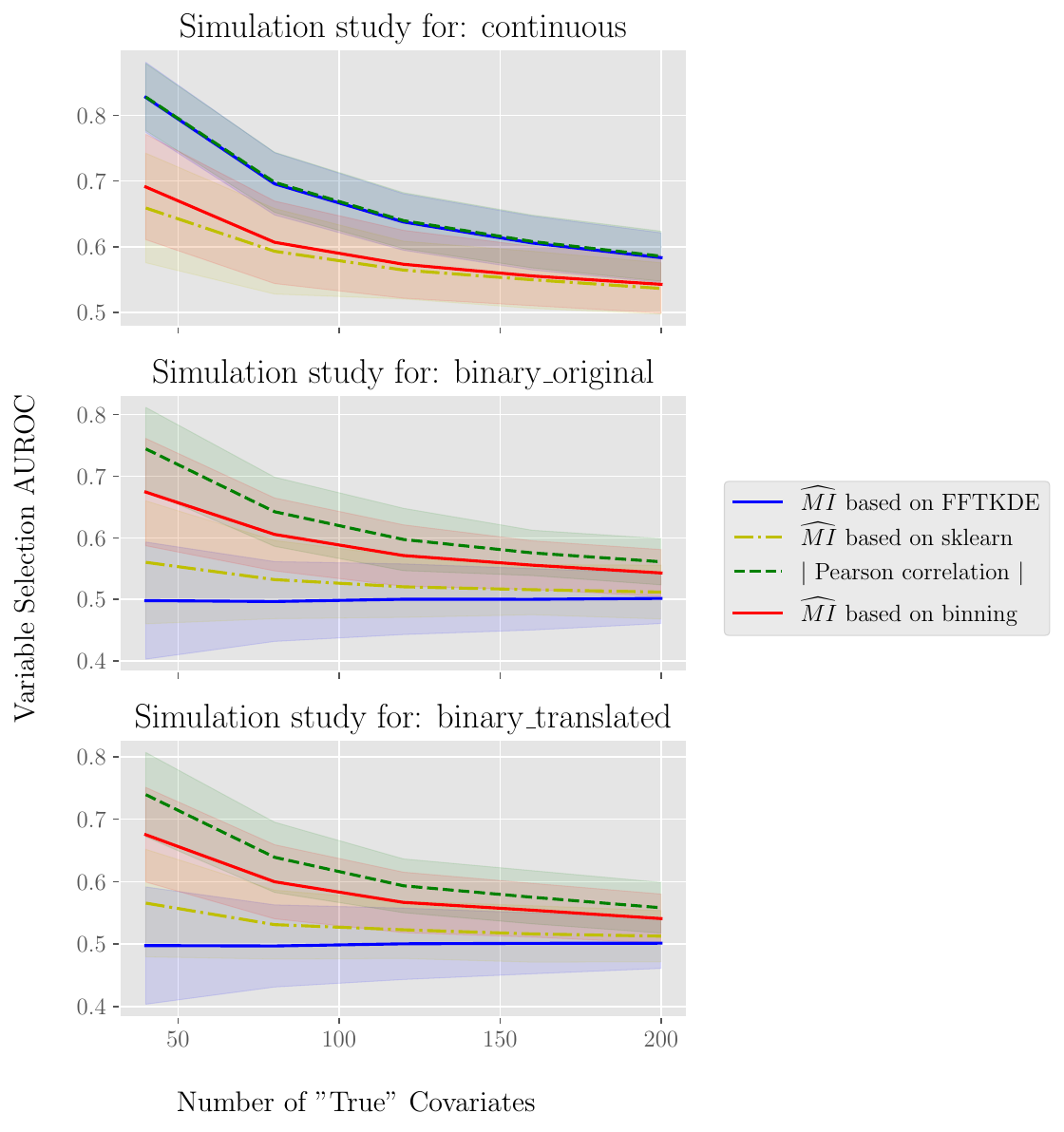}\caption{Variable selection AUROC on the simulated \emph{linear} continuous
and original/translated binary outcomes; the horizontal axis is the
number of \textquotedblleft true\textquotedblright{} covariates used
in the outcome simulation. Means with their $95\%$ confidence intervals were plotted
for $100$ simulation replications. \label{fig:Variable-selection-AUROC-linear}}
\end{figure}

We evaluate the efficacy of our implemented variable screening methods in {\tt fastHDMI} package, including: 1) Mutual information estimation using FFTKDE, 2) Mutual information estimation using \ac{kNN}, 3) Mutual information estimation through binning, and 4) absolute Pearson correlation. Our findings, illustrated in Figures \ref{fig:Variable-selection-AUROC-nonlinear} and \ref{fig:Variable-selection-AUROC-linear}, reveal that for continuous outcomes, the FFTKDE-based mutual information estimator outperforms its counterparts. In scenarios with linear relationships, FFTKDE-based mutual information estimator and absolute Pearson correlation are jointly the most effective. Conversely, for binary outcomes, the binning-based mutual information estimator excels in capturing nonlinear associations, whereas other methodologies display substantially overlapping confidence intervals. In linear association contexts, Pearson correlation emerges as the most effective method for binary outcomes. Interestingly, Pearson correlation, particularly when employed with a balanced number of cases and controls, inherently correlates to a two-sample testing approach, which explains its superior performance for binary outcomes with linearly simulated underlying probability pre-image. 

\begin{figure}[h]
\centering{}\includegraphics[width=1.0\textwidth]{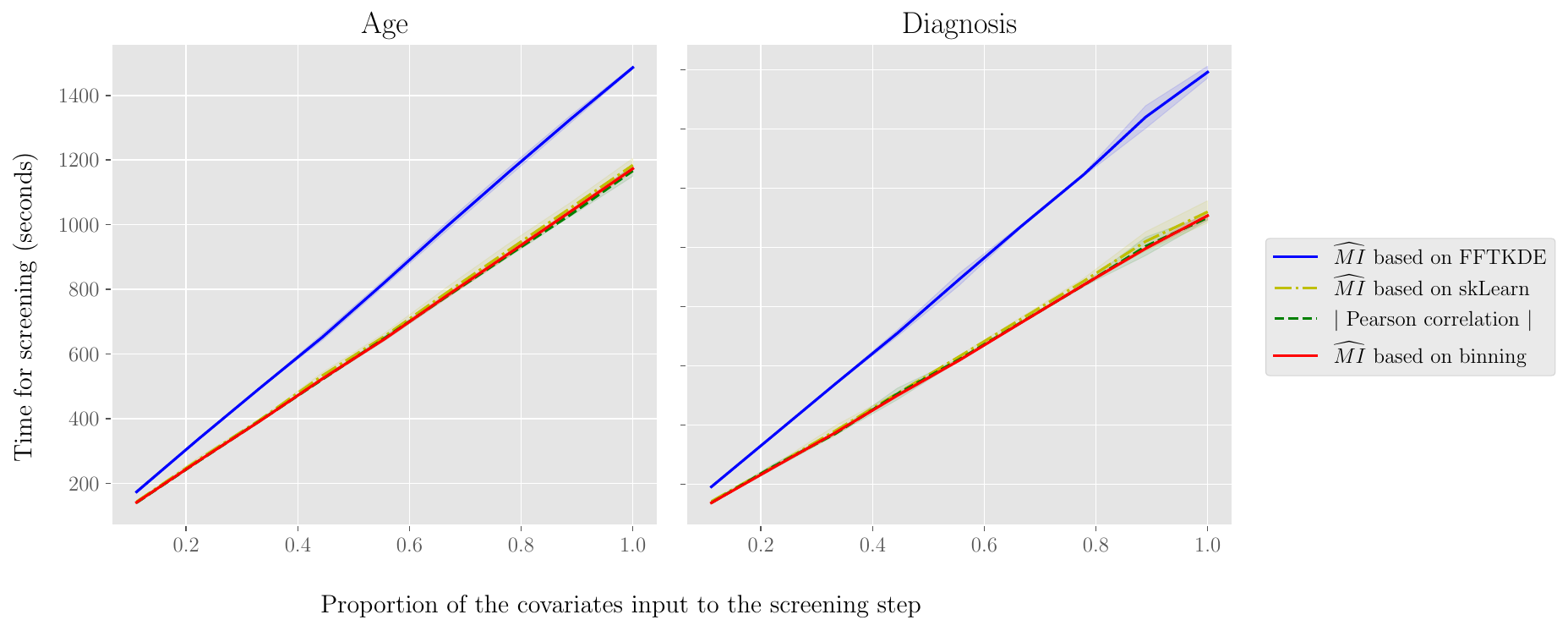}\caption{Running speeds of variable screening for continuous (age) and binary (diagnosis) outcomes utilizing the methods under study. The horizontal axis represents the proportion of features introduced into the screening phase, while the vertical axis measures the time in seconds to complete the screening. The plot displays the mean running times and their corresponding $95\%$ confidence intervals (C.I,), derived from 5 simulation replications. \label{fig:running-speed}}
\end{figure}

All discussed variable screening methods were conducted concurrently on $16$-core CPUs on Compute Canada. The fast Fourier transform (FFT) algorithm is leveraged to significantly enhance the efficiency of the \ac{KDE} estimation process, traditionally viewed as computationally intensive. As depicted in Figure \ref{fig:running-speed}, the execution times to complete the screenings with all the methods implemented in our {\tt fastHDMI} package are assessed. Notably, the \ac{KDE}-based mutual information estimation, often anticipated to be slower, exhibited competitive speed akin to alternative methods, courtesy of the FFT algorithm's effectiveness. This computational efficiency was achieved with the same CPU configuration, while intentionally avoiding multiple data duplications in memory during multiprocessing. Given the substantial size of high--dimensional datasets, duplicating such datasets in memory is generally impractical.

\subsection{Pre-processed ABIDE data case studies \citep{Cameron2013, Barry2020} -- predict age and diagnosis \label{subsec:ABIDE-case-studies} }


In this subsection, we evaluate the performance of various variable screening techniques implemented in the {\tt fastHDMI} package using preprocessed ABIDE data \citep{Cameron2013, Barry2020}. Initially, we deploy the four variable screening methods to identify the features most associated with the outcome. Since we are fitting multiple penalized models, standardization of the selected variables is carried out to achieve a sample mean of $0$ and a standard deviation of $1$. This step is crucial for ensuring consistent penalization across all coefficients of the penalized covariates.

Subsequently, we divide the dataset, stratified by the outcome, into a training set comprising $80\%$ of the observations and a testing set with the remaining $20\%$. This stratification ensures a balanced representation of the outcomes in both sets. For the continuous outcome, age, we employ binning to categorize observations into $30$ bins based on their outcome values, followed by stratification based on the bin labels. This approach allows for the division of the dataset into training and testing sets with similar outcome means, an important factor for reliable prediction performance comparison.

For the continuous outcome variable, age at MRI scan, we fit several models: elastic net, least-angle regression (LARS), least absolute shrinkage and selection operator (LASSO), LASSO-LARS, linear model, Random Forest regressor, and ridge regression. Except for the Random Forest regressor, which utilizes the out-of-bag error scored by $R^2$ for model averaging, all models are tuned using $5$-fold cross-validation with validation set $R^2$ as the scoring function for penalty hyperparameters.

For binary outcomes, diagnosis of autism disorder, we fit both unpenalized and penalized logistic regressions (using $\ell_1$, $\ell_2$, and elastic net penalties), as well as the Random Forest classifier. All models, with the exception of the Random Forest classifier, which uses out-of-bag error scored by Gini impurity for model averaging, are tuned using $5$-fold cross-validation, scored by mean accuracy for the penalty hyperparameters.

Unlike simulation studies in Section \ref{subsec:ABIDE-simulation}, where ``true'' signals are known, case studies lack such definitive benchmarks, necessitating reliance on model-based performance metrics. Hence, we use testing set $R^2$ for continuous outcomes and testing set Area Under the Receiver Operating Characteristic (AUROC) for binary outcomes to evaluate model performance. 

\begin{figure}[h]
\centering{}\includegraphics[width=0.5\textwidth]{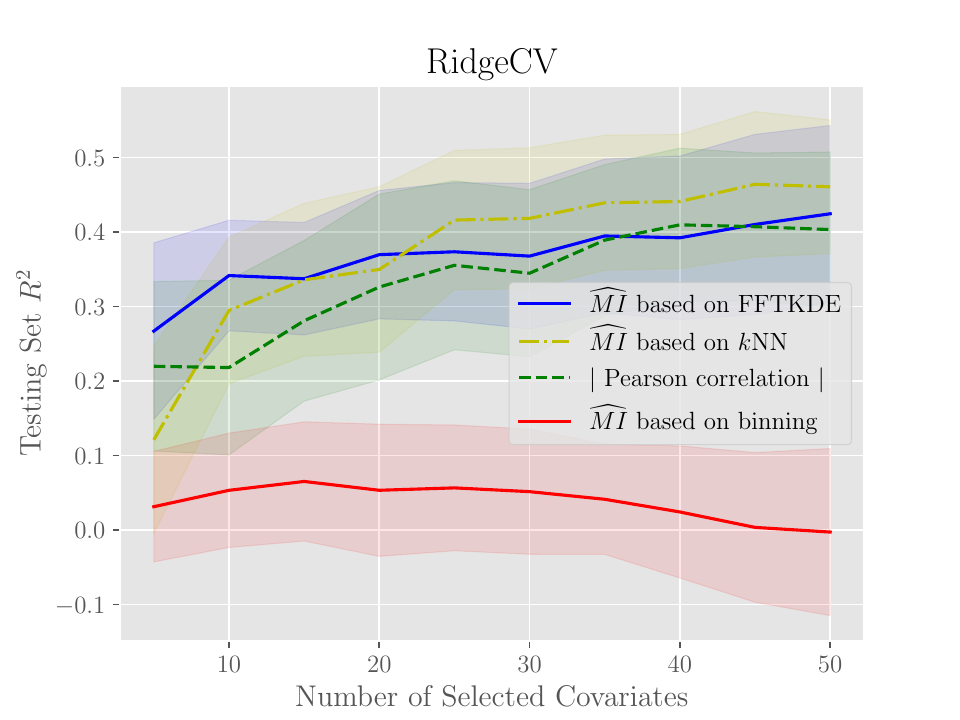}\includegraphics[width=0.5\textwidth]{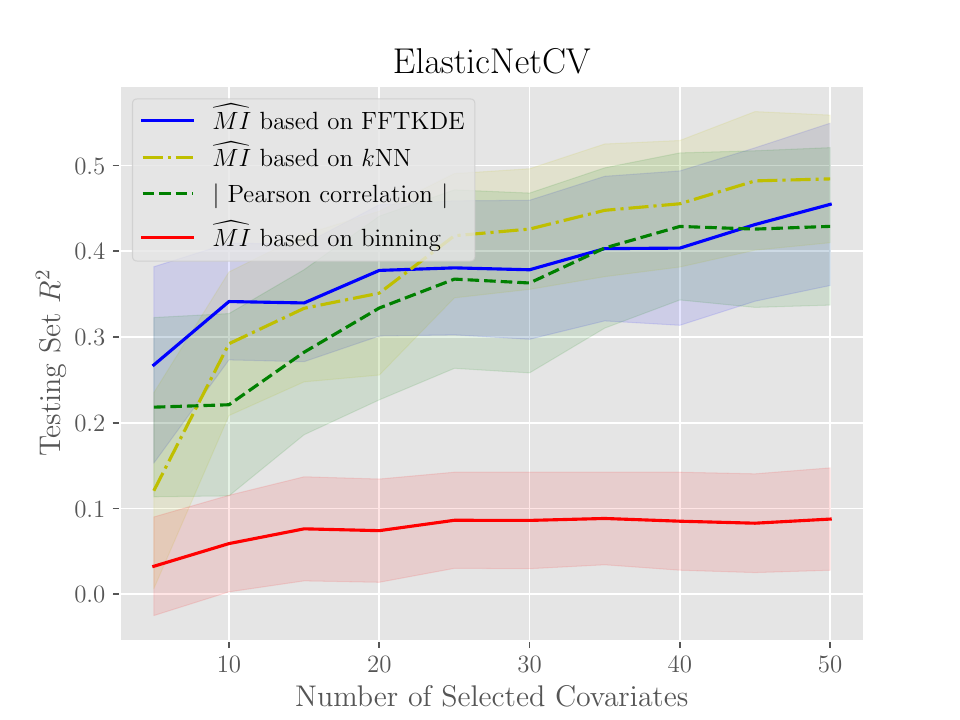} \\ \includegraphics[width=0.5\textwidth]{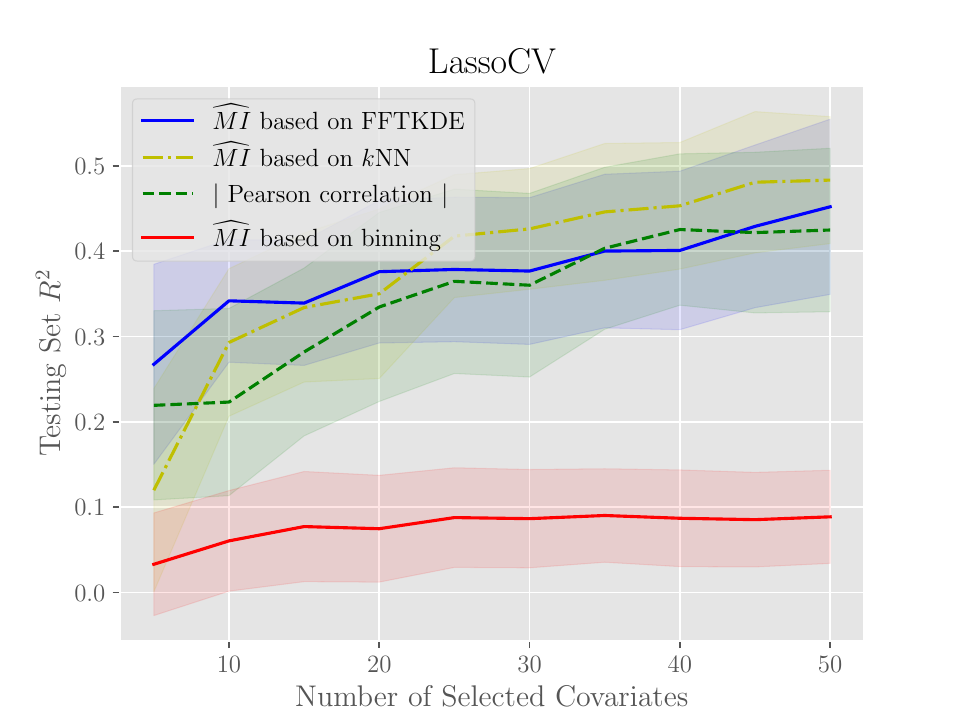}\includegraphics[width=0.5\textwidth]{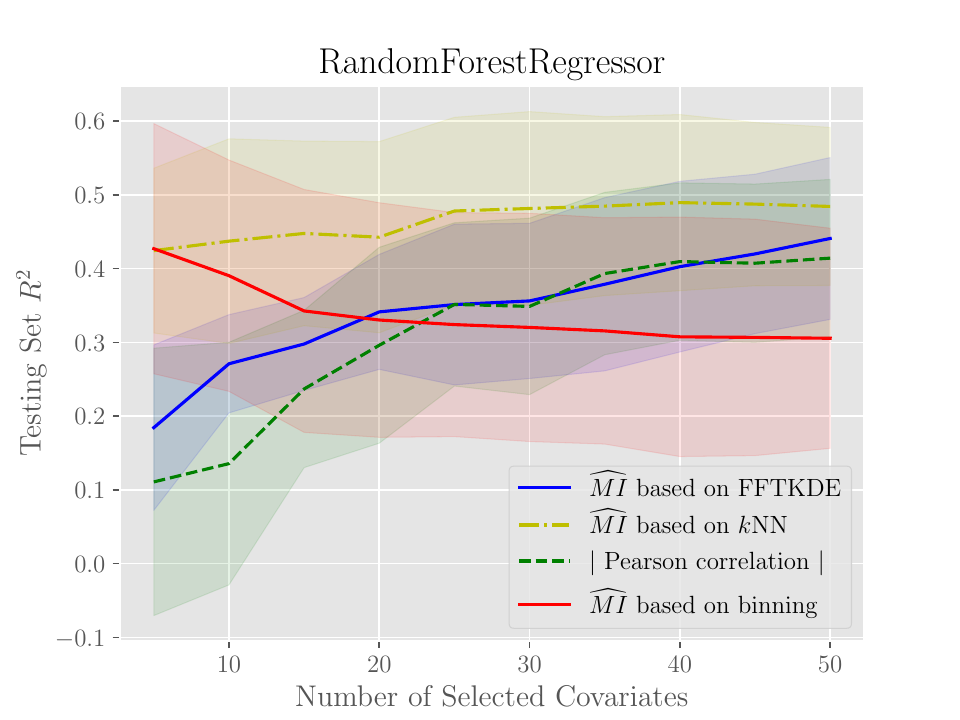}\caption{Testing Set $R^{2}$ for age at the scan outcome v.s. the number of
most associated brain imaging covariates based on the association
measure rankings. Means with their $95\%$ confidence intervals were plotted for $20$
simulation replications. \label{fig:Case-Study-continuous}}
\end{figure}

\begin{figure}[h]
\centering{}\includegraphics[width=0.5\textwidth]{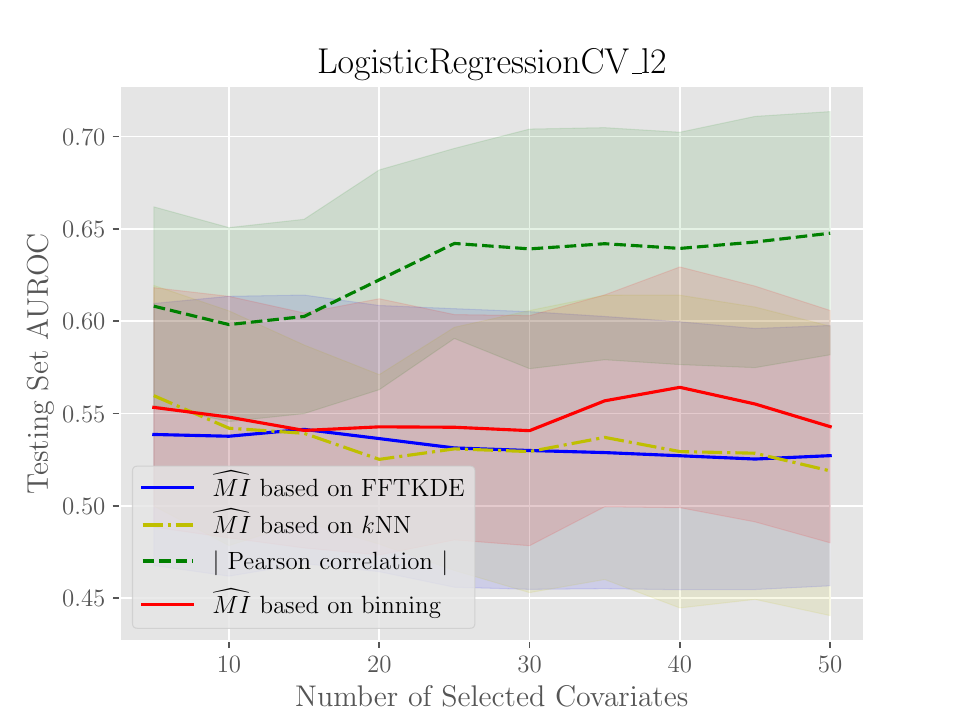}\includegraphics[width=0.5\textwidth]{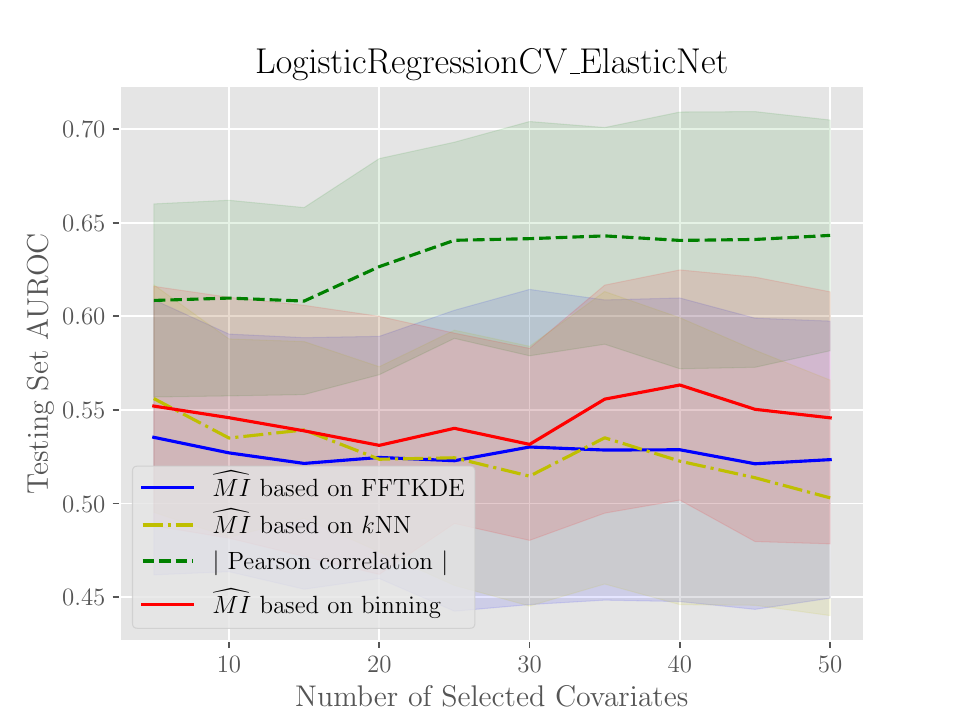}\\ \includegraphics[width=0.5\textwidth]{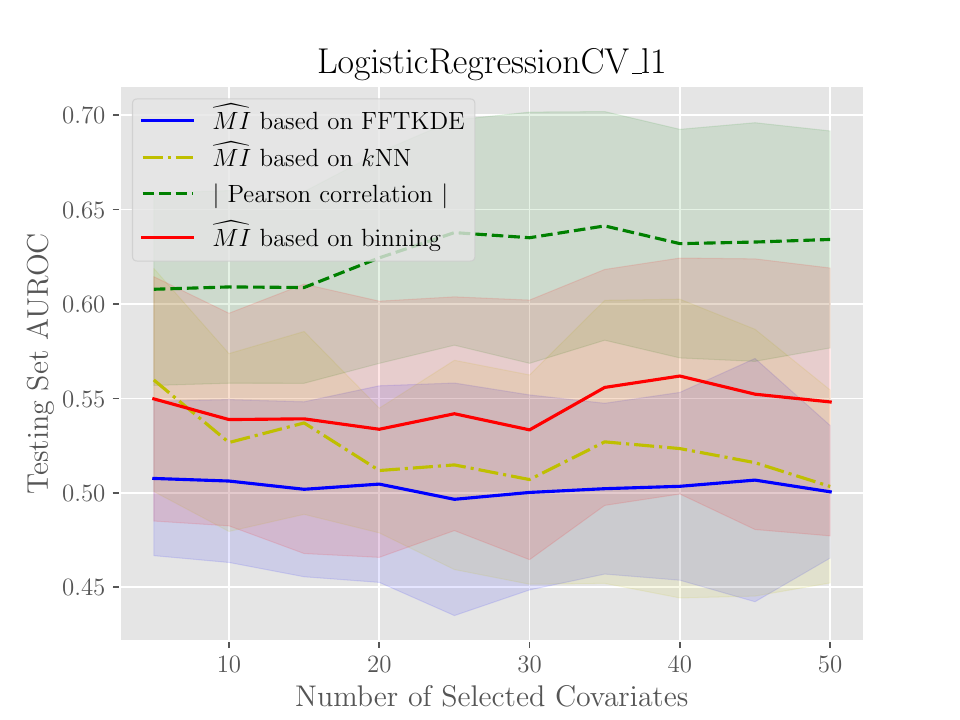}\includegraphics[width=0.5\textwidth]{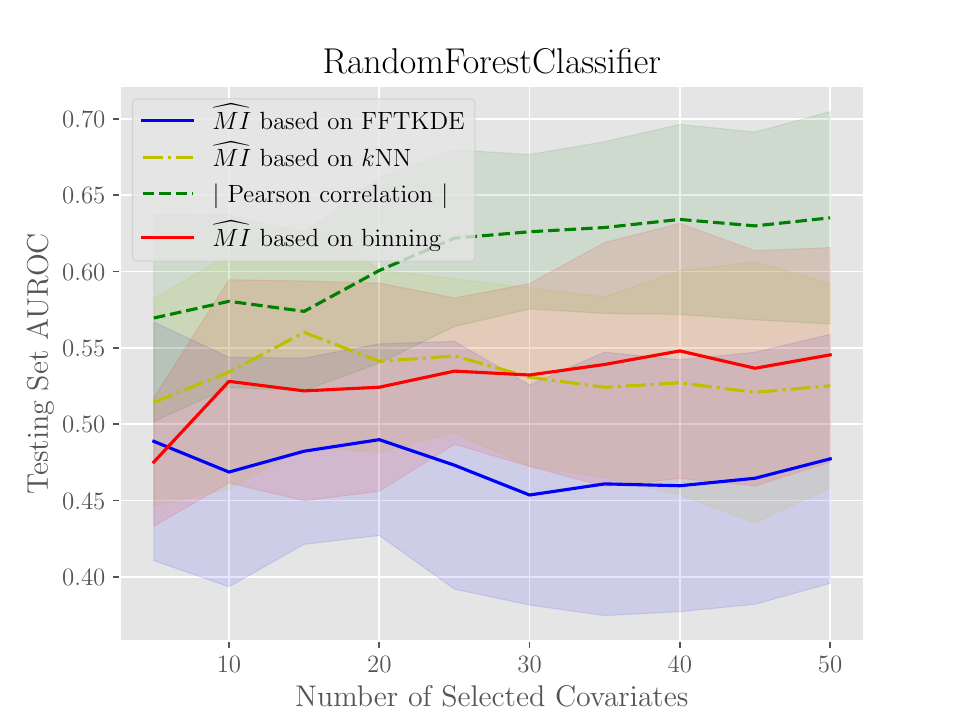}\caption{Testing Set AUROC for autism diagnosis outcome v.s. the number of
most associated brain imaging covariates based on the association
measure rankings. Means with their $95\%$ confidence intervals were plotted for $20$
simulation replications. \label{fig:Case-Study-binary}}
\end{figure}

Figure \ref{fig:Case-Study-continuous} illustrates that in predicting the continuous outcome, age at MRI scan, linear models utilizing brain imaging variables selected using mutual information estimations via FFTKDE or \ac{kNN} emerge as the best-performing. Conversely, models built using variables selected by mutual information estimations based on binning exhibit the least predictive capability. However, within the context of random forest regression, models built using variables chosen through mutual information estimation by \ac{kNN} outperform the rest. Figure \ref{fig:Case-Study-binary} indicates that for the binary outcome of autism diagnosis, models constructed with variables selected via absolute Pearson correlation yield superior predictive performance. This phenomenon could stem from multiple factors, including the linear nature of the assessment model, which favors linear association measures, or a linear relationship between age at MRI scan, the probability of autism diagnosis, and the brain imaging covariates.

\section{Conclusion and Discussion }

In this paper, we introduce the Python package {\tt fastHDMI}, designed to streamline variable screening through three distinct mutual information estimation methods along with absolute Pearson correlation. Our evaluations, conducted on the large, high--dimensional preprocessed ABIDE data \citep{Cameron2013, Barry2020}, affirm {\tt fastHDMI}'s computational efficiency and robustness. Through extensive simulation studies, which encompass both simulations for linear and nonlinear associations, as well as continuous and binary simulated outcomes, we evaluated the performance of each implemented variable screening method. Our findings reveal that for simulated continuous nonlinear outcomes, the FFTKDE-based mutual information estimation method excels in variable selection. Similarly, for simulated binary outcomes with a nonlinear underlying probability preimage, the binning-based mutual information estimation stands out. In the cases of simulated continuous linear outcomes, both absolute Pearson correlation and FFTKDE-based mutual information estimation share the top performance. Furthermore, absolute Pearson correlation is superior for binary outcomes simulated with linear underlying probability preimage. Complementing our simulations, a comprehensive case study on the preprocessed ABIDE data \citep{Cameron2013, Barry2020} showcased the predictive capabilities of models crafted from the most relevant covariates identified by our methods. By pioneering sophisticated variable selection techniques in the domain of high--dimensional neuroimaging data, our work stands as a critical advancement, fostering novel pathways for research exploration and analytical insight within the scientific community. A promising avenue for future research could be to explore variable screening based on non-parametric copula models \citep{Rabhi2019}.

\section{Disclaimer}
All codes to reproduce the simulation and case study results of this paper and outputs from Calcul Quebec/Compute Canada can be found on the following GitHub repository:

\href{https://github.com/Kaiyangshi-Ito/fastHDMI}{https://github.com/Kaiyangshi-Ito/fastHDMI}

\chapter{Accelerated Gradient Methods for Sparse Statistical Learning with Nonconvex Penalties}\label{ch:paper2}

\indent \textbf{Preamble to Manuscript 2.} 
\subsection*{Introduction to the Study and Its Place in the Workflow:}
Manuscript 2 advances the computational frontier by delving into the optimization challenges associated with nonconvex oracle penalties, a critical area when dealing with the complexities of sparse learning of high-dimensional data that follow the initial variable screening. The manuscript adapts Nesterov’s \ac{AG} method, traditionally used for convex objective functions, to handle nonconvexity induced by oracle penalties such as SCAD.

\subsection*{Building on Foundations Established in Manuscript 1:}
The integration of nonconvex optimization techniques is a direct progression from the efficient variable screening method introduced in Manuscript 1. With the relevant variables identified using {\tt fastHDMI}, the need for efficient optimization techniques that can manage the intricacies of the high--dimensional dataset composed of these selected variables becomes apparent. Manuscript 2 addresses this by enhancing the capability of statistical computing methods to converge faster, even when nonconvex penalties are involved, thereby ensuring the computational efficiency of sparse learning.

\subsection*{Innovation and Contribution to Statistical Computing:}
The manuscript's development of a optimization hyperparameter setting based on the complexity upper bound to accelerate convergence represents a significant contribution in statistical computing. By establishing a rate of convergence and providing a new bound for the optimal damping sequence, this study not only enhances the understanding of nonconvex optimization algorithms but also improves the practical application of these methods in high--dimensional settings.

\subsection*{Enhancing Model Performance and Reliability:}
The proposed adaptations allow faster convergence compared to traditional methods, such as the proximal gradient algorithm. This improvement is crucial for handling sophisticated models that emerge from the high--dimensional large datasets prevalent in biostatistics, particularly those involving sparse learning problems. The ability to recover signals more effectively further underscores the practical value of the advances made in this manuscript.

\subsection*{On the Lipschitz--Smooth Constant and More Clarifications:}
Depending on the optimization problem, $L_{\Psi}$ often has a closed form in a context of statistical sparse learning. For example, in the case of a penalized linear model, it is given by $\frac{1}{n}\left\Vert X^{T}X\right\Vert _{2}+L_{\text{SCAD/MCP}}$. In the discussed statistical sparse learning context, nonconvexity arises from the nonconvex penalties. As illustrated in the manuscript, nonconvex penalties typically decompose into a difference of convex form: a convex $\ell_{1}$ component to induce sparsity and a concave component. As discussed in the manuscript, the concave component has a Lipschitz--smooth constant of $L_{\text{SCAD}}=\frac{1}{a-1}$ for SCAD and $L_{\text{MCP}}=\frac{1}{\gamma}$ for MCP, which is often negligible compared to the Lipschitz-smooth constant for the convex smooth component. Previous literature indicates that the greatest eigenvalue of random matrices tend to grow with the number of dimensions. Specifically, studies on the spectral properties of random matrices have shown that the greatest eigenvalue often scales with the dimension of the matrix \citep{Wigner1955, Wigner1958, Marcenko1967, Mehta2004, Bai2010}. This implies that the Lipschitz--smooth constant for the nonconvex smooth component is often negligible compared to that of the convex smooth component in the context of high-dimensional data, where the operator norm of the Hessian typically grows with the greatest eigenvalue of the design matrix. 

Global convergence refers to the property of an algorithm in which, starting from any point in the feasible set, the algorithm will converge to a stationary point.  

The dynamical system interpretation of the momentum methods reveals that the trajectory of the algorithm describes {\em Newtonian particles moving through a viscous medium in a conservative force field} \citep{Qian1999, Su2014, Shi2018a, Attouch2020}. Hence, in the context of statistical computing, when applied to a proper objective function, the trajectory of the optimization algorithm for parameter estimation is bounded, which serves as the rationale behind the boundedness assumption.

\subsection*{Transition to Manuscript 3:}
Having established an efficient framework for optimizing high-dimensional statistical models under nonconvex conditions, Manuscript 3 takes the next logical step by addressing another layer of complexity: robust statistical modeling of correlated observations. The introduction of the $q$Gaussian linear mixed-effects model in Manuscript 3 builds directly on the optimization techniques refined in Manuscript 2, adapting them to models that must also account for correlation among observations. This advancement guarantees that the formulated methods are not just computationally efficient, but also robust to underlying normality assumptions and heavy tails while accounting for correlated observations, assisting in managing the intricacies of biostatistical data.



\newpage
\acresetall

\vspace*{2cm}

\begin{center}
	\Large{Accelerated Gradient Methods for Sparse Statistical Learning with Nonconvex Penalties}
\end{center}

\vspace*{3cm}

\begin{center}
Kai~Yang$^{1}$, Masoud~Asgharian$^{2}$, Sahir~Bhatnagar$^{1}$.
\end{center}

\vspace*{1cm}

\begin{center}
	$^{1}$\textit{Department of Epidemiology, Biostatistics, and Occupational Health, McGill University}\\
	$^{2}$\textit{Department of Mathematics and Statistics, McGill University}\\
\end{center}

\vspace*{1cm}

\begin{center}
	This thesis contains the accepted version of the corresponding paper published in \textit{Statistics and Computing} (\citep{Yang2024}). \\© The Author(s), under exclusive licence to Springer Science+Business Media, LLC, part of Springer Nature 2024
\end{center}

\newpage

\section*{Abstract}
Nesterov's accelerated gradient (AG) is a popular technique to optimize 
objective functions comprising two components: a convex loss and a penalty function. While AG methods perform well for convex penalties, such as the LASSO,  convergence issues may arise when it is applied to nonconvex penalties, such as SCAD. 
A recent proposal generalizes Nesterov's AG method to the nonconvex setting. 
The proposed algorithm requires specification of  several hyperparameters for its practical 
application. 
Aside from some general conditions, there is no explicit rule for selecting the hyperparameters, and how different selection can affect convergence of the algorithm.  
In this article, we propose a hyperparameter setting based on the complexity upper bound to accelerate convergence, and consider the application of this nonconvex AG algorithm to high-dimensional linear and logistic sparse learning problems.
We further establish the rate of convergence and present a simple and useful bound  {to characterize our proposed optimal} damping  sequence. Simulation studies show that convergence can be made, on average, considerably faster than that of the conventional proximal gradient algorithm. Our experiments also show that the proposed method generally outperforms the current state-of-the-art methods in terms of signal recovery.
\newpage

 
\section{Introduction} 

 Sparse learning is an important component of modern data science and is an essential tool for the statistical analysis of high-dimensional data, with significant applications in signal processing and statistical genetics, among others. Penalization
is commonly used to achieve sparsity in parameter estimation. The prototypical optimization problem for obtaining penalized estimators is 
\[
\hat{\boldsymbol{\beta}}\in\arg\min_{\boldsymbol{\beta}\in\mathbb{R}^{q+1}}\left[f\left(\boldsymbol{\beta}\right)+\sum_{j=1}^{q}p_{\lambda}\left(\beta_{j}\right)\right],
\]
where $f:\mathbb{R}^{q+1}\mapsto\mathbb{R}$ is a convex loss function,  $p_\lambda:\mathbb{R}\mapsto\mathbb{R}_{\geq 0}$ constitutes the penalty term, and $\lambda>0$ is the tuning parameter for the penalty. 
Commonly used penalization methods for sparse learning include:
LASSO (Least Absolute Shrinkage and Selection Operator)~\citep{Tibshirani1996}, Elastic Net~\citep{Zou2005}, SCAD (Smoothly Clipped Absolute Deviation)~\citep{Fan2001} and MCP (Minimax Concave Penalty)~\citep{Zhang2010}.
Among these penalties, parameter estimation with
SCAD and MCP leads to a nonconvex objective function. The nonconvexity poses a challenge in statistical computing, as most methods developed for convex objective functions might not converge when applied to the nonconvex counterpart. 


Various approaches have been proposed to carry out
parameter estimation with SCAD or MCP penalties. \citeauthor{Zou2008}~\citep{Zou2008} proposed a local linear approximation, which yields a first-order majorization-minimization (MM) algorithm.  \citeauthor{Kim2008}~\citep{Kim2008} discussed a difference-of-convex
programming (DCP) method for ordinary least square estimators penalized by the SCAD penalty,
which was later generalized by \citeauthor{Wang2013}~\citep{Wang2013} to a general class of nonconvex
penalties to produce a first-order algorithm. These first-order methods belong to the class of proximal gradient descent methods, which are usually inefficient as relaxation is often expensive~\citep{Nesterov2004}. The objective function is often ill-conditioned for sparse learning problems, and gradient descent with constant step size is especially inefficient for high-dimensional problems. Indeed, previous studies have suggested that the condition number of a square random matrix grows linearly with respect to its dimension~\citep{Edelman1988}. Therefore, high-dimensional problems have a large condition number with high probability. Specific to gradient descent with constant step size, the trajectory will oscillate in the directions with a large eigenvalue, moving very slowly toward the directions with a small eigenvalue, making the algorithm inefficient.
\citeauthor{Lee2016}~\citep{Lee2016} developed a modified second-order method
originally designed for the ordinary least square loss function penalized by LASSO with extensions
to SCAD and MCP; this attempt was later extended to generalized linear
models, such as logistic and Poisson regression, and Cox's proportional hazard model. Quasi-Newton methods, or a mixture of first and second-order descent methods, have also been applied on nonconvex penalties~\citep{Ibrahim2012,Ghosh2016}. However, for high-dimensional problems, these second-order methods are slow due to the computational cost of evaluating the secant condition. Concurrently, most first and second-order methods discussed above require a line-search procedure at each step to ensure global convergence, which is prohibitive when the number of parameters to estimate grows large. 
\citeauthor{Breheny2011}~\citep{Breheny2011} implemented a coordinate descent method in the {\tt ncvreg} {\tt R} package to carry out estimation for linear models with least squares loss or logistic regression, penalized by SCAD and MCP. \citeauthor{Mazumder2011}~\citep{Mazumder2011} also implemented a coordinate descent method in the {\tt sparsenet} {\tt R} package, which carries out a closed-form root-finding update in a coordinate-wise manner for penalized linear regression. Similar to how ill-conditioning makes gradient descent inefficient, coordinate descent methods are generally inefficient when the covariate correlations are high~\citep{Friedman2007}. Previous studies have also found that coordinate-wise minimization might not converge for some nonsmooth objective functions~\citep{Spall2012}. 
Furthermore, it is naturally challenging to run coordinate-wise minimization in parallel, as the algorithm must run in a sequential coordinate manner. 


Due to the low computational
cost and adequate memory requirement per iteration, first-order methods
without a line search procedure have become the primary approach for
high-dimensional problems arising from various areas~\citep{Beck2017}. 
For smooth convex objective functions, Nesterov proposed the {\em accelerated
gradient method} (AG) to improve the rate of convergence from $O(1/N)$ for gradient descent
to $O(1/N^{2})$ while achieving global convergence~\citep{Nesterov1983}.
Subsequently, Nesterov extended AG to composite convex problems~\citep{Nesterov2012},
whereas the objective is the sum of a smooth convex function and a simple nonsmooth convex function. With proper step-size choices, Nesterov's
AG was later shown optimal to solve both smooth and nonsmooth
convex programming problems~\citep{Lan2011}. 


Given that sparse learning problems are often high-dimensional,
Nesterov's AG has been frequently used for \emph{convex} problems
in statistical machine learning (e.g.,~\citep{Simon2013, Yang2014, Yu2015, Akyildiz2021}).
However, convergence is questionable if the convexity assumption is violated. 
Recently, \citeauthor{Ghadimi2015}~\citep{Ghadimi2015} generalized the AG method to nonconvex objective functions, hereafter referred to as the nonconvex AG method, and derived the rates of convergence for both smooth and composite objective functions. While this method can be applied to nonconvex sparse learning problems, several hyperparameters must be set prior to running the algorithm and can be difficult to choose in practice. Indeed, the nonconvex AG method has never been applied in the context of sparse statistical learning problems with nonconvex penalties, such as SCAD and MCP. 


This manuscript presents a detailed analysis of the complexity upper bound of the nonconvex AG algorithm and proposes a hyperparameter setting to accelerate convergence (Theorem \ref{thm:convex}). We further establish the rate of convergence (Theorem \ref{thm:1overk}) and present a simple and useful bound to characterize our proposed optimal damping sequence (Theorem \ref{thm:alpha-k-vanishing} and Corollary \ref{cor:alpha-vanishing}). Our simulation studies on penalized linear and logistic models show that the nonconvex AG method with the proposed hyperparameter selector converges considerably faster than other first-order methods. We also compare the signal recovery performance of the algorithm to that of {\tt ncvreg}, the state-of-the-art method based on coordinate descent, showing that the proposed method outperforms the state-of-the-art coordinate descent method. 


The rest of this manuscript is organised as follows. In Sections \ref{sec:Accelerated-Gradient-Method}, \ref{sec:aga}, \ref{sec:theory},
we will present an analysis of the nonconvex AG algorithm by \cite{Ghadimi2015} to illustrate the algorithm as a generalization of Nesterov's AG. We also present formal results about the effect of hyperparameter settings on the complexity upper bound. Section \ref{sec:Simulation-Studies} will include simulation studies for linear and logistic models penalized by SCAD and MCP penalties. The simulation studies show that i) The AG method using our proposed hyperparameter settings converges faster than commonly used first-order methods for data with various $q/n$ and covariate correlation settings; and ii) our method outperforms the current state-of-the-art method, i.e. {\tt ncvreg}, in terms of signal recovery performance, especially when the signal-to-noise ratios are low. The proofs for the theorems are included in the Appendix \ref{sec:proof}.

\section{Motivation and Setup \label{sec:Accelerated-Gradient-Method}} 


Having built on Nesterov's seminal work, \citeauthor{Ghadimi2015}~\citep{Ghadimi2015} 
considered the following composite optimization problem: 
\begin{equation}
    \min_{x\in\mathbb{R}^{q+1}}\Psi\left(x\right)+\chi\left(x\right),\ \Psi\left(x\right)\coloneqq f\left(x\right)+h\left(x\right),
\tag{$\mathcal{P}$}\label{eq:optproblem}
\end{equation}
where $f\in\mathcal{C}_{L_{f}}^{1,1}(\mathbb{R}^{q+1},\mathbb{R})$ is convex, 
$h\in\mathcal{C}_{L_{h}}^{1,1}(\mathbb{R}^{q+1},\mathbb{R})$ is possibly nonconvex, and $\chi$ is a convex function over a bounded domain, and $\mathcal{C}_{L}^{1,1}$ denotes the class of first-order Lipschitz smooth functions with $L$ being the Lipschitz constant. They devised Algorithm \ref{alg:Accelerated-Gradient-Algorithm} discussed in details in next section, and presented a theoretical analysis of their algorithm. 

Some commonly used nonconvex penalties, such as SCAD and MCP, have a form that can naturally be decomposed into summation of a convex and a nonconvex function satisfying the conditions required by \citeauthor{Ghadimi2015}~\citep{Ghadimi2015}. When such penalties are added to a smooth convex deviance measure, such as negative of typical log-likelihoods, the resulting optimization problem follows the form of optimization problem \ref{eq:optproblem}. As we show below this is, in particular, the case when the deviance measure is a quadratic loss and the penalty is either SCAD or MCP. The quadratic loss plays the role of $f$. The other two functions, i.e. $h$ and $\chi$ are specified for both SCAD and MCP penalties. Define 


\begin{equation}p_{\lambda,a,\text{SCAD}}\left(\boldsymbol{\beta}\right)=\chi\left(\boldsymbol{\beta}\right)+h_{\text{SCAD}}\left(\boldsymbol{\beta}\right),
\label{eq:SCAD}
\end{equation}
\begin{equation}p_{\lambda,\gamma,\text{MCP}}\left(\boldsymbol{\beta}\right)=\chi\left(\boldsymbol{\beta}\right)+h_{\text{MCP}}\left(\boldsymbol{\beta}\right);
\label{eq:MCP}
\end{equation}
where $\boldsymbol{\beta}\coloneqq\left[\beta_{0},\beta_{1},\dots,\beta_{q}\right]^{T}$, $\chi\left(\boldsymbol{\beta}\right)=\sum_{j=1}^{q}\lambda\lvert\beta_{j}\rvert$, and 
\begin{align}h_{\text{SCAD}}\left(\boldsymbol{\beta}\right)=\sum_{j=1}^{q}\begin{cases}
0; & \lvert\beta_{j}\rvert\leq\lambda\\
\frac{2\lambda\lvert\beta_{j}\rvert-\beta_{j}^{2}-\lambda^{2}}{2\left(a-1\right)}; & \lambda<\lvert\beta_{j}\rvert<a\lambda\\
\frac{1}{2}\left(a+1\right)\lambda^{2}-\lambda\lvert\beta_{j}\rvert; & \lvert\beta_{j}\rvert\geq a\lambda
\end{cases} & \in\mathcal{C}_{L_{\text{SCAD}}}^{1,1}\label{eq:SCAD-Lsmoothness}\\
h_{\text{MCP}}\left(\boldsymbol{\beta}\right)=\sum_{j=1}^{q}\begin{cases}
-\frac{\beta_{j}^{2}}{2\gamma}; & \lvert\beta_{j}\rvert<\gamma\lambda\\
\frac{1}{2}\gamma\lambda^{2}-\lambda\lvert\beta_{j}\rvert; & \lvert\beta_{j}\rvert\geq\gamma\lambda
\end{cases} & \in\mathcal{C}_{L_{\text{MCP}}}^{1,1}\label{eq:MCP-Lsmoothness}
\end{align}
In the above equations, $\lambda>0,a>2,\gamma>1$ are the penalty tuning parameters. It is trivial that, in \eqref{eq:SCAD} and \eqref{eq:MCP},
$\chi\left(\boldsymbol{\beta}\right)$ is convex and the remaining term is
a first-order smooth concave function. 
In view of the optimization problem \ref{eq:optproblem}, when applying SCAD/MCP on a convex $\mathcal{C}_{L_{\ell}}^{1,1}$ statistical learning objective function, $f=-2\ell$ will be the convex component; $h_{\text{SCAD}},h_{\text{MCP}}$ will be the smooth nonconvex component with $L_{SCAD}=\frac{1}{a-1}$ and $L_{MCP}=\frac{1}{\gamma}$; and $\chi=\sum_{j=1}^{q}\lambda\lvert\beta_j\rvert$ will be the nonsmooth convex component. For high-dimensional statistical learning problems, 
the L-smoothness constant for the smooth nonconvex component, $L_{SCAD}$ and $L_{MCP}$, are often negligible when compared to the greatest singular value of the design matrix~\citep{Meckes2021}. In statistical learning applications, most unconstrained problems can, in fact, be reduced to problems over a bounded domain, as information often suggests the boundedness of the variables.

\section{The Accelerated Gradient Algorithm} \label{sec:aga}
This Section comprises two subsections. Subsection \ref{sec:nagm} includes an algorithm proposed by \citeauthor{Ghadimi2015}~\citep{Ghadimi2015} for solving the composite optimization problem \ref{eq:optproblem}. In Subsection \ref{sec:hnagm} we propose an approach for selecting the hyperparameters of the algorithm by minimizing the complexity upper bound \eqref{eq:complexity-bound}

\subsection{Nonconvex Accelerated Gradient Method} \label{sec:nagm}

Building on Nesterov's AG algorithm, \citeauthor{Ghadimi2015}~\citep{Ghadimi2015} proposed the following algorithm for solving the composite optimization problem \ref{eq:optproblem}. 

\begin{algorithm}[H] 
\caption{Accelerated Gradient Algorithm\label{alg:Accelerated-Gradient-Algorithm}}
\begin{algorithmic}
\Require{starting point $x_0\in \mathbb{R}^{q+1}$, $\{\alpha_k\}$ s.t. $\alpha_1=1$ and $\forall k\geq 2, 0<\alpha_k<1$, $\{\omega_k>0\}$, and $\{\delta_k>0\}$}
\Ensure{Minimizer $x_N^{md}$}
\State 0. Set $x_0^{ag}=x_0$ and $k=1$
\State 1. Set 
\begin{equation}
\label{eqn:x-md}
x_{k}^{md}=\alpha_{k}x_{k-1}^{ag}+\left(1-\alpha_{k}\right)x_{k-1}
\end{equation}
\State 2. Compute $\nabla\Psi\left(x_{k}^{md}\right)$ and set
\begin{align}
x_{k} &= \begin{cases}
x_{k-1}-\delta_{k}\nabla\Psi\left(x_{k}^{md}\right) & \text{(smooth)}\\
\mathcal{P}\left(x_{k-1},\nabla\Psi\left(x_{k}^{md}\right),\delta_{k}\right) & \text{(composite)}
\end{cases}\label{eqn:x-k}\\
x_{k}^{ag} &= \begin{cases}
x_{k}^{md}-\omega_{k}\nabla\Psi\left(x_{k}^{md}\right) & \text{(smooth)}\\
\mathcal{P}\left(x_{k}^{md},\nabla\Psi\left(x_{k}^{md}\right),\omega_{k}\right) & \text{(composite)}
\end{cases}\label{eqn:x-ag}
\end{align}
\State 3. Set $k=k+1$ and go to step 1
\end{algorithmic}
\end{algorithm}

In Algorithm \ref{alg:Accelerated-Gradient-Algorithm}, ``smooth''
represents the updating formulas for smooth problems, and ``composite''
represents the update formulas for composite problems, and $\mathcal{P}$
is the proximal operator defined as: 
\[
\mathcal{P}\left(x,y,c\right)\coloneqq\arg\min_{u\in\mathbb{R}^{q+1}}\left\{ \left\langle y,u\right\rangle +\frac{1}{2c}\left\Vert u-x\right\Vert ^{2}+\chi\left(u\right)\right\} .
\]
It is evident that the composite counter-part of the algorithm is
the Moreau envelope smoothing of the simple nonconvex function; for
this reason, in later analysis of the algorithm, we will use smooth
updating formulas for the sake of parsimony. As an interpretation
of the algorithm, $\left\{ \alpha_{k}\right\} $ controls the damping
of the system, and $\omega_{k}$ controls the step size for the ``gradient
correction'' update for momentum method. In what follows, $\Gamma_{k}$ is defined recursively as:
\[
\Gamma_{k}\coloneqq\begin{cases}
1, & k=1;\\
\left(1-\alpha_{k}\right)\Gamma_{k-1}, & k\geq2.
\end{cases}
\] 
\citeauthor{Ghadimi2015}~\citep{Ghadimi2015} proved that under the following conditions: 
\begin{align}
 & \alpha_{k}\delta_{k}\leq\omega_{k}<\frac{1}{L_{\Psi}},\ \forall k=1,2,\dots N-1\text{ and}\label{eq:convcond1}\\
 & \frac{\alpha_{1}}{\delta_{1}\Gamma_{1}}\geq\frac{\alpha_{2}}{\delta_{2}\Gamma_{2}}\geq\cdots\geq\frac{\alpha_{N}}{\delta_{N}\Gamma_{N}},\label{eq:convcond2}
\end{align}
the rate of convergence for composite optimization problems can be
illustrated by the following complexity upper bound: 
\begin{align}
 & \min_{k=1,\dots,N}\left\Vert \mathcal{G}\left(x_{k}^{md},\nabla\Psi\left(x_{k}^{md}\right),\omega_{k}\right)\right\Vert ^{2}\nonumber \\
 & \hspace{1.2 in} \leq \left[\sum_{k=1}^{N}\Gamma_{k}^{-1}\omega_{k}\left(1-L_{\Psi}\omega_{k}\right)\right]^{-1}\left[\frac{\left\Vert x_{0}-x^{*}\right\Vert ^{2}}{\delta_{1}}+\frac{2L_{h}}{\Gamma_{N}}\left(\left\Vert x^{*}\right\Vert ^{2}+M^{2}\right)\right].\label{eq:complexity-bound}
\end{align}
In the above inequality, $\mathcal{G}\left(x_{k}^{md},\nabla\Psi\left(x_{k}^{md}\right),\omega_{k}\right)$
is the analogue to the gradient for smooth functions defined
by: 
\[
\mathcal{G}\left(x,y,c\right)\coloneqq\frac{1}{c}\left[x-\mathcal{P}\left(x,y,c\right)\right].
\]
In accelerated gradient settings, $x$ corresponds to the past iteration,
$y$ corresponds to the smooth gradient at $x$, and $c$ corresponds
to the step size taken. 

\subsection{Hyperparameters for Nonconvex Accelerated Gradient Method} \label{sec:hnagm}

Here we discuss how hyperparameters, $\alpha_k$, $\omega_k$ and $\delta_k$ can be selected to accelerate convergence of Algorithm \ref{alg:Accelerated-Gradient-Algorithm} by minimizing the complexity upper bound. 
From Lemma \ref{lem:convergence-cond-meaning}, it is clear that the conditions \eqref{eq:convcond1} and \eqref{eq:convcond2} merely present a lower
bound for the vanishing rate of $\left\{ \alpha_{k}\right\} $. We
also observe that the right-hand side of \eqref{eq:bifurcation} is
monotonically increasing with respect to $\alpha_{k}$; thus, to obtain
the maximum values for $\left\{ \alpha_{k}\right\} $, it is sufficient
to maximize $\alpha_{k}$ recursively. 

Using \eqref{eqn:x-md}, \eqref{eqn:x-k}, and \eqref{eqn:x-ag}, we have 
\begin{align*}
\frac{x_{k+1}^{md}-\left(1-\alpha_{k+1}\right)x_{k}^{ag}}{\alpha_{k+1}}= & \frac{x_{k}^{md}-\left(1-\alpha_{k}\right)x_{k-1}^{ag}}{\alpha_{k}}-\delta_{k}\nabla\Psi\left(x_{k}^{md}\right)\text{ and}\\
x_{k}^{ag}= & x_{k}^{md}-\omega_{k}\nabla\Psi\left(x_{k}^{md}\right).
\end{align*}
By sorting out the terms in the above equations, we obtain the following
updating formulas: 
\begin{align}
x_{k}^{ag}= & x_{k}^{md}-\omega_{k}\nabla\Psi\left(x_{k}^{md}\right)\label{eq:grad-correction}\\
x_{k+1}^{md}= & x_{k}^{ag}+\alpha_{k+1}\cdot\left(\frac{1}{\alpha_{k}}-\frac{\delta_{k}}{\omega_{k}}\right)\cdot\left(\omega_{k}\nabla\Psi\left(x_{k}^{md}\right)\right)+\alpha_{k+1}\cdot\left(\frac{1}{\alpha_{k}}-1\right)\left(x_{k}^{ag}-x_{k-1}^{ag}\right)\label{eq:nonconvex-momentum}
\end{align}


Compared to Nesterov's AG, the AG method proposed by Ghadimi and
Lan differs by the convergence conditions \eqref{eq:convcond1} and
\eqref{eq:convcond2}, and the inclusion of the term $\alpha_{k+1}\cdot\left(\frac{1}{\alpha_{k}}-\frac{\delta_{k}}{\omega_{k}}\right)\cdot\left(\omega_{k}\nabla\Psi\left(x_{k}^{md}\right)\right)$
in \eqref{eq:nonconvex-momentum}. Since $\alpha_{k+1}\cdot\left(\frac{1}{\alpha_{k}}-\frac{\delta_{k}}{\omega_{k}}\right)\geq0$
is implied by convergence condition \eqref{eq:convcond1}, this added
term functions as a step to reduce the magnitude of ``gradient correction''
presented in \eqref{eq:grad-correction}: the resulting framework
will keep the same momentum compared to Nesterov's AG, but the momentum
step update will occur at a midpoint between $x_{k}^{ag}$ and $x_{k}^{md}$
to yield $x_{k+1}^{md}$. Such a framework suggests that the proposed
algorithm is merely a midpoint generalization in the gradient correction
step of Nesterov's AG. Therefore, \emph{the acceleration occurs to
the convex component $f$ of the objective function $\Psi$}. Following
this intuition, we proceed to investigate the optimization hyperparameter
settings for the most accelerating effect in Theorem \ref{thm:convex}
based on the idea of minimizing the complexity upper bound \eqref{eq:complexity-bound}
when the objective function is convex; i.e., when $h \equiv 0$. 

It can be deduced from \eqref{eq:bifurcation} that an increasing
sequence of $\left\{ \delta_{k}\right\} $ allows a slower vanishing
rate for $\left\{ \alpha_{k}\right\} $. Specifically, the existence
of $\delta_{1}$ in \eqref{eq:complexity-bound} can be explained
as the following: the momentum initialization step in Algorithm \ref{alg:Accelerated-Gradient-Algorithm}
indicates that $x_{1}^{md}=x_{0}^{ag}=x_{0}$. We also have $x_{1}^{ag}=x_{1}^{md}-\omega_{1}\nabla\Psi\left(x_{1}^{md}\right)=x_{0}^{ag}-\omega_{1}\nabla\Psi\left(x_{0}\right)$
for smooth problems or $x_{1}^{ag}=\mathcal{P}\left(x_{1}^{md},\nabla\Psi\left(x_{1}^{md}\right),\omega_{1}\right)=\mathcal{P}\left(x_{0}^{ag},\nabla\Psi\left(x_{0}\right),\omega_{1}\right)$
for composite problems. In view of \eqref{eq:nonconvex-momentum},
the momentum initializes as $x_{1}^{ag}-x_{0}^{ag}=-\omega_{1}\nabla\Psi\left(x_{0}\right)$
for smooth problems. Thus, should $\delta_{1}<\omega_{1}$
take a smaller value, $\alpha_{2}\cdot\left(\frac{1}{\alpha_{1}}-\frac{\delta_{1}}{\omega_{1}}\right)>0$;
i.e., $x_{2}^{md}$ is a convex combination of $x_{1}^{ag}$ and the
initial point $x_{0}$, and the smaller $\delta_{1}$ is, the closer
$x_{2}^{md}$ is to $x_{0}$. Meanwhile, a smaller $\delta_{1}$
allows a faster increasing sequence $\left\{ \delta_{k}\right\} $;
hence a slower-vanishing sequence $\left\{ \alpha_{k}\right\} $ can
be achieved to incorporate more momentum. This process can be interpreted as follows: when $x_{2}^{md}$ does not retain the full step update
from the initial point $x_{0}$, more initial momentum will be allowed
to accumulate, as the initial momentum is in the same direction as
the update. We therefore choose $\delta_{1}=\omega_{1}$;
i.e., to let $x_{2}^{md}$ retain fully the update from $x_{0}$ in
the direction of $-\omega_{1}\nabla\Psi\left(x_{0}\right)$, such that
no \emph{excess} initial momentum will be needed to account for initial
update deficiency in this direction.

\section{Theoretical Analysis of the Algorithm}\label{sec:theory}


For gradient methods without a line-search procedure, the
step size for the gradient correction is usually set to be a constant.
Based on this convention, we assume $\omega_{k}=\beta$ for $k=1,2,\dots,N$.
Theorem \ref{thm:convex} below presents the optimal choice of hyperparameters under mild conditions. 
\begin{thm}
\label{thm:convex} Assume 
conditions \eqref{eq:convcond1} and \eqref{eq:convcond2} hold. Let $\delta_{1}=\omega_{k}=\omega$ and $h=0$.
Then the complexity upper bound \eqref{eq:complexity-bound} is minimized by: 
\begin{align}
\bar{\alpha}_{k+1}= & \frac{2}{1+\sqrt{1+\frac{4}{\bar{\alpha}_{k}^{2}}}},\ \bar{\alpha}_{1}=1\label{eq:alpha-settings},\\
\bar{\delta}_{k+1}= & \frac{\bar{\omega}}{\bar{\alpha}_{k+1}}\label{eq:lambda-setting},\\
\bar{\omega}= & \frac{2}{3L_{\Psi}}\label{eq:beta-setting}.
\end{align}
\end{thm}
\begin{proof}
See Appendix \ref{sec:proof-thm-1}.
\end{proof}

As illustrated by the proof of the above theorem, the optimization hyperparameter
settings \eqref{eq:alpha-settings}, \eqref{eq:lambda-setting}, and
\eqref{eq:beta-setting} allow for the greatest values of $\left\{ \alpha_{k}\right\} $
under the constant gradient-correction step size and maximum initial
update assumptions; i.e., condition 1. Such settings allow the most
acceleration for the convex component. Although a greater momentum
will result in a much faster convergence at the initial stage of the
algorithm, it will also result in oscillations of larger magnitudes
near the minimizer. Therefore, in the following theorem, we will show
that the complexity upper bound will always maintain $O\left(1/N\right)$
rate of convergence. This observation implies that the accelerated gradient method's worst-case scenario is at least as good as $O\left(1/N\right)$ for gradient descent in terms of the rate of convergence. 

\begin{thm}
\label{thm:1overk} Assume conditions \eqref{eq:convcond1}
and \eqref{eq:convcond2} hold. Then under the assumptions of Theorem
\ref{thm:convex}
, the complexity upper bound is $O\left(1/N\right)$. 
\end{thm}
\begin{proof}
See Appendix \ref{sec:proof-thm-2}.
\end{proof}

The recursive formula for optimal momentum hyperparameter, $\left\{ \alpha_{k}\right\} $, as presented in \eqref{eq:alpha-settings},
is of a rather complicated structure. The next theorem illustrates the vanishing rate of $\left\{ \alpha_{k}\right\} $.  

\begin{thm}
\label{thm:alpha-k-vanishing} 
Let $\bar{\alpha}_{1}=1$ and \eqref{eq:alpha-settings} holds. Then
\begin{equation}
\label{eq:damping-bounds}
\frac{2}{\left(1+a\cdot k^{-b}\right)k+1}<\bar{\alpha}_{k}\leq\frac{2}{k+1}, \quad k=1,\dots,N, 
\end{equation}
for any $a>0,\ 0<b<1$, such that 
\begin{equation}
    \label{eq:k-vanishing-cond} a\left(1-b\right)\cdot2^{2-b}-ab\left(1-b\right)\cdot2^{-b}-1\geq0.
\end{equation}
\end{thm}
\begin{proof}
See Appendix \ref{sec:proof-thm-3}.
\end{proof}

The following corollary establishes a tight bound for the damping sequence, hence providing the speed of convergence of our proposed 
optimal damping sequence $\left\{ \bar{\alpha}_{k}\right\} $
to $\frac{2}{k+1}$. 

\begin{cor} \label{cor:alpha-vanishing}
The lower bound in \eqref{eq:damping-bounds} is maximized at
\[
\bar{a}_k=\frac{2^{\bar{b}_k}}{\left(1-\bar{b}_k\right)\left(4-\bar{b}_k\right)} \quad \text{and} \quad  
\bar{b}_k=\frac{2+5\left(\log\frac{2}{k}\right)+\sqrt{9\left(\log\frac{2}{k}\right)^{2}+4}}{2\left(\log\frac{2}{k}\right)}\label{optimal-b} \quad \text{for} \quad k \geq 8. 
\]
The lower bound in \eqref{eq:damping-bounds} therefore becomes 
\begin{equation}
\frac{k+1}{2}-\bar{\alpha}_{k}^{-1}=O\left(\log k\right)
\end{equation}
\end{cor}


\begin{proof}
See Appendix \ref{sec:proof-cor-1}.
\end{proof}

\begin{figure}[H]
\centering{}\includegraphics[scale=.8]{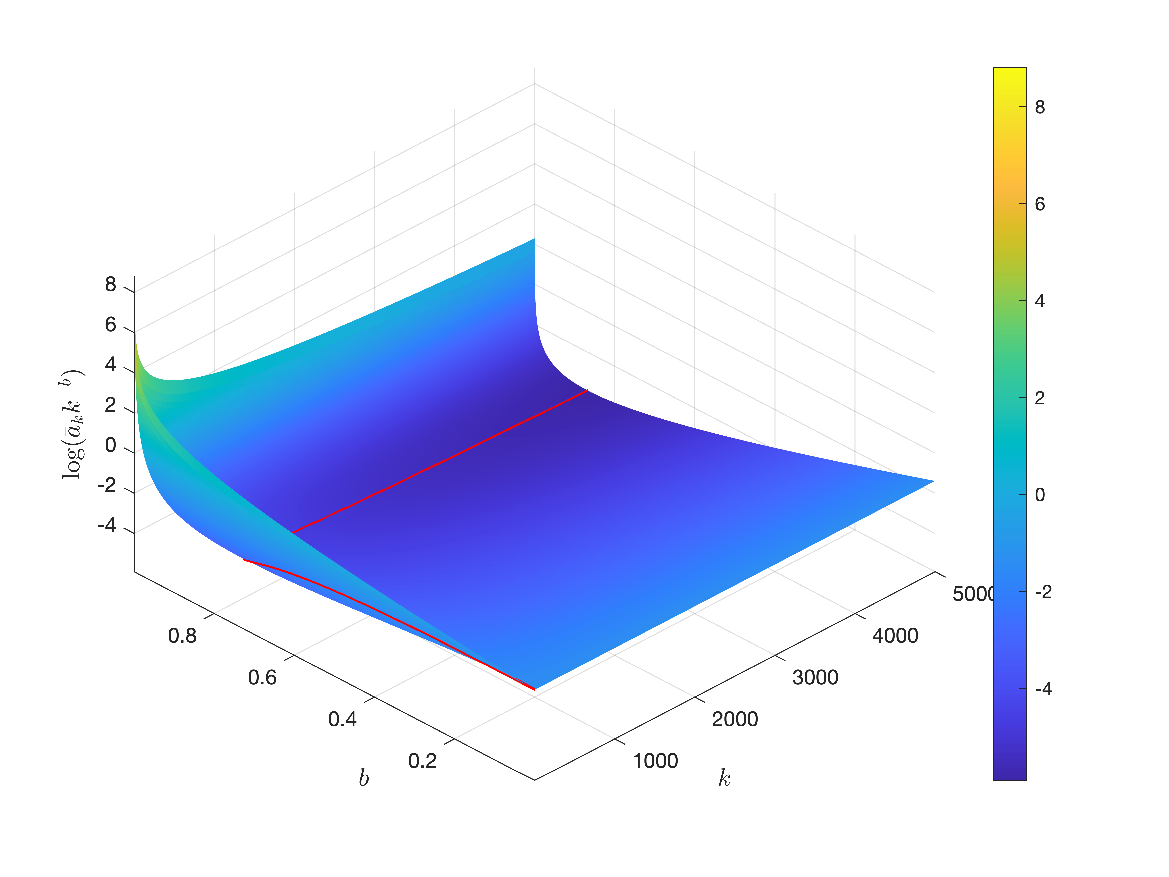}\caption{Numerical plots for Corollary \ref{cor:alpha-vanishing}. The figure plots $\log\left(\bar{a}_kk^{-b}\right)$ v.s. $k$ and $b$; the red line plots its minimizer $\bar{b}_k=\frac{2+5\left(\log\frac{2}{k}\right)+\sqrt{9\left(\log\frac{2}{k}\right)^{2}+4}}{2\left(\log\frac{2}{k}\right)}$ for each $k$. The plot reflects on the speed for the coefficient of $k$ in the denominator of the lower bound in \eqref{eq:damping-bounds} converges to $1$. The red line shows that $\bar{b}_k$ converges to $1$ at an extremely slow rate. \label{fig:numerical-plot}}
\end{figure}

To better illustrate Corollary \ref{cor:alpha-vanishing}, we plot the value of $\log\left(\bar{a}_kk^{-b}\right)$ v.s. $(k,b)$ in Figure \ref{fig:numerical-plot}.
The plot shows that as $k$ grows large, the optimizer $\bar{b}_k$ converges to $1$ at a very slow rate. It also reflects on the speed of $1+\bar{a}_k\cdot k^{-\bar{b}_k}$, the coefficient of $k$ in the denominator of the lower bound in \eqref{eq:damping-bounds}, goes to $1$ as $k$ increases.

\section{Simulation Studies\label{sec:Simulation-Studies} }


In this section, we conduct two sets of simulation studies for nonconvex penalized linear and logistic
models. We first visualize the convergence rates and signal recovery performance for each set of simulation studies using a single simulation replicate. Second, we compare the convergence rates across the first-order methods with varying $q/n$ ratios and covariate correlations for $100$ simulation replications. Lastly, we compare the signal recovery performance using our method to the state-of-the-art method, {\tt ncvreg}~\citep{Breheny2011}, with varying covariate correlations and \ac{SNR} for $100$ simulation replications. Since the iterative complexity differs for the first-order methods and coordinate descent methods, the convergence rates in terms of the number of iterations are not directly comparable. { Thus, we choose to compare the computing time between AG, proximal gradient descent, and coordinate descent.}

\subsection{Simulation Setup} \label{sec:sim-setup}

Linear models with the OLS loss function is a popular method for modelling a continuous response. We aim to achieve signal recovery by solving the following problem for penalized linear models:
\[
\arg\min_{\boldsymbol{\beta}\in\mathbb{R}^{q+1}}\frac{1}{2n}\left\Vert \mathbf{X}\boldsymbol{\beta}-\mathbf{y}\right\Vert _{2}^{2}+\sum_{j=1}^{q}p_{\lambda}\left(\beta_{j}\right),
\]
where $p_\lambda:\mathbb{R}\mapsto\mathbb{R}_{\geq 0}$ is the SCAD or MCP penalty function. 
To compare the convergence rates across the first-order methods, we choose different $q/n$ ratios and the strength of correlation, $\tau$, between the covariates. These two parameters are most likely to impact the convergence rates. Median and corresponding $95\%$ bootstrap confidence intervals from $1000$ bootstrap replications for the number of iterations required for the iterative objective values to make a fixed amount of descent are reported. To compare the signal recovery performance between our AG method and the state-of-the-art package {\tt ncvreg}, we performed $100$ simulation replications with varying SNRs and covariate correlations, as they directly impact the signal recovery performance. The simulation studies we performed adapt the following setups:

\begin{itemize}
\item The total number of observations $n=1000$ for visualization plots and signal recovery performance comparison, and $n=200,500,1000,3000$ for convergence rate { and computing time} comparisons.
\item For visualization purposes, we perform one simulation replicate with the number of covariates $q=2004$, with $4$ nonzero signals being $2,-2,8,-8$. We perform $100$ simulation replications with the number of covariates $q = 2050$, with $5$ blocks of “true” signals equal-spaced with $500$ zeros in-between for convergence rate { and computing time comparison, as well as} signal recovery performance comparison. For each simulation replicate, the blocks of the ``true'' signals are simulated from $N_{10}\left(0.5,1\right)$, $N_{10}\left(5,2\right)$, $N_{10}\left(10,3\right)$, $N_{10}\left(20,4\right)$, $N_{10}\left(50,5\right)$, respectively. 
\item The design matrix, $\mathbf{X}$, is simulated from a multivariate Gaussian distribution
with mean $0$. The covariance matrix $\boldsymbol{\Sigma}$ is a $\tau-$Toeplitz matrix, where $\tau=0.5$ for the visualization plots and $\tau=0.1,0.5,0.9$ for the convergence rate { and computing time comparison, as well as} signal recovery performance comparison. All covariates are standardized; i.e., centered by the sample mean and scaled by the sample standard deviation.
\item The signal-to-noise ratio is set as $\text{SNR}=\frac{\sqrt{\boldsymbol{\beta}_{\text{true}}^{T}\boldsymbol{\Sigma}\boldsymbol{\beta}_{\text{true}}}}{\sigma}$, where $\boldsymbol{\beta}_{\text{true}}$ are the ``true'' coefficient values, and $\sigma$ is used as the residual standard deviation. $\text{SNR}=5$ for visualization plots, $\text{SNR}=3$ for convergence rate comparison, and $\text{SNR}=1,3,7,10$ for signal recovery performance comparison.
\item { For visualization plots, convergence rate and computing time comparisons,} we take $\lambda=0.5,a=3.7$ for SCAD and $\lambda=0.5,\gamma=3$ for MCP, unless otherwise specified. For signal recovery rate comparison, $\lambda$ sequence consists of
$50$ values equal-spaced from $\lambda_{\max}$\footnote{$\lambda_{\max}$ is the minimal value for $\lambda$ such that all
penalized coefficients are estimated as $0$.} to $0$. The tuning parameter $\lambda$ is chosen to minimize the (non-penalized) loss function value on a validation set of the same size as the training set.
\item { For signal recovery performance comparison, we use the same objective function as {\tt ncvreg} to ensure that the same value of penalty tuning parameters results in the same degree of penalization. We also adapt the same strong rule setup as {\tt ncvreg}~\citep{Lee2015}.}
\end{itemize}

{ To compare the gradient-based methods and the coordinate descent method, we compare the computing time when both coded in {\tt Python/CuPy}. The coordinate descent method was coded based on the state-of-the-art pseudo-code~\citep{Breheny2011}. All of the computing was carried out on a NVIDIA A100 GPU with CUDA compute capability of 8.0 on the Narval computing cluster from Calcul Qu\`ebec/Compute Canada. Furthermore, we also excluded the computation of the L-smoothness parameter for the coordinate descent method in our simulations.}

The simulation setups for penalized logistic models are similar to those above for penalized linear models, except that the active coefficients are set differently to account for the exponential scale inherent to the logistic regression. For the single-replicate visualization simulations, we let the $4$ nonzero signals be $0.5,-0.5,0.8,-0.8$. For the simulations with $100$ replications to compare the convergence rate and signal recovery performance, we simulate the $5$ blocks of the ``true'' signals from $N_{10}\left(0.5,1\right)$, $N_{10}\left(0.5,1\right)$,
$N_{10}\left(-0.5,1\right)$, $N_{10}\left(-0.5,1\right)$, $N_{10}\left(1,1\right)$, respectively. The SNR for logistic regression has the same definition as linear models, with Gaussian noise added to the generated continuous predictor $\mathbf{X}\boldsymbol{\beta}_{true}$. The binary outcomes are independent Bernoulli realizations, with probabilities being the logistic transforms of the continuous response.

\subsection{Simulation Results} \label{sec:sim-results}

\subsubsection{Penalized Linear Regression} \label{sec:linear-sim}

\begin{figure}[ht]
\centering{}\includegraphics[scale=0.5]{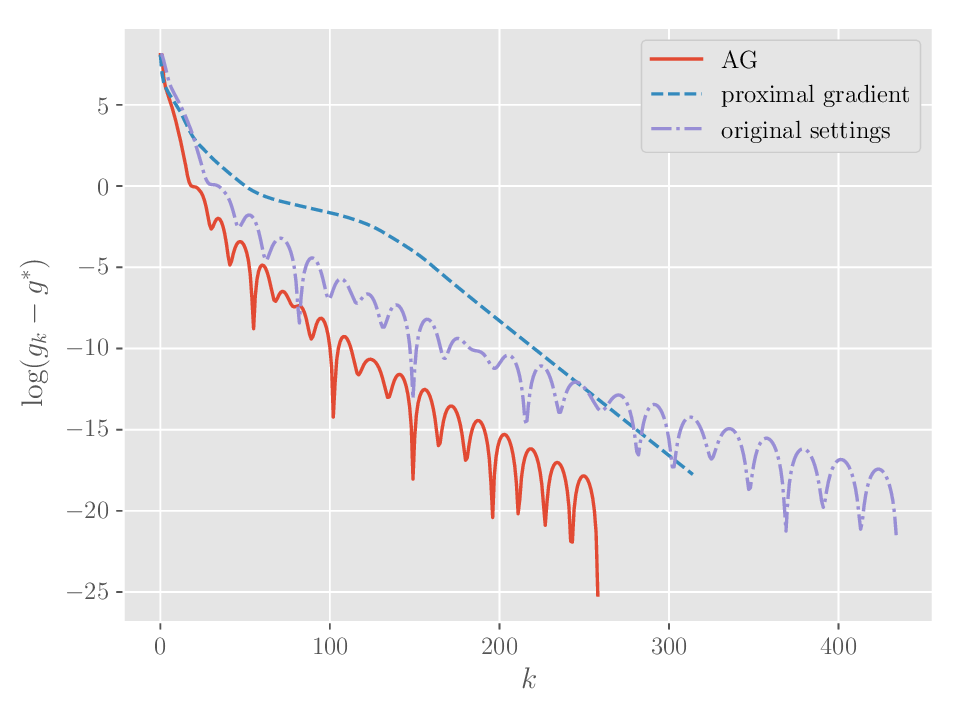}\includegraphics[scale=0.5]{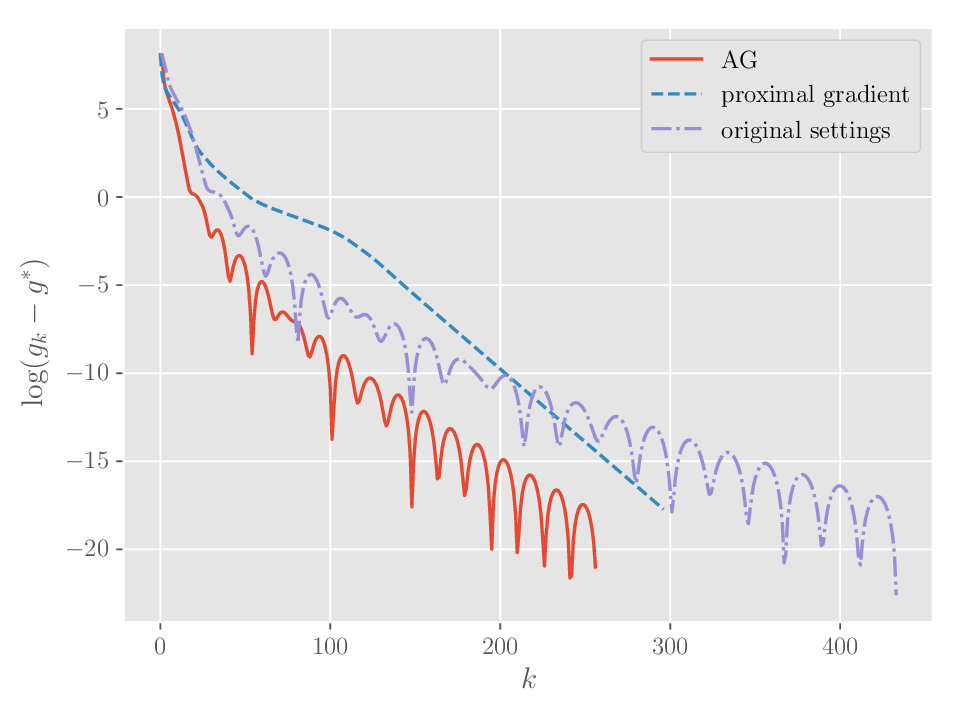}\caption{Convergence rate performance of first-order methods on SCAD (left) and MCP (right)
penalized linear model for a single simulation replicate. $k$ represents the number of iterations, $g_k$ represents the iterative objective function value, and $g^*$ represents the minimum found by the three methods considered. \label{fig:LM-iteration-plot}}
\end{figure}

Figure \ref{fig:LM-iteration-plot} shows the log differences of iterative objective values for a single replicate. This figure visualizes the accelerating effect of the AG method using our proposed hyperparameter settings. Median with the corresponding $95\%$ bootstrap CI of the number of iterations required for the iterative objective function values to make a fixed amount of descent for $100$ simulation replications are reported in Figures \ref{fig:sim-ag-lm-scad}, \ref{fig:sim-ag-lm-mcp} in Appendix \ref{sec:Simulations-LM}. The lack of bars in the reported barplots indicates that the median of $100$ replications breaks down; i.e., the corresponding proximal gradient algorithm fails to converge to the minimizer found by the three algorithms within $2000$ iterations. 
The AG method using our hyperparameter settings converges
much faster than proximal gradient and AG using the original hyperparameter settings
proposed by \citeauthor{Ghadimi2015} for both SCAD and MCP-penalized models discussed
here, as reflected in Figures \ref{fig:LM-iteration-plot}, \ref{fig:sim-ag-lm-scad}, \ref{fig:sim-ag-lm-mcp}. It can also be observed that momentum methods such
as AG are much less likely to be stuck at saddle points or local minimizers
than proximal gradient -- this property is consistent with previous findings~\citep{Jin2017}. Since the proposed AG methods belong to the class
of momentum methods, the AG algorithms do not possess a descent property. As
suggested by a previous study~\citep{Su2014}, oscillation will occur
at the end of the trajectory; the descent property will therefore
vanish. This is also reflected in Figures \ref{fig:LM-iteration-plot},
\ref{fig:logistic-iteration-plot} -- as the trajectory moves close
to the optimizer, the oscillation will start to occur for the AG methods. Among all the first-order methods, the AG method with our proposed hyperparameter settings tends to converge the fastest in all scenarios considered, as illustrated by Figures \ref{fig:sim-ag-lm-scad}, \ref{fig:sim-ag-lm-mcp} in Appendix \ref{sec:Simulations-LM}. The observed standard errors among $100$ simulation replications are rather small, suggesting that the halting time retains predictable for high-dimensional models, which agrees with the recent findings~\citep{Paquette2020}. 

{ Figures \ref{fig:lm-time-scad}, \ref{fig:lm-time-mcp} report median with the corresponding $95\%$ bootstrap CI of the computing time (in seconds) required for the infinity norm of the two consecutive iterations $\left\Vert \boldsymbol{\beta}^{\left(k+1\right)}-\boldsymbol{\beta}^{\left(k\right)}\right\Vert _{\infty}$ to fall below $10^{-4}$ for $100$ simulation replications. It can be observed that the computing time for AG with suggested settings is much shorter than the computing time for coordinate descent.}

\begin{figure}[ht]
\begin{centering}
\includegraphics[scale=0.8]{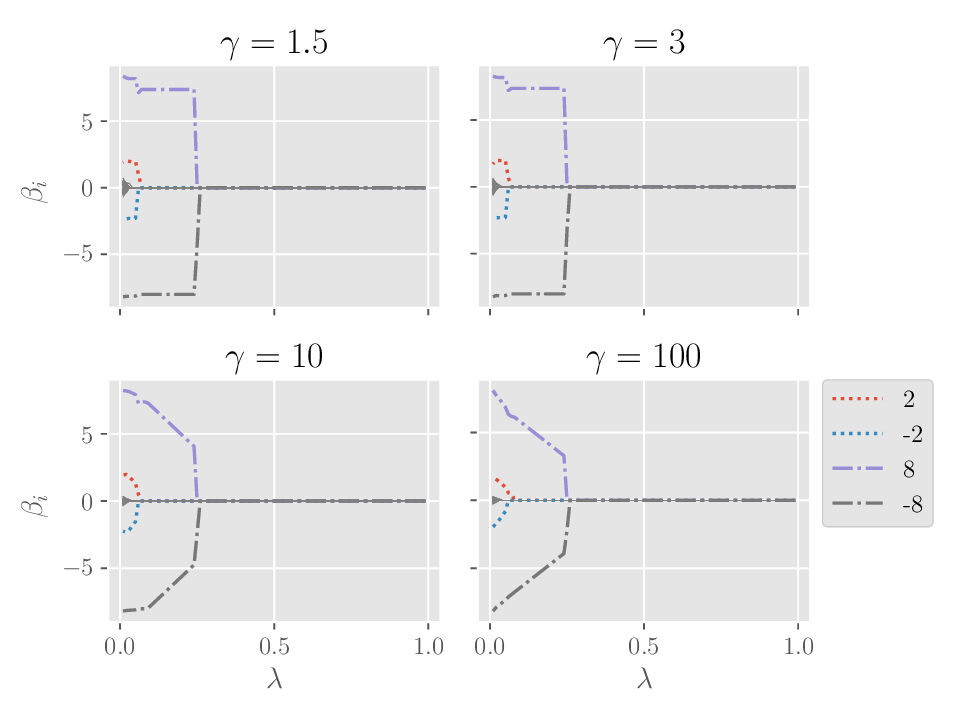}
\par\end{centering}
\caption{Solution paths obtained using the proposed AG method for MCP-penalized linear
model with different values of $\gamma$ for a single simulation replicate. The behaviors of the solution path match the expected from the MCP penalized problems. The solution path behaves similarly to hard-thresholding for a small $\gamma$. As $\gamma$ increases, the solution path will behave more similarly to soft-thresholding. \label{fig:linear-signal}}
\end{figure}

To visualize the signal recovery performance using our proposed method, Figure \ref{fig:linear-signal} plots the solution paths for
the MCP-penalized linear model with different values of $\gamma$. The grey lines in Figure \ref{fig:linear-signal} represent the recovered values for the noise variables. AG method performs very well when applied to signal recovery problems for nonconvex-penalized linear models. Figure \ref{fig:linear-signal} serves as an arbitrary instance that the recovered signals using our method exhibit the expected pattern with MCP -- as $\lambda$ decreases, the degree of penalization decreases, and more false-positive signals will be selected. The stable solution path for the recovered signals suggests that the algorithm does not converge to a point far away from the ``true'' coefficients.

\begin{figure}[ht]
\begin{centering}
\includegraphics[scale=0.8]{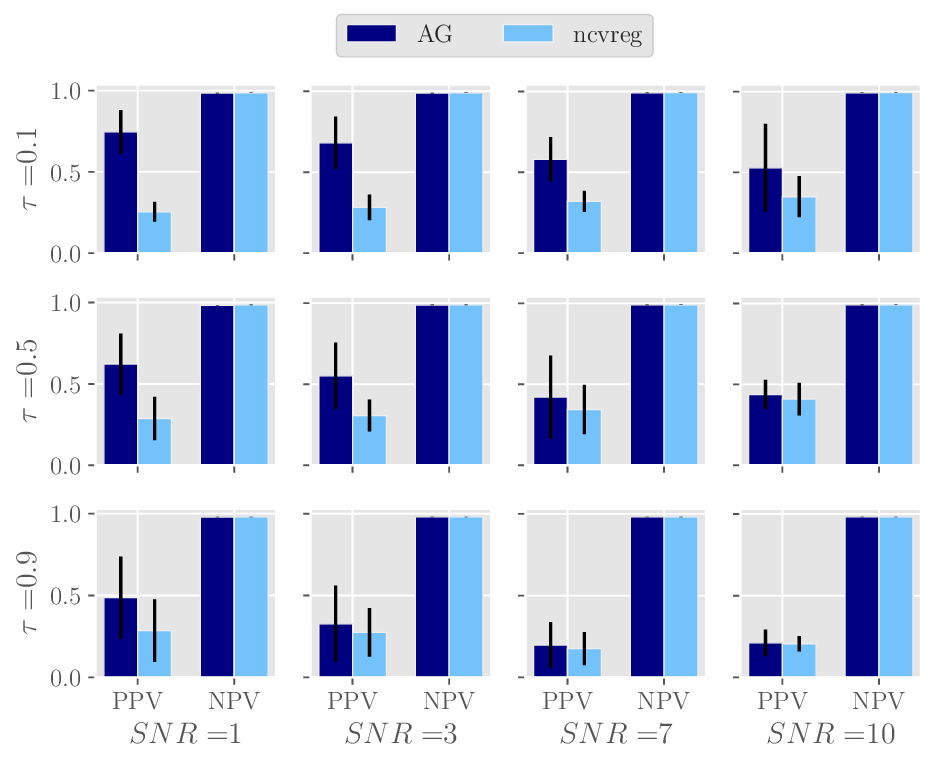}
\par\end{centering}
\caption{Sample means for Positive/Negative Predictive Values (PPV, NPV) of signal detection across different values of covariates correlation ($\tau$) and SNRs for AG with our proposed hyperparameter settings and {\tt ncvreg} on SCAD-penalized linear model over $100$ simulation replications. The error bars represent the standard errors.
\label{fig:LM-SCAD-pv}}
\end{figure}

To further illustrate the signal recovery performance, the means and standard errors for the scaled estimation error $\frac{\left\Vert \boldsymbol{\beta}_{\text{true}}-\hat{\boldsymbol{\beta}}\right\Vert _{2}^{2}}{\left\Vert \boldsymbol{\beta}_{\text{true}}\right\Vert _{2}^{2}}$, positive/negative predictive values (PPV, NPV), and active set cardinality across $100$ replications are reported in Tables \ref{tab:signal-lm-scad} and \ref{tab:signal-lm-mcp} in Appendix \ref{sec:Simulations-LM}. In what follows, $\mathcal{A}$ denotes the set of nonzero ``true'' coefficients and $\hat{\mathcal{A}}$ denotes the set of nonzero coefficients selected by the model. PPV and NPV use the following definitions: 
\begin{equation*}
\text{PPV}\coloneqq \frac{\lvert\mathcal{A}\cap\hat{\mathcal{A}}\rvert}{\lvert\hat{\mathcal{A}}\rvert},\ \text{NPV}\coloneqq \frac{\lvert\mathcal{A}^{C}\cap\hat{\mathcal{A}}^{C}\rvert}{\lvert\hat{\mathcal{A}}^{C}\rvert}.
\end{equation*}
Sample means and standard errors for PPV and NPV from Table \ref{tab:signal-lm-scad} are further visualized in Figure \ref{fig:LM-SCAD-pv}. When applied to sparse
learning problems, the signal recovery performance of our proposed
method often outperforms {\tt ncvreg}, the current state-of-the-art method~\citep{Breheny2011}, particularly in terms of the positive predictive values (PPV). This can be observed from \figref{fig:LM-SCAD-pv} and Tables \ref{tab:signal-lm-scad}, \ref{tab:signal-lm-mcp} from Appendix \ref{sec:Simulations-LM}. This observation is especially evident when the signal-to-noise ratios are low. At the same time, $\ensuremath{\nicefrac{\left\Vert \boldsymbol{\beta}_{\text{true}}-\hat{\boldsymbol{\beta}}\right\Vert _{2}^{2}}{\left\Vert \boldsymbol{\beta}_{\text{true}}\right\Vert _{2}^{2}}}$ for both methods are close. As the SNR increases, the validation set becomes more similar to the training set, causing the chosen model to have a smaller $\lambda$. The model size will therefore increase, which will decrease the value of PPV.

\subsubsection{Penalized Logistic Regression} \label{sec:logistic-sim}

\begin{figure}[H]
\centering{}\includegraphics[scale=0.5]{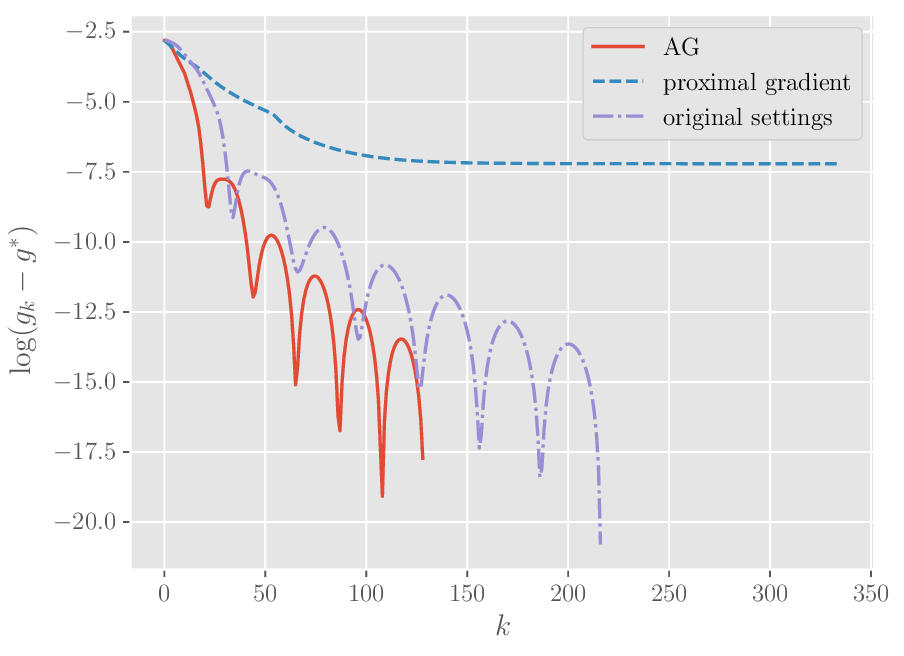}\includegraphics[scale=0.5]{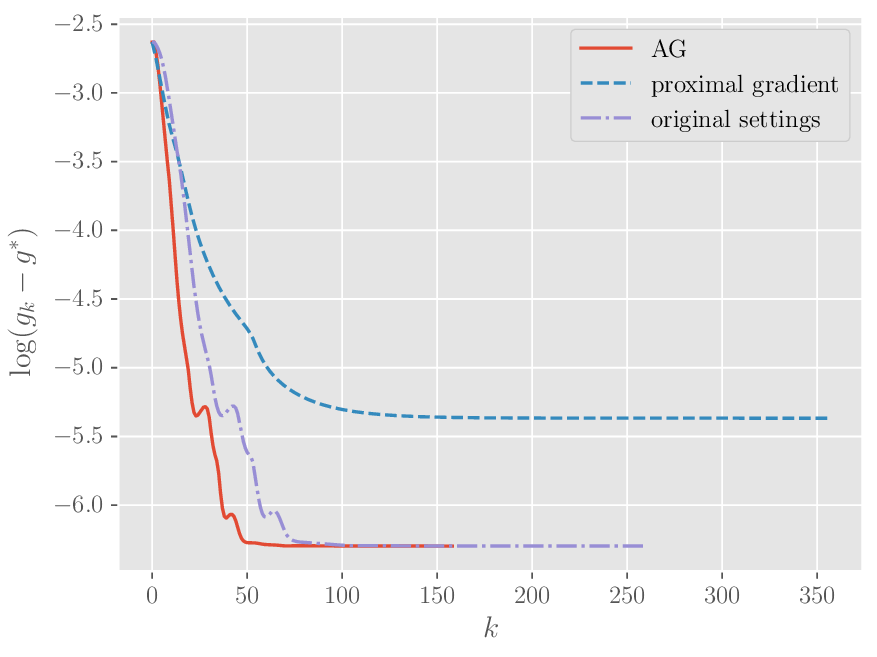}\caption{Convergence rate performance of first-order methods on SCAD (left) and MCP (right)
penalized logistic regression for a single simulation replicate. $k$ represents the number of iterations, $g_k$ represents the iterative objective function value, and $g^*$ represent the minimum found by the three methods considered. \label{fig:logistic-iteration-plot}}
\end{figure}

\begin{figure}[ht]
\begin{centering}
\includegraphics[scale=0.8]{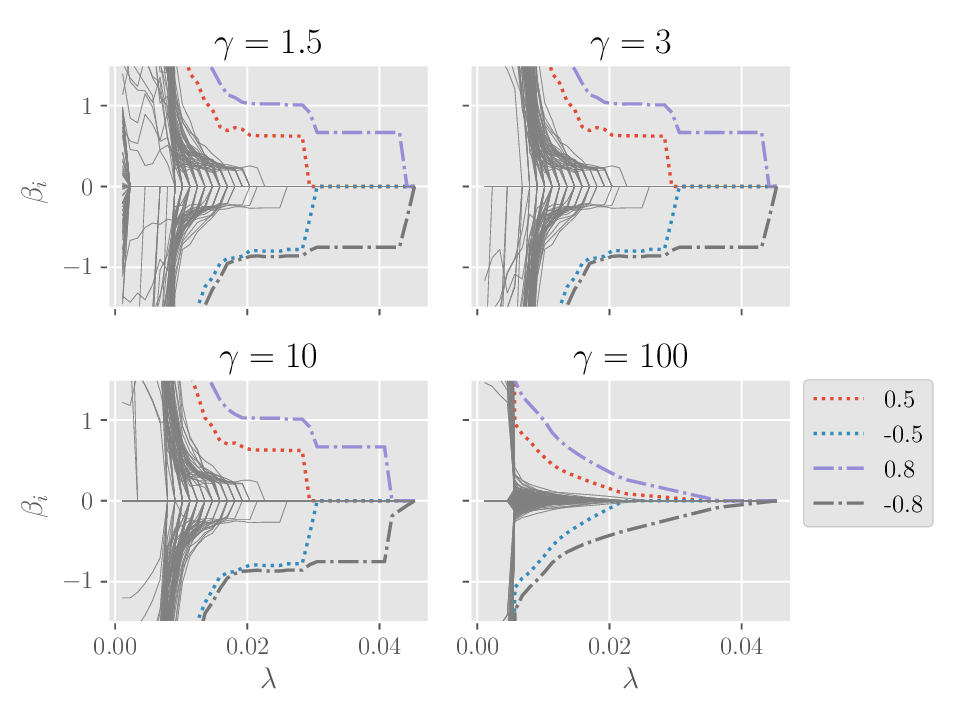}
\par\end{centering}
\caption{Solution paths obtained using the proposed AG method for MCP-penalized logistic
regression with different values of $\gamma$ for a single simulation replicate. The behaviors of the solution path match the expected from the MCP penalized problems. The solution path behaves similarly to hard-thresholding for a small $\gamma$. As $\gamma$ increases, the solution path will behave more similarly to soft-thresholding. \label{fig:logistic-signal}}
\end{figure}

\begin{figure}[ht]
\begin{centering}
\includegraphics[scale=0.8]{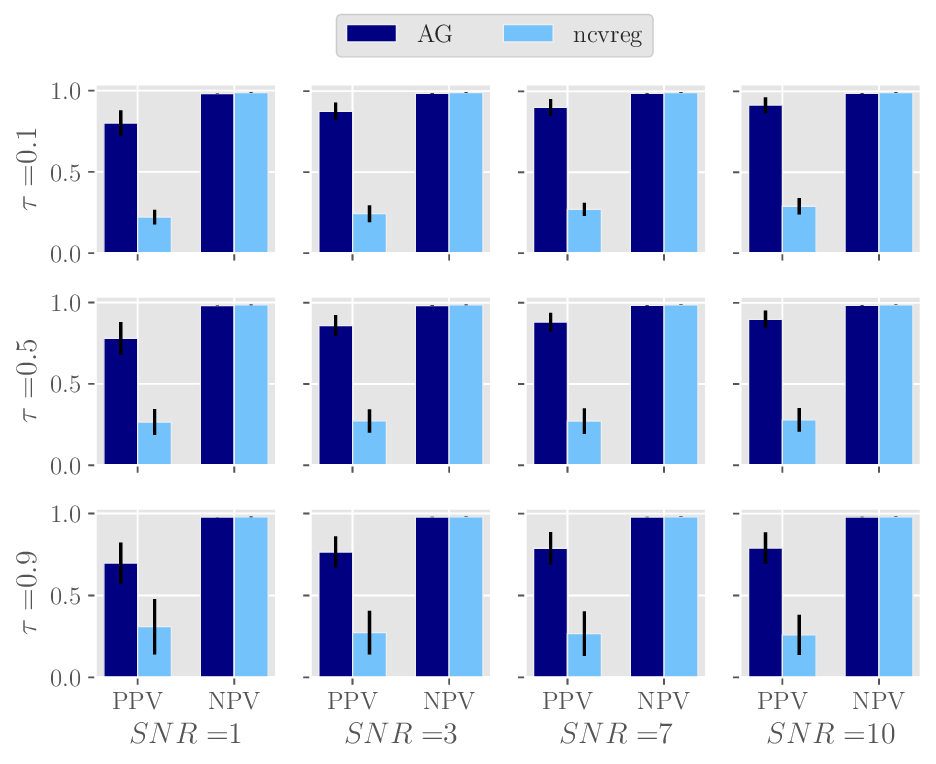}
\par\end{centering}
\caption{Sample means for Positive/Negative Predictive Values (PPV, NPV) of signal detection across different values of covariates correlation ($\tau$) and SNRs for AG with our proposed hyperparameter settings and {\tt ncvreg} on SCAD-penalized logistic model over $100$ simulation replications. The error bars represent the standard error.
\label{fig:logistic-SCAD-pv}}
\end{figure}

The simulation results reflected in Figures \ref{fig:logistic-iteration-plot}, \ref{fig:logistic-signal}, as well as Figures \ref{fig:sim-ag-logistic-scad}, \ref{fig:sim-ag-logistic-mcp} and Tables \ref{tab:signal-logistic-scad}, \ref{tab:signal-logistic-mcp} in Appendix \ref{sec:Simulations-logistic} suggest similar findings for penalized logistic models to our findings for penalized linear models as discussed in Section \ref{sec:linear-sim}. We further note that when applied to penalized logistic models, the coordinate descent method often fails to converge, resulting in overall poor performance in positive predictive values as reflected in Figure \ref{fig:logistic-SCAD-pv} and Tables \ref{tab:signal-logistic-scad}, \ref{tab:signal-logistic-mcp} in Appendix \ref{sec:Simulations-logistic}. When it does converge, the coordinate descent method does so at a very slow rate. In comparison, our proposed method has a convergence guarantee in theory and converges within a reasonable number of iterations in our simulation studies, as shown in Figures \ref{fig:sim-ag-lm-scad}, \ref{fig:sim-ag-lm-mcp} in Appendix \ref{sec:Simulations-logistic}. { In our computing time comparison, we used identical simulation setups and convergence standard for both the AG method and coordinate descent method, running both on a NVIDIA A100 GPU with CUDA compute capability of 8.0 from Compute Canada; the submitted simulation job finished well within $20$ minutes for both SCAD and MCP-penalized logistic models when using the AG method, but exceeded the 7-day computing time limit imposed on the Narval cluster when using the coordinate descent method.}

\section{Discussion}


We considered a recently developed generalization of Nesterov's accelerated gradient method for nonconvex optimization, and we have discussed its potential in sparse statistical learning with nonconvex penalties. An important issue concerning this algorithm is the selection of its sequences of hyperparameters. We present an explicit solution to this problem by minimizing the algorithm's complexity upper bound, hence accelerating convergence of the algorithm. 
Our simulation studies indicate that among first-order methods, the AG method using our proposed hyperparameter settings achieves a convergence rate considerably faster than other first-order methods such as the AG method using the original proposed hyperparameter settings or proximal gradient. Our simulations also show that signal recovery using our proposed method generally outperforms {\tt ncvreg}, the current state-of-the-art method. This performance gain is much more pronounced for penalized linear models when the signal-to-noise ratios are low. For penalized logistic regression, the performance gain observed is consistent across various covariates correlation and signal-to-noise ratio settings. Compared to coordinate-wise minimization methods, our proposed method is less challenged by low signal-to-noise ratios and is feasible to implement in parallel. { Given today's computing facilities, parallel computing is particularly meaningful for large datasets~\citep{Parnell2020}. We also show this gain in parallel computing performance by comparing computing time on a GPU.} Furthermore, our proposed method has weaker convergence conditions and can be applied to a class of problems that do not have an explicit solution to the coordinate-wise objective function. For example, linear mixed models for grouped or longitudinal data involve the inverse of a large covariance matrix. Decomposition of this covariance matrix is necessary to apply the coordinate descent method. However, such decomposition can be computationally costly and numerically unstable~\citep{Quarteroni2007}. On the other hand, matrix decomposition is not needed for first-order methods, as numerically stable yet computationally efficient approaches such as conjugate gradient can be adapted when applying our proposed method. The proposed nonconvex AG method can be applied to a wide range of statistical learning problems, opening various future research opportunities in statistical machine learning and statistical genetics. 

\section{Disclaimer}
All codes to reproduce the simulation results of this paper and outputs from Calcul Quebec/Compute Canada can be found on the following GitHub repository:

\href{https://github.com/Kaiyangshi-Ito/nonconvexAG}{https://github.com/Kaiyangshi-Ito/nonconvexAG}

\chapter{Tsallis Entropy Maximizing Distributions for Robust and Efficient Sparse Learning on Correlated Data}\label{ch:paper3}

\indent \textbf{Preamble to Manuscript 3.}

\subsection*{Introduction to the Study and Its Place in the Workflow:}
Manuscript 3 explores the critical limitations of Gaussian assumptions often made in statistical models, particularly those used to analyze correlated and heterogeneous data typical in biostatistical applications. By proposing the use of the $q$Gaussian distribution, derived from Tsallis entropy maximization, this manuscript introduces a robust alternative capable of accommodating the correlated observations often inherent in biostatistical data, such as genetic and longitudinal studies. This novel approach significantly enhances the flexibility and robustness of statistical models, making it an invaluable addition to the techniques developed in the previous manuscripts.

\subsection*{Building on Optimization and Computational Advances:}
This manuscript extends the efficient computational methodologies refined in Manuscripts 1 and 2 by integrating them into a broader modeling context that includes correlated observations. After establishing efficient variable screening and optimization techniques for sparse learning, the introduction of a new modeling framework that can effectively handle correlations and heterogeneities addresses the next layer of complexity in the analysis of high--dimensional data. The $q$Gaussian model not only offers a solution to the robustness issues posed by the underlying distributional assumptions and heavy tails, but also fits within the computational framework previously developed.

\subsection*{Innovation in Statistical Modeling and Optimization:}
The framework for adapting numerical methods originally designed to find equilibria in flows to tackle composite optimization problems presents a novel approach to address the challenges of statistical computing in sparse learning. This methodology ensures that the models developed are not only theoretically sound but also practically applicable.

\subsection*{Enhancing Data Analysis in Biostatistics:}
By applying the innovative framework mentioned above to the Hager--Zhang conjugate gradient algorithm, Manuscript 3 develops a numerically stable and computationally efficient algorithm for sparse statistical learning. This advancement is crucial for efficiently processing high--dimensional large datasets that biostatistics often deals with. The robustness offered by the $q$Gaussian distribution transforms the landscape of statistical machine learning, making it more robust, hence better suited to the nuanced challenges posed by high--dimensional biostatistical data with correlated observations.

\subsection*{Integration with Previous Manuscripts:}
Manuscript 3 synthesizes and builds on the computational and methodological foundations laid in the first two manuscripts. The variable screening from Manuscript 1 ensures that the most relevant variables are identified for robust modeling, while the optimization techniques from Manuscript 2 provide an efficient optimization algorithm to handle the statistical computing challenges introduced by the $q$Gaussian modeling of high--dimensional data. Together, these manuscripts create a comprehensive workflow for handling high-dimensional biostatistical data, from initial robust screening and modeling of complex data structures to efficient statistical computing algorithms. 



\newpage

\vspace*{2cm}

\begin{center}
	\Large{Tsallis Entropy Maximizing Distributions for Robust and Efficient Sparse Learning on Correlated Data}
\end{center}

\vspace*{3cm}

\begin{center}
Kai~Yang$^{1}$, Masoud~Asgharian$^{2}$, Celia~M.T.~Greenwood$^{1,3}$.
\end{center}

\vspace*{1cm}

\begin{center}
$^{1}$Department of Epidemiology, Biostatistics, and Occupational Health, McGill University 
$^{2}$Department of Mathematics and Statistics, McGill University 
$^{3}$Lady Davis Institute for Medical Research, Montr\'{e}al 
\end{center}
 
\vspace*{1cm}

 
\newpage

\section*{Abstract}

This paper addresses the limitations of Gaussian distribution assumptions
in statistical sparse learning, particularly in modeling correlated
and heterogeneous data. Conventional Gaussian models often lack
robustness towards outliers and underlying distribution assumptions.
To overcome these limitations, we propose the use of the $q$Gaussian
distribution, derived from Tsallis entropy maximization, as a robust
alternative. This is notably relevant in biostatistics, where the presence of correlated observations and heterogeneity, such as in genetic and longitudinal studies, is prevalent.  Our contributions include modeling of correlated data through the re-derived multivariate probability density function from Tsallis entropy maximization, thereby addressing the limitations inherent in conventional Gaussian models. 
Furthermore, we introduce a novel framework that adapts numerical
methods designed for finding equilibria in flows to tackle composite
optimization problems prevalent in statistical sparse learning. Applying
this framework to the Hager-Zhang conjugate gradient algorithm \citep{Hager2005}, we
develop a numerically stable and efficient algorithm for sparse statistical
learning. The $q$Gaussian distribution, informed by the principle
of maximizing Tsallis entropy, presents a viable and flexible alternative
to Gaussian-based methods, potentially transforming the landscape
of statistical machine learning. This paper not only contributes to
the theoretical understanding of statistical distributions and optimization
techniques, but also paves the way for practical data analysis in biostatistics
and related fields.
 
\newpage



\section{Introduction }

In the realm of statistical sparse learning, the pursuit of robust
and efficient methodologies remains paramount, especially when confronted
with the complexities of correlated data. The principle of maximizing
Shannon's entropy stands as a pivotal framework that has led to the
derivation of nearly all frequently utilized statistical distributions
to date \citep{Cover2006}. This principle's application has notably
revealed that the multivariate Gaussian distribution maximizes Shannon's
entropy under first moment and second central moment constraints,
significantly influencing the landscape of statistical sparse learning.
The Gaussian assumption has become a fundamental 
cornerstone of numerous statistical sparse learning problem formulations, and its core assumptions are rarely re-examined or challenged. 

However, the Gaussian distribution's features, particularly its exponential
tail decay and the absence of a shape parameter, can present substantial
limitations. Specifically, the lack of robustness towards outliers
and a limited capacity to accurately represent the distribution's
shape results in violations of the Gaussian assumption in statistical
modeling. Such violations have many practical repercussions, including
the potential for erroneous Type I error rates and the lack of robustness
towards distribution shape when estimating the dispersion parameter,
motivating the development of alternative approaches.

Dispersion or volatility parameters encapsulate critical and often
decisive information about distributions. Their estimations, specifically
in transformations of predicted outcomes, are often indispensable.
For instance, in the context of log-normal distributions, the mean
is directly influenced by the volatility parameter derived from the
underlying Gaussian distribution. Likewise, principles like the Law
of the Unconscious Statistician (LOTUS), which rely on accurate estimation
of volatility and precise understanding of the distribution's shape,
highlight the importance of determining this parameter for dependable
prediction and statistical modeling. 

Volatility estimation is of great importance in finance. Specifically, Ito’s lemma, often used in stochastic calculus for option pricing, explicitly highlights the importance of volatility's contribution. Delving into the realm of stochastic calculus, Ito's lemma provides
a mathematical framework that elegantly captures volatility's impact
on dynamic systems. Ito's lemma states that for a twice-differentiable
function $f$, 
\begin{equation}
df\left(t,X_{t}\right)=\left(\frac{\partial f}{\partial t}+\mu\frac{\partial f}{\partial x}+\frac{1}{2}\sigma^{2}\frac{\partial^{2}f}{\partial x^{2}}\right)dt+\sigma\frac{\partial f}{\partial x}dW_{t},\label{eq:Itos-chain-rule}
\end{equation}
where the term $\sigma^{2}\frac{\partial^{2}f}{\partial x^{2}}$ specifically
denotes the contribution of volatility to changes in the function
$f$. This mathematical representation is pivotal in finance, where
the phenomenon, termed \emph{volatility smile}, challenges the foundational
assumptions of the \ac{BSM} model, signaling empirical
deviations from expected normality in option pricing models. These
deviations have propelled the exploration of alternative distributions
capable of more accurately reflecting market realities \citep{Pena1999}.

In response to these limitations of Gaussian distributions, the $q$Gaussian
distribution, derived from maximizing Tsallis entropy, emerges as a compelling
alternative. The $q$Gaussian distribution is celebrated for its
flexibility in modeling the diverse shapes of bell-curved distributions,
including the ability to account for heavy-tailed distributions ---
a feature crucial for the robust modeling of financial returns. It
provides a more accurate representation of financial returns on platforms
such as the \emph{\ac{NYSE}} and \emph{\ac{NASDAQ}} \citep{Borland2002,Borland2002a,Domingo2017}.
Despite its proven advantages in finance, the incorporation of Tsallis
entropy-maximizing distributions within the domain of statistical
sparse learning and biostatistics remains limited. To the best of
our knowledge, this paper represents the initial endeavor to apply
Tsallis entropy-maximizing distributions for biostatistical data modeling.
Correlated observations, frequently encountered in genetic and longitudinal
studies \citep{Garcia2017,Runcie2019,DandineRoulland2015}, as well
as heterogeneity of the variance, 
will be specifically addressed by our proposed Tsallis entropy-maximizing model for correlated data.

Hence, this paper advocates for the application
of the $q$Gaussian distribution in modeling correlated data within
sparse statistical  learning frameworks. Our approach relaxes the
conventional reliance on normality assumptions;  
We aim to demonstrate that the intricate characteristics of $q$Gaussian
distributions can profoundly enhance the modeling of correlated data,
offering a robust and versatile alternative to conventional Gaussian-based
methods. 

Maximum likelihood estimation is one of the most commonly used estimation
techniques. However, the estimation process encounters notable computational
obstacles when dealing with high--dimensional and extra-large datasets.
Oracle penalties, favored for their efficacy in facilitating variable
selection, present an attractive yet complex solution to sparse learning
problems. However, oracle penalties are notable for their nonconvex and nonsmooth nature \citep{Nikolova2000},
which lead to considerable optimization challenges. Recently, proximal
methods have demonstrated an unmatched speed of convergence, thereby
surpassing most other approaches in efficiently handling estimation
in nonsmooth problems \citep{Hoheisel2020}. Simultaneously, the Krylov
subspace method, recognized among the top ten algorithms for computing in science and engineering 
of the twentieth
century, lays a solid foundation for numerical analysis. The conjugate
gradient method, a prominent member of the Krylov subspace methods,
has been applied extensively in various areas and is a fundamental
numerical tool in solving partial differential equations \citep{Nocedal2000}. While the nonlinear conjugate gradient performs exceptionally well in terms of its convergence speed and numerical stability, much better than accelerated gradient or gradient descent, its global convergence depends on the line--search step, whereas the accelerated gradient and gradient descent methods do not necessarily require the line search step to achieve global convergence \citep{Ghadimi2015, Yang2024}.
Motivated by these methodologies, our paper introduces a proximal conjugate
gradient method that can be applied to solve $q$Gaussian sparse learning problems. This method aims to effectively combine the theoretical strengths of both proximal methods and Krylov subspace techniques. Additionally, our paper addresses the line search step needed for the proximal nonlinear conjugate gradient method to establish global convergence. 

We re-derive the probability density function for the multivariate
$q$Gaussian distribution from a Tsallis entropy maximizing perspective in Lemma \ref{lem:density-normalization}.
This derivation allows for a nuanced understanding and application
of this model in statistical analysis. Furthermore, \emph{our contributions in this paper are as following: }
\begin{enumerate}
\item \emph{We apply the derived 
density to model correlated and heterogeneous data effectively,
while carrying out the sparse statistical learning at the same time. }
\item \emph{Sparse statistical learning involves minimizing a composite
optimization problem, aimed at minimizing a composite objective function
composed of a globally Lipschitz-smooth term, which may be nonconvex,
and a convex nonsmooth term. A variety of numerical methods are available
to find equilibrium points for globally Lipschitz flows. By employing
the Moreau envelope and linearizing the smooth term, we develop a
framework that allows any numerical method designed for finding equilibrium
points in globally Lipschitz flows to be adapted into a numerical
optimization algorithm for minimizing the composite objective function.} 
\item \emph{Leveraging the framework introduced above, we implement it with the state-of-the-art Hager-Zhang conjugate gradient method \citep{Hager2005}. This implementation yields a proximal conjugate gradient algorithm that is not only computationally efficient but also numerically stable, suitable for a wide range of statistical sparse learning challenges. This includes the robust sparse learning approach we devised based on the concept of maximizing the Tsallis entropy distribution. }
\end{enumerate}
The structure of the paper is organized as follows: 

Section \ref{sec:Tsallis-Entropy} delves into the foundational properties
of Tsallis entropy, drawing upon previous literature to establish
a comprehensive background. Following this, Section \ref{sec:Tsallis-Entropy-Maximizing}
introduces the concept of $q-$moments. This section then elaborates
on employing Tsallis entropy maximizing distribution to effectively
model the $q-$correlation structure. In Section \ref{sec:Tsallis-Entropy-Maximizing},
we also re-derive the probability density function maximizing Tsallis
entropy under the first and second central $q-$moment constraints,
incorporating all relevant parameters for a likelihood-based approach
to statistical analysis. 

Our discussion transitions to the challenges and strategies of optimization
in Section \ref{sec:optimization}. This section is twofold; initially,
in Section \ref{subsec:variational-analysis}, we present essential background
knowledge from variational and nonsmooth analysis. This foundation
is critical for our novel contribution: the development of a proximal
framework to transform any first-order numerical optimization algorithm
to a proximal counterpart by leveraging the properties of the Moreau
envelope, detailed in Section \ref{subsec:PCG-framework}. In Section \ref{subsec:Proximal-HZ},
we apply this innovative framework to the state-of-the-art Hager-Zhang
conjugate gradient algorithm. This adaptation produces a proximal
version for tackling sparse statistical learning challenges. The efficacy
of this method is further showcased in Section \ref{sec:optimizing-penalized-mle},
where we outline the application of our proximal Hager-Zhang conjugate
gradient algorithm to optimize a penalized $q$Gaussian likelihood
function. This section also lays out a map from problem formulation
to the practical aspects of prediction using models trained with our
approach.
Finally, Section \ref{sec:Conclusion-and-Discussion}
synthesizes our contributions, offering a reflective conclusion and
proposing avenues for future research. 

\section{\label{sec:Tsallis-Entropy} Tsallis Entropy }

For an arbitrary random variable $X$, Shannon's Entropy \citep{Shannon1948}
poses the definition 
\begin{equation}
H\left(X\right)\coloneqq-\mathbb{E}\log\left(p\left(X\right)\right)=-\int\log\left(p\left(x\right)\right)d\mu_{X},\label{eq:shannon-entropy}
\end{equation}
where $p$ is the likelihood function for $X$. Over a given (likelihood)
function space 
\begin{equation}
    \mathcal{P}\coloneqq\left\{ p\left(x\right)\vert\forall x\in\mathcal{X},\ p\left(x\right)\geq0,\left|\mathbb{E}\log\left(p\left(X\right)\right)\right|<\infty\text{ and }\int_{\mathcal{X}}1d\mu_{X}=1\right\} ,
    \notag
\end{equation}
the Shannon's entropy is a \emph{strictly concave} function, which
implies uniqueness of the maximizer. Many commonly-used distributions
have been shown to maximize Shannon's entropy under certain given
constraints \citep{Cover2006}. For example, uniform distribution,
whether in discrete or continuous case
, are to maximize \eqref{eq:shannon-entropy} over a compact support,
with open sets defined by discrete or Euclidean topology, respectively.
The exponential distribution is defined as maximizing \eqref{eq:shannon-entropy}
over $\mathbb{R}_{\geq0}$ and with a constraint that the first moment
is a constant, $\frac{1}{\lambda}$; where $\lambda$ later turns
out to be the scale parameter. And the Gaussian distribution maximizes
\eqref{eq:shannon-entropy} over $\mathbb{R}$ with given mean and
variance. More examples can be given. For example, the constraint
to obtain a Laplace distribution is a given mean absolute deviation,
etc. 

\emph{Additivity} is a key element of Shannon's entropy. That is,
let $A_{1},A_{2}$ be two independent event sets, then the information
of the intersection $I\left(\mathbb{P}\left(A_{1}\cap A_{2}\right)\right)=I\left(\mathbb{P}\left(A_{1}\right)\cdot\mathbb{P}\left(A_{2}\right)\right)=I\left(\mathbb{P}\left(A_{1}\right)\right)+I\left(\mathbb{P}\left(A_{2}\right)\right)$
--- such homomorphism was considered particularly useful in Shannon's
view \citep{Shannon1948}. 
Later in the 1980s, \cite{Tsallis1988} constructed an entropy
similar to Shannon's entropy but without the additivity property.
To see how Tsallis' entropy was developed, first we look at Tsallis'
$q-$exponential function, which is defined as $\exp_{q}:\mathbb{R}\mapsto\mathbb{R}$,
given by 
\begin{equation}
\exp_{q}x\coloneqq\begin{cases}
\left(\left(1+\left(1-q\right)x\right)\right)^{\frac{1}{1-q}}; & 1+\left(1-q\right)x>0\\
0; & \text{else}
\end{cases}\label{eq:tsalli-q-exp}
\end{equation}
For $q>1$, $\exp_{q}$ is bijective over $\left(0,\frac{1}{q-1}\right)$.
The inverse function, called the $q-$logarithmic function, is given
by 
\begin{equation}
\ln_{q}x\coloneqq\frac{x^{1-q}-1}{1-q}.\label{eq:tsalli-q-log}
\end{equation}
Based on this deformed $q-$exponential function, \cite{Tsallis1988}
developed \emph{Tsallis entropy} by replacing the $\log$ function
in \ref{eq:shannon-entropy} with $q-\log$ function \eqref{eq:tsalli-q-log}
and replacing the expectation with $q-$expectation operator \citep{Tsallis1988}:
\begin{align}
S_{q}\left(X\right) & =-\int_{\mathcal{X}}p^{q}\left(x\right)\ln_{q}p\left(x\right)dx\eqqcolon-\mathbb{E}_{q}\ln_{q}p\left(X\right)\label{eq:tsallis-entropy}\\
 & =\frac{1}{q-1}\left(1-\int_{\mathcal{X}}p^{q}\left(x\right)dx\right)
\end{align}
where $q\in\mathbb{R}\setminus\{1\}$ is a constant, and $\mathbb{E}_{q}f\left(X\right)\coloneqq\int_{\mathcal{X}}f\left(x\right)\cdot p^{q}\left(x\right)dx=\left\langle f\left(x\right),\left(d\mu_{X}\right)^{q}\right\rangle $
is referred to as the \emph{$q-$expectation operator}. Tsallis entropy
is also known as \emph{non-extensive} entropy; namely for arbitrary
independent two random variable $X_{1},X_{2}$: 
\begin{equation}
S_{q}(X_{1},X_{2})=S_{q}(X_{1})\oplus_{q}S_{q}(X_{2}),\label{eq:q-nonextensive}
\end{equation}
where ``$\oplus_{q}$'' is defined as $\forall a,b\in\mathbb{R}$,
\begin{equation}
a\oplus_{q}b\coloneqq a+b+\left(1-q\right)ab.
\end{equation}
Expectation has been used to characterize statistical distributions.
However, one significant drawback of the expectation (linear) operator
is the lack of continuity for some distributions; such as the Cauchy
distribution. Therefore, the $q-$expectation operator, $\mathbb{E}_{q}$,
provides robustness when characterizing the distributions in the real
domain. If the tail of the function vanishes at a rate of $O\left(\left(\log x\right)^{-1}\right)$, the function will not have a proper integral if the support
is unbounded. Thus, for any distribution whose likelihood function
is bounded in uniform norm, $\exists q\in\mathbb{R}_{>0}$ such that
$\mathbb{E}_{q}$ is continuous at the likelihood function in the
function space we are considering. 

\section{\label{sec:Tsallis-Entropy-Maximizing} Tsallis Entropy Maximizing
Distribution to Accommodate the $q-$Correlation Structure }

The Gaussian distribution maximizes the Shannon's entropy in the following
problem: 
\begin{align}
\max\ _{\phi\in \mathcal{P}} & -\int_{\mathbb{R}^{n}}\phi\left(x\right)\log\left(\phi\left(x\right)\right)dx\nonumber \\
\text{s.t. } & \phi\geq0;\nonumber \\
 & \int_{\mathbb{R}^{n}}\phi\left(x\right)dx=1;\nonumber \\
 & \int_{\mathbb{R}^{n}}x\cdot\phi\left(x\right)dx=0;\label{eq:shannon-gaussian-location}\\
 & \int_{\mathbb{R}^{n}}xx^{T}\cdot\phi\left(x\right)dx=\Sigma;\nonumber 
\end{align}
for some $n\in\mathbb{N}_{+}$ and $\Sigma\in\mathbb{R}^{n\times n},\ \Sigma\succ0$. For the sake of parsimony,
in \eqref{eq:shannon-gaussian-location} we assume that the distribution
is centered. To set the central trend parameter, or the mean parameter
in the specific case of the Gaussian distribution, the likelihood function
$\phi$ can be simply translated $x\mapsto x-\mu$ to incorporate
the parameter $\mu$ for the central trend. Entropy functions are invariant
under translation. 

Similarly to how the multivariate Gaussian distribution maximizes Shannon's
entropy in a Euclidean space, the multivariate $q$Gaussian distribution maximizes Tsallis entropy in a Euclidean space. Specifically, the
optimization problem is formulated as: 
\begin{align}
\max\ _{\phi\in L^{q}\left(\mathbb{R}^{n}\right)} & -\int_{\mathbb{R}^{n}}\phi^{q}\left(x\right)dx\label{eq:generalized-entropy-obj}\\
\text{s.t. } & \phi\geq0;\nonumber \\
 & \int_{\mathbb{R}^{n}}\phi\left(x\right)dx=1;\label{eq:normalization-constraint}\\
 & \frac{\int_{\mathbb{R}^{n}}x\cdot\phi^{q}\left(x\right)dx}{\int_{\mathbb{R}^{n}}\phi^{q}\left(x\right)dx}=0;\label{eq:q-first-moment-constraint}\\
 & \frac{\int_{\mathbb{R}^{n}}xx^{T}\cdot\phi^{q}\left(x\right)dx}{\int_{\mathbb{R}^{n}}\phi^{q}\left(x\right)dx}=\Sigma,\label{eq:q-second-moment-constraint}
\end{align}
where $q>1$. The feasible set of Lebesgue space $L^{q}\left(\mathbb{R}^{n}\right)$ is to ensure the well--definedness of Tsallis entropy. The normalization constraint \eqref{eq:normalization-constraint} implies that $\phi\in L^{1}$; however, $\phi\in L^{1}$ does not imply $\phi\in L^{q}$, as the embedding property $L^{1} \subseteq L^{q}$ fails to hold for Lebesgue measure on $\mathbb{R}^{n}$. As an example, consider the one-dimensional example of the probability density function 
\begin{equation}
    \tilde{\phi}\left(x\right)=\begin{cases}
\frac{1}{4}\left|x\right|^{-\frac{1}{2}} & \text{for }x\in\left(-1,1\right)\setminus\left\{ 0\right\} ;\\
0 & \text{else.}
\end{cases}
\end{equation}
Clearly, $\tilde{\phi}\in L^{1}$ but $\tilde{\phi}\not\in L^{2}$. Note that \eqref{eq:q-first-moment-constraint} and
\eqref{eq:q-second-moment-constraint} are the first and second moment
constraints using the $q-$expectation operator $\mathbb{E}_{q}$. As
noted by \cite{Tsukada2005}, maximizing any member of the generalized
class of power-law entropies, including Renyi entropy, Havrda and
Charvat entropy, Arimoto entropy, and Tsallis entropy, all yield the
identical power-law objective function \eqref{eq:generalized-entropy-obj}.
Regarding the constraints, $\int_{\mathbb{R}^{n}}\phi^{q}\left(x\right)dx$
is the normalization factor for the $q-$expectation. \cite{Tsukada2005}
further noted that optimizing the problem formulated above is equivalent
to the following problem: 
\begin{align}
\max\ _{\varphi\in L^{s}\left(\mathbb{R}^{n}\right)} & \int_{\mathbb{R}^{n}}\varphi^{s}\left(x\right)dx\label{eq:generalized-entropy-moment-obj}\\
\text{s.t. } & \varphi\geq0;\nonumber \\
 & \int_{\mathbb{R}^{n}}\varphi\left(x\right)dx=1;\nonumber \\
 & \int_{\mathbb{R}^{n}}x\cdot\varphi\left(x\right)dx=0;\nonumber \\
 & \int_{\mathbb{R}^{n}}xx^{T}\cdot\varphi\left(x\right)dx=\Sigma.\nonumber 
\end{align}
In \eqref{eq:generalized-entropy-moment-obj}, $s\coloneqq q^{-1}\in\left(0,1\right)$, thus $L^{s}\left(\mathbb{R}^{n}\right)$ is a quasi-normed space;
$\varphi\left(x\right)\coloneqq\frac{\phi^{q}\left(x\right)}{\int_{\mathbb{R}^{n}}\phi^{q}\left(x\right)dx}.$
If the maximizer of \eqref{eq:generalized-entropy-moment-obj} is
$\varphi$, then the maximizer of \eqref{eq:generalized-entropy-obj},
$\phi$, will be normalized 
\begin{equation}
\phi\left(x\right)\propto\varphi^{1/q}.
\end{equation}
Several important properties were proposed previously
regarding the $q$Gaussian distributions in previous studies \citep{Vignat2004,Costa2003}.
Notably, 
\begin{enumerate}
\item \label{enu:maximizing-format} Using Bregman information divergence,
Problem \eqref{eq:generalized-entropy-moment-obj} has a unique maximizer
of the form 
\begin{equation}
\varphi\left(x;s\right)=A_{s}\left(1-\left(s-1\right)\beta^{\prime}\left\langle x,\Sigma^{-1}x\right\rangle \right)_{+}^{\frac{1}{s-1}}\label{eq:dual-maximizer}
\end{equation}
for some $s\in\left(\frac{n}{n+2},\infty\right)\setminus\left\{ 1\right\} $,
normalization constant $A_{r}$, and some dispersion parameter $\beta^{\prime}$. 
\item \label{enu:linear-mapping-invariance} If $X\sim q\text{Gaussian}\left(q,\Sigma\right)$,
$H\in\mathbb{R}^{\tilde{n}\times n}$ and $\text{rank}\left(H\right)=\tilde{n}$.
Then $\tilde{X}\sim q\text{Gaussian}\left(\tilde{q},H\Sigma H^{T}\right)$
with 
\begin{equation}
\frac{2}{1-\tilde{q}^{-1}}-\tilde{n}=\frac{2}{1-q^{-1}}-n.\label{eq:linear-mapping-invariance}
\end{equation}
 
\item \label{enu:linear-combination-not-closed} If $X_{1},X_{2}$ are both
$q$Gaussian random vectors but independent, a linear combination
of $H_{1}X_{1}+H_{2}X_{2}$ is not $q$Gaussian. 
\item \label{enu:heavy-tail-bounded-support-duality} The duality property:
if $X\sim q\text{Gaussian}\left(q,\Sigma\right)$ with $1<q<1+\frac{2}{n}$,
let the degree of freedom for $X$ be $m\coloneqq\frac{2}{q-1}-n$
and $\Lambda\coloneqq m\Sigma$, then 
\begin{align*}
\frac{X}{\sqrt{1-\left\langle X,\Lambda^{-1}X\right\rangle }} & \sim q\text{Gaussian}\left(\tilde{q},\frac{m}{m+4}\Sigma\right)\\
\text{with }\frac{1}{\tilde{q}^{-1}-1} & =\frac{1}{1-q^{-1}}-\frac{n}{2}-1,
\end{align*}
and $0<\tilde{q}<1$. 
\end{enumerate}
Property \ref{enu:maximizing-format} will be used in our Lemma \ref{lem:density-normalization}.
Property \ref{enu:linear-mapping-invariance} implies that any components of a $q$Gaussian random vector are also $q$Gaussian, while Property \ref{enu:linear-combination-not-closed}
implies that two independent $q$Gaussian vectors are not jointly
$q$Gaussian. 

By the equivalence of problems \eqref{eq:generalized-entropy-obj}
and \eqref{eq:generalized-entropy-moment-obj} discussed before, \eqref{eq:dual-maximizer}
can be rewritten as 
\begin{equation}
\phi\left(x;q,\Sigma\right)=\left(\alpha-\beta\left\langle x,\Sigma^{-1}x\right\rangle \right)_{+}^{\frac{1}{1-q}}\label{eq:parameter-to-be-found}
\end{equation}
for some constant (parameter) $\alpha,\beta$, $q\in\left(0,1+\frac{2}{n}\right)\setminus\left\{ 1\right\} $;
$x_{+}\coloneqq\max\ \left(0,x\right)$. As shown later in the proof
of Lemma \ref{lem:density-normalization}, the dimension-related
upper bound $1+\frac{2}{n}$ is due to the normalization constraint
\eqref{eq:normalization-constraint}. When $0<q<1$, the density represents
a distribution with bounded support; when $q>1$, the density is a
generalization of the bell curve distributions, and with $q\searrow1$ the
Gaussian distribution is recovered. A higher value of $q$ corresponds
to heavier tails in shape. The duality between the $q$Gaussian random vectors
with $0<q<1$ and $1<q<1+\frac{2}{n}$ was given by \cite{Vignat2004},
which we discussed in Property \ref{enu:heavy-tail-bounded-support-duality} in Section \ref{sec:Tsallis-Entropy-Maximizing}.
Distributions with bounded support correspond to $0<q<1$, and distributions
with heavy tails correspond to $1<q<1+\frac{2}{n}$. For the scope
of this paper, we will focus only on the heavy-tail distributions;
i.e., the case when $q>1$. \cite{Vignat2009} derived the $q$Gaussian
probability density function for $1<q<\frac{n+4}{n+2}$, when the
multivariate $q$Gaussian density becomes the scaled density of the multivariate
student's $t$ distribution. To incorporate the case of $q\in[1+\frac{2}{n+2},1+\frac{2}{n})$,
when the variance does not exist but the $q-$variance can be used
to capture the volatility/dispersion of the data, \cite{Vignat2007b}
also derived the resulting density; however, since a typo was found
in that paper, we re-derive the density in Lemma \ref{lem:density-normalization}.
The parameters presented in the density formula \eqref{eq:$q$Gaussian-heavy-tail-density}
are of particular interest to statisticians, as parameter inference
is the key to statistical analysis and prediction. The case of $q\in[1+\frac{2}{n+2},1+\frac{2}{n})$
will allow the resulting $q$Gaussian distribution to incorporate
the wider class of distributions without finite moments but finite
$q-$moments; such as the Cauchy distribution. Therefore, modeling using the $q$Gaussian distribution with $q$ allowed to take the value in $[1+\frac{2}{n+2},1+\frac{2}{n})$
will be more robust. 
\begin{lem}
\label{lem:density-normalization} When $q\in\left(1,1+\frac{2}{n}\right)$, the unique solution to \eqref{eq:generalized-entropy-obj}
is: 
\begin{equation}
p\left(x;q,\Sigma\right)=\frac{1}{\left|\pi\Sigma\right|^{\nicefrac{1}{2}}}\cdot\frac{\Gamma\left(\frac{1}{q-1}\right)}{\Gamma\left(\frac{1}{q-1}-\frac{n}{2}\right)}\cdot\left(\frac{2}{q-1}-n\right)^{-\frac{n}{2}}\cdot\left(1+\left(\frac{2}{q-1}-n\right)^{-1}\cdot\left\langle x,\Sigma^{-1}x\right\rangle \right)^{\frac{1}{1-q}}.\label{eq:$q$Gaussian-heavy-tail-density}
\end{equation}
\end{lem}

\begin{proof}
By \eqref{eq:dual-maximizer} and the equivalence of the problems \eqref{eq:generalized-entropy-obj}
and \eqref{eq:generalized-entropy-moment-obj}, let the solution to
\eqref{eq:generalized-entropy-obj} be denoted by 
\begin{equation}
p\left(x;q,\Sigma\right)=\frac{1}{Z}\left(\gamma+\left\langle x,\Sigma^{-1}x\right\rangle \right)^{\frac{1}{1-q}}\label{eq:target-maximizer}
\end{equation}
for some $Z,\gamma>0$. Feasibility for problem \ref{eq:generalized-entropy-obj} when $q\in\left(1,1+\frac{2}{n}\right)$
was given in \citep{Vignat2004}. Hence, the strictly concavity of
the objective function \ref{eq:generalized-entropy-obj} implies that the optimal solution is unique. The symmetry of $p\left(x;q,\Sigma\right)$
is implied by \eqref{eq:target-maximizer}; thus, we reformulate the
problem \eqref{eq:generalized-entropy-obj} as the following equivalent
problem: 
\begin{align}
\max\ _{p\in L^{q}\left(\mathbb{R}^{n}\right)} & -\int_{\mathbb{R}^{n}}\left(p\left(x;q,\Sigma\right)\right)^{q}dx\nonumber \\
\text{s.t. } & p\left(x;q,\Sigma\right)\geq0;\nonumber \\
 & \int_{\mathbb{R}_{>0}^{n}}p\left(x;q,\Sigma\right)dx=2^{-n};\nonumber \\
 & \frac{\int_{\mathbb{R}^{n}}x\cdot p^{q}\left(x;q,\Sigma\right)dx}{\int_{\mathbb{R}^{n}}p^{q}\left(x;q,\Sigma\right)dx}=0;\nonumber \\
 & \frac{\int_{\mathbb{R}^{n}}xx^{T}\cdot p^{q}\left(x;q,\Sigma\right)dx}{\int_{\mathbb{R}^{n}}p^{q}\left(x;q,\Sigma\right)dx}=\Sigma.\label{eq:second-q-moment-constraint}
\end{align}
Thus, 
\begin{align}
Z & =2^{n}\int_{\mathbb{R}_{>0}^{n}}\left(\gamma+\left\langle x,\Sigma^{-1}x\right\rangle \right)^{\frac{1}{1-q}}dx\nonumber \\
 & =2^{n}\left|\Sigma^{\nicefrac{1}{2}}\right|\int_{\mathbb{R}_{>0}^{n}}\left(\gamma+\left\langle x,x\right\rangle \right)^{\frac{1}{1-q}}dx\nonumber \\
 & =2^{n}\left|\Sigma\right|^{\nicefrac{1}{2}}\int_{0}^{\infty}r^{n-1}\left(\gamma+r^{2}\right)^{\frac{1}{1-q}}\left(\prod_{i=1}^{n-2}\int_{0}^{\frac{\pi}{2}}\sin^{n-1-i}(\theta_{i})d\theta_{i}\cdot\int_{0}^{\frac{\pi}{2}}1d\theta\right)dr\nonumber \\
 & =2^{n}\left|\Sigma\right|^{\nicefrac{1}{2}}\cdot\left(\prod_{i=1}^{n-2}\int_{0}^{\frac{\pi}{2}}\sin^{n-1-i}(\theta_{i})d\theta_{i}\cdot\int_{0}^{\frac{\pi}{2}}1d\theta\right)\cdot\int_{0}^{\infty}r^{n-1}\left(\gamma+r^{2}\right)^{\frac{1}{1-q}}dr\nonumber \\
 & =\pi2^{n-1}\left|\Sigma\right|^{\nicefrac{1}{2}}\cdot\left(\prod_{i=1}^{n-2}\frac{1}{2}\frac{\Gamma\left(\frac{n-i}{2}\right)\sqrt{\pi}}{\Gamma\left(\frac{n-i+1}{2}\right)}\right)\cdot\int_{0}^{\infty}r^{n-1}\left(\gamma+r^{2}\right)^{\frac{1}{1-q}}dr\nonumber \\
 & =2\pi^{\frac{n}{2}}\left|\Sigma\right|^{\nicefrac{1}{2}}\cdot\left(\Gamma\left(\frac{n}{2}\right)\right)^{-1}\cdot\int_{0}^{\infty}r^{n-1}\left(\gamma+r^{2}\right)^{\frac{1}{1-q}}dr\nonumber \\
 & =2\pi^{\frac{n}{2}}\left|\Sigma\right|^{\nicefrac{1}{2}}\cdot\left(\Gamma\left(\frac{n}{2}\right)\right)^{-1}\cdot\int_{0}^{\infty}\gamma^{\frac{n-1}{2}+\frac{1}{1-q}}\left(\frac{r}{\sqrt{\gamma}}\right)^{n-1}\left(1+\left(\frac{r}{\sqrt{\gamma}}\right)^{2}\right)^{\frac{1}{1-q}}dr\nonumber \\
 & =2\pi^{\frac{n}{2}}\left|\Sigma\right|^{\nicefrac{1}{2}}\cdot\left(\Gamma\left(\frac{n}{2}\right)\right)^{-1}\cdot\gamma^{\frac{n}{2}+\frac{1}{1-q}}\cdot\int_{0}^{\infty}\left(r^{\prime}\right)^{n-1}\left(1+\left(r^{\prime}\right)^{2}\right)^{\frac{1}{1-q}}dr^{\prime}\nonumber \\
 & =2\pi^{\frac{n}{2}}\left|\Sigma\right|^{\nicefrac{1}{2}}\cdot\left(\Gamma\left(\frac{n}{2}\right)\right)^{-1}\cdot\gamma^{\frac{n}{2}+\frac{1}{1-q}}\cdot\int_{0}^{\infty}\left(\left(r^{\prime}\right)^{1-n}\left(1+\left(r^{\prime}\right)^{2}\right)^{\frac{1}{q-1}}\right)^{-1}dr^{\prime}\nonumber \\
 & =2\pi^{\frac{n}{2}}\left|\Sigma\right|^{\nicefrac{1}{2}}\cdot\left(\Gamma\left(\frac{n}{2}\right)\right)^{-1}\cdot\gamma^{\frac{n}{2}+\frac{1}{1-q}}\cdot\frac{1}{2}B\left(\frac{1}{q-1}-\frac{n}{2},\frac{n}{2}\right)\label{eq:beta-function-step}\\
 & =\pi^{\frac{n}{2}}\left|\Sigma\right|^{\nicefrac{1}{2}}\cdot\left(\Gamma\left(\frac{n}{2}\right)\right)^{-1}\cdot\gamma^{\frac{n}{2}+\frac{1}{1-q}}\cdot\frac{\Gamma\left(\frac{1}{q-1}-\frac{n}{2}\right)\Gamma\left(\frac{n}{2}\right)}{\Gamma\left(\frac{1}{q-1}\right)}\nonumber \\
 & =\pi^{\frac{n}{2}}\left|\Sigma\right|^{\nicefrac{1}{2}}\frac{\Gamma\left(\frac{1}{q-1}-\frac{n}{2}\right)}{\Gamma\left(\frac{1}{q-1}\right)}\cdot\gamma^{\frac{n}{2}+\frac{1}{1-q}}.\nonumber 
\end{align}
In step \eqref{eq:beta-function-step}, we use the following formula
for Beta function \citep{Tsukada2005}: 
\begin{equation}
\int_{0}^{\infty}\left(x^{\alpha}\left(1+x^{\lambda}\right)^{\beta}\right)^{-1}dx=\frac{1}{\lambda}B\left(\beta-\frac{1-\alpha}{\lambda},\frac{1-\alpha}{\lambda}\right),\label{eq:beta-function-formula}
\end{equation}
where $\alpha<1,\ \lambda>0,\ \beta>0,\ \lambda\beta>1-\alpha$. Well--definedness of $Z$ and \eqref{eq:beta-function-step} implies that $\frac{1}{q-1}-\frac{n}{2}>0$; i.e.,
$q<1+\frac{2}{n}$, which is the reason for the upper bound for the
choice of $q$. 

\eqref{eq:second-q-moment-constraint} implies that 
\begin{equation}
\text{tr}\left(\Sigma^{-1}\frac{\int_{\mathbb{R}^{n}}xx^{T}\cdot p^{q}\left(x\right)dx}{\int_{\mathbb{R}^{n}}p^{q}\left(x\right)dx}\right)=\text{tr}\left(\Sigma^{-1}\Sigma\right)=n.\label{eq:second-q-moment-trace-trick-1}
\end{equation}
Hence, since $p\left(x\right)$ is symmetric, 
\begin{align}
& \text{tr}\left(\Sigma^{-1}\frac{\int_{\mathbb{R}^{n}}xx^{T}\cdot p^{q}\left(x;q,\Sigma\right)dx}{\int_{\mathbb{R}^{n}}p^{q}\left(x;q,\Sigma\right)dx}\right) \\ & =\text{tr}\left(\Sigma^{-1}\frac{\int_{\mathbb{R}_{>0}^{n}}xx^{T}\cdot p^{q}\left(x;q,\Sigma\right)dx}{\int_{\mathbb{R}_{>0}^{n}}p^{q}\left(x;q,\Sigma\right)dx}\right)\nonumber \\
 & =\frac{\text{tr}\left(\int_{\mathbb{R}_{>0}^{n}}\Sigma^{-1}xx^{T}\cdot p^{q}\left(x;q,\Sigma\right)dx\right)}{\int_{\mathbb{R}_{>0}^{n}}p^{q}\left(x;q,\Sigma\right)dx}\nonumber \\
 & =\frac{\int_{\mathbb{R}_{>0}^{n}}\text{tr}\left(\Sigma^{-1}xx^{T}\cdot p^{q}\left(x;q,\Sigma\right)\right)dx}{\int_{\mathbb{R}_{>0}^{n}}p^{q}\left(x;q,\Sigma\right)dx}\nonumber \\
 & =\frac{\int_{\mathbb{R}_{>0}^{n}}\text{tr}\left(x^{T}\Sigma^{-1}x\right)\cdot p^{q}\left(x;q,\Sigma\right)dx}{\int_{\mathbb{R}_{>0}^{n}}p^{q}\left(x;q,\Sigma\right)dx}\nonumber \\
 & =\frac{\int_{\mathbb{R}_{>0}^{n}}\left\langle x,\Sigma^{-1}x\right\rangle \cdot\left(\frac{1}{Z}\left(\gamma+\left\langle x,\Sigma^{-1}x\right\rangle \right)^{\frac{1}{1-q}}\right)^{q}dx}{\int_{\mathbb{R}_{>0}^{n}}\left(\frac{1}{Z}\left(\gamma+\left\langle x,\Sigma^{-1}x\right\rangle \right)^{\frac{1}{1-q}}\right)^{q}dx}\nonumber \\
 & =\frac{\int_{\mathbb{R}_{>0}^{n}}\left\langle x,\Sigma^{-1}x\right\rangle \cdot\left(\gamma+\left\langle x,\Sigma^{-1}x\right\rangle \right)^{\frac{q}{1-q}}dx}{\int_{\mathbb{R}_{>0}^{n}}\left(\gamma+\left\langle x,\Sigma^{-1}x\right\rangle \right)^{\frac{q}{1-q}}dx}\nonumber \\
 & =\frac{\int_{\mathbb{R}_{>0}^{n}}\left|\Sigma\right|^{\nicefrac{1}{2}}\left\langle x,x\right\rangle \cdot\left(\gamma+\left\langle x,x\right\rangle \right)^{\frac{q}{1-q}}dx}{\int_{\mathbb{R}_{>0}^{n}}\left|\Sigma\right|^{\nicefrac{1}{2}}\left(\gamma+\left\langle x,x\right\rangle \right)^{\frac{q}{1-q}}dx}\nonumber \\
 & =\frac{\int_{\mathbb{R}_{>0}^{n}}\left\langle x,x\right\rangle \cdot\left(\gamma+\left\langle x,x\right\rangle \right)^{\frac{q}{1-q}}dx}{\int_{\mathbb{R}_{>0}^{n}}\left(\gamma+\left\langle x,x\right\rangle \right)^{\frac{q}{1-q}}dx}\nonumber \\
 & =\frac{\int_{0}^{\infty}r^{n-1}\cdot r^{2}\cdot\left(\gamma+r^{2}\right)^{\frac{q}{1-q}}\cdot\left(\prod_{i=1}^{n-2}\int_{0}^{\frac{\pi}{2}}\sin^{n-1-i}(\theta_{i})d\theta_{i}\cdot\int_{0}^{\frac{\pi}{2}}1d\theta\right)dr}{\int_{0}^{\infty}r^{n-1}\cdot\left(\gamma+r^{2}\right)^{\frac{q}{1-q}}\cdot\left(\prod_{i=1}^{n-2}\int_{0}^{\frac{\pi}{2}}\sin^{n-1-i}(\theta_{i})d\theta_{i}\cdot\int_{0}^{\frac{\pi}{2}}1d\theta\right)dr}\nonumber \\
 & =\frac{\int_{0}^{\infty}r^{n+1}\cdot\left(\gamma+r^{2}\right)^{\frac{q}{1-q}}dr}{\int_{0}^{\infty}r^{n-1}\cdot\left(\gamma+r^{2}\right)^{\frac{q}{1-q}}dr}\nonumber \\
 & =\frac{\gamma\int_{0}^{\infty}\left(\left(r^{\prime}\right)^{-n-1}\cdot\left(1+\left(r^{\prime}\right)^{2}\right)^{\frac{q}{q-1}}\right)^{-1}dr^{\prime}}{\int_{0}^{\infty}\left(\left(r^{\prime}\right)^{1-n}\cdot\left(1+\left(r^{\prime}\right)^{2}\right)^{\frac{q}{q-1}}\right)^{-1}dr^{\prime}}\nonumber \\
 & =\frac{\gamma B\left(\frac{q}{q-1}-\frac{n+2}{2},\frac{n+2}{2}\right)}{B\left(\frac{q}{q-1}-\frac{n}{2},\frac{n}{2}\right)}\label{eq:beta-function-usage-2}\\
 & =\gamma\cdot\frac{\Gamma\left(\frac{q}{q-1}-\frac{n}{2}-1\right)\Gamma\left(\frac{n+2}{2}\right)}{\Gamma\left(\frac{q}{q-1}\right)}/\frac{\Gamma\left(\frac{q}{q-1}-\frac{n}{2}\right)\Gamma\left(\frac{n}{2}\right)}{\Gamma\left(\frac{q}{q-1}\right)}\nonumber \\
 & =\gamma\cdot\frac{n}{2}\cdot\left(\frac{q}{q-1}-\frac{n}{2}-1\right)^{-1}.\label{eq:second-q-moment-trace-trick-2}
\end{align}
In step \eqref{eq:beta-function-usage-2}, we used \eqref{eq:beta-function-formula}.
Combining \eqref{eq:second-q-moment-trace-trick-1} and \eqref{eq:second-q-moment-trace-trick-2},
we have 
\begin{equation}
\gamma\cdot\frac{n}{2}\cdot\left(\frac{q}{q-1}-\frac{n}{2}-1\right)^{-1}=n,
\end{equation}
which gives that 
\begin{equation}
\gamma=\frac{2q}{q-1}-n-2=\frac{2}{q-1}-n.
\end{equation}
Thus, the probability density function that maximizes problem \eqref{eq:generalized-entropy-obj}
is: 
\begin{align*}
p\left(x;q,\Sigma\right) & =\left(\pi^{\frac{n}{2}}\left|\Sigma\right|^{\nicefrac{1}{2}}\frac{\Gamma\left(\frac{1}{q-1}-\frac{n}{2}\right)}{\Gamma\left(\frac{1}{q-1}\right)}\cdot\left(\frac{2}{q-1}-n\right)^{\frac{n}{2}+\frac{1}{1-q}}\right)^{-1}\left(\left(\frac{2}{q-1}-n\right)+\left\langle x,\Sigma^{-1}x\right\rangle \right)^{\frac{1}{1-q}}\\
 & =\left(\pi^{\frac{n}{2}}\left|\Sigma\right|^{\nicefrac{1}{2}}\frac{\Gamma\left(\frac{1}{q-1}-\frac{n}{2}\right)}{\Gamma\left(\frac{1}{q-1}\right)}\cdot\left(\frac{2}{q-1}-n\right)^{\frac{n}{2}}\right)^{-1}\left(1+\left(\frac{2}{q-1}-n\right)^{-1}\left\langle x,\Sigma^{-1}x\right\rangle \right)^{\frac{1}{1-q}}\\
 & =\frac{1}{\left|\pi\Sigma\right|^{\nicefrac{1}{2}}}\cdot\frac{\Gamma\left(\frac{1}{q-1}\right)}{\Gamma\left(\frac{1}{q-1}-\frac{n}{2}\right)}\cdot\left(\frac{2}{q-1}-n\right)^{-\frac{n}{2}}\cdot\left(1+\left(\frac{2}{q-1}-n\right)^{-1}\cdot\left\langle x,\Sigma^{-1}x\right\rangle \right)^{\frac{1}{1-q}}.
\end{align*}
\end{proof}
When $q\geq1+\frac{2}{n}$, the solution to problem \eqref{eq:generalized-entropy-obj} does not exist, due to property \ref{enu:maximizing-format} and discussions in the proof. In Lemma \ref{lem:density-normalization}, the presented density \eqref{eq:$q$Gaussian-heavy-tail-density}
outlines a formula for multivariate bell-curve distributions dependent
on the value of $q$. As $q$ shifts from values approaching $1$
from above to values approaching $1+\frac{2}{n}$ form below, the
resulting density transitions from Gaussian through a scaled version
of the multivariate $t-$distribution to Cauchy and beyond. This density
explicitly details all parameters, enabling the application of the
maximum likelihood principle and facilitating the use of maximum likelihood
estimation in modeling correlated data performed in Section \ref{sec:optimizing-penalized-mle}. 

In the context of \eqref{eq:$q$Gaussian-heavy-tail-density}, the \emph{characterization
matrix} \citep{Costa2003}, denoted by $\Sigma$, can undergo modifications
to include the degree of freedom parameter $m\coloneqq\frac{2}{q-1}-n$
\citep{Vignat2005}; specifically,
\begin{align}
p\left(x;q,\Lambda\right) & =\frac{1}{\left|\pi\Lambda\right|^{\nicefrac{1}{2}}}\cdot\frac{\Gamma\left(\frac{m}{2}+\frac{n}{2}\right)}{\Gamma\left(\frac{m}{2}\right)}\cdot\left(1+\left\langle x,\Lambda^{-1}x\right\rangle \right)^{\frac{1}{1-q}},\label{eq:$q$Gaussian-simple-form}\\
\text{where }\Lambda & \coloneqq m\Sigma.\nonumber \\
\text{and }m & \coloneqq\frac{2}{q-1}-n\nonumber 
\end{align}
Below are a few useful remarks related to the $q$Gaussian distribution
and other bell-curve distributions. 
\begin{rem}
\label{rem:incorporate-location} To incorporate the location parameter
$\mu$, \eqref{eq:$q$Gaussian-heavy-tail-density} and \eqref{eq:$q$Gaussian-simple-form}
become 
\begin{align*}
p\left(x;\mu,q,\Sigma\right) & =\frac{1}{\left|\pi\Sigma\right|^{\nicefrac{1}{2}}}\cdot\frac{\Gamma\left(\frac{1}{q-1}\right)}{\Gamma\left(\frac{1}{q-1}-\frac{n}{2}\right)}\cdot\left(\frac{2}{q-1}-n\right)^{-\frac{n}{2}}\\
 & \qquad\cdot\left(1+\left(\frac{2}{q-1}-n\right)^{-1}\cdot\left\langle x-\mu,\Sigma^{-1}\left(x-\mu\right)\right\rangle \right)^{\frac{1}{1-q}};\\
p\left(x;\mu,q,\Lambda\right) & =\frac{1}{\left|\pi\Lambda\right|^{\nicefrac{1}{2}}}\cdot\frac{\Gamma\left(\frac{m}{2}+\frac{n}{2}\right)}{\Gamma\left(\frac{m}{2}\right)}\cdot\left(1+\left\langle x-\mu,\Lambda^{-1}\left(x-\mu\right)\right\rangle \right)^{\frac{1}{1-q}}.
\end{align*}
\end{rem}

\begin{rem}
For random vector $X\sim q\text{Gaussian}\left(q,\Sigma\right)$,
its $q-$covariance is 
\begin{equation}
\mathbb{E}_{q}\left[XX^{T}\right]=\left(\frac{1}{q-1}-\frac{n}{2}\right)^{\frac{1-q}{2}}\left|\pi\Sigma\right|^{\frac{1-q}{2}}\cdot\frac{\Gamma\left(\frac{q}{q-1}-\frac{n}{2}\right)/\left(\Gamma\left(\frac{1}{q-1}-\frac{n}{2}\right)\right)^{q}}{\Gamma\left(\frac{q}{q-1}\right)/\left(\Gamma\left(\frac{1}{q-1}\right)\right)^{q}}\cdot\Sigma.
\end{equation}
\end{rem}

\begin{proof}
We note that 
\begin{align*}
\int_{\mathbb{R}^{n}}p^{q}\left(x;q,\Sigma\right)dx & =\int_{\mathbb{R}^{n}}\left(\frac{1}{\left|\pi\Lambda\right|^{\nicefrac{1}{2}}}\cdot\frac{\Gamma\left(\frac{1}{q-1}\right)}{\Gamma\left(\frac{1}{q-1}-\frac{n}{2}\right)}\cdot\left(1+\left\langle x,\Lambda^{-1}x\right\rangle \right)^{\frac{1}{1-q}}\right)^{q}dx\\
 & =\left(\frac{1}{\left|\pi\Lambda\right|^{\nicefrac{1}{2}}}\cdot\frac{\Gamma\left(\frac{1}{q-1}\right)}{\Gamma\left(\frac{1}{q-1}-\frac{n}{2}\right)}\right)^{q}\cdot\int_{\mathbb{R}^{n}}\left(\left(1+\left\langle x,\Lambda^{-1}x\right\rangle \right)^{\frac{1}{1-q}}\right)^{q}dx\\
 & =\left(\frac{1}{\left|\pi\Lambda\right|^{\nicefrac{1}{2}}}\cdot\frac{\Gamma\left(\frac{1}{q-1}\right)}{\Gamma\left(\frac{1}{q-1}-\frac{n}{2}\right)}\right)^{q}\cdot\int_{\mathbb{R}^{n}}\left(1+\left\langle x,\Lambda^{-1}x\right\rangle \right)^{\frac{q}{1-q}}dx\\
 & =\left(\frac{1}{\left|\pi\Lambda\right|^{\nicefrac{1}{2}}}\cdot\frac{\Gamma\left(\frac{1}{q-1}\right)}{\Gamma\left(\frac{1}{q-1}-\frac{n}{2}\right)}\right)^{q}\left|\Lambda\right|^{\nicefrac{1}{2}}\cdot\int_{\mathbb{R}^{n}}\left(1+\left\langle x,x\right\rangle \right)^{\frac{q}{1-q}}dx\\
 & =2^{n}\left(\frac{1}{\left|\pi\Lambda\right|^{\nicefrac{1}{2}}}\cdot\frac{\Gamma\left(\frac{1}{q-1}\right)}{\Gamma\left(\frac{1}{q-1}-\frac{n}{2}\right)}\right)^{q}\left|\Lambda\right|^{\nicefrac{1}{2}}\cdot\int_{\mathbb{R}_{>0}^{n}}\left(1+\left\langle x,x\right\rangle \right)^{\frac{q}{1-q}}dx\\
 & =2^{n}\left(\frac{1}{\left|\pi\Lambda\right|^{\nicefrac{1}{2}}}\cdot\frac{\Gamma\left(\frac{1}{q-1}\right)}{\Gamma\left(\frac{1}{q-1}-\frac{n}{2}\right)}\right)^{q}\left|\Lambda\right|^{\nicefrac{1}{2}}\cdot\int_{0}^{\infty}r^{n-1}\left(1+r^{2}\right)^{\frac{q}{1-q}}\\
 & \qquad\cdot\left(\prod_{i=1}^{n-2}\int_{0}^{\frac{\pi}{2}}\sin^{n-1-i}(\theta_{i})d\theta_{i}\cdot\int_{0}^{\frac{\pi}{2}}1d\theta\right)dr\\
 & =\pi2^{n-1}\left(\frac{1}{\left|\pi\Lambda\right|^{\nicefrac{1}{2}}}\cdot\frac{\Gamma\left(\frac{1}{q-1}\right)}{\Gamma\left(\frac{1}{q-1}-\frac{n}{2}\right)}\right)^{q}\left|\Lambda\right|^{\nicefrac{1}{2}}\cdot\left(\prod_{i=1}^{n-2}\frac{1}{2}\frac{\Gamma\left(\frac{n-i}{2}\right)\sqrt{\pi}}{\Gamma\left(\frac{n-i+1}{2}\right)}\right)\\
 & \qquad\cdot\int_{0}^{\infty}r^{n-1}\left(1+r^{2}\right)^{\frac{q}{1-q}}dr\\
 & =2\pi^{\frac{n}{2}}\left(\frac{1}{\left|\pi\Lambda\right|^{\nicefrac{1}{2}}}\cdot\frac{\Gamma\left(\frac{1}{q-1}\right)}{\Gamma\left(\frac{1}{q-1}-\frac{n}{2}\right)}\right)^{q}\left|\Lambda\right|^{\nicefrac{1}{2}}\cdot\left(\Gamma\left(\frac{n}{2}\right)\right)^{-1}\cdot\int_{0}^{\infty}\left(r^{1-n}\left(1+r^{2}\right)^{\frac{q}{q-1}}\right)^{-1}dr\\
 & =2\pi^{\frac{n}{2}}\left(\frac{1}{\left|\pi\Lambda\right|^{\nicefrac{1}{2}}}\cdot\frac{\Gamma\left(\frac{1}{q-1}\right)}{\Gamma\left(\frac{1}{q-1}-\frac{n}{2}\right)}\right)^{q}\left|\Lambda\right|^{\nicefrac{1}{2}}\cdot\left(\Gamma\left(\frac{n}{2}\right)\right)^{-1}\cdot\frac{1}{2}B\left(\frac{q}{q-1}-\frac{n}{2},\frac{n}{2}\right)\\
 & =\pi^{\frac{n}{2}}\left(\frac{1}{\left|\pi\Lambda\right|^{\nicefrac{1}{2}}}\cdot\frac{\Gamma\left(\frac{1}{q-1}\right)}{\Gamma\left(\frac{1}{q-1}-\frac{n}{2}\right)}\right)^{q}\left|\Lambda\right|^{\nicefrac{1}{2}}\cdot\left(\Gamma\left(\frac{n}{2}\right)\right)^{-1}\cdot\frac{\Gamma\left(\frac{q}{q-1}-\frac{n}{2}\right)\Gamma\left(\frac{n}{2}\right)}{\Gamma\left(\frac{q}{q-1}\right)}\\
 & =\left(\frac{1}{\left|\pi\Lambda\right|^{\nicefrac{1}{2}}}\cdot\frac{\Gamma\left(\frac{1}{q-1}\right)}{\Gamma\left(\frac{1}{q-1}-\frac{n}{2}\right)}\right)^{q}\left|\pi\Lambda\right|^{\nicefrac{1}{2}}\cdot\frac{\Gamma\left(\frac{q}{q-1}-\frac{n}{2}\right)}{\Gamma\left(\frac{q}{q-1}\right)}\\
 & =\left|\pi\Lambda\right|^{\frac{1-q}{2}}\cdot\frac{\Gamma\left(\frac{q}{q-1}-\frac{n}{2}\right)/\left(\Gamma\left(\frac{1}{q-1}-\frac{n}{2}\right)\right)^{q}}{\Gamma\left(\frac{q}{q-1}\right)/\left(\Gamma\left(\frac{1}{q-1}\right)\right)^{q}}.
\end{align*}
From \eqref{eq:q-second-moment-constraint}, we then have the following
expression for the $q-$variance-covariance matrix 
\begin{align}
\mathbb{E}_{q}\left[XX^{T}\right] & =\int_{\mathbb{R}^{n}}xx^{T}\cdot p^{q}\left(x;q,\Sigma\right)dx\nonumber \\
 & =\int_{\mathbb{R}^{n}}p^{q}\left(x;q,\Sigma\right)dx\cdot\Sigma\nonumber \\
 & =\left|\pi\Lambda\right|^{\frac{1-q}{2}}\cdot\frac{\Gamma\left(\frac{q}{q-1}-\frac{n}{2}\right)/\left(\Gamma\left(\frac{1}{q-1}-\frac{n}{2}\right)\right)^{q}}{\Gamma\left(\frac{q}{q-1}\right)/\left(\Gamma\left(\frac{1}{q-1}\right)\right)^{q}}\cdot\Sigma\nonumber \\
 & =m^{\frac{1-q}{2}}\left|\pi\Sigma\right|^{\frac{1-q}{2}}\cdot\frac{\Gamma\left(\frac{q}{q-1}-\frac{n}{2}\right)/\left(\Gamma\left(\frac{1}{q-1}-\frac{n}{2}\right)\right)^{q}}{\Gamma\left(\frac{q}{q-1}\right)/\left(\Gamma\left(\frac{1}{q-1}\right)\right)^{q}}\cdot\Sigma\nonumber \\
 & =\left(\frac{1}{q-1}-\frac{n}{2}\right)^{\frac{1-q}{2}}\left|\pi\Sigma\right|^{\frac{1-q}{2}}\Sigma\cdot\frac{\Gamma\left(\frac{q}{q-1}-\frac{n}{2}\right)/\left(\Gamma\left(\frac{1}{q-1}-\frac{n}{2}\right)\right)^{q}}{\Gamma\left(\frac{q}{q-1}\right)/\left(\Gamma\left(\frac{1}{q-1}\right)\right)^{q}}.\label{eq:q-variance-covariance}
\end{align}
\end{proof}
\begin{rem}
Multivariate $t$ distributions with degree of freedom of $\frac{2}{q-1}-n$
and the scale matrix of $\Sigma$ are $q$Gaussian with shape parameter
$q$ and scale matrix $\Sigma$. 
\end{rem}

\begin{rem}
The multivariate Cauchy distributions \citep{Lee2014a} in $\mathbb{R}^{n}$
with scale matrix $\frac{1}{2}\Sigma$ are $q$Gaussian with shape
parameter $q=1+\frac{2}{n+1}$ and scale matrix $\Sigma$. 
\end{rem}

\begin{rem}
For a random vector $X\sim q\text{Gaussian}\left(q,\Sigma\right)$,
the variance-covariance matrix exists if and only if $q<1+\frac{2}{n+2}$;
following the same procedure to derive the variance-covariance matrix
for multivariate $t$ distribution yields that,\emph{ if existing,}
\begin{equation}
\mathbb{E}\left[XX^{T}\right]=\frac{m}{m-2}\Sigma,\label{eq:variance-covariance}
\end{equation}
where $m=\frac{2}{q-1}-n$. 
\end{rem}

The remarks above on $q$Gaussian distributions reveal their flexibility
in incorporating a location parameter, $\mu$, and adapting to multivariate
contexts through detailed formulas. Remarkably, these distributions
bridge with the class of multivariate bell curve distributions including
Gaussian, scaled $t$, and Cauchy distributions under certain conditions
on the shape parameter $q$. The elaboration of $q-$correlation and the $q-$variable-covariance matrix underscores the capability of
these distributions to model and understand the intricacies of correlated
data effectively. 

\section{\label{sec:optimization} Proximal Conjugate Gradient Algorithm }

As delineated in Section \ref{subsec:optimization-motivation}, our
optimization scenario is predominantly quadratic in nature; therefore,
the conjugate gradient approach has potential for fast convergence
and numerical stability. This insight forms the basis for our introduction
of a proximal conjugate gradient algorithm framework, tailored to
navigate the complexities introduced by the nonconvex penalized $q$Gaussian
likelihood function for sparse statistical learning. To lay the groundwork
for this discussion, we begin with a overview of relevant concepts
in variational and nonsmooth analysis, presented in Section \ref{subsec:variational-analysis}. The results presented in Section \ref{subsec:variational-analysis} can be found in recent textbooks on variational and nonsmooth analysis, such as \citep{Rockafellar2010, Clarke1990, Morduchovic2018, Mordukhovich2006, Mordukhovich2006a, Bauschke2011}. 

\subsection{\label{subsec:variational-analysis} A Review on Variational and
Nonsmooth Analysis }

Let $\mathcal{C}^{k,\alpha_{H}}$ with $k\in\mathbb{N}_{\geq0}$ and $\alpha_{H}\in\left[0,1\right]$ denote the function space such that $\forall F\in\mathcal{C}^{k,\alpha_{H}}$, $F$ is $k$th continuously differentiable, and $D^{k}F$ is globally H\"{o}lder continuous with exponent $\alpha_{H}$; clearly, when $\alpha_{H}=1$, $D^{k}F$ is globally Lipschitz continuous. In this subsection, we will state the results from variational and
nonsmooth analysis related to the following optimization problem:
\begin{equation}
\min\ _{x\in\mathbb{R}^{p+1}}f\left(x\right)\coloneqq g\left(x\right)+h\left(x\right),\label{eq:general-opt-problem}
\end{equation}
where $f\in\mathcal{C}^{0,0}\left(\mathbb{R}^{p+1},\mathbb{R}\right)$
is a locally-Lipschitz proper function, $g\in\mathcal{C}^{1,1}\left(\mathbb{R}^{p+1},\mathbb{R}\right)$
is globally $L_{\nabla g}-$smooth and possibly nonconvex, and $h\in\mathcal{C}^{0,0}\left(\mathbb{R}^{p+1},\mathbb{R}\right)$
is a convex locally-Lipschitz function, possibly nonsmooth.   The globally Lipschitz property
of $\nabla g$ can be alternatively addressed by carrying out the optimization
over a compact set. In such scenarios, given that $\nabla g$ is locally
Lipschitz, it inherently becomes globally Lipschitz when restricted
to a compact set. 

Results from convex analysis suggest that $g,h$ are Clarke regular;
thus, $f$ is Clarke regular. The Clarke's directional derivative,
defined by 
\begin{align*}
f^{\circ}\left(x;d\right) & \coloneqq\lim_{y\rightarrow x}\sup_{t\searrow0}\frac{f\left(y+td\right)-f\left(y\right)}{t}\\
 & =\inf_{\delta>0}\sup_{\left\Vert y-x\right\Vert \leq\delta,0<t<\delta}\frac{f\left(y+td\right)-f\left(y\right)}{t},
\end{align*}
exists for all $x\in\mathbb{R}^{p+1}$ since $f$ is Clarke regular.
The Clarke subdifferential, denoted by $\partial_{\circ}$, is a set-valued
mapping defined by 
\begin{equation}
\partial_{\circ}f\left(x\right)\coloneqq\left\{ \phi\in\mathbb{R}^{p+1}\vert\forall d\in\mathbb{R}^{p+1},\ \left\langle \phi,d\right\rangle \leq f^{\circ}\left(x;d\right)\right\} .
\end{equation}
Since $f$ is a locally Lipschitz function, $\forall x\in\mathbb{R}^{p+1},\ \partial_{\circ}f\left(x\right)\neq\emptyset$.
Fundamental convex analysis results show that $\forall x\in\mathbb{R}^{p+1},\ \partial_{\circ}f\left(x\right)$
is compact, convex, and upper-semicontinuous. $\forall x,d\in\mathbb{R}^{p+1}$,
and we also have 
\begin{equation}
f^{\circ}\left(x;d\right)=\max\ _{u\in\partial_{\circ}f\left(x\right)}\left\langle u,\frac{d}{\left\Vert d\right\Vert }\right\rangle .\label{eq:max-inner-product-def}
\end{equation}
Furthermore, \eqref{eq:max-inner-product-def} is upper-semicontinuous
with respect to $x$. Simple convex geometry results conclude that 
\begin{equation}
    \left\{ \left(v,-1\right)\vert v\in\partial_{\circ}f\left(x\right)\right\} =N_{\text{epi }f}\left(x,f\left(x\right)\right),
\end{equation}
where $N_{\text{epi }f}\left(x,f\left(x\right)\right)$ denotes the \emph{normal cone} to $\text{epi }f$ at the point $\left(x,f\left(x\right)\right)$. 

Since $g$ is smooth, $\partial_{\circ}g\left(x\right)=\left\{ \nabla g\left(x\right)\right\} $
is a singleton. Then $\partial_{\circ}f\left(x\right)=\partial_{\circ}g\left(x\right)+\partial_{\circ}h\left(x\right)$,
and 
\begin{align}
f^{\circ}\left(x;d\right) & =\max\ _{u\in\partial_{\circ}f\left(x\right)}\left\langle u,\frac{d}{\left\Vert d\right\Vert }\right\rangle \nonumber \\
 & =\max\ _{u\in\left(\nabla g\left(x\right)+\partial_{\circ}h\left(x\right)\right)}\left\langle u,\frac{d}{\left\Vert d\right\Vert }\right\rangle \nonumber \\
 & =\left\langle \nabla g\left(x\right),\frac{d}{\left\Vert d\right\Vert }\right\rangle +\max\ _{v\in\partial_{\circ}h\left(x\right)}\left\langle v,\frac{d}{\left\Vert d\right\Vert }\right\rangle \label{eq:clarke-subdifferential-decomposition}\\
 & =g^{\circ}\left(x;d\right)+h^{\circ}\left(x;d\right) \notag.
\end{align}

Let 
\begin{equation}
M_{\rho}t\left(x\right)\coloneqq\left(t\square\left(\frac{1}{2\rho}\left\Vert \cdot\right\Vert ^{2}\right)\right)\left(x\right)=\inf_{y\in\mathbb{R}^{p+1}}t\left(y\right)+\frac{1}{2\rho}\left\Vert y-x\right\Vert ^{2} \label{eq:moreau-envelope-defn}
\end{equation}
denote the Moreau envelope operator parameterized by $\rho\in\mathbb{R}_{>0}$
applied on an arbitrary proper, lower semi-continuous, locally Lipschitz function $t\in\mathcal{C}^{0,0}\left(\mathbb{R}^{p+1},\mathbb{R}\right)$,
where ``$\square$'' denotes the infimal convolution operator. We
have that the Moreau envelope is a smoothing operator, specifically,
\begin{equation}
\text{epi }t+\text{epi }\frac{1}{2\rho}\left\Vert \cdot\right\Vert ^{2}\subseteq \text{epi }M_{\rho}t,\label{eq:epi-decomposition}
\end{equation}
where ``$\text{epi}$'' denotes the epigraph. Clearly, $M_{\rho}t\left(x\right)\leq t\left(x\right)$, since $\left(0,0\right)\in\text{epi }\frac{1}{2\rho}\left\Vert \cdot\right\Vert ^{2}$ implies that $\text{epi }t=\text{epi }t+\left(0,0\right)\subseteq\text{epi }t+\text{epi }\frac{1}{2\rho}\left\Vert \cdot\right\Vert ^{2}\subseteq\text{epi }M_{\rho}t$. When $t$ is convex, \eqref{eq:epi-decomposition} takes the equal sign; i.e., the infimal convolution becomes the exact infimal convolution. 

Consider the affine function 
\begin{equation}
A\left(x\right)\coloneqq\left\langle a,x\right\rangle +b,
\end{equation}
simple algebra shows that the Moreau envelope applied on $A$ is 
\begin{equation}
M_{\rho}A\left(x\right)=\left\langle a,x\right\rangle +b+\frac{\rho}{2}\left\Vert a\right\Vert ^{2}=A\left(x\right)+\frac{\rho}{2}\left\Vert a\right\Vert ^{2}
\end{equation}
for some $a,b\in\mathbb{R}^{p+1}$. Moreover, the following affine
addition property is often used in proximal algorithms, mainly due
to the fact that the epigraph of an affine function is a half-space that
the Moreau envelope applied on: 
\begin{align}
M_{\rho}\left(t+A\right)\left(x\right) & =M_{\rho}t\left(x-\rho a\right)+\left\langle a,x\right\rangle +b-\frac{\rho}{2}\left\Vert a\right\Vert ^{2}\label{eq:Moreau-envelope-affine-addition}
\end{align}
Let 
\begin{equation}
\text{prox}_{\rho t}\left(x\right)\coloneqq\arg M_{\rho}t\left(x\right)=\underset{y\in\mathbb{R}^{p+1}}{\arg\min\ }t\left(y\right)+\frac{1}{2\rho}\left\Vert y-x\right\Vert ^{2}
\end{equation}
denote the proximal operator, a set-valued mapping; we have 
\begin{equation}
\text{prox}_{\rho t}=\left(I+\rho\partial_{\circ}t\right)^{-1}\label{eq:proximal-resolvency}
\end{equation}
is the resolvent of the Clarke's subdifferential operator $\rho\partial_{\circ}t$. 

For nonsmooth problems, proximal methods are often used. Fundamental
convex analysis results show that: 
\begin{enumerate}
\item the Moreau envelope $M_{\rho}t\left(x\right)$ is twice differentiable;
thus, its gradient $\nabla M_{\rho}t\left(x\right)$ is well-defined. 
\item If $t$ is convex, $\text{prox}_{\rho t}\left(x\right)$ is a singleton.
\emph{For the sake of parsimony, with a slight abuse of notation,
we use $\text{prox}_{\rho t}$ to represent a function in this case.}
It follows that both $\text{prox}_{\rho t}$ and $\nabla M_{\rho}t$
are firmly non-expansive, and that 
\begin{equation}
\nabla M_{\rho}t\left(x\right)=\rho^{-1}\left(x-\text{prox}_{\rho t}\left(x\right)\right).
\end{equation}
\end{enumerate}
The results from variational and nonsmooth analysis in this subsection
have laid the foundation for proving the properties discussed
in Section \ref{subsec:PCG-framework}. 

\subsection{\label{subsec:PCG-framework} Proximal Conjugate Gradient Framework }

Proximal methods are powerful optimization techniques and are particularly
adept at handling problems characterized by sparsity, which usually
leads to an optimization problem that is nonsmooth \citep{Nikolova2000}.
Proximal algorithm tends to outperform other methods by far for nonsmooth
problems \citep{Yu2017,Li2016}. On another ground, Krylov subspace
methods represent a cornerstone of numerical analysis, providing a
powerful framework for solving large-scale optimization problems efficiently
\citep{Saad2003}. Krylov subspace methods exhibit a remarkable property
of convergence acceleration and vastly improved numerical stability,
making them indispensable tools in the numerical analyst's toolkit. 

Having reviewed the related results from variational and non-smooth analysis in Section \ref{subsec:variational-analysis}, we are
ready to introduce our main optimization framework to combine proximal
methods and conjugate gradient together. The essence of proximal algorithms
lies upon the Moreau envelope's smoothing on the objective function. Indeed,
proximal methods minimize $M_{\rho}f$ instead of $f$, thus avoiding nonsmoothness since $M_{\rho}f$ is a smooth function. In this view,
proximal algorithms are, in fact, minimizing the Moreau envelope of
the objective function. Thus, a wide class of numerical optimization
algorithms can easily have their proximal version. Among those, conjugate
gradients, a type of Krylov subspace method, are the state-of-the-art
methods in smooth optimization due to their computational and memory
efficiency, scalability, and numerical stability. 

Prior to introducing our proximal conjugate gradient update framework,
we will first show the equivalency of the optimization problem to
minimize \eqref{eq:general-opt-problem} and its the Moreau envelope.
In nonconvex optimization, the main task for numerical optimization
is to find a Clarke stationary point of the objective function, for
which we show in Theorem \ref{lem:equivalent-clarke-stationary-point}
that the set of Clarke stationary point of $f$ is identical to that
of $M_{\rho}f$ for $\rho\in\left(0,L_{\nabla g}^{-1}\right)$. 
\begin{lem}
\label{lem:equivalent-clarke-stationary-point} $\forall\bar{x}\in\mathbb{R}^{p+1},\rho\in\left(0,L_{\nabla g}^{-1}\right)$,
\begin{equation}
0\in\partial_{\circ}f\left(\bar{x}\right)\Leftrightarrow\nabla M_{\rho}f\left(\bar{x}\right)=0,\label{eq:equiv-equilibrium}
\end{equation}
\end{lem}

\begin{proof}
Consider arbitrary $x\in\mathbb{R}^{p+1},\ \rho\in\left(0,L_{\nabla g}^{-1}\right)$.
As discussed previously, the gradient of the Moreau envelope $\nabla M_{\rho}f\left(x\right)=\rho^{-1}\left(x-\text{prox}_{\rho f}\left(x\right)\right)$
implies that 
\begin{equation}
\text{prox}_{\rho f}\left(x\right)=x-\rho\nabla M_{\rho}f\left(x\right),
\end{equation}
which implies the following first-order (necessary) optimality condition
for Clarke's stationary point: 
\begin{equation}
0\in\rho^{-1}\left(x-\rho\nabla M_{\rho}f\left(x\right)-x\right)+\partial_{\circ}f\left(x-\rho\nabla M_{\rho}f\left(x\right)\right).
\end{equation}
The relation above is simplified to 
\begin{equation}
\nabla M_{\rho}f\left(x\right)\in\partial_{\circ}f\left(x-\rho\nabla M_{\rho}f\left(x\right)\right)=\nabla g\left(x-\rho\nabla M_{\rho}f\left(x\right)\right)+\partial_{\circ}h\left(x-\rho\nabla M_{\rho}f\left(x\right)\right).\label{eq:moreau-gradient-subdifferential}
\end{equation}
Consider arbitrary $\bar{x}\in\mathbb{R}^{p+1},\rho\in\left(0,L_{\nabla g}^{-1}\right)$. 

``$\Rightarrow$'' of \eqref{eq:equiv-equilibrium}: 

Let $0\in\partial_{\circ}f\left(\bar{x}\right)=\nabla g\left(\bar{x}\right)+\partial_{\circ}h\left(\bar{x}\right)$;
i.e., $\bar{x}$ is a Clarke stationary point of $f$. Then $-\nabla g\left(\bar{x}\right)\in\partial_{\circ}h\left(\bar{x}\right)$.
Since $h$ is convex, \eqref{eq:moreau-gradient-subdifferential}
implies that 
\begin{equation}
\left\langle -\nabla g\left(\bar{x}\right)-\left(\nabla M_{\rho}f\left(\bar{x}\right)-\nabla g\left(\bar{x}-\rho\nabla M_{\rho}f\left(\bar{x}\right)\right)\right),\rho\nabla M_{\rho}f\left(\bar{x}\right)\right\rangle \geq0.
\end{equation}
Simplification gives 
\begin{equation}
\left\langle \nabla g\left(\bar{x}-\rho\nabla M_{\rho}f\left(\bar{x}\right)\right)-\nabla g\left(\bar{x}\right),\nabla M_{\rho}f\left(\bar{x}\right)\right\rangle \geq\left\Vert \nabla M_{\rho}f\left(\bar{x}\right)\right\Vert ^{2}.\label{eq:equiv-inequ-1}
\end{equation}
By Cauchy-Schwartz inequality, 
\begin{equation}
\left\langle \nabla g\left(\bar{x}-\rho\nabla M_{\rho}f\left(\bar{x}\right)\right)-\nabla g\left(\bar{x}\right),\nabla M_{\rho}f\left(\bar{x}\right)\right\rangle \leq L_{\nabla g}\cdot\rho\left\Vert \nabla M_{\rho}f\left(\bar{x}\right)\right\Vert ^{2}.\label{eq:equiv-inequ-2}
\end{equation}
Since $\rho<L_{\nabla g}^{-1}$, \eqref{eq:equiv-inequ-1} and \eqref{eq:equiv-inequ-2}
imply that 
\begin{equation}
\left\Vert \nabla M_{\rho}f\left(\bar{x}\right)\right\Vert ^{2}\leq\left\langle \nabla g\left(\bar{x}-\rho\nabla M_{\rho}f\left(\bar{x}\right)\right)-\nabla g\left(\bar{x}\right),\nabla M_{\rho}f\left(\bar{x}\right)\right\rangle <\left\Vert \nabla M_{\rho}f\left(\bar{x}\right)\right\Vert ^{2},
\end{equation}
which implies that 
\begin{equation}
\nabla M_{\rho}f\left(\bar{x}\right)=0;
\end{equation}
i.e., $\bar{x}$ is the stationary point of $M_{\rho}f$, hence a
Clarke's stationary point. 

``$\Leftarrow$'' of \eqref{eq:equiv-equilibrium}: 

Let $\nabla f_{\rho}\left(\bar{x}\right)=0$; i.e. $\bar{x}$ is a
stationary point of $M_{\rho}f$. It follows directly from \eqref{eq:moreau-gradient-subdifferential}
that 
\begin{equation}
0=\nabla M_{\rho}f\left(\bar{x}\right)\in\partial_{\circ}f\left(\bar{x}-\rho\nabla M_{\rho}f\left(\bar{x}\right)\right)=\partial_{\circ}f\left(\bar{x}\right);
\end{equation}
i.e., $\bar{x}$ is a Clarke stationary point of $f$. 
\end{proof}
The vast majority of optimization algorithms for smooth objective functions
require Lipschitz continuity of the objective function. Thus, we are
to propose the following Lemma to show the Lipschitz continuity of
the gradient of the Moreau envelope of $f$. 
\begin{lem}
\label{lem:moreau-Lipschitz} $\forall\rho\in\left(0,L_{\nabla g}^{-1}\right)$,
$\exists L_{\nabla M_{\rho}f}\in\mathbb{R}_{>0}$ such that 
\begin{equation}
\forall x,y\in\mathbb{R}^{p+1},\ \left\Vert \nabla M_{\rho}f\left(x\right)-\nabla M_{\rho}f\left(y\right)\right\Vert \leq L_{\nabla M_{\rho}f}\left\Vert x-y\right\Vert .
\end{equation}
\end{lem}

\begin{proof}
Consider arbitrary $x,y\in\mathbb{R}^{p+1}$. From \eqref{eq:moreau-gradient-subdifferential},
since $h$ is convex, 
\begin{dmath}
\left\langle \nabla M_{\rho}f\left(x\right)-\nabla g\left(x-\rho\nabla M_{\rho}f\left(x\right)\right)-\left(\nabla M_{\rho}f\left(y\right)-\nabla g\left(y-\rho\nabla M_{\rho}f\left(y\right)\right)\right),x-\rho\nabla M_{\rho}f\left(x\right)-\left(y-\rho\nabla M_{\rho}f\left(y\right)\right)\right\rangle \geq0.
\end{dmath}
Simplification gives 
\begin{dmath}
\left\langle \nabla M_{\rho}f\left(x\right)-\nabla M_{\rho}f\left(y\right)-\left(\nabla g\left(x-\rho\nabla M_{\rho}f\left(x\right)\right)-\nabla g\left(y-\rho\nabla M_{\rho}f\left(y\right)\right)\right),x-y-\rho\left(\nabla M_{\rho}f\left(x\right)-\nabla M_{\rho}f\left(y\right)\right)\right\rangle \geq0.
\end{dmath}
Let $\delta_{\nabla M_{\rho}f}\coloneqq\nabla M_{\rho}f\left(x\right)-\nabla M_{\rho}f\left(y\right)$,
$\delta_{\nabla g}\coloneqq\nabla g\left(x-\rho\nabla f_{\rho}\left(x\right)\right)-\nabla g\left(y-\rho\nabla f_{\rho}\left(y\right)\right)$,
and $\delta_{x,y}\coloneqq x-y$, then 
\begin{align*}
0 & \leq\left\langle \delta_{\nabla M_{\rho}f}-\delta_{\nabla g},\delta_{x,y}-\rho\delta_{\nabla M_{\rho}f}\right\rangle \\
 & =-\rho\left\Vert \delta_{\nabla M_{\rho}f}\right\Vert ^{2}+\rho\left\langle \delta_{\nabla g},\delta_{\nabla M_{\rho}f}\right\rangle +\left\langle \delta_{\nabla M_{\rho}f},\delta_{x,y}\right\rangle -\left\langle \delta_{\nabla g},\delta_{x,y}\right\rangle \\
 & \leq-\rho\left\Vert \delta_{\nabla M_{\rho}f}\right\Vert ^{2}+\rho\left\Vert \delta_{\nabla g}\right\Vert \cdot\left\Vert \delta_{\nabla M_{\rho}f}\right\Vert +\left\Vert \delta_{\nabla M_{\rho}f}\right\Vert \cdot\left\Vert \delta_{x,y}\right\Vert +\left\Vert \delta_{\nabla g}\right\Vert \cdot\left\Vert \delta_{x,y}\right\Vert \\
 & \leq-\rho\left\Vert \delta_{\nabla M_{\rho}f}\right\Vert ^{2}+\rho L_{\nabla g}\left(\left\Vert \delta_{x,y}\right\Vert +\rho\left\Vert \delta_{\nabla M_{\rho}f}\right\Vert \right)\cdot\left\Vert \delta_{\nabla M_{\rho}f}\right\Vert \\
 & \qquad+\left\Vert \delta_{\nabla M_{\rho}f}\right\Vert \cdot\left\Vert \delta_{x,y}\right\Vert +L_{\nabla g}\left(\left\Vert \delta_{x,y}\right\Vert +\rho\left\Vert \delta_{\nabla M_{\rho}f}\right\Vert \right)\cdot\left\Vert \delta_{x,y}\right\Vert 
\end{align*}
Simplification of the above inequality gives 
\begin{equation}
\left\Vert \delta_{\nabla M_{\rho}f}\right\Vert \leq\frac{2L_{g}\rho+1+\sqrt{8L_{g}\rho+1}}{2\rho\left(1-L_{\nabla g}\rho\right)}\left\Vert \delta_{x,y}\right\Vert ;
\end{equation}
i.e., 
\begin{equation}
\left\Vert \nabla M_{\rho}f\left(x\right)-\nabla M_{\rho}f\left(y\right)\right\Vert \leq L_{\nabla M_{\rho}f}\left\Vert x-y\right\Vert ,
\end{equation}
where 
\begin{equation}
L_{\nabla M_{\rho}f}\coloneqq\frac{2L_{g}\rho+1+\sqrt{8L_{g}\rho+1}}{2\rho\left(1-L_{\nabla g}\rho\right)}>0.
\end{equation}
Following this idea, we introduce our proximal conjugate gradient
framework in Algorithm \ref{alg:proximal-general-moreau}. 
\end{proof}
\begin{algorithm}[H]
\begin{algorithmic}[1]   \State Input: A fixed value of $\rho\in\left(0,\rho^{-1}\right)$ \State Calculate the gradient of the Moreau envelope: $s^{\left( k\right)} \coloneqq \nabla M_{\rho}f\left(x^{\left( k\right)}\right)$   \State $d^{\left( k\right)} \coloneqq -s^{\left( k\right)}+\beta^{\left( k\right)}\cdot d^{\left(k-1\right)}$   \State Line search to find $\alpha^{\left( k\right)}$ for the update $x^{\left(k+1\right)}\coloneqq x^{\left( k\right)}+\alpha^{\left(k\right)} d^{\left(k\right)}$ \State Update $x^{\left(k+1\right)} \coloneqq x^{\left( k\right)}+\alpha^{\left(k\right)} d^{\left(k\right)}$ \end{algorithmic} 

\caption{Proximal Point Algorithm \label{alg:proximal-general-moreau} }
\end{algorithm}

In the above algorithm, $\beta^{\left(k\right)}$ is the conjugate parameter. The
significant meaning of Algorithm \ref{alg:proximal-general-moreau}
is that for any global convergent numerical method to find the equilibria of a globally Lipschitz flow, which generally include the global convergent first-order methods, Algorithm \ref{alg:proximal-general-moreau} can transform
such a method to a proximal counterpart. 

For some objective functions, the gradient of the Moreau envelope
can be calculated directly. However, calculation for the Moreau envelope's
gradient is not tractable for many objective functions whose smooth
component $g$ is of complicated form. Motivated by this, we further
consider the following the Moreau envelope of the objective function
with linearized $g$, such linearization step is frequently used in
proximal algorithms for statistical sparse learning problems (e.g.,
\citep{Nesterov2004,Ghadimi2013,Yang2024}). 

Consider the linearized surrogate of \eqref{eq:general-opt-problem}, the locally Lipschitz function 
$\tilde{f}\in\mathcal{C}^{0,0}\left(\mathbb{R}^{p+1},\mathbb{R}\right)$,
defined by 
\begin{align}
\tilde{f}\left(x;u\right) & \coloneqq\left\langle u,x\right\rangle +h\left(x\right)\label{eq:linearizing-nonconvex}\\
\text{prox}_{\rho\tilde{f}}\left(x;u\right) & =\underset{y\in\mathbb{R}^{p+1}}{\arg\min\ }\left\{ \left\langle u,y\right\rangle +\frac{1}{2\rho}\left\Vert y-x\right\Vert ^{2}+h\left(y\right)\right\} \label{eq:practical-proximal}\\
\nabla_{x}M_{\rho}\tilde{f}\left(x;u\right) & =\rho^{-1}\left(x-\text{prox}_{\rho\tilde{f}}\left(x;u\right)\right)\label{eq:practical-moreau-gradient}
\end{align}
$\text{prox}_{\rho\tilde{f}}\left(x;u\right)$ is the proximal operator
applied on $\tilde{f}$, and $\nabla_{x}M_{\rho}\tilde{f}\left(x;u\right)$
is the gradient of the Moreau envelope of $\tilde{f}$. The linearization
term $\left\langle u,x\right\rangle $ in \eqref{eq:linearizing-nonconvex}
depends on $u$. Recognize that $\tilde{f}\left(x;u\right)$ is linearizing
the nonconvex smooth component $g$ in \eqref{eq:opt-problem} when
$u=\nabla g\left(x\right)$. 

We establish several definitions for subsequent utilization. Define
the mapping $\tilde{g}_{\rho}=I-\rho\nabla g\in \mathcal{C}^{0,0}\left(\mathbb{R}^{p+1},\mathbb{R}^{p+1}\right)$
\emph{for some $\rho\in\left(0,L_{\nabla g}^{-1}\right)$}, the locally Lipschitz property of $\tilde{g}_{\rho}$ follows from $g\in \mathcal{C}^{1,1}$; i.e.,
$\tilde{g}_{\rho}\left(x\right)\coloneqq x-\rho\nabla g\left(x\right).$
The following Lemma identifies some fundamental property of $\tilde{g}_{\rho}$. 
\begin{lem}
\label{lem:bijectivity} $\tilde{g}_{\rho}$ is a bijective from $\mathbb{R}^{p+1}$
to $\mathbb{R}^{p+1}$, and $\tilde{g}_{\rho}^{-1}$ is globally Lipschitz
with constant $\left(1-\rho L_{\nabla g}\right)^{-1}$. 
\end{lem}

\begin{proof}
\emph{Injectivity proof: }

Consider arbitrary $x_{1},x_{2}\in\mathbb{R}^{p+1}$. Since $\rho\in\left(0,L_{\nabla g}^{-1}\right)$,
$x_{1}-\rho\nabla g\left(x_{1}\right)=x_{2}-\rho\nabla g\left(x_{2}\right)$
implies that 
\begin{equation}
\left\Vert x_{1}-x_{2}\right\Vert =\rho\left\Vert \nabla g\left(x_{1}\right)-\nabla g\left(x_{2}\right)\right\Vert \leq\rho L_{\nabla g}\left\Vert x_{1}-x_{2}\right\Vert <\left\Vert x_{1}-x_{2}\right\Vert ,
\end{equation}
hence $x_{1}=x_{2}$. This shows that $\tilde{g}_{\rho}$ is a injective
mapping. 

\emph{Surjectivity proof: }

Consider arbitrary $y_{1},y_{2}\in\mathbb{R}^{p+1}$. Consider arbitrary
$z\in\mathbb{R}^{p+1}$. Define mapping $\mathcal{T}\left(y\right)\coloneqq z+\rho\nabla g\left(y\right)$,
then 
\begin{align*}
\left\Vert \mathcal{T}\left(y_{1}\right)-\mathcal{T}\left(y_{2}\right)\right\Vert  & =\left\Vert z+\rho\nabla g\left(y_{1}\right)-\left(z+\rho\nabla g\left(y_{2}\right)\right)\right\Vert \\
 & =\rho\left\Vert \nabla g\left(y_{1}\right)-\nabla g\left(y_{2}\right)\right\Vert \\
 & \leq\rho L_{\nabla g}\left\Vert y_{1}-y_{2}\right\Vert \\
 & <\left\Vert y_{1}-y_{2}\right\Vert .
\end{align*}
Thus, $\mathcal{T}$ is a contraction mapping, since $\mathbb{R}^{p+1}$
equipped with Euclidean topology is a Banach space, by Banach fixed
point theorem, $\mathcal{T}$ has a fixed point; i.e., $\exists y\in\mathbb{R}^{p+1}$
such that $y=z+\rho\nabla g\left(y\right)$, or equivalently, $\tilde{g}_{\rho}\left(y\right)=y-\rho\nabla g\left(y\right)=z$.
Thus, $\mathbb{R}^{p+1}\subseteq\tilde{g}_{\rho}\left(\mathbb{R}^{p+1}\right)$. 

\emph{Globally Lipschitz constant derivation for inverse map: }

Since $\nabla g$ is globally $L_{\nabla g}-$Lipschitz, 
\begin{align}
\left\Vert \tilde{g}_{\rho}\left(y_{1}\right)-\tilde{g}_{\rho}\left(y_{2}\right)\right\Vert  & =\left\Vert y_{1}-\rho\nabla g\left(y_{1}\right)-\left(y_{2}-\rho\nabla g\left(y_{2}\right)\right)\right\Vert \nonumber \\
 & =\left\Vert y_{1}-y_{2}-\rho\left(\nabla g\left(y_{1}\right)-\nabla g\left(y_{2}\right)\right)\right\Vert \nonumber \\
 & \geq\left|\left\Vert y_{1}-y_{2}\right\Vert -\left\Vert \rho\left(\nabla g\left(y_{1}\right)-\nabla g\left(y_{2}\right)\right)\right\Vert \right|\nonumber \\
 & =\left|\left\Vert y_{1}-y_{2}\right\Vert -\rho\left\Vert \nabla g\left(y_{1}\right)-\nabla g\left(y_{2}\right)\right\Vert \right|\nonumber \\
 & =\left\Vert y_{1}-y_{2}\right\Vert -\rho\left\Vert \nabla g\left(y_{1}\right)-\nabla g\left(y_{2}\right)\right\Vert \label{eq:triangle-inequ}\\
 & \geq\left(1-\rho L_{\nabla g}\right)\left\Vert y_{1}-y_{2}\right\Vert \nonumber 
\end{align}
where \eqref{eq:triangle-inequ} is due to the fact that 
\begin{equation}
\rho\left\Vert \nabla g\left(y_{1}\right)-\nabla g\left(y_{2}\right)\right\Vert \leq\rho L_{\nabla g}\left\Vert y_{1}-y_{2}\right\Vert <\left\Vert y_{1}-y_{2}\right\Vert .
\end{equation}
Since $\tilde{g}_{\rho}$ is surjective, consider arbitrary $z_{1},z_{2}\in\mathbb{R}^{p+1}$
let $y_{1}\coloneqq\tilde{g}_{\rho}^{-1}\left(z_{1}\right)$ and $y_{2}\coloneqq\tilde{g}_{\rho}^{-1}\left(z_{2}\right)$,
then 
\begin{equation}
\left\Vert \tilde{g}_{\rho}^{-1}\left(z_{1}\right)-\tilde{g}_{\rho}^{-1}\left(z_{2}\right)\right\Vert \leq\left(1-\rho L_{\nabla g}\right)^{-1}\left\Vert z_{1}-z_{2}\right\Vert .
\end{equation}
\end{proof}
Define 
\begin{equation}
\mathcal{G}_{\rho\tilde{f}}\left(x\right)\coloneqq\nabla_{x}M_{\rho}\tilde{f}\left(x;u\right)\label{eq:target-gradient-flow}
\end{equation}
with $u=\nabla g\left(x\right)$; i.e., $\mathcal{G}_{\rho\tilde{f}}\left(x\right)$
is the gradient of the Moreau envelope of $\tilde{f}$. 

Similarly to Lemma \ref{lem:equivalent-clarke-stationary-point} and
\ref{lem:moreau-Lipschitz}, we are to prove that the set of Clarke's
stationary of \eqref{eq:opt-problem} is identical to the set $\left\{ \bar{x}\in\mathbb{R}^{p+1}\vert\mathcal{G}_{\rho\tilde{f}}\left(\bar{x}\right)=0\right\} $
in Lemma \ref{lem:equivalent-clarke-stationary-point-pratical},
and then we are to show that \eqref{eq:target-gradient-flow} is globally
Lipschitz in Lemma \ref{lem:moreau-Lipschitz-practical}. 
\begin{lem}
\label{lem:equivalent-clarke-stationary-point-pratical} $\forall\bar{x}\in\mathbb{R}^{p+1},\rho\in\mathbb{R}_{>0}$,
\begin{equation}
0\in\partial_{\circ}f\left(\bar{x}\right)\Leftrightarrow\mathcal{G}_{\rho\tilde{f}}\left(\bar{x}\right)=0.\label{eq:equiv-equilibrium-practical}
\end{equation}
\end{lem}

\begin{proof}
Consider arbitrary $x\in \mathbb{R}^{p+1}$. The $\tilde{f}$ is convex since it is a sum of convex function $h$
and a linear mapping of $x$, which is convex. 
\begin{align}
\mathcal{G}_{\rho\tilde{f}}\left(x\right) & =\rho^{-1}\left(x-\text{prox}_{\rho\tilde{f}}\left(x;\nabla g\left(x\right)\right)\right)\label{eq:convex-moreau-envelope-grad}\\
 & =\rho^{-1}\left(x-\text{prox}_{\rho h}\left(x-\rho\nabla g\left(x\right)\right)\right)\label{eq:moreau-envelope-grad-decomposition}\\
 & =\nabla g\left(x\right)+\rho^{-1}\left(x-\rho\nabla g\left(x\right)-\text{prox}_{\rho h}\left(x-\rho\nabla g\left(x\right)\right)\right)\nonumber \\
 & =\nabla g\left(x\right)+\left(\nabla M_{\rho}h\right)\circ\tilde{g}_{\rho}\left(x\right)\label{eq:moreau-envelope-grad-decomposition-final-form}
\end{align}
\eqref{eq:moreau-envelope-grad-decomposition} is due to the affine
addition property of proximal mapping. From \eqref{eq:linearizing-nonconvex}
and \eqref{eq:convex-moreau-envelope-grad}, 
\begin{align}
\mathcal{G}_{\rho\tilde{f}}\left(x\right) & =\rho^{-1}\left(x-\text{prox}_{\rho\tilde{f}}\left(x,\nabla g\left(x\right)\right)\right)\nonumber \\
\implies\text{prox}_{\rho\tilde{f}}\left(x,\nabla g\left(x\right)\right) & =x-\rho\cdot\mathcal{G}_{\rho\tilde{f}}\left(x\right)\nonumber \\
\implies0 & \in\rho^{-1}\left(x-\rho\cdot\mathcal{G}_{\rho\tilde{f}}\left(x\right)-x\right)+\partial_{\circ}\tilde{f}\left(x-\rho\cdot\mathcal{G}_{\rho\tilde{f}}\left(x\right)\right)\nonumber \\
\implies0 & \in-\mathcal{G}_{\rho\tilde{f}}\left(x\right)+\nabla g\left(x\right)+\partial_{\circ}h\left(x-\rho\cdot\mathcal{G}_{\rho\tilde{f}}\left(x\right)\right)\nonumber \\
\implies\mathcal{G}_{\rho\tilde{f}}\left(x\right) & \in\nabla g\left(x\right)+\partial_{\circ}h\left(x-\rho\cdot\mathcal{G}_{\rho\tilde{f}}\left(x\right)\right)\label{eq:practical-subdifferential}\\
\implies\mathcal{G}_{\rho\tilde{f}}\left(x\right)-\nabla g\left(x\right) & \in\partial_{\circ}h\left(x-\rho\cdot\mathcal{G}_{\rho\tilde{f}}\left(x\right)\right).\nonumber 
\end{align}
Thus, since $h$ is convex, $\forall v\in\partial_{\circ}h\left(x\right)$,
\begin{align}
\left\langle \mathcal{G}_{\rho\tilde{f}}\left(x\right)-\nabla g\left(x\right)-v,x-\rho\cdot\mathcal{G}_{\rho\tilde{f}}\left(x\right)-x\right\rangle  & \geq0\nonumber \\
\implies\left\langle \mathcal{G}_{\rho\tilde{f}}\left(x\right)-\nabla g\left(x\right)-v,\mathcal{G}_{\rho\tilde{f}}\left(x\right)\right\rangle  & \leq0\nonumber \\
\implies\left\Vert \mathcal{G}_{\rho\tilde{f}}\left(x\right)\right\Vert ^{2} & \leq\left\langle \nabla g\left(x\right)+v,\mathcal{G}_{\rho\tilde{f}}\left(x\right)\right\rangle \label{eq:moreau-gradient-CS}\\
 & \leq\left\Vert \nabla g\left(x\right)+v\right\Vert \cdot\left\Vert \mathcal{G}_{\rho\tilde{f}}\left(x\right)\right\Vert \nonumber \\
\implies\left\Vert \mathcal{G}_{\rho\tilde{f}}\left(x\right)\right\Vert  & \leq\left\Vert \nabla g\left(x\right)+v\right\Vert ,\label{eq:moreau-gradient-minimal}
\end{align}
provided that $\left\Vert \mathcal{G}_{\rho\tilde{f}}\left(x\right)\right\Vert \neq0$.
Basic results on the Moreau envelope shows that $\mathcal{G}_{\rho\tilde{f}}\left(x\right)=0$
implies that $x$ is a Clarke stationary point of $\tilde{f}\left(x\right)$. 

Now we are proceed to prove \eqref{eq:equiv-equilibrium-practical}: 

``$\Rightarrow$'':

Consider arbitrary $\bar{x}\in\mathbb{R}^{p+1}$ and $\rho\in\mathbb{R}_{>0}$.
Let $0\in\partial_{\circ}f\left(\bar{x}\right)=\nabla g\left(\bar{x}\right)+\partial_{\circ}h\left(\bar{x}\right)$;
i.e., $\bar{x}$ is a Clarke stationary point of $f$. Then $\exists v\in\partial_{\circ}h\left(\bar{x}\right)$
such that $\nabla g\left(\bar{x}\right)+v=0$. \eqref{eq:moreau-gradient-minimal}
implies that 
\begin{equation}
\left\Vert \mathcal{G}_{\rho\tilde{f}}\left(\bar{x}\right)\right\Vert \leq\left\Vert \nabla g\left(\bar{x}\right)+v\right\Vert =0.
\end{equation}
Thus, $\mathcal{G}_{\rho\tilde{f}}\left(\bar{x}\right)=0$. 

``$\Leftarrow$'': 

Consider arbitrary $\bar{x}\in\mathbb{R}^{p+1}$. Let $\mathcal{G}_{\rho\tilde{f}}\left(\bar{x}\right)=0$;
i.e., $\mathcal{G}_{\rho\tilde{f}}\left(\bar{x}\right)=0$ is stationary.
\eqref{eq:practical-subdifferential} implies that 
\begin{equation}
0=\mathcal{G}_{\rho\tilde{f}}\left(\bar{x}\right)\in\nabla g\left(\bar{x}\right)+\partial_{\circ}h\left(\bar{x}-\rho\cdot\mathcal{G}_{\rho\tilde{f}}\left(\bar{x}\right)\right)=\nabla g\left(\bar{x}\right)+\partial_{\circ}h\left(\bar{x}\right)=\partial_{\circ}f\left(\bar{x}\right).
\end{equation}
Thus, $\bar{x}$ is a Clarke stationary point of $f$. 
\end{proof}
\begin{lem}
\label{lem:moreau-Lipschitz-practical} $\forall\rho\in\mathbb{R}_{>0}$,
$\exists L_{\mathcal{G}_{\rho\tilde{f}}}\in\mathbb{R}_{>0}$ such
that 
\begin{equation}
\forall x,y\in\mathbb{R}^{p+1},\ \left\Vert \mathcal{G}_{\rho\tilde{f}}\left(x\right)-\mathcal{G}_{\rho\tilde{f}}\left(y\right)\right\Vert \leq L_{\mathcal{G}_{\rho\tilde{f}}}\left\Vert x-y\right\Vert .
\end{equation}
\end{lem}

\begin{proof}
Consider arbitrary $x,y\in\mathbb{R}^{p+1}$ and $\rho\in\mathbb{R}_{>0}$.
Let $u\coloneqq\nabla g\left(x\right)$ and $v\coloneqq\nabla g\left(y\right)$,
\begin{align}
\left\Vert \mathcal{G}_{\rho\tilde{f}}\left(x\right)-\mathcal{G}_{\rho\tilde{f}}\left(y\right)\right\Vert  & =\left\Vert \nabla_{x}M_{\rho}\tilde{f}\left(x;u\right)-\nabla_{y}M_{\rho}\tilde{f}\left(y;v\right)\right\Vert \nonumber \\
 & =\left\Vert \nabla_{x}M_{\rho}\tilde{f}\left(x;u\right)-\nabla_{x}M_{\rho}\tilde{f}\left(x;v\right)+\nabla_{x}M_{\rho}\tilde{f}\left(x;v\right)-\nabla_{y}M_{\rho}\tilde{f}\left(y;v\right)\right\Vert \nonumber \\
 & \leq\left\Vert \nabla_{x}M_{\rho}\tilde{f}\left(x;u\right)-\nabla_{x}M_{\rho}\tilde{f}\left(x;v\right)\right\Vert +\left\Vert \nabla_{x}M_{\rho}\tilde{f}\left(x;v\right)-\nabla_{y}M_{\rho}\tilde{f}\left(y;v\right)\right\Vert \nonumber \\
 & \leq\left\Vert u-v\right\Vert +\left\Vert \nabla_{x}M_{\rho}\tilde{f}\left(x;v\right)-\nabla_{y}M_{\rho}\tilde{f}\left(y;v\right)\right\Vert \label{eq:lemma-4-result}\\
 & \leq\left\Vert u-v\right\Vert + \rho^{-1}\left\Vert x-y\right\Vert \label{eq:convex-moreau-grad-nonexpansive}\\
 & \leq L_{\nabla g}\left\Vert x-y\right\Vert + \rho^{-1}\left\Vert x-y\right\Vert \nonumber \\
 & =\left(L_{\nabla g}+\rho^{-1}\right)\left\Vert x-y\right\Vert \nonumber 
\end{align}
\eqref{eq:lemma-4-result} is due to Lemma 4 in \citep{Ghadimi2013},
and \eqref{eq:convex-moreau-grad-nonexpansive} is due to $\tilde{f}\left(\cdot;v\right)$
is convex and the fact that the gradient of a convex function's Moreau
envelope is $\rho^{-1}-$Lipschitz. Therefore, let 
\begin{equation}
L_{\mathcal{G}_{\rho\tilde{f}}}\coloneqq L_{\nabla g}+\rho^{-1}
\end{equation}
and we have 
\begin{equation}
\left\Vert \mathcal{G}_{\rho\tilde{f}}\left(x\right)-\mathcal{G}_{\rho\tilde{f}}\left(y\right)\right\Vert \leq L_{\mathcal{G}_{\rho\tilde{f}}}\left\Vert x-y\right\Vert .
\end{equation}
\end{proof}

Furthermore, \eqref{eq:moreau-envelope-grad-decomposition-final-form}
suggests that 
\begin{equation}
\mathcal{G}_{\rho\tilde{f}}=\nabla g+\left(\nabla M_{\rho}h\right)\circ\tilde{g}_{\rho}=\nabla g+\left(\nabla M_{\rho}h\right)\circ\left(Id-\rho\nabla g\right).
\end{equation}
Hence, 
\begin{align}
Id-\rho\mathcal{G}_{\rho\tilde{f}} & =Id-\rho\nabla g-\rho\left(\nabla M_{\rho}h\right)\circ\left(Id-\rho\nabla g\right)\nonumber \\
 & =\tilde{g}_{\rho}-\rho\left(\nabla M_{\rho}h\right)\circ\tilde{g}_{\rho}\nonumber \\
 & =\left(Id-\rho\left(\nabla M_{\rho}h\right)\right)\circ\tilde{g}_{\rho}\label{eq:flow-composition}\\
 & =\tilde{g}_{\rho}^{-1}\circ\left(\tilde{g}_{\rho}\circ\left(Id-\rho\left(\nabla M_{\rho}h\right)\right)\right)\circ\tilde{g}_{\rho}\label{eq:topologically-conjugate}
\end{align}
shows that $\tilde{g}_{\rho}^{-1}\circ\left(\tilde{g}_{\rho}\circ\left(Id-\rho\left(\nabla M_{\rho}h\right)\right)\right)\circ\tilde{g}_{\rho}:\mathbb{R}^{p+1}\mapsto\mathbb{R}^{p+1}$
equals to $Id-\rho\mathcal{G}_{\rho\tilde{f}}:\mathbb{R}^{p+1}\mapsto\mathbb{R}^{p+1}$.
Since $\tilde{g}_{\rho}$ ibijectivee, and that $\tilde{g}_{\rho}$
and $\tilde{g}_{\rho}^{-1}$ are continuous due to the globally Lipschitz
property from Lemma \ref{lem:bijectivity}, $\tilde{g}_{\rho}$ is
a homeomorphism. Hence, $Id-\rho\mathcal{G}_{\rho\tilde{f}}$ and
$\tilde{g}_{\rho}\circ\left(Id-\rho\left(\nabla M_{\rho}h\right)\right)$
are topologically equivalent  mappings via the homeomorphism $\tilde{g}_{\rho}$.
Lemma \ref{lem:moreau-Lipschitz-practical} implies that $\mathcal{G}_{\rho\tilde{f}}$
is globally Lipschitz, which sufficiently implies by the Cauchy-Lipschitz
theorem that the differential equation 
\begin{equation}
\dot{x}\coloneqq\frac{dx}{dt}=\mathcal{G}_{\rho\tilde{f}}\left(x\right)
\end{equation}
has a unique solution for any given initial value condition. Thus,
$\mathcal{G}_{\rho\tilde{f}}$ generates a unique flow under a given
initial value condition. 

The operator equations presented above can be understood as demonstrating
how $\mathcal{G}_{\rho\tilde{f}}$ functions analogously to a gradient
operator. Specifically, $Id-\rho\mathcal{G}_{\rho\tilde{f}}$ represents
executing a descent operation in the $-\mathcal{G}_{\rho\tilde{f}}$
direction with a step size $\rho$. Similarly, $Id-\rho\left(\nabla M_{\rho}h\right)$
represents a single gradient descent step with $\rho$ as the step
size with objective function $M_{\rho}h$, the Moreau envelope of
$h$; while $\tilde{g}_{\rho}=Id-\rho\nabla g$ reflects a gradient
descent step with objective function $g$, again with $\rho$ as the
step size. Equation \eqref{eq:flow-composition} elucidates that a
descent in the $-\mathcal{G}_{\rho\tilde{f}}$ direction is identical
to first performing a one-step gradient descent on $g$, followed
by $M_{\rho}h$; or performing gradient descents in a converse order
yields a topologically equivalence via the homeomorphism $\tilde{g}_{\rho}$,
as shown in \eqref{eq:topologically-conjugate}. 

In short summary, the approach based on linearization of the smooth term and the Moreau envelope
enables us to build equivalence between identifying Clarke stationary
points of the original nonsmooth objective function \eqref{eq:opt-problem}
and finding equilibria of the (unique) flow generated by $\mathcal{G}_{\rho\tilde{f}}$,
as demonstrated in Lemma \ref{lem:equivalent-clarke-stationary-point-pratical}.
The task of finding equilibria within a globally Lipschitz continuous
flow, such as the $\mathcal{G}_{\rho\tilde{f}}$ flow, is well explored
within mathematics, particularly in the realms of dynamical systems
and numerical analysis (see, for example, \citep{Quarteroni2007,Atkinson1989,Lubich2006,Hubbard1995,Helmke1994}).
Cauchy-Lipschitz theorem establishes the uniqueness of solutions to
initial value problems for globally Lipschitz continuous flows; while the existence of equilibria is a direct result of Brouwer fixed-point theorem. Numerical
methods for dynamical systems, including methods for finding equilibria
of the flow, are largely based on this uniqueness result. This is
one reason that the vast majority of numerical methods in the context
of dynamical systems require the flow to be globally Lipschitz. It is
important to note that these numerical strategies, widely applied
across dynamical systems, do not hinge on the flow being derived from
a conservative field. As such, the process of formulating a potential
function for $\mathcal{G}_{\rho\tilde{f}}$ is not a prerequisite
for employing numerical techniques to determine its equilibria. This
perspective underscores the versatility of numerical methods in dynamical
systems in finding the equilibria of flows, regardless of the explicit
existence of a potential function, a stance corroborated by various
sources in the literature \citep{Quarteroni2007,Atkinson1989,Lubich2006,Hubbard1995,Helmke1994,Ross2019,Riahi2018}.
In this view, the construction of a potential function for $\mathcal{G}_{\rho\tilde{f}}$
is generally not necessary when deploying numerical analysis methods
to find its equilibria. 

In the context of nonlinear conjugate gradient algorithms for optimization,
achieving global convergence on nonconvex objective functions that
are globally Lipschitz-smooth implies that such methods can reliably
find equilibria within the corresponding flow dynamics \citep{Ross2019,Riahi2018}.
These algorithms typically incorporate a line search step, which may
use a surrogate objective function instead of the original. This surrogate
can be a constructed potential, Lyapunov, or energy function, offering
flexibility when finding the potential function for $\mathcal{G}_{\rho\tilde{f}}$
poses challenges \citep{Ross2019,Clarke2004,Sontag1998}. 

When it is feasible to construct a potential function whose gradient with respect
to $x$ is $\mathcal{G}_{\rho\tilde{f}}$, the associated
objective function and its gradient become more manageable, allowing
for direct global convergence arguments. If constructing a potential
function with respect to $x$ for the $\left(\nabla M_{\rho}h\right)\circ\tilde{g}_{\rho}\left(x\right)$
term in \eqref{eq:moreau-envelope-grad-decomposition-final-form}
or $\nabla g\circ\left(Id-\rho\left(\nabla M_{\rho}h\right)\right)$
in \eqref{eq:topologically-conjugate} is tractable, the objective
function with gradient being \eqref{eq:moreau-envelope-grad-decomposition-final-form}
or $\tilde{g}_{\rho}\circ\left(Id-\rho\left(\nabla M_{\rho}h\right)\right)$
can hence be easily constructed. Thus, arguments for global convergence for methods
based on the objective function and its gradient directly follow
to prove the global convergence of the numerical optimization algorithm when
applied to the constructed potential function for $\mathcal{G}_{\rho\tilde{f}}$.
We remark that $Id-\rho\mathcal{G}_{\rho\tilde{f}}$ and $\tilde{g}_{\rho}\circ\left(Id-\rho\left(\nabla M_{\rho}h\right)\right)$
generate two topologically equivalent flows via homeomorphism $\tilde{g}_{\rho}$;
thus, their equilibria can be transformed by $\tilde{g}_{\rho}$ and
share the same stability. In the context of numerical optimization, this implies
that a fixed point $\bar{x}$ for the mapping $Id-\rho\mathcal{G}_{\rho\tilde{f}}$
corresponds bijectively to a fixed point $\tilde{g}_{\rho}\left(\bar{x}\right)$
for $\tilde{g}_{\rho}\circ\left(Id-\rho\left(\nabla M_{\rho}h\right)\right)$.
Characterized by the first-order optimality condition in optimization
of smooth functions, or equivalently, the stationary condition in dynamical
system, 
\begin{equation}
\mathcal{G}_{\rho\tilde{f}}\left(\bar{x}\right)=\nabla g\left(\bar{x}\right)+\left(\nabla M_{\rho}h\right)\circ\tilde{g}_{\rho}\left(\bar{x}\right)=0\Leftrightarrow\nabla M_{\rho}h\left(\tilde{g}_{\rho}\left(\bar{x}\right)\right)+\nabla g\left(\tilde{g}_{\rho}\left(\bar{x}\right)-\rho\nabla M_{\rho}h\left(\tilde{g}_{\rho}\left(\bar{x}\right)\right)\right)=0.
\end{equation}
This approach is practical because the literature on first-order numerical
optimization techniques frequently includes proofs of global convergence
for methods that depend on the objective function and its gradient
(for example, see \citep{Fletcher1964,Polak1969,Hestenes1952,Dai1999,Hager2005}).
Alternatively, construction of a potential function for $\mathcal{G}_{\rho\tilde{f}}$ is often not necessary due to the fact that fixed-point methods finding equilibria for a flow mostly establish convergence properties based on Banach fixed point theorem. This theorem guarantees convergence through intrinsic flow properties, obviating the need for a potential function \citep{Burden2016, Atkinson1989, Agarwal2009}. Conventionally, the use of line search based on the objective function and its gradient has been applied in some numerical methods to ensure global convergence. However, with the rapid growth of research in high--dimensional statistical machine learning and large-scale optimization, evaluations of the objective function often proven to be inefficient. Consequently, recent years have seen the exploration of two main alternatives. For instance, two different types of approaches for global convergent nonlinear conjugate gradient methods have been proposed without the conventional objective function-based line search procedure. One type of approach ensures global convergence by utilizing a line search mechanism that depends only on the nonlinear equation that generates the flow \citep{Feng2017, Snyman1985, Snyman2004, Kafka2019}; that is, the gradient function for smooth optimization, or $\mathcal{G}_{\rho\tilde{f}}$ in our case. As an example, under the smoothness assumption, the first-order optimality condition for an exact line search often solves for $\alpha$ with the current value $x^{\left(k\right)}$ and the search direction $d^{\left(k\right)}$ from $\left\langle \mathcal{G}_{\rho\tilde{f}}\left(x^{\left(k\right)}+\alpha\cdot d^{\left(k\right)}\right),d^{\left(k\right)}\right\rangle =0$, an equation dependent only on $\mathcal{G}_{\rho\tilde{f}}$ but not any surrogate objective function. From a practical perspective, this one-dimensional root finding problem can be carried out efficiently using the Brent root finding algorithm \citep{Brent1971}. The other approach suggests achieving global convergence either without the need for line search \citep{Shi2005, Chen2018, Sun2001, Wu2011a, Wang2006, Zhou2009} or by meeting a condition related to the Zoutendijk condition to replace the Wolfe-Powell conditions of sufficient descent (Armijo) and curvature \citep{Neumaier2024}. Additionally, in scenarios where the fulfillment of a sufficient descent (Armijo) condition is imperative, the formulation of a surrogate objective function becomes essential. Considering \eqref{eq:moreau-envelope-grad-decomposition-final-form}, where a surrogate objective is required for the line search phase, it could be formulated as: 
\begin{align}
 & g\left(x\right)+\left(M_{\rho}h\right)\circ\tilde{g}_{\rho}\left(x\right)\nonumber \\
 & =g\left(x\right)+\left(M_{\rho}h\right)\circ\tilde{g}_{\rho}\left(x\right)+\text{constant}\nonumber \\
 & =g\left(x\right)+\left\langle \nabla g\left(x\right),\text{prox}_{\rho h}\left(x-\rho\nabla g\left(x\right)\right)-x\right\rangle +\frac{1}{2\rho}\left\Vert \text{prox}_{\rho h}\left(x-\rho\nabla g\left(x\right)\right)-x\right\Vert ^{2}\label{eq:evaluate-obj}\\
 & \qquad+h\left(\text{prox}_{\rho h}\left(x-\rho\nabla g\left(x\right)\right)\right)+\text{constant}\nonumber 
\end{align}
This formulation, denoted as \eqref{eq:evaluate-obj}, represents
a quadratic approximation of $g$ plus the nonsmooth term $h$, evaluated
at $\text{prox}_{\rho h}\left(x-\rho\nabla g\left(x\right)\right)$. This type of formulation has often been used for the line search step in previous studies \citep{Beck2009, Kanzow2020}. 
The addition of the term $-\left\langle \nabla g\left(x\right),x\right\rangle $
acts as a constant in \eqref{eq:linearizing-nonconvex}, analogous
to fixing the value of $u$ as $\nabla g\left(x\right)$ for linearization.
This constant term, $-\left\langle \nabla g\left(x\right),x\right\rangle $,
doesn't alter the gradient of the Moreau envelope \eqref{eq:practical-moreau-gradient}
or the proximal point \eqref{eq:practical-proximal}, serving to frame
the quadratic approximation of $g\left(\text{prox}_{\rho h}\left(x-\rho\nabla g\left(x\right)\right)\right)$. 

Evaluation of $\text{prox}_{\rho h}$ in \eqref{eq:evaluate-obj}
is tractable and efficient for many functions, such as the $\ell_{1}$
norm commonly encountered in sparse statistical learning can be efficiently
computed via the soft-thresholding function. Given that line search
rules such as the Wolfe-Powell or Armijo-Goldstein conditions require
only the difference in the value of the objective function at two points to decide
on the step size, the constant term in \eqref{eq:evaluate-obj} can
be disregarded. Subsequent global convergence arguments stem from
the fixed-point theory analysis of the numerical methods deployed to find the equilibria of the $\mathcal{G}_{\rho\tilde{f}}$
flow. Another possible surrogate objective function inspired by the
quadratic Lyapunov function for the $\mathcal{G}_{\rho\tilde{f}}$
flow could be $\frac{1}{2}\left\Vert \mathcal{G}_{\rho\tilde{f}}\right\Vert ^{2}$,
attains its minimal value $0$ exactly at the $\mathcal{G}_{\rho\tilde{f}}$
flow's equilibria. This quadratic approach simplifies evaluation,
but it may not offer insights into the potential function's
landscape, potentially limiting the numerical algorithm's acceleration
capabilities if such an algorithm uses the landscape information to ensure the sufficient descent (Armijo) condition. Therefore, formulating the surrogate objective function
preserving the landscape of the original objective function as
outlined in \eqref{eq:evaluate-obj} is preferable. 

Building on the above discussion, we introduce our practical proximal
conjugate gradient framework in Algorithm \ref{alg:proximal-practical-moreau}. 

\begin{algorithm}[H]
\begin{algorithmic}[1]   \State Input: A fixed value of $\rho\in\left(0,\rho^{-1}\right)$ \State Calculate the proximal value $p^{\left( k\right)} \coloneqq \text{prox}_{{\rho}^{-1}h}\left(x^{\left( k\right)}-{\rho}^{-1}\cdot \nabla g\left(x^{\left( k\right)}\right)\right)$   \State Calculate $\mathcal{G}_{\rho\tilde{f}}\left(x^{\left( k\right)}\right)$: $s^{\left( k\right)} \coloneqq {\rho}\left(x^{\left( k\right)}-p^{\left( k\right)}\right)$   \State $d^{\left( k\right)} \coloneqq -s^{\left( k\right)}+\beta^{\left( k\right)}\cdot d^{\left(k-1\right)}$   \State Line search to find $\alpha^{\left( k\right)}$ for the update $x^{\left(k+1\right)}\coloneqq x^{\left( k\right)}+\alpha^{\left(k\right)} d^{\left(k\right)}$, if needed.  \State Update $x^{\left(k+1\right)} \coloneqq x^{\left( k\right)}+\alpha^{\left(k\right)} d^{\left(k\right)}$ \end{algorithmic} 

\caption{Computationally Tractable Proximal Conjugate Gradient Update Scheme \label{alg:proximal-practical-moreau} }
\end{algorithm}

In Algorithm \ref{alg:proximal-practical-moreau}, $\beta^{\left(k\right)}$ functions
as the conjugate parameter. Unlike Algorithm \ref{alg:proximal-general-moreau},
Algorithm \ref{alg:proximal-practical-moreau} facilitates the update
process without the need to compute $\nabla M_{\rho}f\left(x^{\left(k\right)}\right)$.
This adaptation is significantly valuable in practical scenarios,
especially in statistical sparse learning challenges characterized
by a complicated smooth component $g$ alongside a simple nonsmooth
convex component $h$. In such cases, computing $\text{prox}_{\rho h}$
is markedly more tractable and efficient than $\text{prox}_{\rho f}$.
This approach is particularly beneficial for sparse statistical learning
issues, where sparsity is commonly induced by an $\ell_{1}$ penalty
term. 

\subsection{\label{subsec:Proximal-HZ} Proximal Hager-Zhang \citep{Hager2005}
Conjugate Gradient }

The nonlinear conjugate gradient method represents the pinnacle of
first-order techniques for addressing smooth optimization challenges.
Various versions of nonlinear conjugate gradient methods have been
introduced, including the Fletcher-Reeves (FR) method \citep{Fletcher1964},
the modified Polak-Ribiere-Polyak (PRP+) method \citep{Polak1969,Gilbert1992},
the Hestenes-Stiefel (HS) method \citep{Hestenes1952}, the Dai-Yuan
(DY) method \citep{Dai1999}, and the Hager-Zhang (HZ) method \citep{Hager2005}.
These versions have all demonstrated global convergence with nonconvex
globally Lipschitz-smooth objective functions. Among these, the Hager-Zhang
conjugate gradient method is notable for delivering the best numerical
performance on large-scale datasets, as indicated in previous research
\citep{Hager2006}. Building on this, having introduced our practical
proximal conjugate gradient update mechanism in Algorithm \ref{alg:proximal-practical-moreau},
we aim to extend this approach by adapting the smooth Hager-Zhang
nonlinear conjugate gradient method to its proximal version
in Algorithm \ref{alg:proximal-HZ-CG}. 

\begin{algorithm}
\begin{algorithmic}[1]  \State \textbf{Input:} Initial point $x^{\left( 0\right)}$; $g\in\mathcal{C}^{1,1}\left(\mathbb{R}^{p+1},\mathbb{R}\right)$; locally-Lipschitz, convex $h\in\mathcal{C}^{0,0}\left(\mathbb{R}^{p+1},\mathbb{R}\right)$; the smoothing parameter for the Moreau envelope $\rho\in\left(0,\rho^{-1}\right)$; $k \coloneqq 0$  \State \textbf{Output:} $p$  \State $k += 1$  \State Calculate the gradient for $g$: $g^{\left( 0\right)} \coloneqq \nabla g\left(x^{\left( 0\right)}\right)$  \State Calculate the proximal value $p^{\left( 0\right)} \coloneqq \text{prox}_{\rho,h}\left(x^{\left( 0\right)}-\rho\cdot g^{\left( 0\right)}\right)$  \State Calculate the gradient analog: $s^{\left( 0\right)} \coloneqq x^{\left( 0\right)}-p^{\left( 0\right)}$  \State $d^{\left( 0\right)} \coloneqq -s^{\left( 0\right)}$  \State Perform the line search with $d^{\left( 0\right)}$ with step size $\alpha^{\left( 0\right)}$ \State Update $x_1 \coloneqq x^{\left( 0\right)} + \alpha^{\left( 0\right)} d^{\left( 0\right)}$  \While{not converged}  \State $k += 1$ \State Calculate the gradient for $g$: $g^{\left( k\right)} \coloneqq \nabla g\left(x^{\left( k\right)}\right)$  \State Calculate the proximal value $p^{\left( k\right)} \coloneqq \text{prox}_{\rho,h}\left(x^{\left( k\right)}-\rho\cdot g^{\left( k\right)}\right)$  \State Calculate the gradient analog: $s^{\left( k\right)} \coloneqq x^{\left( k\right)}-p^{\left( k\right)}$  \State $d^{\left( k\right)} \coloneqq -s^{\left( k\right)}+\bar{\beta}^{\left( k\right)}\cdot d^{\left(k-1\right)}$  \State Perform the line search with $d^{\left( k\right)}$ with step size $\alpha^{\left( k\right)}$ based on Wolfe-Powell conditions \label{eqn:d-k} \State Update $x^{\left(k+1\right)} \coloneqq x^{\left( k\right)} + \alpha^{\left( k\right)} d^{\left( k\right)}$  \State Check for convergence  \EndWhile  \State \textbf{return} $p^{\left( k\right)}$  \end{algorithmic} 

\caption{Proximal Hager-Zhang \citep{Hager2005} Conjugate Gradient \label{alg:proximal-HZ-CG}}
\end{algorithm}

In Algorithm \ref{alg:proximal-HZ-CG}, Hager-Zhang's conjugate parameter
$\bar{\beta}^{\left(k\right)}$ is defined as \citep{Hager2005}: 
\begin{align*}
y^{\left(k\right)} & \coloneqq s^{\left(k+1\right)}-s^{\left(k\right)}\\
\beta^{\left(k\right)} & \coloneqq\frac{1}{\left\langle d^{\left(k\right)},y^{\left(k\right)}\right\rangle }\cdot\left\langle y^{\left(k\right)}-2\frac{\left\Vert y^{\left(k\right)}\right\Vert ^{2}}{\left\langle d^{\left(k\right)},y^{\left(k\right)}\right\rangle }d^{\left(k\right)},s^{\left(k+1\right)}\right\rangle \\
\eta^{\left(k\right)} & \coloneqq-\frac{1}{\left\Vert d^{\left(k\right)}\right\Vert \min\ \left\{ \eta,\left\Vert s^{\left(k\right)}\right\Vert \right\} }\\
\bar{\beta}^{\left(k\right)} & \coloneqq\max\ \left\{ \beta^{\left(k\right)},\eta^{\left(k\right)}\right\} 
\end{align*}

It was proven that if the line search step in Algorithm \ref{alg:proximal-HZ-CG}
satisfies Wolfe-Powell conditions and the gradient is globally Lipschitz,
Hager-Zhang conjugate gradient achieves global convergence finding a stationary point for
a smooth nonconvex objective function. In a dynamical system view, this corresponds to the global attraction property of the trajectory of the numerical algorithm to find equilibria for globally Lipschitz flows.
Lemma \ref{lem:moreau-Lipschitz-practical} implies that $\mathcal{G}_{\rho\tilde{f}}$,
or $s^{\left(k\right)}$ in Algorithm \ref{alg:proximal-HZ-CG}, are globally Lipschitz.
Thus, by Lemma \ref{lem:equivalent-clarke-stationary-point-pratical},
Algorithm \ref{alg:proximal-HZ-CG} yields the Clarke stationary point
of $f$. Based on the arguments in Section \ref{subsec:PCG-framework}, if the potential
function for $\mathcal{G}_{\rho\tilde{f}}$ is tractable to construct, the Wolfe-Powell line search in Algorithm \ref{alg:proximal-HZ-CG} can be carried out using the potential function of $\mathcal{G}_{\rho\tilde{f}}$ as the surrogate objective function; alternatively, an exact line search can be carried out by finding $\alpha$ that satisfies $\left\langle \mathcal{G}_{\rho\tilde{f}}\left(x^{\left(k\right)}+\alpha\cdot d^{\left(k\right)}\right),d^{\left(k\right)}\right\rangle =0$ --- such an exact line search can usually be carried out efficiently using Brent's method to find a root of a one-dimensional equation in $\mathbb{R}_{>0}$ \citep{Brent1971}. Furthermore, the descent property of $d^{\left(k\right)}$ was shown by \cite{Hager2005} independent of the line searches, which guarantees that $\left\langle \mathcal{G}_{\rho\tilde{f}}\left(x^{\left(k\right)}+\alpha\cdot d^{\left(k\right)}\right),d^{\left(k\right)}\right\rangle =0$ has a positive root. Moreover, another line search to ensure global convergence can be carried out by backtracking to find $\alpha^{\left(k\right)}$ satisfying 
\begin{equation}
    -\left\langle \mathcal{G}_{\rho\tilde{f}}\left(x^{\left(k\right)}+c_{1}\cdot\alpha^{\left(k\right)}d^{\left(k\right)}\right),d^{\left(k\right)}\right\rangle \geq c_{1}c_{2}\cdot\alpha^{\left(k\right)}\left\Vert d^{\left(k\right)}\right\Vert ^{2} \label{eq:feng-backtracking-line-search},
\end{equation}
where $c_{1},c_{2}\in\mathbb{R}_{>0}$ are constant to be chosen. When $\mathcal{G}_{\rho\tilde{f}}$ is pseudo-monotone in the sense of Karamardian \citep{Karamardian1976}, since the global Lipschitz property was established for $\mathcal{G}_{\rho\tilde{f}}$ in Lemma \ref{lem:moreau-Lipschitz-practical}, global convergence was proven for this backtracking line search method \citep{Feng2017}. We conclude this section with the observation that certain conjugate gradient methods obviate the need for line search procedures by determining the step size directly from $s^{(k)}$ and $d^{(k)}$, as exemplified in \citep{Chen2018}.

\section{\label{sec:optimizing-penalized-mle} Optimizing Algorithm and Prediction for Penalized $q$Gaussian Likelihood Problems }

\subsection{Problem Formulation }

Using the $q$Gaussian distribution to model the data will undoubtedly
enhance the robustness towards the underlying distributional assumption
and outliers. However, unlike the Gaussian distribution, two independent
$q$Gaussian random vectors are not jointly $q$Gaussian. Thus, we
take the following approach to model the data. Let $\mathbf{X}_{\text{train}}\in\mathbb{R}^{n_{\text{train}}\times\left(p+1\right)}$, $\ y_{\text{train}}\in\mathbb{R}^{n_{\text{train}}}$
denote the training design matrix and outcome, $\mathbf{X}_{\text{val}}\in\mathbb{R}^{n_{\text{val}}\times\left(p+1\right)},\ y_{\text{val}}\in\mathbb{R}^{n_{\text{val}}}$
denote the validation design matrix and outcome, and $\mathbf{X}_{\text{test}}\in\mathbb{R}^{n_{\text{test}}\times\left(p+1\right)},\ y_{\text{test}}\in\mathbb{R}^{n_{\text{test}}}$
denote the testing design matrix and outcome. Let 
\begin{equation}
\mathbf{X}\coloneqq\left[\mathbf{X}_{\text{train}}^{T},\mathbf{X}_{\text{val}}^{T},\mathbf{X}_{\text{test}}^{T}\right]^{T}\in\mathbb{R}^{n\times\left(p+1\right)}
\end{equation}
denote the design matrix for the entire dataset, and let 
\begin{equation}
y\coloneqq\left[y_{\text{train}}^{T},y_{\text{val}}^{T},y_{\text{test}}^{T}\right]^{T}\in\mathbb{R}^{n}
\end{equation}
denote the outcome for the entire dataset. Instead of assuming the $q$Gaussian
distribution for the training, validation and testing set separately,
we assume that
\begin{equation}
y\sim q\text{Gaussian}\left(q,\mathbf{X}\theta,\Sigma\right),
\end{equation}
where $\theta\in\mathbb{R}^{p+1}$ denotes the coefficients for regression,
and $\Sigma$ denotes the characteristic/scale matrix for the entire
data. Clearly, 
\begin{align*}
\mathbf{X}_{\text{train}} & =\left[I_{n_{\text{train}}\times n_{\text{train}}},0_{n_{\text{train}}\times n_{\text{val}}},0_{n_{\text{train}}\times n_{\text{test}}}\right]\mathbf{X}\\
y_{\text{train}} & =\left[I_{n_{\text{train}}\times n_{\text{train}}},0_{n_{\text{train}}\times n_{\text{val}}},0_{n_{\text{train}}\times n_{\text{test}}}\right]y
\end{align*}
implies that 
\begin{equation}
y_{\text{train}}\sim q\text{Gaussian}\left(q_{\text{train}},\mathbf{X}_{\text{train}}\theta,\Sigma_{\text{train}}\right)\label{eq:$q$Gaussian-distributed}
\end{equation}
where by the linear mapping closeness property \ref{enu:linear-mapping-invariance},
\begin{equation}
\Sigma_{\text{train}}=\left[I_{n_{\text{train}}\times n_{\text{train}}},0_{n_{\text{train}}\times n_{\text{val}}},0_{n_{\text{train}}\times n_{\text{test}}}\right]\Sigma\left[I_{n_{\text{train}}\times n_{\text{train}}},0_{n_{\text{train}}\times n_{\text{val}}},0_{n_{\text{train}}\times n_{\text{test}}}\right]^{T}
\end{equation}
is the $n_{\text{train}}\times n_{\text{train}}$ block diagonal matrix
of $\Sigma$ corresponding to the training data. By \eqref{enu:linear-mapping-invariance},
\begin{equation}
\frac{2}{1-q_{\text{train}}^{-1}}-n_{\text{train}}=\frac{2}{1-q^{-1}}-n,
\end{equation}
which implies that 
\begin{equation}
\frac{1}{q_{\text{train}}-1}-n_{\text{train}}=\frac{1}{q-1}-n.\label{eq:recover-q}
\end{equation}
\eqref{eq:recover-q} allows us to recover $q$ from the training
procedure. Above formulas for the training data and parameters can
trivially be applied to the validation and the testing data and parameters;
thus, validation and test can carried out easily from the model build
from the training data. 

For $q-$correlated data, often times, the $q-$correlation structure
is inferred or given prior to the model fitting; thus, we assume that
the $q-$correlation structure is given as $\Psi$ and we estimate the volatility / dispersion / scale parameter $\sigma^{2}>0$ such that
\begin{equation}
\Sigma=\sigma^{2}\Psi.
\end{equation}
Trivially, $\Psi_{\text{train}}$ is the block diagonal matrix of
$\Psi$ corresponding to the training data and $\Sigma_{\text{train}}=\sigma^{2}\Psi_{\text{train}}$. 

We are now ready to formulate our likelihood loss function. To utilize
$q$Gaussian distribution to model the $q-$correlated observations, we estimate the value of $q$ such that $q$ is allowed to
vary, and the model will thus be more robust towards a wide class of distributions.
Therefore, we choose to build the model using \eqref{eq:$q$Gaussian-heavy-tail-density},
since the dispersion matrix $\Lambda$ \eqref{eq:$q$Gaussian-simple-form}
depends on $q$. We formulate our maximization of our log-likelihood function
as the following from \eqref{eq:$q$Gaussian-distributed} and \eqref{eq:$q$Gaussian-heavy-tail-density}:
\begin{dgroup*}
  \begin{dmath*}
\underset{q_{\text{train}}\in(1,1+\frac{2}{n_{\text{train}}}),\theta\in\mathbb{R}^{p+1},\sigma^{2}\in\mathbb{R}_{>0}}{\arg\max\ }  \log(\frac{1}{\left|\sigma^{2}\Psi_{\text{train}}\right|^{\nicefrac{1}{2}}}\cdot\frac{\Gamma\left(\frac{1}{q_{\text{train}}-1}\right)}{\Gamma\left(\frac{1}{q_{\text{train}}-1}-\frac{n_{\text{train}}}{2}\right)}\cdot\left(\frac{2}{q_{\text{train}}-1}-n_{\text{train}}\right)^{-\frac{n_{\text{train}}}{2}}
  \cdot\left(1+\left(\frac{2}{q_{\text{train}}-1}-n_{\text{train}}\right)^{-1}\cdot\left\langle y_{\text{train}}-\mathbf{X}_{\text{train}}\theta,\left(\sigma^{2}\Psi_{\text{train}}\right)^{-1}\left(y_{\text{train}}-\mathbf{X}_{\text{train}}\theta\right)\right\rangle \right)^{\frac{1}{1-q_{\text{train}}}}).
\end{dmath*}
\end{dgroup*}
To address the high--dimensional data concerns, Oracle penalties are
incorporated to carry out variable selection. To penalize the log-likelihood
loss function to achieve variable selection, we formulate the following
problem: 
\begin{dgroup*}
  \begin{dmath*}
\underset{q_{\text{train}}\in(1,1+\frac{2}{n_{\text{train}}}),\theta\in\mathbb{R}^{p+1},\sigma^{2}\in\mathbb{R}_{>0}}{\arg\min\ }-\log(\frac{1}{\left|\sigma^{2}\Psi_{\text{train}}\right|^{\nicefrac{1}{2}}}\cdot\frac{\Gamma\left(\frac{1}{q_{\text{train}}-1}\right)}{\Gamma\left(\frac{1}{q_{\text{train}}-1}-\frac{n_{\text{train}}}{2}\right)}\cdot\left(\frac{2}{q_{\text{train}}-1}-n_{\text{train}}\right)^{-\frac{n_{\text{train}}}{2}}
 \qquad\cdot\left(1+\left(\frac{2}{q_{\text{train}}-1}-n_{\text{train}}\right)^{-1}\cdot\sigma^{-2}\cdot\left(\left\langle y_{\text{train}}-\mathbf{X}_{\text{train}}\theta,\Psi_{\text{train}}^{-1}\left(y_{\text{train}}-\mathbf{X}_{\text{train}}\theta\right)\right\rangle +2n_{\text{train}}\sum_{j=2}^{p+1}w\left(\theta_{j}\right)\right)\right)^{\frac{1}{1-q_{\text{train}}}})\nonumber 
 \end{dmath*}
\begin{dmath}
\Leftrightarrow  \underset{q_{\text{train}}\in(1,1+\frac{2}{n_{\text{train}}}),\theta\in\mathbb{R}^{p+1},\sigma^{2}\in\mathbb{R}_{>0}}{\arg\min\ }\frac{n}{2}\log\sigma^{2}-\log\Gamma\left(\frac{1}{q_{\text{train}}-1}\right)+\log\Gamma\left(\frac{1}{q_{\text{train}}-1}-\frac{n_{\text{train}}}{2}\right)+\frac{n_{\text{train}}}{2}\log\left(\frac{2}{q_{\text{train}}-1}-n_{\text{train}}\right)
 \quad+\frac{1}{q_{\text{train}}-1}\log\left(1+\left(\frac{2}{q_{\text{train}}-1}-n_{\text{train}}\right)^{-1}\cdot\sigma^{-2}\cdot\left(\left\langle y_{\text{train}}-\mathbf{X}_{\text{train}}\theta,\Psi_{\text{train}}^{-1}\left(y_{\text{train}}-\mathbf{X}_{\text{train}}\theta\right)\right\rangle +2n_{\text{train}}\sum_{j=2}^{p+1}w\left(\theta_{j}\right)\right)\right)\label{eq:opt-problem}
\end{dmath}
\end{dgroup*}
In the above formulated problem, $w$ is the Oracle penalty function,
and we are not to penalize the intercept term. The $2n_{\text{train}}$
multiplier is to ensure that the penalization effect is consistent with
the number of training observations. We choose to put the penalty
term together with the quadratic term without the variance scale parameter
$\sigma^{2}$ for two reasons: first, the optimization problem is
more tractable under such problem formulation; second, we do not wish
to let the value of $\sigma^{2}$ perturb the degree of penalization.
Comparing to penalized the log-likelihood directly, we choose to penalize
the quadratic component directly as it is more tractable. It was shown
that doing so will preserve Oracle properties \citep{Nikolova2000}
of penalized estimators. 

For the optimization procedure, we will proceed in a blockwise manner;
i.e., we will optimize $q_{\text{train}},\theta,\sigma^{2}$ separately
in each iteration. More details will be given in the following
subsections. 

\subsection{Minimizing with respect to $q_{\text{train}}$ and $\sigma^{2}$ }

With all the other parameters fixed, the sub-problem to minimize with
respect to $\sigma^{2}$ is 
\begin{align}
\underset{\sigma^{2}\in\mathbb{R}_{>0}}{\arg\min\ } & \frac{n}{2}\log\sigma^{2}+\frac{1}{q_{\text{train}}-1}\log(1+\left(\frac{2}{q_{\text{train}}-1}-n_{\text{train}}\right)^{-1}\cdot\sigma^{-2}\label{eq:subproblem-sigma2}\\
 & \quad\cdot\left(\left\langle y_{\text{train}}-\mathbf{X}_{\text{train}}\theta,\Psi_{\text{train}}^{-1}\left(y_{\text{train}}-\mathbf{X}_{\text{train}}\theta\right)\right\rangle +2n_{\text{train}}\sum_{j=2}^{p+1}w\left(\theta_{j}\right)\right))\nonumber 
\end{align}
which has a smooth objective function with respect to $\sigma^{2}$.
The first-order optimality condition 
\begin{align*}
\frac{n}{2}= & \frac{1}{q_{\text{train}}-1}\\
 & \cdot\frac{\left(\frac{2}{q_{\text{train}}-1}-n_{\text{train}}\right)^{-1}\cdot\left(\left\langle y_{\text{train}}-\mathbf{X}_{\text{train}}\theta,\Psi_{\text{train}}^{-1}\left(y_{\text{train}}-\mathbf{X}_{\text{train}}\theta\right)\right\rangle +2n_{\text{train}}\sum_{j=2}^{p+1}w\left(\theta_{j}\right)\right)}{\sigma^{2}+\left(\frac{2}{q_{\text{train}}-1}-n_{\text{train}}\right)^{-1}\cdot\left(\left\langle y_{\text{train}}-\mathbf{X}_{\text{train}}\theta,\Psi_{\text{train}}^{-1}\left(y_{\text{train}}-\mathbf{X}_{\text{train}}\theta\right)\right\rangle +2n_{\text{train}}\sum_{j=2}^{p+1}w\left(\theta_{j}\right)\right)}
\end{align*}
implies that the optimal value for the subproblem \eqref{eq:subproblem-sigma2} takes
minimizer 
\begin{align}
\overline{\sigma^{2}} & =\left(\frac{1}{q_{\text{train}}-1}/\frac{n}{2}-1\right)\cdot\left(\frac{2}{q_{\text{train}}-1}-n_{\text{train}}\right)^{-1}\label{eq:sigma2-closed-form-minimizer}\\
 & \quad\cdot\left(\left\langle y_{\text{train}}-\mathbf{X}_{\text{train}}\theta,\Psi_{\text{train}}^{-1}\left(y_{\text{train}}-\mathbf{X}_{\text{train}}\theta\right)\right\rangle +2n_{\text{train}}\sum_{j=2}^{p+1}w\left(\theta_{j}\right)\right)>0,\nonumber 
\end{align}
which is feasible. The feasible set for $q_{\text{train}}$ is $\left(1,1+\frac{2}{n_{\text{train}}}\right)$,
in this view, when $n_{\text{train}}$ is large, the numerical stability
will be an issue if minimization is carried out with respect
to $q_{\text{train}}$ directly. Thus, we choose to minimize with
respect to $\frac{1}{q_{\text{train}}-1}\in\left(\frac{n_{\text{train}}}{2},\infty\right)$. 

First of all, we are to prove that such minimization is feasible. 
\begin{lem}
The objective function \eqref{eq:opt-problem} has a local minimizer
in $\left(\frac{n_{\text{train}}}{2},\infty\right)$ with respect
to $\frac{1}{q_{\text{train}}-1}$. 
\end{lem}

\begin{proof}
Since the objective function \eqref{eq:opt-problem} is continuous
and smooth with respect to $\frac{1}{q_{\text{train}}-1}$, we only
need to analyze the derivative when $\frac{1}{q_{\text{train}}-1}\searrow0$
and $\frac{1}{q_{\text{train}}-1}\rightarrow\infty$.

$\frac{1}{q_{\text{train}}-1}\rightarrow\infty$:

Stirling's formula states that 
\begin{equation}
\lim_{x\rightarrow\infty}\frac{\Gamma\left(x\right)}{\sqrt{\frac{2\pi}{x}}\left(\frac{x}{e}\right)^{x}\left(1+O\left(x^{-1}\right)\right)}=1.
\end{equation}
Thus, 
\begin{equation}
\lim_{\frac{1}{q_{\text{train}}-1}\rightarrow\infty}\frac{\Gamma\left(\frac{1}{q_{\text{train}}-1}\right)}{\Gamma\left(\frac{1}{q_{\text{train}}-1}-\frac{n_{\text{train}}}{2}\right)}\cdot\left(\frac{2}{q_{\text{train}}-1}-n_{\text{train}}\right)^{-\frac{n_{\text{train}}}{2}}=1
\end{equation}
then 
\begin{equation}
\lim_{\frac{1}{q_{\text{train}}-1}\rightarrow\infty}-\log\left(\frac{\Gamma\left(\frac{1}{q_{\text{train}}-1}\right)}{\Gamma\left(\frac{1}{q_{\text{train}}-1}-\frac{n_{\text{train}}}{2}\right)}\cdot\left(\frac{2}{q_{\text{train}}-1}-n_{\text{train}}\right)^{-\frac{n_{\text{train}}}{2}}\right)=0.
\end{equation}
We also have 
\begin{align*}
 & \frac{1}{q_{\text{train}}-1}\log(1+\left(\frac{2}{q_{\text{train}}-1}-n_{\text{train}}\right)^{-1}\cdot\sigma^{-2}\\
 & \quad\cdot\left(\left\langle y_{\text{train}}-\mathbf{X}_{\text{train}}\theta,\Psi_{\text{train}}^{-1}\left(y_{\text{train}}-\mathbf{X}_{\text{train}}\theta\right)\right\rangle +2n_{\text{train}}\sum_{j=2}^{p+1}w\left(\theta_{j}\right)\right))\\
 & =O\left(\left(\frac{1}{q_{\text{train}}-1}\right)/\log\left(\frac{1}{q_{\text{train}}-1}\right)\right),
\end{align*}
which implies that this term will goes to infinity as $\frac{1}{q_{\text{train}}-1}\rightarrow\infty$.
Thus, the objective function \eqref{eq:opt-problem} goes to infinity
as $\frac{1}{q_{\text{train}}-1}\rightarrow\infty$. 

$\frac{1}{q_{\text{train}}-1}\searrow\frac{n_{\text{train}}}{2}$: 

Since $\Gamma\left(\frac{1}{q_{\text{train}}-1}-\frac{n_{\text{train}}}{2}\right)\rightarrow\infty$
as $\frac{1}{q_{\text{train}}-1}\searrow\frac{n_{\text{train}}}{2}$. 

The penalized log-likelihood involving $\frac{1}{q_{\text{train}}-1}$
can be simplified as 
\begin{align*}
 & -\log\frac{\Gamma\left(\frac{1}{q_{\text{train}}-1}\right)}{\Gamma\left(\frac{1}{q_{\text{train}}-1}-\frac{n_{\text{train}}}{2}\right)}\cdot(\left(\frac{2}{q_{\text{train}}-1}-n_{\text{train}}\right)\\
 & \qquad+\left\langle y_{\text{train}}-\mathbf{X}_{\text{train}}\theta,\left(\sigma^{2}\Psi_{\text{train}}\right)^{-1}\left(y_{\text{train}}-\mathbf{X}_{\text{train}}\theta\right)\right\rangle +2n_{\text{train}}\sum_{j=2}^{p+1}w\left(\theta_{j}\right))^{\frac{1}{1-q_{\text{train}}}}\rightarrow\infty
\end{align*}
as $\frac{1}{q_{\text{train}}-1}\searrow\frac{n_{\text{train}}}{2}$.
Thus, the subproblem to minimize with respect to $\frac{1}{q_{\text{train}}-1}$
is coercive on $\left(\frac{n_{\text{train}}}{2},\infty\right)$.
Coercivity implies that any minimizing sequence $\left\{ \left(\frac{1}{q_{\text{train}}-1}\right)_{j}\right\} $
must be contained within a bounded subset of $\left(\frac{n_{\text{train}}}{2},\infty\right)$.
Thus, Bolzano--Weierstrass theorem implies the existence of a convergent
subsequence. Let $\left\{ \left(\frac{1}{q_{\text{train}}-1}\right)_{j_{k}}\right\} $
be one such subsequence, and let $\overline{\left(\frac{1}{q_{\text{train}}-1}\right)}$
be its limit. Since the subproblem has a continuous objective function
with respect to $\frac{1}{q_{\text{train}}-1}$, the objective function
is lower--semicontinuous and the value of the objective function at
$\overline{\left(\frac{1}{q_{\text{train}}-1}\right)}$ is less than or equal to
the value of the objective function at $\left(\frac{1}{q_{\text{train}}-1}\right)_{j_{k}}$
for all $k=1,2,\dots,\infty$. Thus, since $\left\{ \left(\frac{1}{q_{\text{train}}-1}\right)_{j}\right\} $
is a minimizing sequence, the value of the objective function at $\overline{\left(\frac{1}{q_{\text{train}}-1}\right)}$
is less than or equal to the infimum of the objective function on
$\left(\frac{n_{\text{train}}}{2},\infty\right)$. Hence, since the
entire minimizing sequence is contained in $\left(\frac{n_{\text{train}}}{2},\infty\right)$,
$\overline{\left(\frac{1}{q_{\text{train}}-1}\right)}\in\left(\frac{n_{\text{train}}}{2},\infty\right)$
solves the subproblem of minimizing with respect to $\frac{1}{q_{\text{train}}-1}$. 

Unlike $\sigma^{2}$, the minimizer for $\frac{1}{q_{\text{train}}-1}$
is not in closed form, and the evaluation of the derivative with respect
to $\frac{1}{q_{\text{train}}-1}$ can not be carried out efficiently.
Thus, we apply Brent's line-search method to optimize the $\frac{1}{q_{\text{train}}-1}$
subproblem \citep{Brent1971}. 
\end{proof}

\subsection{Minimizing with respect to $\theta$ \label{subsec:optimization-motivation}}

Minimizing with respect to $\theta$ involves a nonconvex smooth function
and a convex nonsmooth function, which is termed a composite problem. In 
Section \ref{subsec:Proximal-HZ}, we developed a proximal conjugate
gradient algorithm for such composite optimization. 

In this part, we will establish an important remark regarding the
Oracle penalty. The subproblem we are to minimize with respect to
$\theta$ is 
\begin{align}
 & \underset{\theta\in\mathbb{R}^{p+1}}{\arg\min\ }\left\langle y_{\text{train}}-\mathbf{X}_{\text{train}}\theta,\Psi_{\text{train}}^{-1}\left(y_{\text{train}}-\mathbf{X}_{\text{train}}\theta\right)\right\rangle +2n_{\text{train}}\sum_{j=2}^{p+1}w\left(\theta_{j}\right)\nonumber \\
\Leftrightarrow & \underset{\theta\in\mathbb{R}^{p+1}}{\arg\min\ }\frac{1}{2n_{\text{train}}}\left\langle y_{\text{train}}-\mathbf{X}_{\text{train}}\theta,\Psi_{\text{train}}^{-1}\left(y_{\text{train}}-\mathbf{X}_{\text{train}}\theta\right)\right\rangle +\sum_{j=2}^{p+1}w\left(\theta_{j}\right)\nonumber \\
\Leftrightarrow & \underset{\theta\in\mathbb{R}^{p+1}}{\arg\min\ }\frac{1}{2n_{\text{train}}}\left\langle \theta,\mathbf{X}_{\text{train}}^{T}\Psi_{\text{train}}^{-1}\mathbf{X}_{\text{train}}\theta\right\rangle -2\left\langle y_{\text{train}},\Psi_{\text{train}}^{-1}\mathbf{X}_{\text{train}}\theta\right\rangle +\sum_{j=2}^{p+1}w\left(\theta_{j}\right)\label{eq:central-trend-subproblem}
\end{align}
$w$ can be chosen as oracle penalties such as SCAD/MCP penalties.
And it has been shown that both SCAD / MCP penalties admit a difference-of-convex
decomposition to a first-order smooth concave term plus $\lambda$
times $\ell_{1}$ penalty. The quadratic loss function is clearly
convex and smooth. This justifies our assumption for the objective
function. To carry out the proximal Hager-Zhang conjugate gradient
method proposed in Section \ref{subsec:Proximal-HZ}, we need to calculate
$L_{\nabla g}$, the $L-$smoothness constant for the smooth component.
Previous work suggests $L_{\nabla g}=\max\ \left\{ \text{max eigenvalue of }\frac{1}{n_{\text{train}}}\mathbf{X}_{\text{train}}^{T}\Psi_{\text{train}}^{-1}\mathbf{X}_{\text{train}},c_{\text{penalty}}\right\} $,
where $c_{\text{penalty}}$ is the $L-$smoothness constant for the
smooth component of the penalty, which will be $\frac{1}{a-1}$ for
SCAD and $\frac{1}{\gamma}$ for MCP \citep{Yang2024}. 
\begin{rem}
For high dimensional data, often times, the number of covariates exceeds
the number of observations; i.e, $\text{null}\left(\mathbf{X}_{\text{train}}\right)\neq\emptyset$.
Both SCAD/MCP penalties take constant values in $B_{\infty}\left(0,c\right)$\footnote{$B_{\infty}\left(0,c\right)$ denotes the open ball in uniform norm
\emph{in the corresponding space}, centered at the origin with radius
$c$. }; where $c=a\lambda$ for SCAD and $c=\gamma\lambda$ for MCP. Given
any stationary point $\bar{\theta}$ in the nonempty solution set
defined by $\mathbf{X}_{\text{train}}^{T}\Psi_{\text{train}}^{-1}\mathbf{X}_{\text{train}}-\mathbf{X}_{\text{train}}^{T}y_{\text{train}}=0$.
For the set $\bar{\theta}+\text{null}\left(\mathbf{X}_{\text{train}}\right)\setminus B_{\infty}\left(0,c\right)$,
each point in the relative interior (which is nonempty) of this set
is a Clarke stationary point, which implies that any algorithm with
a starting point in this set will converge in $0$ steps. This might
pose an issue for signal recovery, since $\text{null}\left(\mathbf{X}_{\text{train}}\right)$
is a vector subspace and some points can be very far from the origin. 
\end{rem}

\begin{rem}
In view of the subproblem with respect to $\theta$, it is
trivial that the minimizer for \eqref{eq:central-trend-subproblem}
does not depend on the other parameters, which are $q$ and $\sigma^{2}$.
Since the $q$Gaussian distribution is a generalization for all bell curve distributions, the estimation of the central trend
using the maximum likelihood principle for bell curve distributions is
equivalent to minimize a quadratic function, which has a breakdown point
of $0$. 

Taking into account the optimization subproblem with respect to $\theta$, it is evident that
the solution to \eqref{eq:central-trend-subproblem} remains unaffected
by the other parameters, namely $q$ and $\sigma^{2}$. Given that the $q$Gaussian distribution
extends the framework of bell curve distributions, the problem \eqref{eq:central-trend-subproblem} implies that estimating the central trend through
the maximum likelihood principle for all bell curve distributions is equivalent
to minimizing a quadratic function of the central trend, therefore characterized by a breakdown point
of $0$. 
\end{rem}

\subsection{Prediction for $y_{\text{test}}$ }

To show how prediction can be made, we will show the methods to
predict $y_{\text{test}}$ using the trained model in this subsection.
The same method applies for validation when predictions on $y_{\text{val}}$
are needed or to predict any new data based on the trained model. Since
the data are mutually $q$Gaussian, \eqref{eq:recover-q} implies
that 
\begin{equation}
\frac{1}{q_{\text{train}}-1}-n_{\text{train}}=\frac{1}{q_{\text{val}}-1}-n_{\text{val}}=\frac{1}{q_{\text{test}}-1}-n_{\text{test}}=\frac{1}{q-1}-n,
\end{equation}
which will be used to recover the value of the shape parameter $q_{\text{val}},q_{\text{test}}$.
Note that when $n_{\text{new}}$ data points are introduced, the total
number of observations $n$ changes from $n_{\text{train}}+n_{\text{val}}+n_{\text{test}}$
to $n_{\text{train}}+n_{\text{val}}+n_{\text{test}}+n_{\text{new}}$,
thus, the value of $q_{\cdot}$ will change for the entire dataset.
However, $q_{\text{train}}$ stays the same; thus, we suggest inferring the shape parameter for each data set based on $q_{\text{train}}$
directly using the equation above. With $q_{\cdot}$ calculated, it
is straightforward to estimate the $q-$variance-covariance matrix
$\mathbb{E}_{q}\left[\left(y_{\cdot}-\mathbf{X}\theta\right)\left(y_{\cdot}-\mathbf{X}\theta\right)^{T}\right]$
based on \eqref{eq:q-variance-covariance}; or, \emph{if existing,}
the variance-covariance matrix $\mathbb{E}\left[\left(y_{\cdot}-\mathbf{X}\theta\right)\left(y_{\cdot}-\mathbf{X}\theta\right)^{T}\right]$
based on \eqref{eq:variance-covariance}. 



\section{\label{sec:Conclusion-and-Discussion} Conclusion and Discussion }

This paper explores the field of statistical sparse learning, focusing
on modeling correlated data through the lens of maximizing Tsallis
entropy. It addresses the limitations inherent in the conventional
Gaussian distribution, notably its lack of robustness towards outliers
and underlying shape assumptions, by advocating for the $q$Gaussian
distribution. This distribution, derived from Tsallis entropy maximization,
represents a novel approach to handling correlated data and heterogeneity
--- elements frequently encountered in biostatistical contexts involving
genetic and longitudinal studies. 

This paper encompasses a re-derived probability density function
for the multivariate $q$Gaussian distribution based on Tsallis entropy
maximization. Statistical modeling based on the derived density paves the
way for the analysis of correlated data and heterogeneity and enables
variable selection. Furthermore, we have developed an innovative framework
capable of converting any numerical method, originally designed to
identify equilibria in flows, into a tool for tackling composite
optimization problems that are prevalent in statistical sparse learning.
By applying this framework to the Hager-Zhang conjugate gradient algorithm,
we have crafted an effective and stable algorithm tailored to the
challenges of sparse statistical learning. Given the abundance of
methods for numerically identifying equilibria for globally Lipschitz
flows, our approach significantly broadens the arsenal of techniques
available to address sparse statistical learning optimization challenges. 

In conclusion, our research positions the $q$Gaussian distribution,
underpinned by maximizing Tsallis entropy, as a robust and adaptable
alternative to Gaussian-based methodologies in statistical sparse
learning on correlated data. This breakthrough not only confronts
the traditional limitation of Gaussian assumptions, but also paves
the way for expanded investigation into Tsallis entropy-maximizing
distributions, particularly within the domain of biostatistics and
allied disciplines. 

Future directions for research include the exploration of the log-linear
model through the lens of Tsallis entropy maximization, akin to approaches
previously based on Shannon's entropy. Moreover, the study of the
phenomenon called \emph{volatility smirk} in financial return data may benefit
from employing the log-$q$Gaussian distribution ---- a transformation
of the $q$Gaussian distribution, which can provide deeper insights
into the nuances of financial markets. Additionally, in the field
of statistical computing research, our framework that transforms numerical
methods for identifying flow equilibria into algorithms for solving
composite optimization problems opens numerous avenues for future
research, especially in a sparse learning context.

\chapter{Discussion}\label{ch:discussion}





In the field of statistical analysis and supervised statistical machine learning applied to high--dimensional, extra-large datasets, which are prevalent in genetics and neuroimaging, the typical workflow initiates with the screening of variables to pinpoint those most relevant to the outcome, subsequently constructing the model based on these selected variables from the screening step. These procedures typically follow an unsupervised preprocessing step, such as genotyping quality control or Hardy-Weinberg equilibrium filtering in genetic data. For example, simulation studies and case studies utilizing preprocessed ABIDE data \citep{Cameron2013, Barry2020}, as presented in my first manuscript, exemplify this workflow in a real--world context. Within this structured approach of high--dimensional biostatistical analysis, the variable screening methods presented in my first manuscript emerge as a pivotal tool, effectively handling the nonlinear association between the outcome and covariates in the screening step. Furthermore, the $q$Gaussian modeling introduced in my third manuscript offers a flexible modeling framework that extends beyond the Gaussian assumption, enabling the adaptation of distributional shapes to encompass heavier tails, thus enhancing outlier accommodation and providing a more robust estimate of the volatility parameter. Furthermore, the optimization techniques developed in manuscripts 2 and 3 lay the groundwork for efficiently undertaking computationally intensive tasks from sparse learning in high--dimensional large datasets, thereby facilitating deeper insights from various statistical modeling approaches applied to high--dimensional biomedical data.


The mutual information-based screening tool developed in my first manuscript is adaptable for use in any high--dimensional dataset encountered in the field of biostatistics, including those in genetics and neuroimaging. This method's efficacy, as demonstrated in my first manuscript, stems from leveraging the computational speed using the Fast Fourier Transform (FFT), ensuring that the running time of FFT-based Kernel Density Estimation (KDE) remains competitive with alternative methods. It is important to underscore the inherent trade-offs in statistical analysis: between statistical efficiency and the breadth of underlying assumptions. Mutual information and copula methods, free from the constraints of linearity, excel with nonlinearly associated data for screening tasks but at the expense of statistical efficiency. Nonparametric methods like KDE, by relaxing distributional assumptions, similarly trade some statistical efficiency for robustness towards the underlying distributional assumptions. However, for univariate variable screening, where only two variables are considered at each iteration, this compromise on statistical efficiency is of less concern. Consider the idea behind the ``curse of dimensionality'': for two continuous variables, a dataset with as few as $900$ data points will allow $30$ data points per dimension -- following this rule of thumb, the $2D$ surface for the bivariate density can be estimated fairly well using the information encompassed in the dataset, making nonparametric estimation of the measure of association particularly suited for variable screening tasks. Furthermore, the existence of nonparametric methods for estimating copulas, as highlighted in \citep{Rabhi2019}, opens up intriguing avenues for future research. Specifically, comparing the efficacy of nonparametric copula estimation with that of nonparametric mutual information estimation in the context of variable screening presents a promising direction for further research.

The first manuscript focuses on variable screening using marginal association. In some cases, when certain variables exhibit significantly higher association measures than others, they should be selected. The number of variables to include is often based on external knowledge. Generally, the asymptotic distribution of the mutual information estimator can help guide the decision on the number of variables to include in the model. Therefore, further exploration of this topic presents an interesting direction for future research. Additionally, for variable selection based on a joint model, the number of covariates often depends on when it gives the best predictive performance. However, for ultra high--dimensional data, it is usually infeasible to find the set of variables that will give the best performance.

\begin{figure}[H]
\centering{}\includegraphics[width=0.5\textwidth]{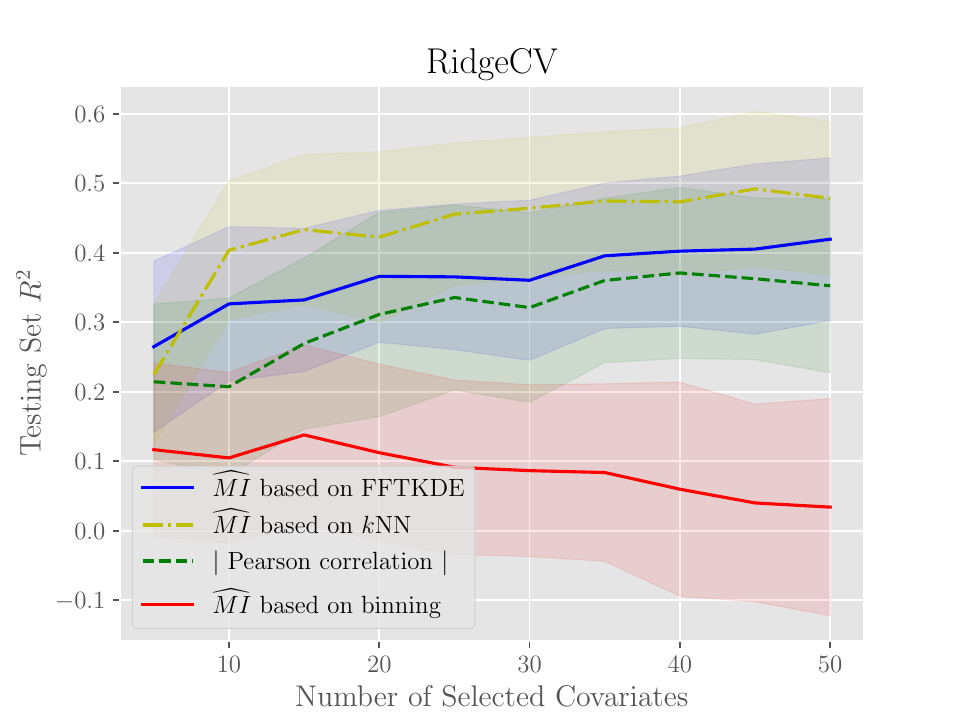}\includegraphics[width=0.5\textwidth]{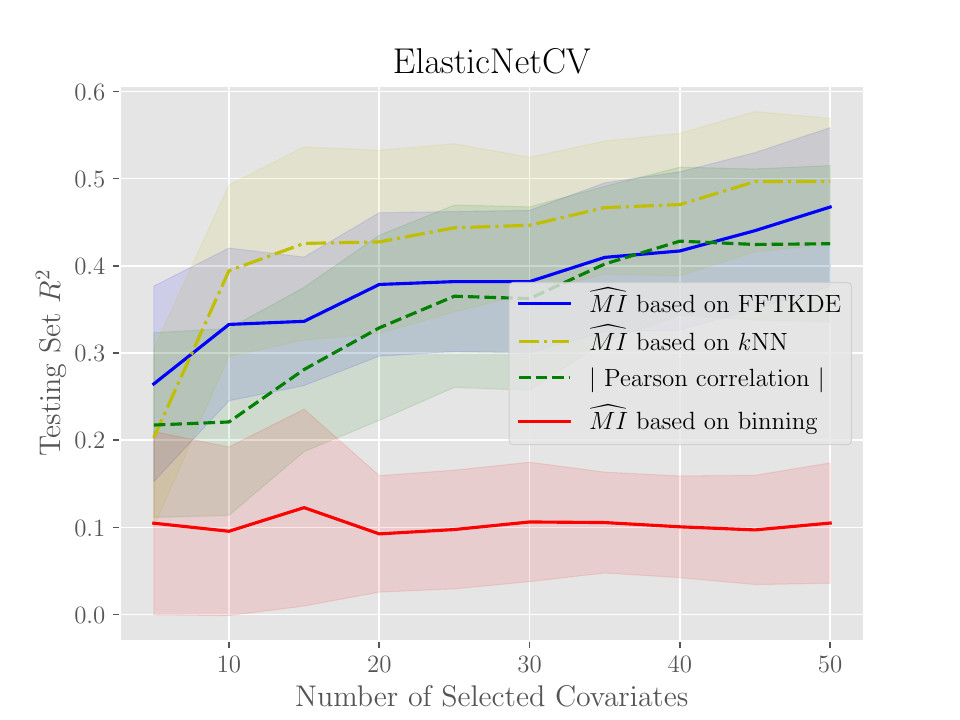} \\ \includegraphics[width=0.5\textwidth]{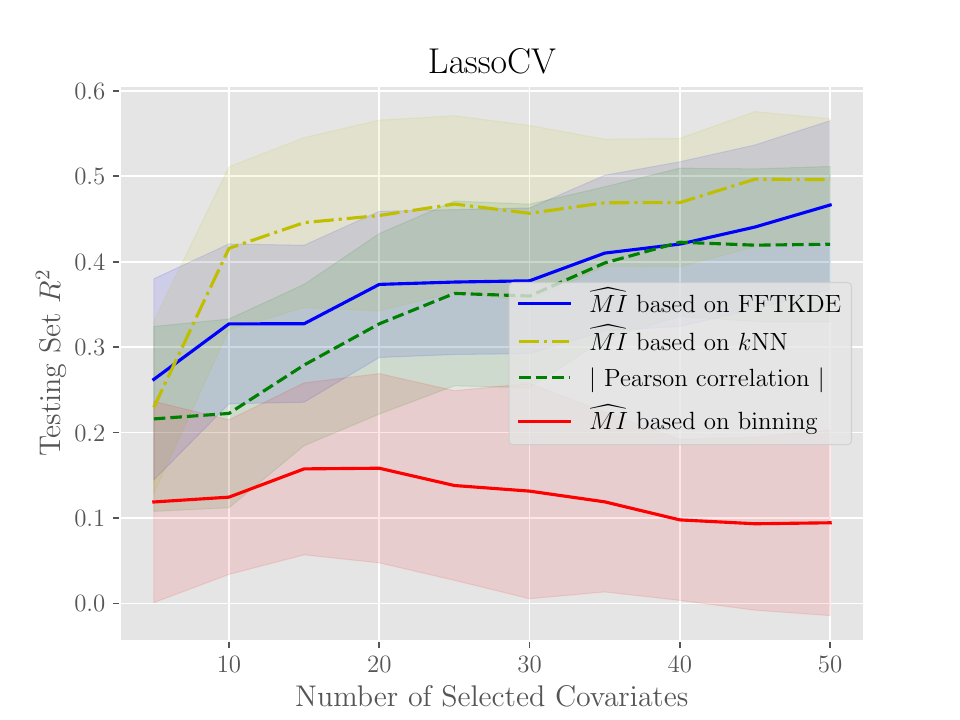}\includegraphics[width=0.5\textwidth]{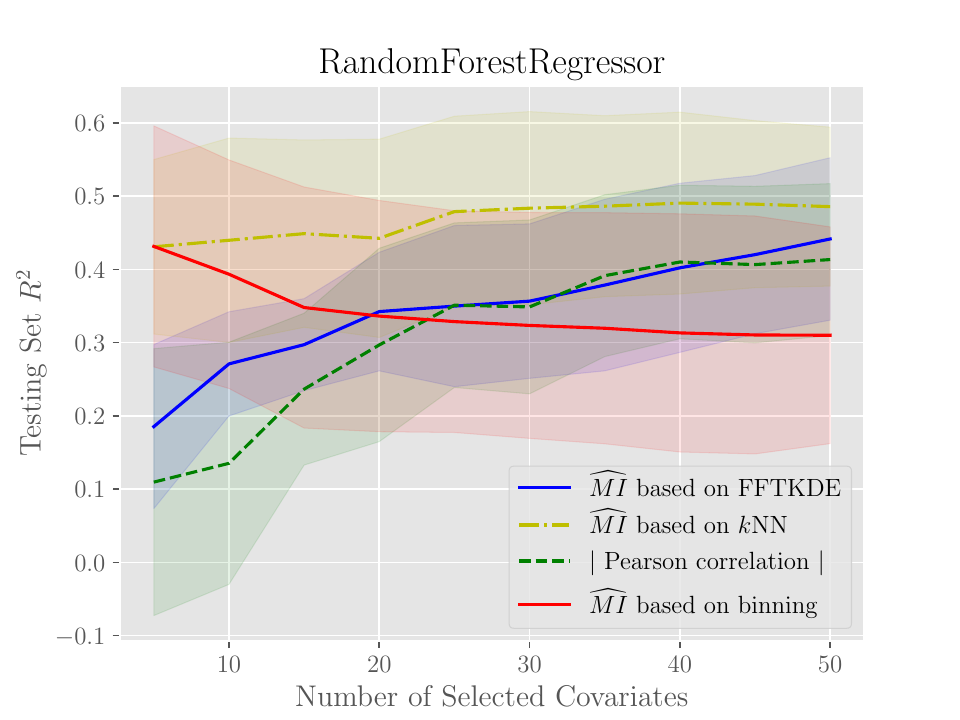}\caption{Testing Set $R^{2}$ for age at the scan outcome v.s. the number of
most associated brain imaging covariates based on the association
measure rankings. \emph{The most associated brain imaging covariates are then input to the spline transformer using Bernstein polynomial of degree $3$ to produce the data for model--fitting.} Means with their $95\%$ confidence intervals were plotted for $20$
simulation replications. \label{fig:poly3-case-Study-continuous}}
\end{figure}

\begin{figure}[H]
\centering{}\includegraphics[width=0.5\textwidth]{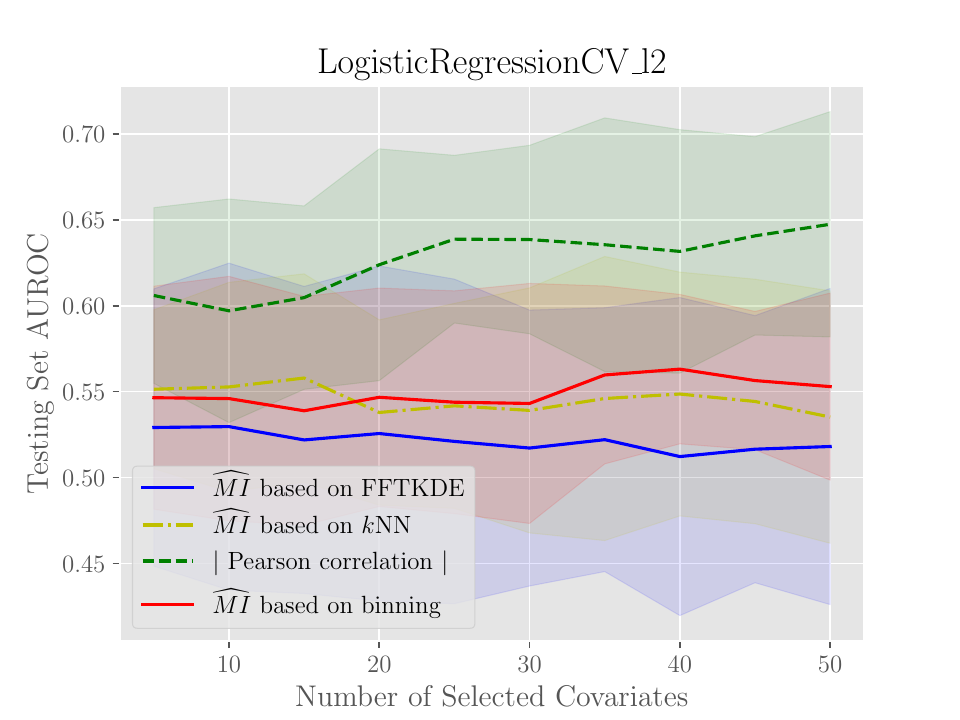}\includegraphics[width=0.5\textwidth]{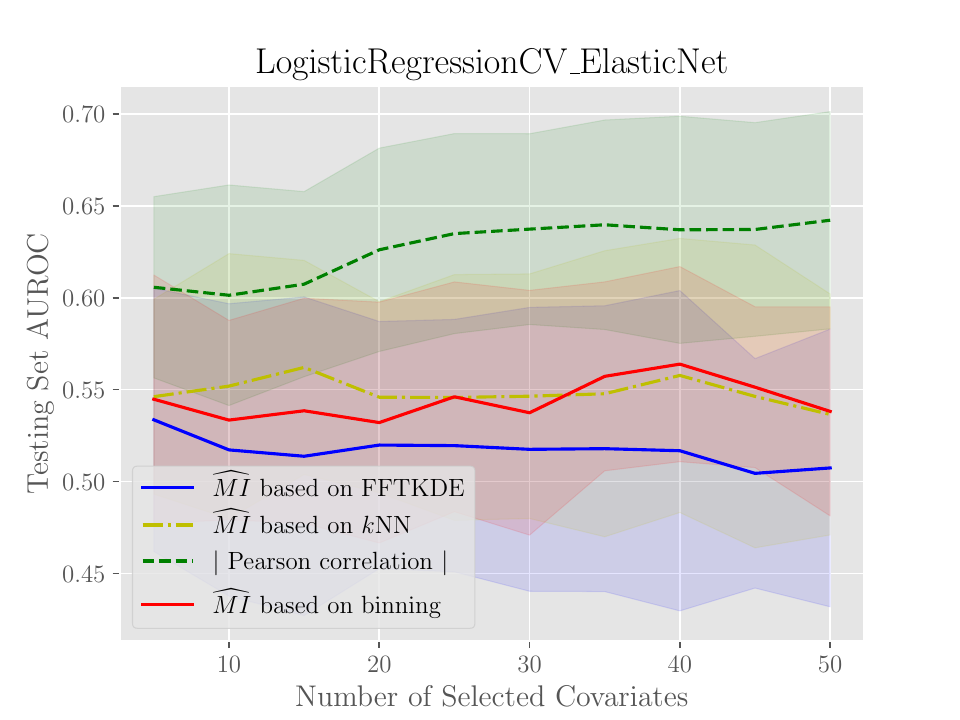}\\ \includegraphics[width=0.5\textwidth]{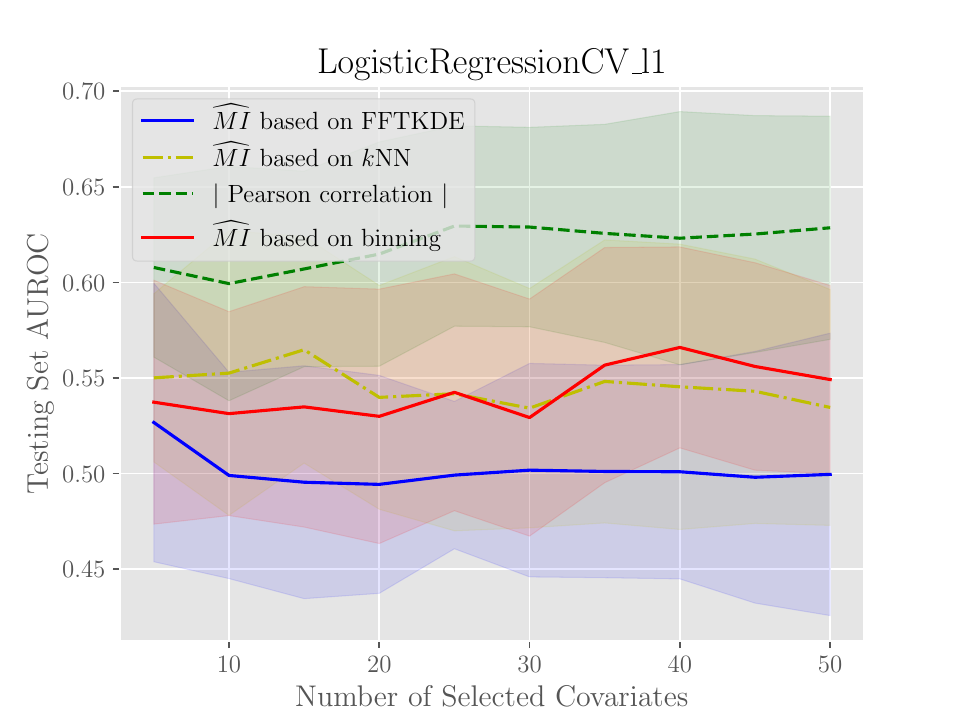}\includegraphics[width=0.5\textwidth]{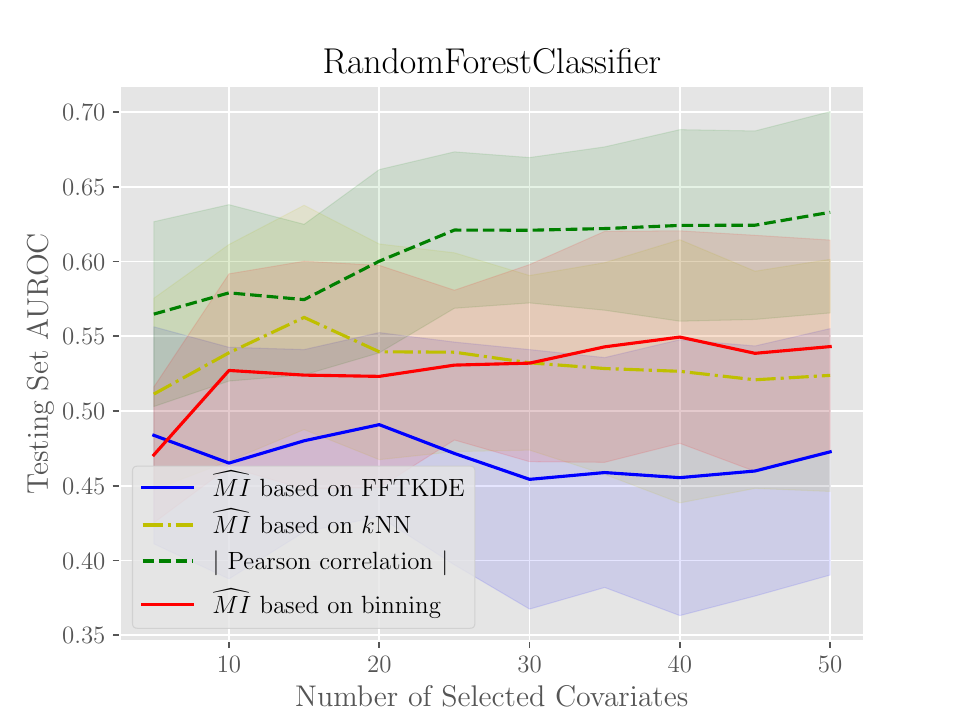}\caption{Testing Set AUROC for autism diagnosis outcome v.s. the number of
most associated brain imaging covariates based on the association
measure rankings. \emph{The most associated brain imaging covariates are then input to the spline transformer using Bernstein polynomial of degree $3$ to produce the data for model--fitting.} Means with their $95\%$ confidence intervals were plotted for $20$
simulation replications. \label{fig:poly3-case-Study-binary}}
\end{figure}

The case studies aimed at predicting age and autism diagnosis in my first manuscript predominantly utilized penalized (generalized) linear models, largely because these models are currently considered state-of-the-art. However, to address the potential limitations of linearity, I extended the scope of these studies by also fitting the models on the splines produced by Bernstein polynomials of degree $3$ on the selected covariates \citep{Racine2022} and repeating the same model fitting process. This spline transformation approach allows for a nuanced exploration beyond the underlying linearity assumptions of the (generalized) linear models. The results, illustrated in Figures \ref{fig:poly3-case-Study-continuous} and \ref{fig:poly3-case-Study-binary}, in fact, corroborate our observations from the first manuscript as shown in Figures \ref{fig:Case-Study-continuous} and \ref{fig:Case-Study-binary}.

In manuscripts 2 and 3, substantial advances are made in statistical computing, particularly in the optimization of composite problems that are commonly encountered in sparse learning -- encompassing scenarios like penalized least squares, robust objective functions (for example, Huber loss), log-likelihood, partial log-likelihood, the generalized method of moments (GMM), and more. A common characteristic of these (unpenalized) objective functions is their Lipschitz smoothness; hence, when penalized with sparse penalties, the smooth component of the penalized objective function will retain the Lipschitz smoothness. 
Consider, for example, different forms of generalized linear models (GLM) that are often used in the field of biostatistics, for the link function $g$, the general form of GLM is to model 

\begin{equation}
g\left(\mathbb{E}\left(\mathbf{y}\vert\mathbf{X}\right)\right)=\mathbf{X}\boldsymbol{\beta}.
\notag
\end{equation}

A result from the fact that all norms are equivalent in finite-dimensional spaces is that all linear operators mapping from a finite-dimensional normed linear space to any normed linear space are bounded. Thus, the finite-dimensional design matrix $\mathbf{X}$ clearly has a bounded operator norm, thus $\mathbf{X}\boldsymbol{\beta}$ is globally Lipschitz smooth with respect to $\boldsymbol{\beta}$. The link function, in fact, often satisfies $g\in \mathcal{C}^{2}$ and has a bounded second derivative, thus implies $g\in \mathcal{C}^{1,1}$. When the link function does not have a bounded second derivative over the Euclidean space such as exponential function for Poisson regression, restricting the regression on a compact set will make any locally Lipschitz-smooth optimization problem globally Lipschitz-smooth over the restricted compact set. As argued in our second manuscript, the vast majority of statistical learning problems can be considered as optimizing over a closed ball centered at the origin with a large but finite radius in practice. Subsequently, the log-likelihood function is usually globally Lipschitz-smooth with respect to the parameters of interest. Since continuity is invariant under function composition, the objective function is globally Lipschitz-smooth as long as differentiability conditions allow, which is mostly the case for (unpenalized) statistical loss functions. 

The adaptation to mixed-effects models, frequently used in biological data modeling and longitudinal studies, involves incorporating finite-dimensional design matrices for mixed effects. As discussed previously, finite-dimensional design matrices have a bounded operator norm, thus their linear mapping is globally Lipschitz-smooth. This fact guarantees the preservation of Lipschitz smoothness for the objective functions of the mixed effects model variants stemming from the objective function of the GLMs discussed before, enabling the optimization methods developed in manuscripts 2 and 3 to be effectively applied to them. Consequently, these algorithms are computationally efficient, and their first-order nature ensures a low memory consumption. This is vital for analysis of high-dimensional large biological datasets whose size by far exceeds the memory bottleneck. 

Sometimes, the challenge of deriving closed-form expressions for the gradient necessitates the use of numerical tools such as auto-differentiation. This approach is feasible and practical, thanks to advances in various auto-differentiation technologies due to rapid research in training deep neural networks. Nevertheless, conducting error analysis in the application of auto-differentiation, particularly for gradient calculations in first-order optimization algorithms, presents a rich area for further investigation. This exploration is also pertinent given the explosive growth in neural network research, where insights from numerical analysis can significantly contribute to advances in both statistical computing and deep learning.

An important part of manuscripts 2 and 3 is established based on Moreau envelope, also known as Moreau-Yosida regularization, which was originally established as a crucial concept within functional analysis in Hilbert spaces \citep{Moreau1965}, before it was recognized for its extensive applicability in optimization and variational analysis in finite-dimensional settings. This concept also facilitates a nuanced discussion on the trade-off between statistical efficiency and robustness towards outliers. Traditionally, the mean, minimizing the $L_2$ norm, is considered statistically efficient but vulnerable to outliers as it has a breakdown point of $0$; while the median, minimizing the $L_1$ norm, offers robustness towards outliers with a breakdown point of $0.5$ at the expense of statistical efficiency. The Huber M-estimator balances between mean and median, is the result of minimizing Huber loss function, defined by: 

\begin{equation}
    L_{\delta}\left(a\right)=\begin{cases}
\frac{1}{2}\theta^{2}, & \left|\theta\right|\leq\delta;\\
\delta\left(\left|\theta\right|-\frac{1}{2}\delta\right), & \left|\theta\right|>\delta.
\end{cases}
\notag
\end{equation}

\begin{figure}[H]
\centering{}\includegraphics[width=0.8\textwidth]{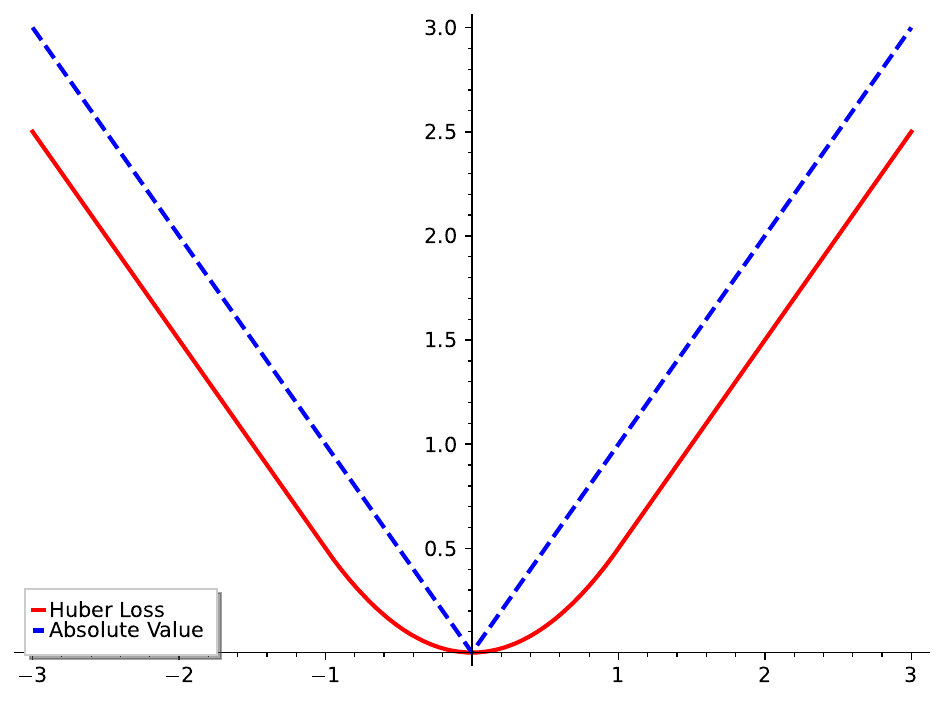}\caption{(scaled) Huber loss function and the absolute value function} \label{fig:huber-loss}
\end{figure}

Immediately, 

\begin{equation}
    \text{epi }\frac{L_{\delta}}{\delta}=\text{epi }\left|\cdot\right|+\text{epi }\frac{1}{2\delta}\left(\cdot\right)^{2}, 
    \notag
\end{equation}

which, based on our discussion of variational and nonsmooth analysis in the third manuscript, gives the (exact) infimal convolution equality  

\begin{equation}
    \frac{L_{\delta}}{\delta}=\left|\cdot\right|\square\frac{1}{2\delta}\left(\cdot\right)^{2}=M_{\delta}\left|\cdot\right|;
    \notag
\end{equation}

that is, the Huber loss function scaled by $\frac{1}{\delta}$ is the smoothing of the absolute value function by Moreau envelope parametrized by smoothing parameter $\delta$. Figure \ref{fig:huber-loss} visualizes how the scaled Huber loss function acts as Moreau envelope smoothing the absolute value function. This fact connects the famous Huber M-estimator, which maintains a balance between statistical efficiency and robustness to outliers for central trend estimation, to the famous Moreau envelope, which is the foundational work for our discussion of proximal methods in manuscripts 2 and 3. 

Further on the trade-off between statistical efficiency and robustness towards outliers and distributional assumption, the $q$Gaussian distribution from my third manuscript, when compared to the conventional Gaussian distribution, allows adjustment to the shape parameter. Recall that when $q\searrow 1$, $q$Gaussian random variable becomes Gaussian. The trade-off here is that the estimation of $q$ costs one degree of freedom. Nevertheless, we deem the trade-off of one degree of freedom usually a worthwhile exchange for the capability to account for the shape of the underlying distribution. 

Moreau envelope, explored in detail within manuscripts 2 and 3, represents just one of the many concepts from functional analysis and operator theory that extends the relevance of statistics beyond statistical computing. Dynamical system also plays a crucial role in a broader statistical landscape, transcending research in statistical computing. The principle of maximum likelihood estimation, which seeks to minimize the negative log-likelihood, is intricately linked to dynamical systems through the concept of observed Fisher information. Observed Fisher information measures the curvature of this minimizer to reflect its stability, which can be used to estimate the expected Fisher information, and provides a lower bound on the variance of any unbiased estimator as dictated by the \ac{CRLB}. To sum up in one sentence, the stability of the minimizer of the negative log-likelihood function directly reflects on the variance of any unbiased estimator. 

Much of the content of this thesis content in manuscripts 2 and 3 focuses on sparse learning, as known as variable selection, using sparse penalties. A challenge in this area of statistical learning involves the determination of hyper-parameters, such as those for sparse penalties in sparse learning or the bandwidth matrix for kernel density estimation. Conventionally, selecting these hyperparameters has relied on a data-driven methodology employing zero-order techniques such as grid search, supplemented by bootstrap validation or cross-validation. This approach, while effective, can be computationally intensive, particularly with large, high--dimensional datasets. Recent developments in implicit differentiation offer a promising avenue for a more computationally efficient choice of hyperparameters \citep{Blondel2022, Bertrand2020, Bertrand2022} -- this research suggests that leveraging implicit differentiation to speed up hyperparameter optimization presents a compelling direction for future research in statistical computing, particularly in the realms of sparse learning and variable selection. The work presented in this thesis revolves around nonconvex penalties, which result in non--unique local minimizers. This poses a significant challenge for implicit differentiation in bi--level optimization when tuning penalty hyper--parameters. Unlike \ac{LASSO}, where the set of minimizers can be proven to be a singleton under certain conditions, the presence of multiple local minimizers in nonconvex penalties complicates the adaptation of such methods. 

The third manuscript extensively builds on the concept of Tsallis entropy maximization, leading to the formulation of the $q$Gaussian likelihood. This approach enables the modeling of distributions with power-law decay for their tails, thus allowing heavier tails compared to the exponential decay observed in Gaussian distributions. Drawing parallels with Shannon entropy's application in solving likelihood equations, as noted by \cite{Calcagni2019}, leveraging the principle of entropy maximization often leads to the formulation of a dual problem alongside the primal likelihood maximization problem. This duality sometimes can reveal unique characteristics and computational strategies applicable to a wide array of models frequently used in biostatistics, for example, log-linear models. From this perspective, investigating the role of Tsallis entropy maximization in addressing likelihood maximization issues, particularly for cases related to the log-$q$Gaussian distribution, emerges as a compelling research pathway situated at the confluence of biostatistics and physics.

In addition to the discussion of proximal methods established based on Moreau envelope smoothing, detailed in manuscripts 2 and 3, there exists an alternative technique that applies smoothing directly to the nonsmooth objective function \citep{Chen2010}. This technique is akin to the concept of mollification, a term often encountered in discussions of partial differential equations or functional analysis. In manuscript 2, we showed that SCAD and MCP penalties can be represented in difference-of-convex form, combining a convex $\ell_1$ norm with a smooth, concave term -- refer to equations \eqref{eq:SCAD}, \eqref{eq:MCP}, \eqref{eq:SCAD-Lsmoothness}, and \eqref{eq:MCP-Lsmoothness} in manuscript 2 for an in-depth explanation. Each term, including the $\ell_1$ norm, undergoes independent smoothing if needed. Smoothing can be applied directly to both terms separately. To smooth out the $\ell_1$ norm term, since the absolute value of the term takes value $-1$ on $\mathbb{R}_{<0}$ and $1$ on $\mathbb{R}_{>0}$, the first-order derivative can be made continuous by any sigmoid function. For example: scaled and translated logistic function 
\begin{equation}
    \frac{1.0 \, e^{\left(\delta_{\mathit{moll}} \theta\right)} - 1.0}{1.0 \, e^{\left(\delta_{\mathit{moll}} \theta\right)} + 1.0},
    \notag
\end{equation}
with an integral serving as a smoothed out (``mollified'') version of $\ell_1$ penalty 
\begin{equation}
    -\frac{\delta_{\mathit{moll}} \theta + 2 \, \log\left(2\right) - 2 \, \log\left(e^{\left(\delta_{\mathit{moll}} \theta\right)} + 1\right)}{\delta_{\mathit{moll}}}.
    \notag
\end{equation}

\begin{figure}[H]
\centering{}\includegraphics[width=0.5\textwidth]{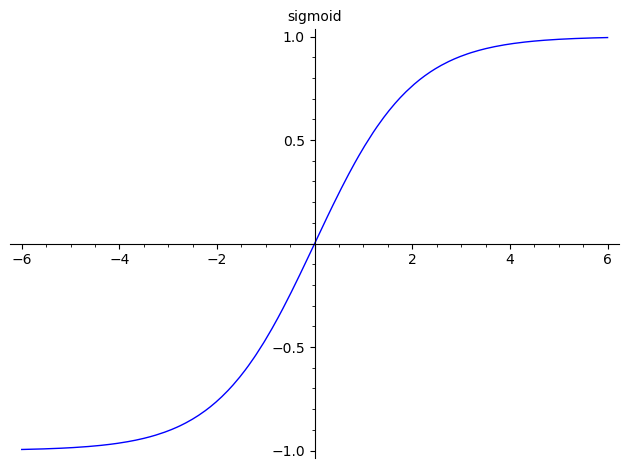}\includegraphics[width=0.5\textwidth]{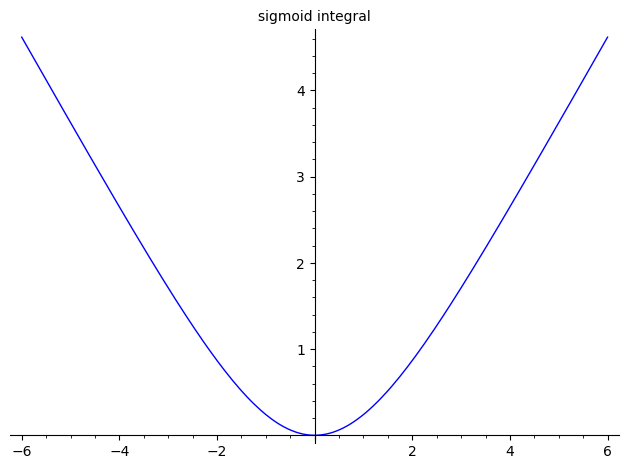}
\caption{Scaled and Translated Logistic Function to Make the Discontinuous Derivative Continuous and Its Integral to ``Mollify'' $\ell_1$} \label{fig:l1-mollification}
\end{figure}

In the above equations, $\delta_{moll}\geq 0$ denotes the smoothing parameter and $\theta$ denotes the penalized coefficient; $\delta_{moll}= 0$ recovers the nonsmooth $\ell_1$ penalty component. The figure shown in Figure \ref{fig:l1-mollification} illustrates the smoothing (referred to as ``mollification'') that is similar to, but not identical to, the role of a mollifier, which is defined as a function that is infinitely differentiable and has a compact support. A similar smoothing effect can be achieved for the $\ell_1$ function using other sigmoid functions, such as $\text{arc}\tan$, which ensure the derivative converges to $-1$ as $\theta\rightarrow -\infty$ and $1$ as $\theta\rightarrow \infty$. Notably, most sigmoid functions possess infinite-order smoothness, rendering the smoothed out $\ell_1$ infinitely differentiable. This approach to smoothing enables the application of smoothing techniques to nonsmooth objective functions \citep{Chen2010}.

For the concave term within the oracle penalty, already first-order smooth, smoothing to achieve second-order smoothness $\mathcal{C}^2$ is feasible. The equation for the smoothed MCP concave component serves as an illustration:
\begin{dmath*}
p_{MMCP,\lambda,\gamma,\delta_{MCP}}\left(\theta\right)=\begin{cases}
-\frac{\theta^{2}}{2\gamma}; & \left|\theta\right|<\gamma\lambda-\delta_{MCP}\\
-(\frac{\gamma^{3}\lambda^{3}-3\delta_{MCP}\gamma^{2}\lambda^{2}+3\delta_{MCP}^{2}\gamma\lambda-\delta_{MCP}^{3}+3\left(\gamma\lambda+\delta_{MCP}\right)\theta^{2}}{12\delta_{MCP}\gamma} & \gamma\lambda-\delta_{MCP}\leq\left|\theta\right|<\gamma\lambda+\delta_{MCP}\\
\qquad-\frac{\left(3\gamma^{2}\lambda^{2}-6\delta_{MCP}\gamma\lambda+3\delta_{MCP}^{2}+\theta^{2}\right)\left|\theta\right|}{12\delta_{MCP}\gamma});\\
-\lambda\left|\theta\right|+\frac{3\gamma^{2}\lambda^{2}+\delta_{MCP}^{2}}{6\gamma}; & \left|\theta\right|\geq\gamma\lambda+\delta_{MCP}
\end{cases}
\end{dmath*}

In these expressions, $\delta_{MCP}\geq 0$ denotes the smoothing parameter for the MCP concave component, and $\theta$ the penalized coefficient, with the remaining parameters as defined in manuscript 2, aligning with the original MCP formulation \citep{Zhang2010a}; $\delta_{MCP}= 0$ recovers the original first-order smooth MCP concave component. Smoothing the first-order smooth concave component further to achieve $\mathcal{C}^2$ or higher smoothness potentially allows for establishing $\mathcal{C}^2$-diffeomorphism, augmenting discussions in manuscript 3 on the dynamical system and optimization. This advancement could be pivotal for employing statistical computing algorithms within Morse theory, suggesting a unified framework for a broader spectrum of optimization problems discussed in this thesis. Mentioning this further smoothing of the smooth concave component, if leaving the $\ell_1$ nonsmooth component unaltered, preserves the Oracle property of the estimator \citep{Nikolova2000}.

Throughout this thesis, the Bayesian methodology for analyzing high--dimensional data was not extensively covered, primarily due to the significant computational cost of posterior calculations. Such computations tend to be significantly more resource-intensive than those required by frequentist or likelihood-based approaches, making them less practical for large, high--dimensional datasets. However, recent advances have discussed posterior computations in a Hilbert space setting \citep{RiutortMayol2022, Sprungk2017}. The work presented in manuscripts 2 and 3 primarily engages with concepts in a Euclidean space, which can be extended to an infinite-dimensional space. For example, the accelerated gradient technique discussed in the second manuscript is adaptable to a Hilbert space context, aligned with prior analyses of first-order optimization algorithms analyzed using a dynamical system reflecting a Hessian-driven damping mechanism \citep{Attouch2020}. Consequently, adapting optimization strategies to a Hilbert space setting may offer more efficient computational alternatives for Bayesian analysis, considering that posterior computations can be viewed as optimization problems within a function space.

\chapter{Conclusion}\label{ch:conc}





This thesis met its objectives by systematically addressing several challenges in statistical computing and modeling in the analysis of high-dimensional biological data that are frequently encountered in neuroimaging and genetics. Each manuscript within this thesis contributes to a cohesive workflow that enhances our ability to draw meaningful insights from complex datasets.

In the first manuscript, I introduced {\tt fastHDMI}, a Python package specifically designed for efficient variable screening within high--dimensional contexts. This tool is robust to nonlinear associations, which is essential for the statistical analysis of many datasets. The application of {\tt fastHDMI} to the preprocessed {\ac{ABIDE}} \citep{Cameron2013, Barry2020} dataset exemplifies its practical utility and transformative potential in real--world scenarios. This manuscript establishes the groundwork for the typical workflow in biostatistical analysis, which starts with variable screening to identify those variables most relevant to the outcome, pivotal for the subsequent modeling steps.

Building on this, the second manuscript advanced our computational capabilities by developing efficient statistical computing methods for sparse learning that utilize nonconvex penalties, addressing significant computational challenges. The manuscript's focus on optimizing hyperparameter settings based on complexity bounds significantly enhances the efficiency of statistical computing, particularly in handling large--scale high--dimensional data. This manuscript has been published on Statistics and Computing \citep{Yang2024}.

The third manuscript further refines our approach by introducing the $q$Gaussian linear mixed-effects model. This innovative model provides a robust alternative to conventional Gaussian models by accommodating broader distributional shapes and heavier tails. This advancement is critical for modeling correlated and heterogeneous data often encountered in biostatistics, such as in genetic and longitudinal studies, thus enhancing the robustness and flexibility of statistical analyses. In addition, the innovative framework for converting numerical methods for finding equilibria of dynamical systems into optimization algorithms for composite objective functions, often found in sparse penalized objectives, offers significant insights. As a result, various researches could leverage this using the proposed framework to adapt a numerical method for dynamical systems to a numerical algorithm for composite optimization problem that is prevalent in the statistical computation of sparse learning.

Collectively, these manuscripts create a comprehensive approach for robust and efficient analysis of high-dimensional data by seamlessly integrating solutions to variable screening robust to nonlinearity and distributional assumption, optimize objective functions with nonconvexity and nonsmoothness induced by Oracle sparse penalties, and model correlated data structures robust to distributional assumptions and heavy tails. The implications of this research are substantial, providing robust, scalable, and computationally efficient methodologies that improve our capacity to analyze and interpret high--dimensional large datasets. By improving the efficiency and robustness of these statistical learning processes, this thesis supports significant advancements in personalized medicine, enhances our understanding of complex genetic interactions and brain functions, and fosters the development of better diagnostic and therapeutic strategies. The methodologies developed here set a new standard in statistical computing for high-dimensional data analysis, paving the way for future research that will expand their applications in diverse fields in science and medicine.


\begin{appendices}
	\chapter{Appendix to Manuscript 1} 


\section{Methodology Consideration \label{apdx:method-contribution}}

For a function $f$ defined over an Euclidean space $\mathbb{R}^{n}$, its (continuous) Fourier transform is defined as 
\begin{equation}    \left(\mathcal{F}f\right)\left(\xi\right)\coloneqq\int_{\mathbb{R}^{n}}f\left(x\right)\exp\left(-2\pi i\cdot\left\langle x,\xi\right\rangle \right)dx,
\end{equation}
a linear operator. The Fourier series is then the synthesis formula. Consider a square-integrable function space $L^{2}\left(\left[-\pi,\pi\right]\right)$,
the fundamental results of Fourier analysis \citep{Stein2003} conclude that $\left\{ \phi_{k}\coloneqq\exp\left(ikx\right)\vert k\in\mathbb{Z}\right\} $
is an orthonormal and complete basis for this Hilbert space with the
inner product being defined by 
\begin{equation}
\forall f,g\in L^{2}\left(\left[-\pi,\pi\right]\right),\ \left\langle f,g\right\rangle \coloneqq\frac{1}{2\pi}\int_{-\pi}^{\pi}f\left(x\right)\bar{g}\left(x\right)dx.
\notag
\end{equation}
We remark the that inner product for a complex Hilbert space is linear for the first argument and anti-linear for the second argument. The Fourier series that represents any function $f\in L^{2}\left(\left[-\pi,\pi\right]\right)$
is then 
\begin{equation}
f=\sum_{k=-\infty}^{\infty}\left\langle f,\phi_{k}\right\rangle \phi_{k}.
\notag
\end{equation}
Clearly, ($1D$ continuous) Fourier transform is to extend the idea of decomposing functions on the interval $\left[-\pi,\pi\right]$ to analyzing them across $\mathbb{R}$ by scaling the frequency domain. This approach applies analogously to higher-dimensional situations. The completeness of the Fourier basis is given by the Fourier theorem, while the uniqueness of continuous Fourier transform and the inverse Fourier transform under certain conditions is a key result in Fourier analysis \citep{Stein2003}. An important property of the Fourier series/continuous Fourier transform is the convolution property:
\begin{equation}
\forall f,g\in L^{2}\left(\left[-\pi,\pi\right]\right),\ \mathcal{F}\left(f*g\right)=\left(\mathcal{F}f\right)\cdot\left(\mathcal{F}g\right),
\notag
\end{equation}
where $\mathcal{F}$ denotes the Fourier transform. 

For a finite number of data points, \emph{discrete Fourier transform
(\ac{DFT})} can be used to approximate a function using the Fourier basis
$\left\{ \phi_{k}\right\} $ mentioned above. In the context of our discussion of \ac{DFT}, for a slight abuse of notions, let $\mathcal{F}$ also represent the Fourier series. In physical space, the equispaced grid of points is usually scaled first to match the domain of the \ac{DFT} transform, often chosen as $\left[-\pi,\pi\right]$ for $1D$ data or $\left[-\pi,\pi\right]\times\left[-\pi,\pi\right]$ for $2D$ data. \ac{DFT} then transforms
the function values evaluated at the equispaced data points in the
physical space to Fourier coefficients in the frequency space by multiplication
of the following matrix, called \ac{DFT} matrix: 
\begin{equation}
\Psi\coloneqq N^{-\frac{1}{2}}\left[\begin{array}{ccccc}
\psi^{0} & \psi^{0} & \psi^{0} & \cdots & \psi^{0}\\
\psi^{0} & \psi & \psi^{2} &  & \psi^{N-1}\\
\psi^{0} & \psi^{2} & \psi^{4} &  & \psi^{2\left(N-1\right)}\\
 & \vdots &  & \ddots & \vdots\\
\psi^{0} & \psi^{N-1} & \psi^{2\left(N-1\right)} & \cdots & \psi^{\left(N-1\right)\left(N-1\right)}
\end{array}\right],
\notag
\end{equation}
where $\psi\coloneqq\exp\left(-\frac{1}{N}2\pi i\right)$. \ac{FFT} is an algorithm to efficiently perform the \ac{DFT} for
a finite number of data points, reducing the complexity from $O\left(N^{2}\right)$
to $O\left(N\log N\right)$ \citep{Cooley1965}. Inverse \ac{FFT} can
be done similarly. 

In a two-dimensional space, the \ac{DFT} of the function $f$ is based on the
projection on a $2D$ Fourier basis $\left\{ \phi_{k}\coloneqq\exp\left(ikx+ijy\right)\vert k,j\in\mathbb{Z}\right\} $.
The convolution property and \ac{FFT} in a $2D$ space is then similar
to that of the $1D$ space \citep{Stein2003, Cooley1965}. 

Based on above, kernel density estimation can be computed efficiently
using the convolution property of Fourier transform and \ac{FFT} \citep{Silverman1982}. \cite{Silverman1982} further demonstrated the outstanding numerical performance of \ac{FFTKDE}.
Specifically, the kernel density estimation for $N$ data points is
\begin{equation}
\hat{f}\left(x;\Omega\right)\coloneqq N^{-1}\sum_{j=1}^{N}K\left(x-x_{j};\Omega\right),
\notag
\end{equation}
where $K$ denotes the kernel and $\Omega$ denotes the bandwidth
matrix. Thus, KDE can be carried out efficiently by
\begin{equation}
\hat{f}\left(x;\Omega\right)=N^{-1}\sum_{j=1}^{N}K\left(x;\Omega\right)*\delta\left(x-x_{j}\right),
\notag
\end{equation}
where $\delta$ is Dirac delta, which functions as a ``spike'' and
has Fourier transform being a constant function depending only on the chosen normalization constant of the Fourier transform. This allows $\hat{f}$
to be calculated efficiently, since the convolution property of Fourier transform implies
that 
\begin{equation}
\mathcal{F}\left(\hat{f}\right)\left(x;\Omega\right)=\mathcal{F}\left(K\right)\left(x;\Omega\right)\cdot\mathcal{F}\left(\delta\right)\left(x-x_{j}\right).
\notag
\end{equation}
Then, $\hat{f}\left(x;\Omega\right)$ evaluated on a $2D$ equispaced grid
can be calculated using I\ac{FFT}. Therefore, the evaluated density value on the $2D$ equispaced grid can
be used to calculate the mutual information estimation, specifically, 
\begin{equation}
\widehat{MI}\left(Y,X_{j}\right)=\int_{\text{supp}\left(Y\right)}\int_{\text{supp}\left(X_{j}\right)}\hat{f}_{Y,X_{j}}\left(y,x_{j}\right)\cdot\log\frac{\hat{f}_{Y,X_{j}}\left(y,x_{j}\right)}{\hat{f}_{Y}\left(y\right)\cdot\hat{f}_{X_{j}}\left(x_{j}\right)}dx_{j}dy
\label{eq:numerical-FFTKDE-MI-estimator}
\end{equation}

In \eqref{eq:numerical-FFTKDE-MI-estimator}, $\hat{f}_{Y}\left(y\right)$, $\hat{f}_{X_{j}}\left(x_{j}\right)$, and the expectation estimator itself can be numerically computed using the forward Euler method. Notably, employing the \ac{FFT} for the integration of density functions often fails to deliver satisfactory numerical results, primarily attributed to the inherent periodic characteristics of the method. \eqref{eq:numerical-FFTKDE-MI-estimator} is the equation that we use to calculate the \ac{FFTKDE} mutual information estimator. 

The estimation of mutual information using another nonparametric method, $k$NN \citep{Faivishevsky2008,Kraskov2004,Victor2002,Pal2010,Lord2018,Gao2014}, was also discussed in the paper.
The estimation of mutual information based on $k$NN can be viewed through
the lens of $k$NN density estimator. The bivariate $k$NN density
estimator can be given by 
\[
\hat{f}\left(x;k\right)\coloneqq\frac{k}{N}\cdot\left(\pi\cdot R^{2}\left(x;k\right)\right)^{-1},
\]
where $R\left(x;k\right)$ denotes the Euclidean distance from $x$
to its $k$-nearest-neighbor. In the context of a bivariate density estimator,
$\pi\cdot R^{2}\left(x;k\right)$ represents the area of the Euclidean-normed
closed ball centered at $x$ that includes the $k$-nearest-neighbors
of $x$. Following the idea of empirical CDF, the probability that
a data point is included in this closed ball is $\frac{k}{N}$;
assuming that the density inside the closed ball remains constant, the estimate of such density will be the probability of being included
in the closed ball divided by the area of the closed ball, which is
the bivariate density estimator described above. The multivariate case
with more than two variables can be established in a similar way. 

	\chapter{Appendix to Manuscript 2}  \label{appchap2}


 \section{Proofs} \label{sec:proof}
We first establish the following Lemma needed for the proof of Theorem \ref{thm:convex}.
\subsection{Proof of Theorem \ref{thm:convex} }
\label{sec:proof-thm-1}

The following lemma is needed in the proof of Theorem 1. 
\begin{lem}
\label{lem:convergence-cond-meaning} Assume that $\forall k=1,2,\dots,N$,
the convergence conditions \eqref{eq:convcond1} and \eqref{eq:convcond2}
hold, then we have the following recursive relation: 
\begin{equation}
\alpha_{k+1}\leq\frac{1}{1+\frac{\delta_{k}/\delta_{k+1}}{\alpha_{k}}}\label{eq:bifurcation}.
\end{equation}
\end{lem}
\begin{proof}
The convergence conditions \eqref{eq:convcond1} and \eqref{eq:convcond2}
gives that $\forall k=1,2,\dots,N-1$, 
\[
\alpha_{k+1}\delta_{k+1}\leq\omega_{k+1}\Leftrightarrow\alpha_{k+1}\leq\frac{\omega_{k+1}}{\delta_{k+1}},\text{ and}
\]
\[
\frac{\alpha_{k}}{\delta_{k}\Gamma_{k}}\geq\frac{\alpha_{k+1}}{\delta_{k+1}\Gamma_{k+1}}\Leftrightarrow\frac{\alpha_{k}}{\delta_{k}}\geq\frac{\alpha_{k+1}}{\delta_{k+1}\left(1-\alpha_{k+1}\right)}\Leftrightarrow\alpha_{k+1}\leq\frac{\alpha_{k}\delta_{k+1}}{\alpha_{k}\delta_{k+1}+\delta_{k}}.
\]
Following above two inequalities, we have that 
\begin{equation}
\alpha_{k+1}\leq\min\left\{ \frac{\omega_{k+1}}{\delta_{k+1}},\frac{\alpha_{k}\delta_{k+1}}{\alpha_{k}\delta_{k+1}+\delta_{k}}\right\} \label{eq:alpha-upper-bound}.
\end{equation}
We observe that in \eqref{eq:alpha-upper-bound}, $\frac{\omega_{k+1}}{\delta_{k+1}}$
is monotonically decreasing with respect to $\delta_{k+1}$ on $\mathbb{R}_{+}$;
while $\frac{\alpha_{k}\delta_{k+1}}{\alpha_{k}\delta_{k+1}+\delta_{k}}$
is monotonically increasing with respect to $\delta_{k+1}$ on $\mathbb{R}_{+}$.
This suggests: 
\begin{equation}
\arg\max_{\delta_{k+1}>0}\left(\min\left\{ \frac{\omega_{k+1}}{\delta_{k+1}},\frac{\alpha_{k}\delta_{k+1}}{\alpha_{k}\delta_{k+1}+\delta_{k}}\right\} \right)=\left\{ \frac{\omega_{k+1}+\sqrt{\omega_{k+1}^{2}+\frac{4\omega_{k+1}\delta_{k}}{\alpha_{k}}}}{2}\right\} .\label{eq:get-lambda}
\end{equation}
That is, the inequality constraints conditions \eqref{eq:convcond1}
and \eqref{eq:convcond2} for convergence are merely a lower bound
on the\emph{ vanishing rate} of $\left\{ \alpha_{k}\right\} $. Therefore
it follows from \eqref{eq:convcond1} and the (necessary) optimality
condition for \eqref{eq:get-lambda} that 
\begin{equation}
\alpha_{k+1}\leq\frac{2\omega_{k+1}}{\omega_{k+1}+\sqrt{\omega_{k+1}^{2}+\frac{4\omega_{k+1}\delta_{k}}{\alpha_{k}}}}\leq\frac{2}{1+\sqrt{1+\frac{4\delta_{k}}{\alpha_{k}\omega_{k+1}}}}=\frac{2}{1+\sqrt{1+\frac{4\delta_{k}/\delta_{k+1}}{\alpha_{k}\alpha_{k+1}}}}.\label{eq:alpha-bound}
\end{equation}
By simplifying \eqref{eq:bifurcation}, we have: 
\[
\alpha_{k+1}\leq\frac{1}{1+\frac{\delta_{k}/\delta_{k+1}}{\alpha_{k}}}.
\]
\end{proof}

We now proceed with the proof of Theorem 1. 
\begin{proof}
The complexity upper bound \eqref{eq:complexity-bound} under the
given conditions can be simplified as: 
\begin{align}
 & \left[\sum_{k=1}^{N}\Gamma_{k}^{-1}\omega_{k}\left(1-L_{\Psi}\omega_{k}\right)\right]^{-1}\left[\frac{\left\Vert x_{0}-x^{*}\right\Vert ^{2}}{\delta_{1}}+\frac{2L_{f}}{\Gamma_{N}}\left(\left\Vert x^{*}\right\Vert ^{2}+M^{2}\right)\right]\nonumber \\
= & \left[\sum_{k=1}^{N}\Gamma_{k}^{-1}\omega_{k}\left(1-L_{\Psi}\omega_{k}\right)\right]^{-1}\cdot\frac{\left\Vert x_{0}-x^{*}\right\Vert ^{2}}{\delta_{1}}\nonumber \\
= & \frac{1}{\omega\left(1-L_{\Psi}\omega\right)}\left(\sum_{k=1}^{N}\Gamma_{k}^{-1}\right)^{-1}\cdot\frac{\left\Vert x_{0}-x^{*}\right\Vert ^{2}}{\omega}\nonumber \\
= & \left(\sum_{k=1}^{N}\Gamma_{k}^{-1}\right)^{-1}\cdot\frac{\left\Vert x_{0}-x^{*}\right\Vert ^{2}}{\omega^{2}\left(1-L_{\Psi}\omega\right)}\label{eq:convex-complex-bound}.
\end{align}
Observe that $\left(\sum_{k=1}^{N}\Gamma_{k}^{-1}\right)^{-1}$ is
monotonically decreasing with respect to $\alpha_{k}$ for all $k=1,2,\dots,N$.
This property implies that \eqref{eq:convex-complex-bound} is minimized when
$\alpha_{k}$ attains its greatest value for $k=1,2,\dots,N$. 

Condition $\delta_{1}=\omega_{k}=\omega$ gives that 
\[
\omega_{1}=\delta_{1}=\alpha_{1}\delta_{1}.
\]
Since the upper bound for $\alpha_{k+1}$ presented in \eqref{eq:bifurcation}
is monotonically increasing with respect to $\alpha_{k}$, it then
follows inductively from the (necessary) optimality condition of \eqref{eq:alpha-upper-bound}
that 
\[
\alpha_{k+1}\leq\frac{1}{1+\frac{\delta_{k}/\delta_{k+1}}{\alpha_{k}}}=\frac{1}{1+\frac{\alpha_{k+1}}{\alpha_{k}^{2}}},
\]
which simplifies to 
\[
\alpha_{k+1}\leq\frac{2}{1+\sqrt{1+\frac{4}{\alpha_{k}^{2}}}}.
\]
While $\omega^{2}\left(1-L_{\Psi}\omega\right)$ should be maximized
to minimize the value of \eqref{eq:convex-complex-bound}, which implies
the minimizer for $\omega$ is 
\[
\bar{\omega}=\frac{2}{3L_{\Psi}}.
\]
And $\bar{\lambda}_{k+1}=\frac{\bar{\omega}}{\bar{\alpha}_{k+1}}$
follows directly form the necessary optimality condition for \eqref{eq:alpha-upper-bound}. It is trivial to check that $\left(\left\{ \bar{\alpha}_{k}\right\} ,\left\{ \bar{\delta}_{k}\right\} ,\bar{\omega}\right)$ is feasible under given constraints \eqref{eq:convcond1} and \eqref{eq:convcond2}.
\end{proof}

\subsection{Proof of Theorem \ref{thm:1overk}}
\label{sec:proof-thm-2}
\begin{proof}
Consider arbitrary $k=2,\dots,N$, then $\alpha_{k}\in\left(0,1\right)$
by definition. In the convergence conditions \eqref{eq:convcond1}
and \eqref{eq:convcond2}, this gives us that 
\[
\frac{\alpha_{k+1}}{\alpha_{k}}\leq\frac{2}{\alpha_{k}+\sqrt{\alpha_{k}^{2}+4}}\in\left(\frac{\sqrt{5}-1}{2},1\right).
\]
Thus, $\left\{ \alpha_{k}\right\} $ is a bounded monotonically decreasing
sequence, and $\alpha_{2}\leq\frac{2}{1+\sqrt{1+\frac{4}{1^{2}}}}=\frac{\sqrt{5}-1}{2}$
further implies that $\forall k\geq2,\ \alpha_{k}\in(0,\frac{\sqrt{5}-1}{2}]$. 

For all $k\geq2$, $\alpha_{k}\in\left(0,1\right)$ implies that $1-\alpha_{k}\in\left(0,1\right)$.
Therefore, $\Gamma_{k}^{-1}=\frac{1}{\left(1-\alpha_{2}\right)\left(1-\alpha_{3}\right)\cdots\left(1-\alpha_{k}\right)}$
is monotonically increasing with respect to $k$. Thus, $\sum_{k=1}^{N}\Gamma_{k}^{-1}=O\left(N\right)$,
which implies that $\left(\sum_{k=1}^{N}\Gamma_{k}^{-1}\right)^{-1}\cdot C_{1}=O\left(1/N\right)$. 

Observe that 
\begin{align}
 & 0<\left(\Gamma_{N}\sum_{k=1}^{N}\frac{1}{\Gamma_{k}}\right)^{-1}=\frac{1}{N\cdot\Gamma_{N}}\cdot\frac{N}{\sum_{k=1}^{N}\frac{1}{\Gamma_{k}}}\nonumber \\
 & \leq\frac{1}{N\cdot\Gamma_{N}}\cdot\left(\prod_{k=1}^{N}\Gamma_{k}\right)^{\frac{1}{N}}=\frac{1}{N}\cdot\left(\prod_{k=1}^{N}\frac{\Gamma_{k}}{\Gamma_{N}}\right)^{\frac{1}{N}}\label{eq:HM-GM}\\
 & =\frac{1}{N}\cdot\left(\prod_{k=1}^{N}\frac{\Gamma_{N}}{\Gamma_{k}}\right)^{-\frac{1}{N}}=\frac{1}{N}\cdot\left(\prod_{k=2}^{N}\left(1-\alpha_{k}\right)^{k}\right)^{-\frac{1}{N}}\nonumber \\
 & =\frac{1}{N}\cdot\prod_{k=2}^{N}\left(1-\alpha_{k}\right)^{-\frac{k}{N}},\nonumber 
\end{align}
where the inequality in \eqref{eq:HM-GM} follows from the harmonic
mean-geometric mean inequality. 

Consider arbitrary $N\in\mathbb{N}$, now we are to prove that $\forall k=1,2,\dots,N,\ \alpha_{k}\leq\frac{2}{k+1}$.
By definition, $\alpha_{1}=1\leq1$. Assume that $\alpha_{k}\leq\frac{2}{k+1}$,
then by the convergence conditions, 
\begin{align*}
\alpha_{k+1}\leq & \frac{2}{1+\sqrt{1+\frac{4}{\alpha_{k}^{2}}}}\\
\leq & \frac{2}{1+\sqrt{1+4/\left(\frac{2}{k+1}\right)^{2}}}\\
= & \frac{2}{1+\sqrt{2+2k+k^{2}}}\\
< & \frac{2}{k+2}.
\end{align*}
Thus, by mathematical induction, $\forall k=1,2,\dots,N,\ \alpha_{k}\leq\frac{2}{k+1}$.
Hence, $\sum_{k=1}^{N}\frac{k}{N}\alpha_{k}<\sum_{k=1}^{N}\frac{k}{N}\cdot\frac{2}{k}=\sum_{k=1}^{N}\frac{2}{N}=2<\infty$
as $N\rightarrow\infty$. 

Furthermore, we have that $\forall x\in(0,\frac{\sqrt{5}-1}{2}],\ -\log\left(1-x\right)<x$.
Combined with the fact that $\forall k\geq2,\ \alpha_{k}\in(0,\frac{\sqrt{5}-1}{2}]$,
we have that $\forall k\geq2,\ -\log\left(1-\alpha_{k}\right)<\alpha_{k}$.
Thus, 
\[
\log\left(\prod_{k=2}^{N}\left(1-\alpha_{k}\right)^{-\frac{k}{N}}\right)=-\sum_{k=2}^{N}\frac{k}{N}\log\left(1-\alpha_{k}\right)<\sum_{k=2}^{N}\frac{k}{N}\alpha_{k}\leq2<\infty.
\]

Therefore, $\prod_{k=2}^{N}\left(1-\alpha_{k}\right)^{-\frac{k}{N}}$
is also upper bounded as $N\rightarrow\infty$, which implies that
\[
\left(\sum_{k=1}^{N}\frac{\Gamma_{N}}{\Gamma_{k}}\right)^{-1}\leq\frac{1}{N}\cdot\prod_{k=2}^{N}\left(1-\alpha_{k}\right)^{-\frac{k}{N}}=O\left(1/N\right).
\]
Hence, $\left(\sum_{k=1}^{N}\frac{\Gamma_{N}}{\Gamma_{k}}\right)^{-1}\cdot C_{2}=O\left(1/N\right)$.
Therefore, $\left(\sum_{k=1}^{N}\Gamma_{k}^{-1}\right)^{-1}\cdot C_{1}+\left(\sum_{k=1}^{N}\frac{\Gamma_{N}}{\Gamma_{k}}\right)^{-1}\cdot C_{2}=O\left(1/N\right)$. 
\end{proof}

\subsection{Proof of Theorem \ref{thm:alpha-k-vanishing}}
\label{sec:proof-thm-3}
\begin{proof}
$\bar{\alpha}_{k}\leq\frac{2}{k+1}$ for $k=1,2,\dots,N$ has already been proved
in the proof of Theorem \ref{thm:1overk}. For the left inequality,
note that $\bar{\alpha}_{1}=1\geq\frac{2}{2+a}$ for $a>0$; for $k\geq2$,
we are to prove a stronger inequality: 
\begin{equation}
\bar{\alpha}_{k}\geq\frac{2}{\sqrt{\left(1+a\cdot k^{-b}\right)k\left[\left(1+a\cdot k^{-b}\right)k+2\right]}}\label{eq:alpha-close-lower}.
\end{equation}
For $k=2$, condition \eqref{eq:k-vanishing-cond} implies that
\begin{equation}
\label{eq:k-vanising-form2}
a\cdot2^{-b}\geq\frac{1}{\left(1-b\right)\left(4-b\right)}>\frac{1}{4}>\sqrt{5}-2 \text{ for }0<b<1,
\end{equation}
which suggests $\bar{\alpha}_{2}=\frac{2}{1+\sqrt{5}}\geq\frac{2}{\sqrt{\left(1+a\cdot2^{-b}\right)\cdot2\left[\left(1+a\cdot2^{-b}\right)\cdot2+2\right]}}$
by simple algebra. Assume \eqref{eq:alpha-close-lower} holds for
$k=t$, then 
\begin{align}
\bar{\alpha}_{t+1}= & \frac{2}{1+\sqrt{1+\frac{4}{\bar{\alpha}_{t}^{2}}}}\nonumber \\
\geq & \frac{2}{1+\sqrt{1+4/\left(2/\sqrt{\left(1+a\cdot t^{-b}\right)t\left[\left(1+a\cdot t^{-b}\right)t+2\right]}\right)^{2}}}\nonumber \\
= & \frac{2}{1+\sqrt{1+\left(1+a\cdot t^{-b}\right)t\left[\left(1+a\cdot t^{-b}\right)t+2\right]}}\nonumber \\
= & \frac{2}{\left(1+a\cdot t^{-b}\right)t+2}\nonumber \\
\geq & \frac{2}{\sqrt{\left(1+a\cdot\left(t+1\right)^{-b}\right)\left(t+1\right)\left[\left(1+a\cdot\left(t+1\right)^{-b}\right)\left(t+1\right)+2\right]}}\label{eq:many-terms};
\end{align}
and \eqref{eq:many-terms} follows from 
\begin{align}
 & \left(1+a\cdot\left(t+1\right)^{-b}\right)\left(t+1\right)\left[\left(1+a\cdot\left(t+1\right)^{-b}\right)\left(t+1\right)+2\right]-\left[\left(1+a\cdot t^{-b}\right)t+2\right]^{2}\nonumber \\
= & a^{2}\left[\left(t+1\right)^{2-2b}-t^{2-2b}\right]+2at\left[\left(t+1\right)^{1-b}-t^{1-b}\right]+4a\left[\left(t+1\right)^{1-b}-t^{1-b}\right]-1\nonumber \\
\geq & 2at\left[\left(t+1\right)^{1-b}-t^{1-b}\right]-1\nonumber \\
= & 2at^{2-b}\left[\left(1+\frac{1}{t}\right)^{1-b}-1\right]-1\nonumber \\
\geq & 2at^{2-b}\left[1+\left(1+b\right)t^{-1}-\frac{1}{2}b\left(1-b\right)t^{-2}-1\right]-1\label{eq:binom-approx}\\
= & 2a\left(1-b\right)t^{1-b}-ab\left(1-b\right)t^{-b}-1\geq0\label{eq:cond-2-works}.
\end{align}
\eqref{eq:binom-approx} follows from binomial approximation
inequality; $a>0$ and $0<b<1$ suggest that $2a\left(1-b\right)k^{1-b}-ab\left(1-b\right)k^{-b}-1$
is monotonically increasing with respect to $k$ for $k>0$, condition
\eqref{eq:k-vanishing-cond} therefore implies that $2a\left(1-b\right)k^{1-b}-ab\left(1-b\right)k^{-b}-1\geq0$
for all $k\geq2$, which is \eqref{eq:cond-2-works}. 

And proof of the left inequality for $k\geq2$ proceeds as the following:
\begin{align*}
\bar{\alpha}_{k}\geq & \frac{2}{\sqrt{\left(1+a\cdot k^{-b}\right)k\left[\left(1+a\cdot k^{-b}\right)k+2\right]}}\\
> & \frac{2}{\sqrt{\left(1+a\cdot k^{-b}\right)k\left[\left(1+a\cdot k^{-b}\right)k+2\right]+1}}\\
= & \frac{2}{\left(1+a\cdot k^{-b}\right)k+1}.
\end{align*}
\end{proof}

\subsection{Proof of Corollary \ref{cor:alpha-vanishing}}
\label{sec:proof-cor-1}
\begin{proof}
Observe that the lower bound of \eqref{eq:damping-bounds} is monotonically decreasing with respect to $a$ under given conditions. Constraint \eqref{eq:k-vanishing-cond} implies \eqref{eq:k-vanising-form2}, which further suggests that $$a\geq\frac{2^{b}}{\left(1-b\right)\left(4-b\right)}>0\text{ for }0<b<1;$$
i.e., $\bar{a}_k=\frac{\left(2/k\right)^{\bar{b}_k}}{\left(1-\bar{b}_k\right)\left(4-\bar{b}_k\right)}$. Thus, maximizing the lower bound of \eqref{eq:damping-bounds} is equivalent to minimize the convex function $\log\frac{\left(2/k\right)^{b}}{\left(1-b\right)\left(4-b\right)}$ with respect to $b$ over a open set $\left(0,1\right)$. First--order sufficient optimality condition gives the unique optimizer 
\[
\bar{b}_k=\frac{2+5\left(\log\frac{2}{k}\right)+\sqrt{9\left(\log\frac{2}{k}\right)^{2}+4}}{2\left(\log\frac{2}{k}\right)}\in\left(0,1\right)
\]
for $k\geq 8$. 
Simple algebra shows that $\lim_{k\rightarrow\infty}\frac{\bar{a}_kk^{1-\bar{b}_k}}{\log k}=\frac{2}{3}e$. Thus, the lower bound in Theorem \ref{thm:alpha-k-vanishing} becomes $\frac{k+1}{2}-\bar{\alpha}_{k}^{-1}=O\left(\log k\right)$.
\end{proof}

\section{Further Simulations}
\subsection{Penalized Linear Model \label{sec:Simulations-LM} }

\begin{figure}[H]
\begin{centering}
\includegraphics[width=\textwidth]{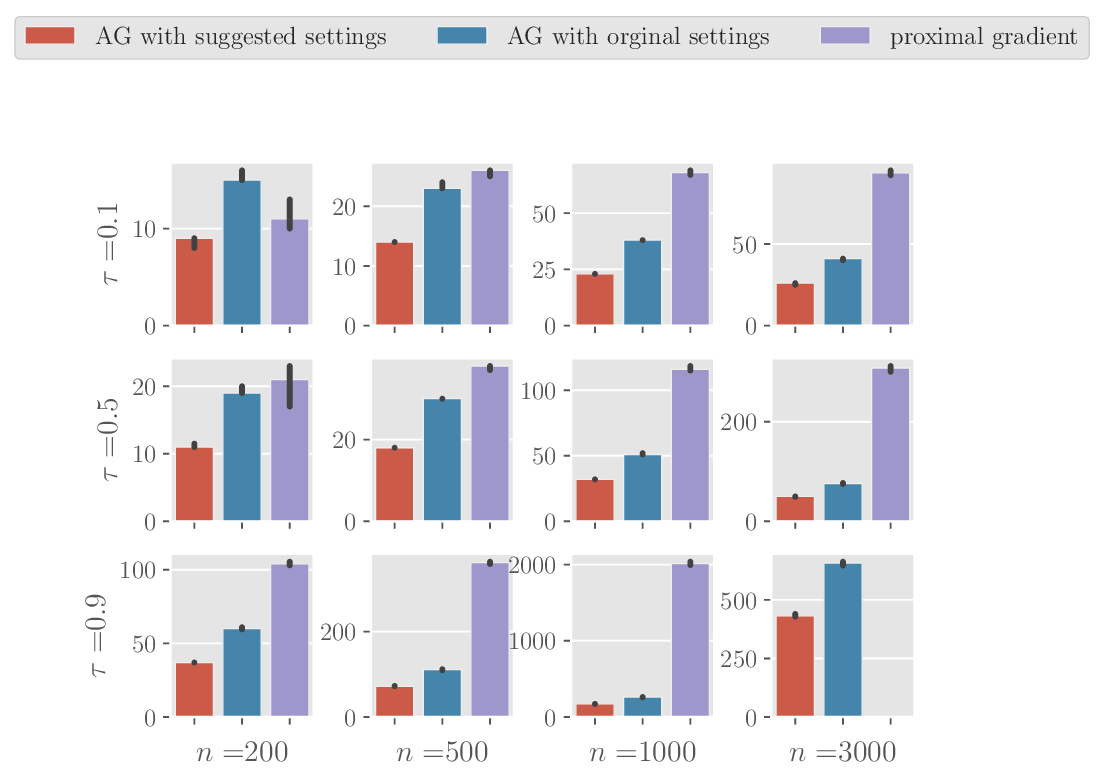}
\par\end{centering}
\caption{Median for the number of iterations required for the iterative objective value to reach $g^{*}+e^{3}$
on SCAD-penalized linear model for AG with our proposed hyperparameter
settings, AG with original settings, and proximal gradient over $100$ simulation replications, across varying covariates correlation ($\tau$) and $q/n$ values. The error bars represent the $95\%$ CIs from $1000$ bootstrap replications, $g^*$ represents the minimum per iterate found by the three methods considered.
\label{fig:sim-ag-lm-scad}}
\end{figure}

\begin{figure}[H]
\begin{centering}
\includegraphics[width=\textwidth]{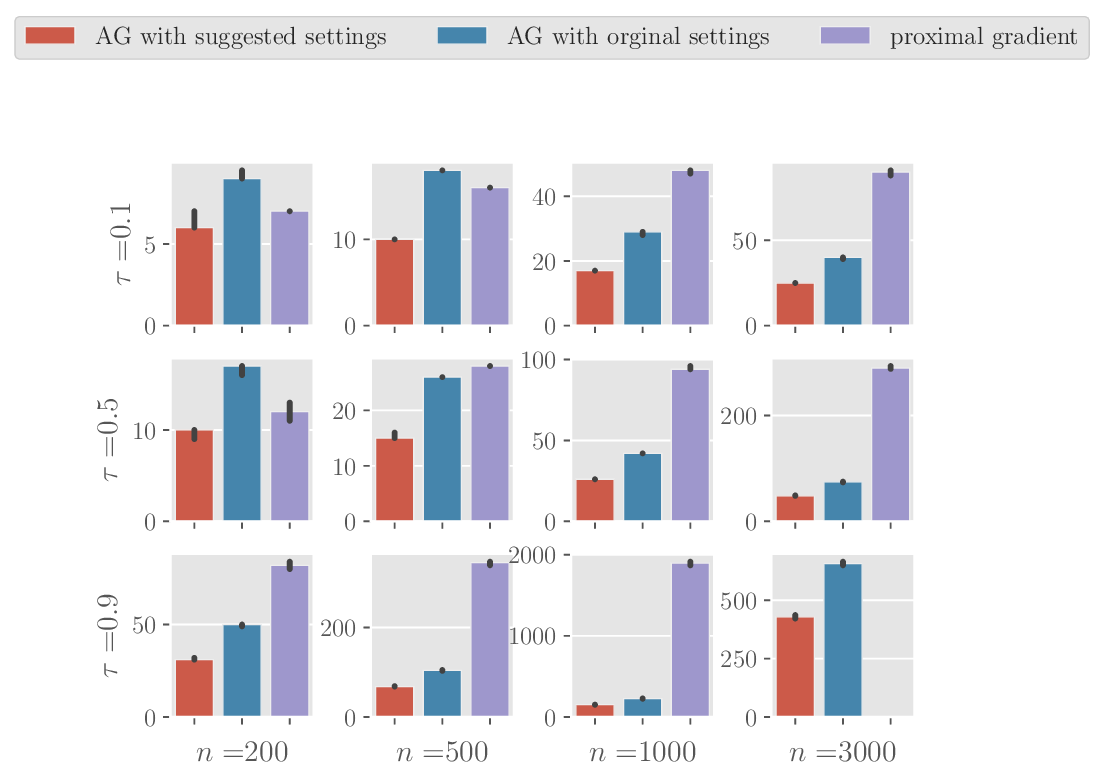}
\par\end{centering}
\caption{Median for the number of iterations required for iterative objective values to reach $g^{*}+e^{3}$
on MCP-penalized linear model for AG with our proposed hyperparameter settings, AG with original settings, and proximal gradient over $100$ simulation replications, across varying covariates correlation ($\tau$) and $q/n$ values. The error bars represent the $95\%$ CIs from $1000$ bootstrap replications, $g^*$ represents the minimum per iterate found by the three methods considered. \label{fig:sim-ag-lm-mcp}}
\end{figure}

In \figref{fig:sim-ag-lm-scad} and \ref{fig:sim-ag-lm-mcp},
the red bar represents AG using our proposed hyperparameter settings,
blue bar represents proximal gradient, and the purple bar represents AG using
the original hyperparameter settings~\citep{Ghadimi2015}. It is evident
that for penalized linear models, AG using our hyperparameter settings
outperforms proximal gradient or AG using the original proposed hyperparameter settings considerably.

\begin{figure}[H]
\begin{centering}
\includegraphics[width=\textwidth]{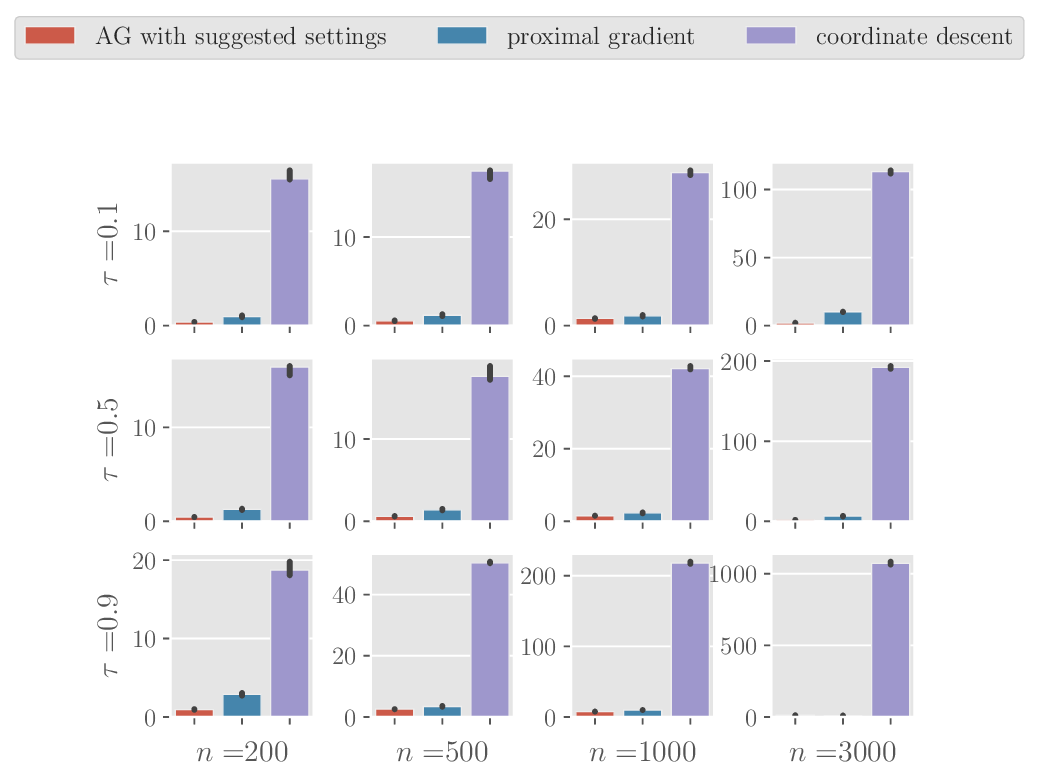}
\par\end{centering}
\caption{ Median for the computing time (in seconds) required for $\left\Vert \boldsymbol{\beta}^{\left(k+1\right)}-\boldsymbol{\beta}^{\left(k\right)}\right\Vert _{\infty}$ to fall below $10^{-4}$ on SCAD-penalized linear model for AG with our proposed hyperparameter
settings, proximal gradient, and coordinate descent over $100$ simulation replications, across varying covariates correlation ($\tau$) and $q/n$ values. The error bars represent the $95\%$ CIs from $1000$ bootstrap replications, $g^*$ represents the minimum per iterate found by the three methods considered.
\label{fig:lm-time-scad}}
\end{figure}

\begin{figure}[H]
\begin{centering}
\includegraphics[width=\textwidth]{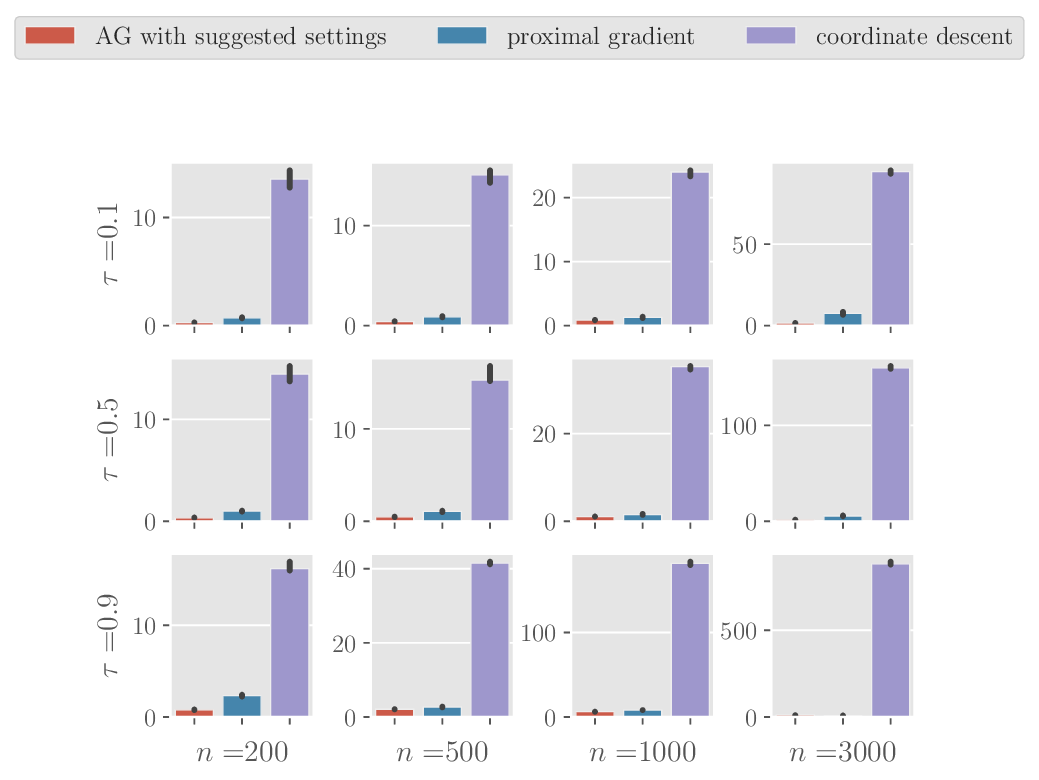}
\par\end{centering}
\caption{ Median for the computing time (in seconds) required for $\left\Vert \boldsymbol{\beta}^{\left(k+1\right)}-\boldsymbol{\beta}^{\left(k\right)}\right\Vert _{\infty}$ to fall below $10^{-4}$ on MCP-penalized linear model for AG with our proposed hyperparameter
settings, proximal gradient, and coordinate descent over $100$ simulation replications, across varying covariates correlation ($\tau$) and $q/n$ values. The error bars represent the $95\%$ CIs from $1000$ bootstrap replications, $g^*$ represents the minimum per iterate found by the three methods considered.
\label{fig:lm-time-mcp}}
\end{figure}

\newpage
{ In \figref{fig:lm-time-scad} and \ref{fig:lm-time-mcp},
the red bar represents AG using our proposed hyperparameter settings,
blue bar represents proximal gradient, and the purple bar represents coordinate descent. It is evident
that for penalized linear models, AG using our hyperparameter settings
outperforms coordinate descent significantly in terms of computing time.}

\newpage

\begin{table}[H]
\begin{center}

\caption{Signal recovery performance (sample mean and standard error of $\nicefrac{\left\Vert \boldsymbol{\beta}_{\text{true}}-\hat{\boldsymbol{\beta}}\right\Vert _{2}^{2}}{\left\Vert \boldsymbol{\beta}_{\text{true}}\right\Vert _{2}^{2}}$, Positive/Negative Predictive Values (PPV, NPV) for signal detection, and active set cardinality $\lvert\hat{\mathcal{A}}\rvert$) for {\tt ncvreg} and AG with our proposed hyperparameter settings on SCAD-penalized linear model over $100$ simulation replications, across varying values of SNRs and covariates correlations ($\tau$). \label{tab:signal-lm-scad}}
\begin{tabular}{c|ccc}
\hline 
\multicolumn{1}{c|}{$\nicefrac{\left\Vert \boldsymbol{\beta}_{\text{true}}-\hat{\boldsymbol{\beta}}\right\Vert _{2}^{2}}{\left\Vert \boldsymbol{\beta}_{\text{true}}\right\Vert _{2}^{2}}$} & $\tau=0.1$ & $0.5$ & $0.9$\\
\hline 
$\text{SNR}=1$, AG &  $0.128 (0.021)$ &  $0.521 (0.114)$ &  $2.839 (0.497)$ \\
$\text{SNR}=1$, \texttt{ncvreg} &   $0.131 (0.02)$ &  $0.485 (0.102)$ &  $2.929 (0.525)$ \\
$\text{SNR}=3$, AG &   $0.05 (0.009)$ &  $0.156 (0.035)$ &  $2.075 (0.339)$ \\
$\text{SNR}=3$, \texttt{ncvreg} &  $0.052 (0.009)$ &  $0.156 (0.028)$ &  $2.087 (0.357)$ \\
$\text{SNR}=7$, AG &  $0.022 (0.004)$ &  $0.085 (0.014)$ &  $1.278 (0.262)$ \\
$\text{SNR}=7$, \texttt{ncvreg} &  $0.021 (0.004)$ &  $0.083 (0.015)$ &    $1.3 (0.262)$ \\
$\text{SNR}=10$, AG &  $0.016 (0.003)$ &  $0.065 (0.011)$ &  $1.163 (0.207)$ \\
$\text{SNR}=10$, \texttt{ncvreg} &  $0.015 (0.003)$ &  $0.063 (0.013)$ &   $1.167 (0.22)$ \\
\hline 
\hline 
\multicolumn{1}{c|}{PPV} & $\tau=0.1$ & $0.5$ & $0.9$\\
\hline 
$\text{SNR}=1$, AG &  $0.747 (0.134)$ &  $0.622 (0.188)$ &   $0.488 (0.25)$ \\
$\text{SNR}=1$, \texttt{ncvreg} &  $0.255 (0.061)$ &  $0.287 (0.132)$ &   $0.286 (0.19)$ \\
$\text{SNR}=3$, AG &  $0.681 (0.162)$ &  $0.551 (0.206)$ &  $0.327 (0.234)$ \\
$\text{SNR}=3$, \texttt{ncvreg} &  $0.282 (0.079)$ &  $0.307 (0.098)$ &  $0.275 (0.148)$ \\
$\text{SNR}=7$, AG &   $0.58 (0.138)$ &   $0.42 (0.257)$ &  $0.197 (0.141)$ \\
$\text{SNR}=7$, \texttt{ncvreg} &   $0.32 (0.065)$ &  $0.344 (0.152)$ &  $0.175 (0.101)$ \\
$\text{SNR}=10$, AG &  $0.528 (0.272)$ &   $0.437 (0.09)$ &  $0.211 (0.081)$ \\
$\text{SNR}=10$, \texttt{ncvreg} &  $0.349 (0.127)$ &    $0.409 (0.1)$ &  $0.206 (0.047)$ \\
\hline 
\hline 
\multicolumn{1}{c|}{NPV} & $\tau=0.1$ & $0.5$ & $0.9$\\
\hline 
$\text{SNR}=1$, AG &  $0.984 (0.001)$ &  $0.984 (0.001)$ &  $0.979 (0.001)$ \\
$\text{SNR}=1$, \texttt{ncvreg} &  $0.987 (0.001)$ &  $0.986 (0.001)$ &   $0.98 (0.001)$ \\
$\text{SNR}=3$, AG &  $0.989 (0.001)$ &  $0.988 (0.002)$ &   $0.98 (0.001)$ \\
$\text{SNR}=3$, \texttt{ncvreg} &   $0.99 (0.001)$ &  $0.989 (0.001)$ &   $0.98 (0.001)$ \\
$\text{SNR}=7$, AG &  $0.992 (0.001)$ &  $0.991 (0.001)$ &  $0.981 (0.001)$ \\
$\text{SNR}=7$, \texttt{ncvreg} &  $0.993 (0.001)$ &  $0.991 (0.001)$ &  $0.981 (0.001)$ \\
$\text{SNR}=10$, AG &  $0.993 (0.001)$ &  $0.992 (0.001)$ &  $0.982 (0.001)$ \\
$\text{SNR}=10$, \texttt{ncvreg} &  $0.993 (0.001)$ &  $0.992 (0.001)$ &  $0.982 (0.001)$ \\
\hline 
\hline 
\multicolumn{1}{c|}{$\lvert\hat{\mathcal{A}}\rvert$} & $\tau=0.1$ & $0.5$ & $0.9$\\
\hline 
$\text{SNR}=1$, AG &     $25.82 (8.08)$ &   $31.58 (17.056)$ &  $23.11 (15.166)$ \\
$\text{SNR}=1$, \texttt{ncvreg} &  $100.88 (25.582)$ &   $94.32 (41.572)$ &  $42.01 (20.592)$ \\
$\text{SNR}=3$, AG &   $42.78 (14.003)$ &   $55.48 (20.653)$ &  $42.83 (16.308)$ \\
$\text{SNR}=3$, \texttt{ncvreg} &  $120.17 (33.554)$ &  $101.75 (29.498)$ &  $46.72 (16.252)$ \\
$\text{SNR}=7$, AG &   $61.89 (21.881)$ &   $97.88 (36.736)$ &  $86.71 (26.567)$ \\
$\text{SNR}=7$, \texttt{ncvreg} &   $115.4 (23.845)$ &  $107.19 (31.445)$ &    $89.74 (23.1)$ \\
$\text{SNR}=10$, AG &  $101.21 (66.968)$ &   $81.17 (25.325)$ &   $70.8 (11.642)$ \\
$\text{SNR}=10$, \texttt{ncvreg} &   $123.5 (52.077)$ &   $90.58 (40.419)$ &  $71.47 (10.954)$ \\
\hline 
\end{tabular}

\end{center}
\end{table}

\newpage

\begin{table}[H]
\begin{center}

\caption{Signal recovery performance (sample mean and standard error of $\nicefrac{\left\Vert \boldsymbol{\beta}_{\text{true}}-\hat{\boldsymbol{\beta}}\right\Vert _{2}^{2}}{\left\Vert \boldsymbol{\beta}_{\text{true}}\right\Vert _{2}^{2}}$, Positive/Negative Predictive Values (PPV, NPV), and active set cardinality $\lvert\hat{\mathcal{A}}\rvert$ for signal detection) for {\tt ncvreg} and AG with our proposed hyperparameter settings on MCP-penalized linear model over $100$ simulation replications, across varying values of SNRs and covariates correlations ($\tau$). \label{tab:signal-lm-mcp}}
\begin{tabular}{c|ccc}
\hline 
\multicolumn{1}{c|}{$\nicefrac{\left\Vert \boldsymbol{\beta}_{\text{true}}-\hat{\boldsymbol{\beta}}\right\Vert _{2}^{2}}{\left\Vert \boldsymbol{\beta}_{\text{true}}\right\Vert _{2}^{2}}$} & $\tau=0.1$ & $0.5$ & $0.9$\\
\hline 
$\text{SNR}=1$, AG &  $0.133 (0.022)$ &  $0.563 (0.124)$ &   $2.839 (0.39)$ \\
$\text{SNR}=1$, \texttt{ncvreg} &  $0.126 (0.019)$ &  $0.494 (0.112)$ &   $2.86 (0.427)$ \\
$\text{SNR}=3$, AG &   $0.049 (0.01)$ &  $0.169 (0.034)$ &  $1.997 (0.329)$ \\
$\text{SNR}=3$, \texttt{ncvreg} &  $0.048 (0.009)$ &  $0.161 (0.032)$ &    $1.92 (0.34)$ \\
$\text{SNR}=7$, AG &  $0.021 (0.004)$ &  $0.088 (0.016)$ &  $1.503 (0.329)$ \\
$\text{SNR}=7$, \texttt{ncvreg} &   $0.02 (0.004)$ &  $0.086 (0.017)$ &  $1.416 (0.302)$ \\
$\text{SNR}=10$, AG &  $0.014 (0.003)$ &  $0.059 (0.011)$ &  $1.084 (0.272)$ \\
$\text{SNR}=10$, \texttt{ncvreg} &  $0.014 (0.003)$ &  $0.059 (0.013)$ &  $1.134 (0.248)$ \\
\hline 
\hline 
\multicolumn{1}{c|}{PPV} & $\tau=0.1$ & $0.5$ & $0.9$\\
\hline 
$\text{SNR}=1$, AG &   $0.85 (0.081)$ &  $0.744 (0.161)$ &  $0.616 (0.208)$ \\
$\text{SNR}=1$, \texttt{ncvreg} &  $0.435 (0.085)$ &  $0.407 (0.135)$ &  $0.387 (0.154)$ \\
$\text{SNR}=3$, AG &  $0.842 (0.119)$ &   $0.732 (0.21)$ &  $0.506 (0.286)$ \\
$\text{SNR}=3$, \texttt{ncvreg} &  $0.505 (0.112)$ &  $0.514 (0.121)$ &   $0.366 (0.18)$ \\
$\text{SNR}=7$, AG &  $0.761 (0.175)$ &  $0.646 (0.293)$ &  $0.505 (0.218)$ \\
$\text{SNR}=7$, \texttt{ncvreg} &  $0.541 (0.128)$ &  $0.547 (0.173)$ &  $0.483 (0.201)$ \\
$\text{SNR}=10$, AG &  $0.801 (0.099)$ &  $0.489 (0.134)$ &  $0.375 (0.225)$ \\
$\text{SNR}=10$, \texttt{ncvreg} &  $0.559 (0.107)$ &  $0.476 (0.135)$ &  $0.377 (0.225)$ \\
\hline 
\hline 
\multicolumn{1}{c|}{NPV} & $\tau=0.1$ & $0.5$ & $0.9$\\
\hline 
$\text{SNR}=1$, AG &  $0.983 (0.001)$ &  $0.982 (0.001)$ &  $0.979 (0.001)$ \\
$\text{SNR}=1$, \texttt{ncvreg} &  $0.986 (0.001)$ &  $0.984 (0.001)$ &    $0.979 (0.0)$ \\
$\text{SNR}=3$, AG &  $0.988 (0.001)$ &  $0.986 (0.001)$ &   $0.98 (0.001)$ \\
$\text{SNR}=3$, \texttt{ncvreg} &  $0.989 (0.001)$ &  $0.987 (0.001)$ &   $0.98 (0.001)$ \\
$\text{SNR}=7$, AG &  $0.991 (0.001)$ &  $0.989 (0.001)$ &  $0.981 (0.001)$ \\
$\text{SNR}=7$, \texttt{ncvreg} &  $0.992 (0.001)$ &  $0.989 (0.001)$ &  $0.981 (0.001)$ \\
$\text{SNR}=10$, AG &  $0.992 (0.001)$ &   $0.99 (0.001)$ &  $0.982 (0.001)$ \\
$\text{SNR}=10$, \texttt{ncvreg} &  $0.993 (0.001)$ &   $0.99 (0.001)$ &  $0.982 (0.001)$ \\
\hline 
\hline 
\multicolumn{1}{c|}{$\lvert\hat{\mathcal{A}}\rvert$} & $\tau=0.1$ & $0.5$ & $0.9$\\
\hline 
$\text{SNR}=1$, AG &    $19.7 (4.584)$ &     $20.6 (9.45)$ &    $12.5 (8.163)$ \\
$\text{SNR}=1$, \texttt{ncvreg} &  $51.61 (13.612)$ &  $47.32 (16.093)$ &  $20.25 (11.411)$ \\
$\text{SNR}=3$, AG &   $30.55 (8.437)$ &   $34.52 (16.44)$ &  $25.37 (14.373)$ \\
$\text{SNR}=3$, \texttt{ncvreg} &  $60.14 (15.873)$ &  $48.08 (13.783)$ &   $31.0 (13.981)$ \\
$\text{SNR}=7$, AG &  $44.45 (14.273)$ &  $56.95 (32.804)$ &  $31.96 (25.048)$ \\
$\text{SNR}=7$, \texttt{ncvreg} &   $66.7 (20.364)$ &  $58.36 (24.633)$ &  $33.38 (25.617)$ \\
$\text{SNR}=10$, AG &   $43.23 (11.26)$ &  $64.65 (12.923)$ &  $46.58 (18.186)$ \\
$\text{SNR}=10$, \texttt{ncvreg} &   $65.36 (13.06)$ &  $67.16 (15.483)$ &  $46.07 (19.223)$ \\
\hline 
\end{tabular}

\end{center}
\end{table}

\subsection{Penalized Logistic Regression
\label{sec:Simulations-logistic} }

\begin{figure}[H]
\begin{centering}
\includegraphics[width=\textwidth]{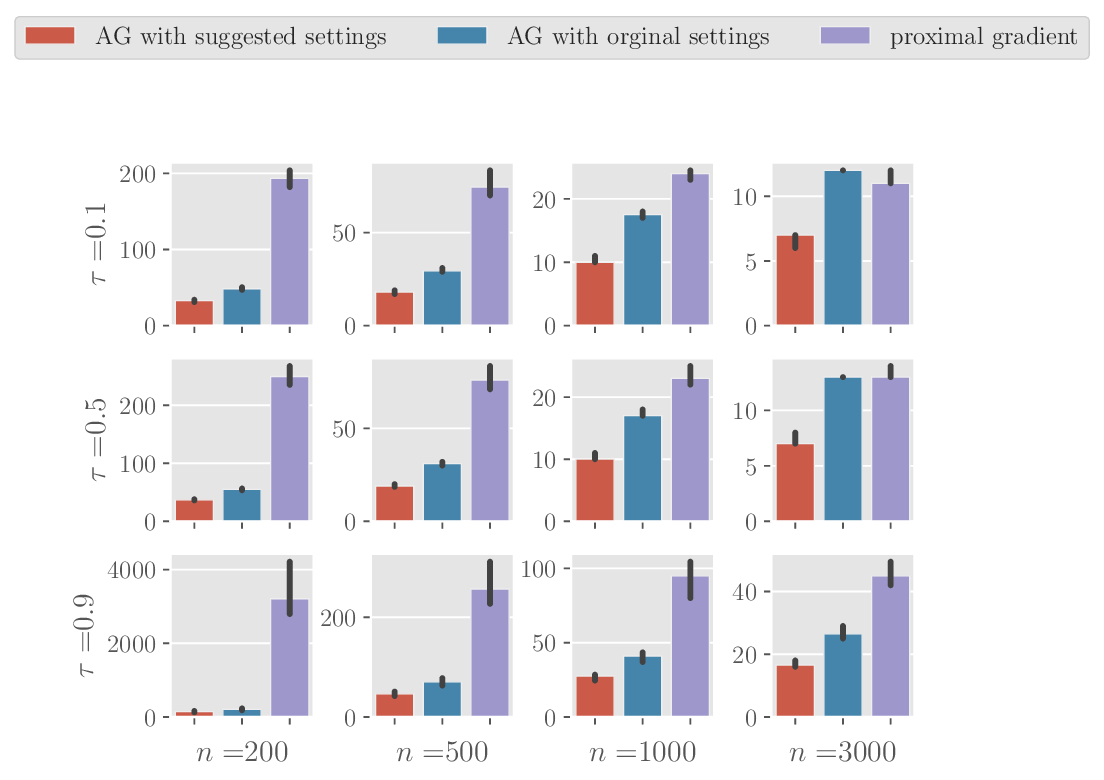}
\par\end{centering}
\caption{Median for the number of iterations required for the iterative objective values to reach $g^{*}+e^{2}$
on SCAD-penalized logistic regression for AG with our proposed hyperparameter
settings, AG with original settings, and proximal gradient over $100$ simulation replications, across varying covariates correlation ($\tau$) and $q/n$ values. The error bars represent the $95\%$ CIs from $1000$ bootstrap replications, $g^*$ represents the minimum per iterate found by the three methods considered.
\label{fig:sim-ag-logistic-scad}}
\end{figure}

\begin{figure}[H]
\begin{centering}
\includegraphics[width=\textwidth]{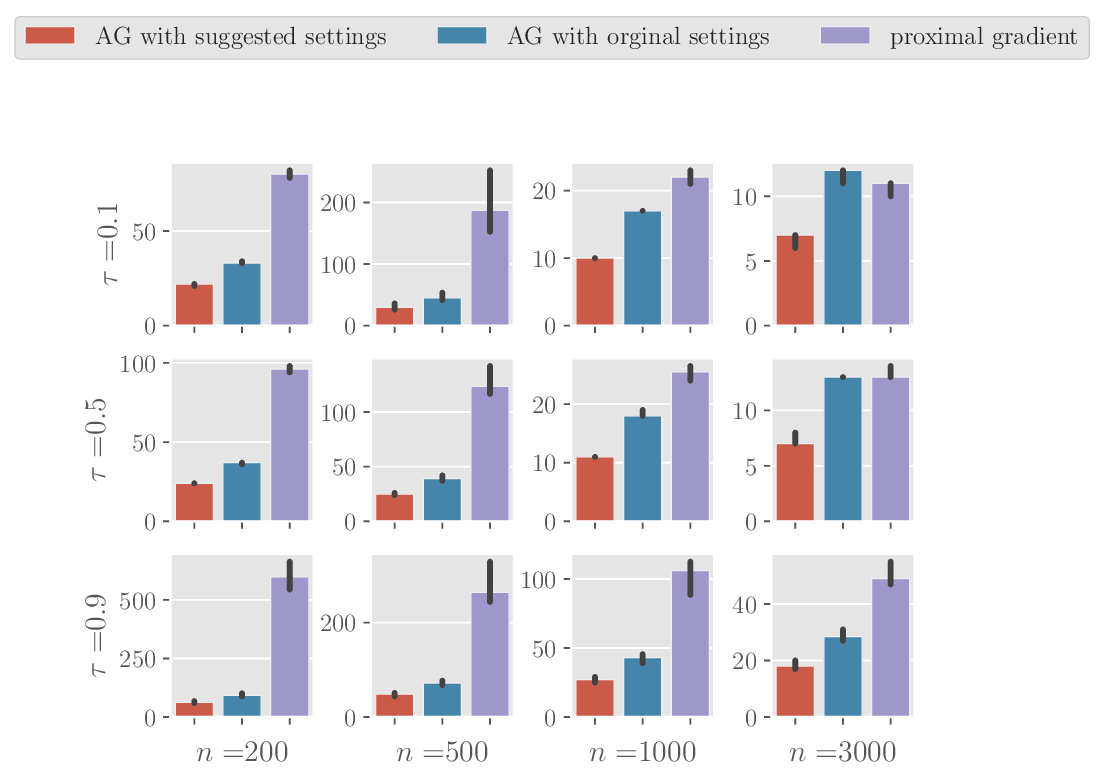}\caption{Median for the number of iterations required for iterative objective values to reach $g^{*}+e^{2}$
on MCP-penalized logistic regression for AG with our proposed hyperparameter
settings, AG with original settings, and proximal gradient over $100$ simulation replications, across varying covariates correlation ($\tau$) and $q/n$ values. The error bars represent the $95\%$ CIs from $1000$ bootstrap replications, $g^*$ represents the minimum per iterate found by the three methods considered.
\label{fig:sim-ag-logistic-mcp}}
\par\end{centering}
\end{figure}

Figure \eqref{fig:sim-ag-logistic-scad} and \eqref{fig:sim-ag-logistic-mcp}
suggest that much less iterations are needed for our method to achieve
the same amount of descent in comparison of AG with original proposed
settings for penalized logistic models.  

\newpage

\begin{table}[H]
\begin{center}

\caption{Signal recovery performance (sample mean and standard error of $\nicefrac{\left\Vert \boldsymbol{\beta}_{\text{true}}-\hat{\boldsymbol{\beta}}\right\Vert _{2}^{2}}{\left\Vert \boldsymbol{\beta}_{\text{true}}\right\Vert _{2}^{2}}$, Positive/Negative Predictive Values (PPV, NPV), and active set cardinality $\lvert\hat{\mathcal{A}}\rvert$ for signal detection) for {\tt ncvreg} and AG with our proposed hyperparameter settings on SCAD-penalized logistic model over $100$ simulation replications, across varying values of SNRs and covariates correlations ($\tau$). \label{tab:signal-logistic-scad}}
\begin{tabular}{c|ccc}
\hline 
\multicolumn{1}{c|}{$\nicefrac{\left\Vert \boldsymbol{\beta}_{\text{true}}-\hat{\boldsymbol{\beta}}\right\Vert _{2}^{2}}{\left\Vert \boldsymbol{\beta}_{\text{true}}\right\Vert _{2}^{2}}$} & $\tau=0.1$ & $0.5$ & $0.9$\\
\hline 
$\text{SNR}=1$, AG &  $0.768 (0.047)$ &   $0.81 (0.041)$ &   $0.896 (0.04)$ \\
$\text{SNR}=1$, \texttt{ncvreg} &  $0.803 (0.033)$ &   $0.84 (0.033)$ &  $0.903 (0.037)$ \\
$\text{SNR}=3$, AG &  $0.556 (0.057)$ &  $0.656 (0.054)$ &  $0.839 (0.056)$ \\
$\text{SNR}=3$, \texttt{ncvreg} &  $0.603 (0.053)$ &  $0.682 (0.055)$ &  $0.813 (0.053)$ \\
$\text{SNR}=7$, AG &  $0.377 (0.076)$ &  $0.521 (0.073)$ &  $0.779 (0.072)$ \\
$\text{SNR}=7$, \texttt{ncvreg} &  $0.438 (0.054)$ &  $0.537 (0.074)$ &  $0.735 (0.074)$ \\
$\text{SNR}=10$, AG &  $0.311 (0.077)$ &  $0.474 (0.073)$ &  $0.757 (0.079)$ \\
$\text{SNR}=10$, \texttt{ncvreg} &  $0.377 (0.064)$ &  $0.481 (0.079)$ &  $0.712 (0.078)$ \\
\hline 
\hline 
\multicolumn{1}{c|}{PPV} & $\tau=0.1$ & $0.5$ & $0.9$\\
\hline 
$\text{SNR}=1$, AG &    $0.8 (0.079)$ &    $0.779 (0.1)$ &  $0.697 (0.126)$ \\
$\text{SNR}=1$, \texttt{ncvreg} &  $0.221 (0.045)$ &  $0.265 (0.079)$ &  $0.309 (0.169)$ \\
$\text{SNR}=3$, AG &  $0.875 (0.054)$ &  $0.859 (0.065)$ &  $0.765 (0.096)$ \\
$\text{SNR}=3$, \texttt{ncvreg} &  $0.244 (0.052)$ &  $0.273 (0.072)$ &  $0.273 (0.133)$ \\
$\text{SNR}=7$, AG &  $0.901 (0.052)$ &  $0.881 (0.057)$ &  $0.788 (0.098)$ \\
$\text{SNR}=7$, \texttt{ncvreg} &    $0.27 (0.04)$ &  $0.271 (0.079)$ &  $0.267 (0.136)$ \\
$\text{SNR}=10$, AG &  $0.915 (0.048)$ &  $0.899 (0.054)$ &  $0.789 (0.097)$ \\
$\text{SNR}=10$, \texttt{ncvreg} &    $0.29 (0.05)$ &  $0.279 (0.072)$ &   $0.26 (0.123)$ \\
\hline 
\hline 
\multicolumn{1}{c|}{NPV} & $\tau=0.1$ & $0.5$ & $0.9$\\
\hline 
$\text{SNR}=1$, AG &  $0.982 (0.001)$ &   $0.98 (0.001)$ &  $0.978 (0.001)$ \\
$\text{SNR}=1$, \texttt{ncvreg} &  $0.987 (0.002)$ &  $0.985 (0.002)$ &   $0.98 (0.001)$ \\
$\text{SNR}=3$, AG &  $0.985 (0.002)$ &  $0.982 (0.001)$ &  $0.979 (0.001)$ \\
$\text{SNR}=3$, \texttt{ncvreg} &   $0.99 (0.002)$ &  $0.987 (0.002)$ &   $0.98 (0.001)$ \\
$\text{SNR}=7$, AG &  $0.987 (0.002)$ &  $0.984 (0.001)$ &  $0.979 (0.001)$ \\
$\text{SNR}=7$, \texttt{ncvreg} &  $0.992 (0.001)$ &  $0.988 (0.001)$ &   $0.98 (0.001)$ \\
$\text{SNR}=10$, AG &  $0.988 (0.002)$ &  $0.984 (0.001)$ &  $0.979 (0.001)$ \\
$\text{SNR}=10$, \texttt{ncvreg} &  $0.992 (0.001)$ &  $0.988 (0.001)$ &   $0.98 (0.001)$ \\
\hline 
\hline 
\multicolumn{1}{c|}{$\lvert\hat{\mathcal{A}}\rvert$} & $\tau=0.1$ & $0.5$ & $0.9$\\
\hline 
$\text{SNR}=1$, AG &     $17.07 (3.91)$ &     $13.4 (3.365)$ &    $7.62 (2.134)$ \\
$\text{SNR}=1$, \texttt{ncvreg} &  $120.14 (28.882)$ &   $86.49 (24.421)$ &  $39.41 (19.448)$ \\
$\text{SNR}=3$, AG &    $23.34 (4.203)$ &    $16.59 (3.459)$ &    $8.69 (2.082)$ \\
$\text{SNR}=3$, \texttt{ncvreg} &   $134.85 (29.96)$ &   $98.48 (28.434)$ &  $42.47 (15.014)$ \\
$\text{SNR}=7$, AG &     $26.98 (4.58)$ &    $19.46 (3.659)$ &    $9.79 (2.246)$ \\
$\text{SNR}=7$, \texttt{ncvreg} &  $130.33 (22.255)$ &  $105.03 (28.123)$ &  $48.81 (19.059)$ \\
$\text{SNR}=10$, AG &    $27.95 (4.462)$ &    $19.57 (3.141)$ &   $10.24 (2.346)$ \\
$\text{SNR}=10$, \texttt{ncvreg} &  $124.58 (23.016)$ &   $103.49 (27.66)$ &  $50.64 (21.138)$ \\
\hline 
\end{tabular}

\end{center}
\end{table}

\newpage

\begin{table}[H]
\begin{center}

\caption{Signal recovery performance (sample mean and standard error of $\nicefrac{\left\Vert \boldsymbol{\beta}_{\text{true}}-\hat{\boldsymbol{\beta}}\right\Vert _{2}^{2}}{\left\Vert \boldsymbol{\beta}_{\text{true}}\right\Vert _{2}^{2}}$, Positive/Negative Predictive Values (PPV, NPV), and active set cardinality $\lvert\hat{\mathcal{A}}\rvert$ for signal detection) for {\tt ncvreg} and AG with our proposed hyperparameter settings on MCP-penalized logistic model over $100$ simulation replications, across varying values of SNRs and covariates correlations ($\tau$). \label{tab:signal-logistic-mcp}}
\begin{tabular}{c|ccc}
\hline 
\multicolumn{1}{c|}{$\nicefrac{\left\Vert \boldsymbol{\beta}_{\text{true}}-\hat{\boldsymbol{\beta}}\right\Vert _{2}^{2}}{\left\Vert \boldsymbol{\beta}_{\text{true}}\right\Vert _{2}^{2}}$} & $\tau=0.1$ & $0.5$ & $0.9$\\
\hline 
$\text{SNR}=1$, AG &  $0.769 (0.044)$ &  $0.808 (0.041)$ &  $0.897 (0.043)$ \\
$\text{SNR}=1$, \texttt{ncvreg} &  $0.795 (0.036)$ &  $0.829 (0.032)$ &  $0.903 (0.038)$ \\
$\text{SNR}=3$, AG &  $0.555 (0.058)$ &  $0.654 (0.053)$ &  $0.834 (0.054)$ \\
$\text{SNR}=3$, \texttt{ncvreg} &  $0.605 (0.049)$ &  $0.674 (0.054)$ &  $0.825 (0.057)$ \\
$\text{SNR}=7$, AG &   $0.383 (0.08)$ &  $0.521 (0.069)$ &   $0.779 (0.07)$ \\
$\text{SNR}=7$, \texttt{ncvreg} &  $0.438 (0.057)$ &   $0.533 (0.07)$ &  $0.761 (0.071)$ \\
$\text{SNR}=10$, AG &   $0.31 (0.079)$ &  $0.469 (0.073)$ &  $0.753 (0.076)$ \\
$\text{SNR}=10$, \texttt{ncvreg} &  $0.381 (0.061)$ &   $0.48 (0.082)$ &  $0.737 (0.077)$ \\
\hline 
\hline 
\multicolumn{1}{c|}{PPV} & $\tau=0.1$ & $0.5$ & $0.9$\\
\hline 
$\text{SNR}=1$, AG &   $0.879 (0.06)$ &  $0.859 (0.058)$ &  $0.779 (0.087)$ \\
$\text{SNR}=1$, \texttt{ncvreg} &  $0.372 (0.068)$ &  $0.401 (0.106)$ &  $0.375 (0.157)$ \\
$\text{SNR}=3$, AG &   $0.906 (0.05)$ &   $0.889 (0.05)$ &  $0.805 (0.086)$ \\
$\text{SNR}=3$, \texttt{ncvreg} &   $0.43 (0.065)$ &  $0.445 (0.106)$ &  $0.395 (0.126)$ \\
$\text{SNR}=7$, AG &  $0.919 (0.044)$ &   $0.903 (0.05)$ &  $0.809 (0.102)$ \\
$\text{SNR}=7$, \texttt{ncvreg} &  $0.463 (0.063)$ &   $0.45 (0.104)$ &  $0.417 (0.145)$ \\
$\text{SNR}=10$, AG &  $0.918 (0.045)$ &  $0.911 (0.038)$ &  $0.804 (0.111)$ \\
$\text{SNR}=10$, \texttt{ncvreg} &  $0.502 (0.069)$ &  $0.468 (0.095)$ &  $0.412 (0.137)$ \\
\hline 
\hline 
\multicolumn{1}{c|}{NPV} & $\tau=0.1$ & $0.5$ & $0.9$\\
\hline 
$\text{SNR}=1$, AG &  $0.981 (0.001)$ &   $0.98 (0.001)$ &  $0.978 (0.001)$ \\
$\text{SNR}=1$, \texttt{ncvreg} &  $0.986 (0.002)$ &  $0.983 (0.001)$ &  $0.978 (0.001)$ \\
$\text{SNR}=3$, AG &  $0.985 (0.002)$ &  $0.982 (0.001)$ &  $0.979 (0.001)$ \\
$\text{SNR}=3$, \texttt{ncvreg} &  $0.989 (0.002)$ &  $0.985 (0.001)$ &  $0.979 (0.001)$ \\
$\text{SNR}=7$, AG &  $0.987 (0.002)$ &  $0.984 (0.001)$ &   $0.98 (0.001)$ \\
$\text{SNR}=7$, \texttt{ncvreg} &  $0.991 (0.002)$ &  $0.986 (0.001)$ &   $0.98 (0.001)$ \\
$\text{SNR}=10$, AG &  $0.988 (0.002)$ &  $0.984 (0.001)$ &   $0.98 (0.001)$ \\
$\text{SNR}=10$, \texttt{ncvreg} &  $0.991 (0.001)$ &  $0.987 (0.001)$ &   $0.98 (0.001)$ \\
\hline 
\end{tabular}
\begin{tabular}{c|ccc}
\hline 
\multicolumn{1}{c|}{$\lvert\hat{\mathcal{A}}\rvert$} & $\tau=0.1$ & $0.5$ & $0.9$\\
\hline 
$\text{SNR}=1$, AG &   $13.86 (3.082)$ &   $11.42 (2.776)$ &   $6.72 (1.744)$ \\
$\text{SNR}=1$, \texttt{ncvreg} &  $59.83 (14.138)$ &   $42.1 (12.546)$ &  $19.72 (8.393)$ \\
$\text{SNR}=3$, AG &   $21.86 (4.313)$ &   $15.84 (3.036)$ &   $8.84 (1.938)$ \\
$\text{SNR}=3$, \texttt{ncvreg} &  $66.57 (13.203)$ &    $48.28 (14.5)$ &  $22.81 (9.784)$ \\
$\text{SNR}=7$, AG &   $25.75 (4.776)$ &   $18.78 (3.189)$ &  $10.33 (2.565)$ \\
$\text{SNR}=7$, \texttt{ncvreg} &  $69.44 (11.876)$ &  $52.54 (13.638)$ &  $24.63 (8.741)$ \\
$\text{SNR}=10$, AG &   $27.53 (4.649)$ &   $19.55 (3.093)$ &  $11.06 (2.877)$ \\
$\text{SNR}=10$, \texttt{ncvreg} &  $65.38 (10.776)$ &  $51.66 (12.785)$ &  $25.59 (9.428)$ \\
\hline 
\end{tabular}

\end{center}
\end{table}

\end{appendices}

\clearpage


\addcontentsline{toc}{chapter}{References}
\bibHeading{References}
\bibliographystyle{plainnat}
\bibliography{KYthesisbib.bib}

\begin{thebibliography}{194}
\providecommand{\natexlab}[1]{#1}
\providecommand{\url}[1]{\texttt{#1}}
\expandafter\ifx\csname urlstyle\endcsname\relax
  \providecommand{\doi}[1]{doi: #1}\else
  \providecommand{\doi}{doi: \begingroup \urlstyle{rm}\Url}\fi

\bibitem[Adeli et~al.(2017)Adeli, Wu, Saghafi, An, Shi, and Shen]{Adeli2017}
Ehsan Adeli, Guorong Wu, Behrouz Saghafi, Le~An, Feng Shi, and Dinggang Shen.
\newblock Kernel-based joint feature selection and max-margin classification
  for early diagnosis of parkinson's disease.
\newblock \emph{Scientific Reports}, 7\penalty0 (1), January 2017.
\newblock \doi{10.1038/srep41069}.

\bibitem[Agarwal et~al.(2009)Agarwal, Meehan, and O'Regan]{Agarwal2009}
Ravi~P. Agarwal, Maria Meehan, and Donal O'Regan.
\newblock \emph{Fixed point theory and applications}.
\newblock Number 141 in Cambridge tracts in mathematics. Cambridge Univ. Press,
  Cambridge [u.a.], digitally printed version, paperpack re-issue edition,
  2009.
\newblock ISBN 9780521802505.

\bibitem[Akyildiz and M\'{i}guez(2021)]{Akyildiz2021}
{\"O}mer~Deniz Akyildiz and Joaqu\'{i}n M\'{i}guez.
\newblock Convergence rates for optimised adaptive importance samplers.
\newblock \emph{Statistics and Computing}, 31\penalty0 (2), January 2021.
\newblock \doi{10.1007/s11222-020-09983-1}.

\bibitem[Atkinson(1989)]{Atkinson1989}
Kendall~E. Atkinson.
\newblock \emph{An Introduction to Numerical Analysis}.
\newblock Wiley, New York [u.a.], 2. ed., [14. print] edition, 1989.
\newblock ISBN 0471624896.
\newblock Bibliogr. S. 665.

\bibitem[Attouch et~al.(2020)Attouch, Chbani, Fadili, and Riahi]{Attouch2020}
Hedy Attouch, Zaki Chbani, Jalal Fadili, and Hassan Riahi.
\newblock First-order optimization algorithms via inertial systems with hessian
  driven damping.
\newblock \emph{Mathematical Programming}, November 2020.
\newblock \doi{10.1007/s10107-020-01591-1}.

\bibitem[Bai and Silverstein(2010)]{Bai2010}
Zhidong Bai and Jack~W. Silverstein.
\newblock \emph{Spectral Analysis of Large Dimensional Random Matrices}.
\newblock Springer New York, 2010.
\newblock ISBN 9781441906618.
\newblock \doi{10.1007/978-1-4419-0661-8}.

\bibitem[Barry et~al.(2020)Barry, Bhagwat, Misic, Poline, and
  Greenwood]{Barry2020}
Amadou Barry, Nikhil Bhagwat, Bratislav Misic, Jean-Baptiste Poline, and Celia
  M.~T. Greenwood.
\newblock Asymmetric influence measure for high dimensional regression.
\newblock \emph{Communications in Statistics - Theory and Methods}, 51\penalty0
  (16):\penalty0 5461--5487, November 2020.
\newblock \doi{10.1080/03610926.2020.1841793}.

\bibitem[Bauschke and Combettes(2011)]{Bauschke2011}
Heinz~H. Bauschke and Patrick~L. Combettes.
\newblock \emph{Convex Analysis and Monotone Operator Theory in Hilbert
  Spaces}.
\newblock SpringerLink. Springer New York, New York, NY, 2011.
\newblock ISBN 9781441994677.

\bibitem[Beck(2017)]{Beck2017}
Amir Beck.
\newblock \emph{First-order methods in optimization}.
\newblock Society for Industrial and Applied Mathematics Mathematical
  Optimization Society, Philadelphia Philadelphia, 2017.
\newblock ISBN 9781611974997.

\bibitem[Beck and Teboulle(2009)]{Beck2009}
Amir Beck and Marc Teboulle.
\newblock A fast iterative shrinkage-thresholding algorithm for linear inverse
  problems.
\newblock \emph{SIAM Journal on Imaging Sciences}, 2\penalty0 (1):\penalty0
  183--20, 2009.
\newblock URL
  \url{https://proxy.library.mcgill.ca/login?url=https://search.proquest.com/docview/925336974?accountid=12339}.
\newblock Copyright - Copyright] {\copyright} 2009 Society for Industrial and
  Applied Mathematics; Last updated - 2012-07-02.

\bibitem[Bell and Drew(2018)]{Bell2018}
Daniel Bell and Zach Drew.
\newblock Voxel size, September 2018.

\bibitem[Bertrand et~al.(2020)Bertrand, Klopfenstein, Blondel, Vaiter,
  Gramfort, and Salmon]{Bertrand2020}
Quentin Bertrand, Quentin Klopfenstein, Mathieu Blondel, Samuel Vaiter,
  Alexandre Gramfort, and Joseph Salmon.
\newblock Implicit differentiation of lasso-type models for hyperparameter
  optimization.
\newblock In Hal~Daum{\'{e}} III and Aarti Singh, editors, \emph{Proceedings of
  the 37\textsuperscript{th} International Conference on Machine Learning},
  volume 119 of \emph{Proceedings of Machine Learning Research}, pages
  810--821. PMLR, 13--18 Jul 2020.
\newblock URL \url{https://proceedings.mlr.press/v119/bertrand20a.html}.

\bibitem[Bertrand et~al.(2022)Bertrand, Klopfenstein, Massias, Blondel, Vaiter,
  Gramfort, and Salmon]{Bertrand2022}
Quentin Bertrand, Quentin Klopfenstein, Mathurin Massias, Mathieu Blondel,
  Samuel Vaiter, Alexandre Gramfort, and Joseph Salmon.
\newblock Implicit differentiation for fast hyperparameter selection in
  non-smooth convex learning.
\newblock \emph{Journal of Machine Learning Research}, 23\penalty0
  (149):\penalty0 1--43, 2022.
\newblock URL \url{http://jmlr.org/papers/v23/21-0486.html}.

\bibitem[Bhatnagar et~al.(2019)Bhatnagar, Yang, Lu, Schurr, Loredo-Osti,
  Forest, Oualkacha, and Greenwood]{Bhatnagar2019}
Sahir~R Bhatnagar, Yi~Yang, Tianyuan Lu, Erwin Schurr, JC~Loredo-Osti, Marie
  Forest, Karim Oualkacha, and Celia~MT Greenwood.
\newblock Simultaneous snp selection and adjustment for population structure in
  high dimensional prediction models.
\newblock \emph{bioRxiv}, 2019.
\newblock \doi{10.1101/408484}.
\newblock URL \url{https://www.biorxiv.org/content/early/2019/07/15/408484}.

\bibitem[Birg{\'{e}} and Rozenholc(2006)]{Birge2006}
Lucien Birg{\'{e}} and Yves Rozenholc.
\newblock How many bins should be put in a regular histogram.
\newblock \emph{ESAIM: Probability and Statistics}, 10:\penalty0 24--45,
  January 2006.
\newblock ISSN 1262-3318.
\newblock \doi{10.1051/ps:2006001}.

\bibitem[Blondel et~al.(2022)Blondel, Berthet, Cuturi, Frostig, Hoyer,
  Llinares-Lopez, Pedregosa, and Vert]{Blondel2022}
Mathieu Blondel, Quentin Berthet, Marco Cuturi, Roy Frostig, Stephan Hoyer,
  Felipe Llinares-Lopez, Fabian Pedregosa, and Jean-Philippe Vert.
\newblock Efficient and modular implicit differentiation.
\newblock In S.~Koyejo, S.~Mohamed, A.~Agarwal, D.~Belgrave, K.~Cho, and A.~Oh,
  editors, \emph{Advances in Neural Information Processing Systems}, volume~35,
  pages 5230--5242. Curran Associates, Inc., 2022.
\newblock URL
  \url{https://proceedings.neurips.cc/paper_files/paper/2022/file/228b9279ecf9bbafe582406850c57115-Paper-Conference.pdf}.

\bibitem[Borland(2002{\natexlab{a}})]{Borland2002}
Lisa Borland.
\newblock A theory of non-gaussian option pricing.
\newblock \emph{Quantitative Finance}, 2\penalty0 (6):\penalty0 415--431,
  December 2002{\natexlab{a}}.
\newblock \doi{10.1080/14697688.2002.0000009}.

\bibitem[Borland(2002{\natexlab{b}})]{Borland2002a}
Lisa Borland.
\newblock Option pricing formulas based on a non-gaussian stock price model.
\newblock \emph{Physical Review Letters}, 89\penalty0 (9):\penalty0 098701,
  August 2002{\natexlab{b}}.
\newblock \doi{10.1103/physrevlett.89.098701}.

\bibitem[Botev et~al.(2010)Botev, Grotowski, and Kroese]{Botev2010}
Z.~I. Botev, J.~F. Grotowski, and D.~P. Kroese.
\newblock Kernel density estimation via diffusion.
\newblock \emph{The Annals of Statistics}, 38\penalty0 (5), October 2010.
\newblock \doi{10.1214/10-aos799}.

\bibitem[Breheny and Huang(2011)]{Breheny2011}
Patrick Breheny and Jian Huang.
\newblock Coordinate descent algorithms for nonconvex penalized regression,
  with applications to biological feature selection.
\newblock \emph{Annals of Applied Statistics 2011, Vol. 5, No. 1, 232-253},
  April 2011.
\newblock \doi{10.1214/10-AOAS388}.

\bibitem[Brent(1971)]{Brent1971}
R.~P. Brent.
\newblock An algorithm with guaranteed convergence for finding a zero of a
  function.
\newblock \emph{The Computer Journal}, 14\penalty0 (4):\penalty0 422--425,
  April 1971.
\newblock ISSN 1460-2067.
\newblock \doi{10.1093/comjnl/14.4.422}.

\bibitem[B{\"{u}}hlmann et~al.(2014)B{\"{u}}hlmann, Kalisch, and
  Meier]{Buehlmann2014}
Peter B{\"{u}}hlmann, Markus Kalisch, and Lukas Meier.
\newblock High-dimensional statistics with a view toward applications in
  biology.
\newblock \emph{Annual Review of Statistics and Its Application}, 1\penalty0
  (1):\penalty0 255--278, January 2014.
\newblock \doi{10.1146/annurev-statistics-022513-115545}.

\bibitem[Burden(2016)]{Burden2016}
Richard~L. Burden.
\newblock \emph{Numerical analysis}.
\newblock Cengage Learning, Boston, MA, tenth edition edition, 2016.
\newblock ISBN 1305253663.
\newblock Includes bibliographical references and index.

\bibitem[Bycroft et~al.(2018)Bycroft, Freeman, Petkova, Band, Elliott, Sharp,
  Motyer, Vukcevic, Delaneau, O'Connell, Cortes, Welsh, Young, Effingham,
  McVean, Leslie, Allen, Donnelly, and Marchini]{Bycroft2018}
Clare Bycroft, Colin Freeman, Desislava Petkova, Gavin Band, Lloyd~T. Elliott,
  Kevin Sharp, Allan Motyer, Damjan Vukcevic, Olivier Delaneau, Jared
  O'Connell, Adrian Cortes, Samantha Welsh, Alan Young, Mark Effingham, Gil
  McVean, Stephen Leslie, Naomi Allen, Peter Donnelly, and Jonathan Marchini.
\newblock The {UK} biobank resource with deep phenotyping and genomic data.
\newblock \emph{Nature}, 562\penalty0 (7726):\penalty0 203--209, October 2018.
\newblock \doi{10.1038/s41586-018-0579-z}.

\bibitem[Calcagn{\`{i}} et~al.(2019)Calcagn{\`{i}}, Finos, Alto{\'{e}}, and
  Pastore]{Calcagni2019}
Antonio Calcagn{\`{i}}, Livio Finos, Gianmarco Alto{\'{e}}, and Massimiliano
  Pastore.
\newblock A maximum entropy procedure to solve likelihood equations.
\newblock \emph{Entropy}, 21\penalty0 (6):\penalty0 596, June 2019.
\newblock ISSN 1099-4300.
\newblock \doi{10.3390/e21060596}.

\bibitem[Cameron et~al.(2013)Cameron, Yassine, Carlton, Francois, Alan,
  Andr{\'{a}}s, Budhachandra, John, Qingyang, Michael, Chaogan, and
  Pierre]{Cameron2013}
Craddock Cameron, Benhajali Yassine, Chu Carlton, Chouinard Francois, Evans
  Alan, Jakab Andr{\'{a}}s, Khundrakpam Budhachandra, Lewis John, Li~Qingyang,
  Milham Michael, Yan Chaogan, and Bellec Pierre.
\newblock The neuro bureau preprocessing initiative: open sharing of
  preprocessed neuroimaging data and derivatives.
\newblock \emph{Frontiers in Neuroinformatics}, 7, 2013.
\newblock \doi{10.3389/conf.fninf.2013.09.00041}.

\bibitem[Castellan(2000)]{Castellan2000}
Gw{\'{e}}na{\"{e}}lle Castellan.
\newblock S{\'{e}}lection d{\textquoteright}histogrammes {\`{a}}
  l{\textquoteright}aide d{\textquoteright}un crit{\`{e}}re de type akaike.
\newblock \emph{Comptes Rendus de l{\textquoteright}Acad{\'{e}}mie des Sciences
  - Series I - Mathematics}, 330\penalty0 (8):\penalty0 729--732, April 2000.
\newblock ISSN 0764-4442.
\newblock \doi{10.1016/s0764-4442(00)00250-0}.

\bibitem[Chai et~al.(2009)Chai, Walther, Beck, and Fei-Fei]{Chai2009}
Barry Chai, Dirk~B. Walther, Diane~M. Beck, and Li~Fei-Fei.
\newblock Exploring functional connectivity of the human brain using
  multivariate information analysis.
\newblock In \emph{Proceedings of the 22\textsuperscript{nd} International
  Conference on Neural Information Processing Systems}, NIPS'09, pages
  270--278, Red Hook, NY, USA, 2009. Curran Associates Inc.
\newblock ISBN 9781615679119.

\bibitem[Chandrashekar and Sahin(2014)]{Chandrashekar2014}
Girish Chandrashekar and Ferat Sahin.
\newblock A survey on feature selection methods.
\newblock \emph{Computers {\&}amp$\mathsemicolon$ Electrical Engineering},
  40\penalty0 (1):\penalty0 16--28, January 2014.
\newblock \doi{10.1016/j.compeleceng.2013.11.024}.

\bibitem[Chatterjee et~al.(2016)Chatterjee, Shi, and
  Garc{\'{\i}}a-Closas]{Chatterjee2016}
Nilanjan Chatterjee, Jianxin Shi, and Montserrat Garc{\'{\i}}a-Closas.
\newblock Developing and evaluating polygenic risk prediction models for
  stratified disease prevention.
\newblock \emph{Nature Reviews Genetics}, 17\penalty0 (7):\penalty0 392--406,
  May 2016.
\newblock \doi{10.1038/nrg.2016.27}.

\bibitem[Chen et~al.(2018)Chen, Luo, Han, and Chen]{Chen2018}
Cuiling Chen, Liling Luo, Caihong Han, and Yu~Chen.
\newblock Global convergence of an extended descent algorithm without line
  search for unconstrained optimization.
\newblock \emph{Journal of Applied Mathematics and Physics}, 06\penalty0
  (01):\penalty0 130--137, 2018.
\newblock ISSN 2327-4379.
\newblock \doi{10.4236/jamp.2018.61013}.

\bibitem[Chen and Zhou(2010)]{Chen2010}
Xiaojun Chen and Weijun Zhou.
\newblock Smoothing nonlinear conjugate gradient method for image restoration
  using nonsmooth nonconvex minimization.
\newblock \emph{{SIAM} Journal on Imaging Sciences}, 3\penalty0 (4):\penalty0
  765--790, January 2010.
\newblock \doi{10.1137/080740167}.

\bibitem[Clarke(2004)]{Clarke2004}
Francis Clarke.
\newblock \emph{Lyapunov Functions and Feedback in Nonlinear Control}, pages
  267--282.
\newblock Springer Berlin Heidelberg, May 2004.
\newblock ISBN 9783540399834.
\newblock \doi{10.1007/978-3-540-39983-4_17}.

\bibitem[Clarke(1990)]{Clarke1990}
Francis~H. Clarke.
\newblock \emph{Optimization and nonsmooth analysis}.
\newblock Number~5 in Classics in applied mathematics. Society for Industrial
  and Applied Mathematics (SIAM, 3600 Market Street, Floor 6, Philadelphia, PA
  19104), Philadelphia, Pa, 1990.
\newblock ISBN 9781611971309.
\newblock Reprint. Originally published: New York : Wiley, 1983.

\bibitem[Cole et~al.(2017)Cole, Ritchie, Bastin, Hern{\'{a}}ndez, Maniega,
  Royle, Corley, Pattie, Harris, Zhang, Wray, Redmond, Marioni, Starr, Cox,
  Wardlaw, Sharp, and Deary]{Cole2017}
J~H Cole, S~J Ritchie, M~E Bastin, M~C~Vald{\'{e}}s Hern{\'{a}}ndez,
  S~Mu{\~{n}}oz Maniega, N~Royle, J~Corley, A~Pattie, S~E Harris, Q~Zhang, N~R
  Wray, P~Redmond, R~E Marioni, J~M Starr, S~R Cox, J~M Wardlaw, D~J Sharp, and
  I~J Deary.
\newblock Brain age predicts mortality.
\newblock \emph{Molecular Psychiatry}, 23\penalty0 (5):\penalty0 1385--1392,
  April 2017.
\newblock \doi{10.1038/mp.2017.62}.

\bibitem[Collins et~al.(1994)Collins, Neelin, Peters, and Evans]{Collins1994}
D.~Louis Collins, Peter Neelin, Terrence Peters, and Alan~C. Evans.
\newblock Automatic 3d intersubject registration of mr volumetric data in
  standardized talairach space.
\newblock \emph{Journal of Computer Assisted Tomography}, 18:\penalty0
  192--205, 1994.
\newblock URL \url{https://api.semanticscholar.org/CorpusID:8026836}.

\bibitem[Combrisson et~al.(2022)Combrisson, Allegra, Basanisi, Ince, Giordano,
  Bastin, and Brovelli]{Combrisson2022}
Etienne Combrisson, Michele Allegra, Ruggero Basanisi, Robin~A.A. Ince,
  Bruno~L. Giordano, Julien Bastin, and Andrea Brovelli.
\newblock Group-level inference of information-based measures for the analyses
  of cognitive brain networks from neurophysiological data.
\newblock \emph{{NeuroImage}}, 258:\penalty0 119347, September 2022.
\newblock \doi{10.1016/j.neuroimage.2022.119347}.

\bibitem[Cooley and Tukey(1965)]{Cooley1965}
James~W. Cooley and John~W. Tukey.
\newblock An algorithm for the machine calculation of complex fourier series.
\newblock \emph{Mathematics of Computation}, 19\penalty0 (90):\penalty0
  297--301, 1965.
\newblock ISSN 1088-6842.
\newblock \doi{10.1090/s0025-5718-1965-0178586-1}.

\bibitem[Costa et~al.(2003)Costa, Hero, and Vignat]{Costa2003}
Jose Costa, Alfred Hero, and Christophe Vignat.
\newblock On solutions to multivariate maximum $\alpha$-entropy problems.
\newblock In \emph{Lecture Notes in Computer Science}, pages 211--226. Springer
  Berlin Heidelberg, 2003.
\newblock \doi{10.1007/978-3-540-45063-4_14}.

\bibitem[Cover and Thomas(2006)]{Cover2006}
T.~M. Cover and Joy~A. Thomas.
\newblock \emph{Elements of Information Theory}.
\newblock John Wiley \& Sons, Inc., 2006.
\newblock ISBN 9780471241959.

\bibitem[Dai and Yuan(1999)]{Dai1999}
Y.~H. Dai and Y.~Yuan.
\newblock A nonlinear conjugate gradient method with a strong global
  convergence property.
\newblock \emph{SIAM Journal on Optimization}, 10\penalty0 (1):\penalty0
  177--182, January 1999.
\newblock ISSN 1095-7189.
\newblock \doi{10.1137/s1052623497318992}.

\bibitem[Dale et~al.(1999)Dale, Fischl, and Sereno]{Dale1999}
Anders~M. Dale, Bruce Fischl, and Martin~I. Sereno.
\newblock Cortical surface-based analysis.
\newblock \emph{NeuroImage}, 9\penalty0 (2):\penalty0 179--194, February 1999.
\newblock ISSN 1053-8119.
\newblock \doi{10.1006/nimg.1998.0395}.

\bibitem[Dandine-Roulland and Perdry(2015)]{DandineRoulland2015}
Claire Dandine-Roulland and Herv{\'{e}} Perdry.
\newblock The use of the linear mixed model in human genetics.
\newblock \emph{Human Heredity}, 80\penalty0 (4):\penalty0 196--206, 2015.
\newblock ISSN 1423-0062.
\newblock \doi{10.1159/000447634}.

\bibitem[Domingo et~al.(2017)Domingo, d{\textquoteright}Onofrio, and
  Flandoli]{Domingo2017}
Dario Domingo, Alberto d{\textquoteright}Onofrio, and Franco Flandoli.
\newblock Boundedness vs unboundedness of a noise linked to tsallis
  q-statistics: The role of the overdamped approximation.
\newblock \emph{Journal of Mathematical Physics}, 58\penalty0 (3), March 2017.
\newblock ISSN 1089-7658.
\newblock \doi{10.1063/1.4977081}.

\bibitem[Edelman(1988)]{Edelman1988}
Alan Edelman.
\newblock Eigenvalues and condition numbers of random matrices.
\newblock \emph{{SIAM} Journal on Matrix Analysis and Applications}, 9\penalty0
  (4):\penalty0 543--560, October 1988.
\newblock \doi{10.1137/0609045}.

\bibitem[Epanechnikov(1969)]{Epanechnikov1969}
V.~A. Epanechnikov.
\newblock Non-parametric estimation of a multivariate probability density.
\newblock \emph{Theory of Probability \& Its Applications}, 14\penalty0
  (1):\penalty0 153--158, January 1969.
\newblock ISSN 1095-7219.
\newblock \doi{10.1137/1114019}.

\bibitem[Faivishevsky and Goldberger(2008)]{Faivishevsky2008}
Lev Faivishevsky and Jacob Goldberger.
\newblock Ica based on a smooth estimation of the differential entropy.
\newblock In \emph{Proceedings of the 21\textsuperscript{st} International
  Conference on Neural Information Processing Systems}, NIPS'08, pages
  433--440, Red Hook, NY, USA, 2008. Curran Associates Inc.
\newblock ISBN 9781605609492.

\bibitem[Fan and Li(2001)]{Fan2001}
Jianqing Fan and Runze Li.
\newblock Variable selection via nonconcave penalized likelihood and its oracle
  properties.
\newblock \emph{Journal of the American Statistical Association}, 96\penalty0
  (456):\penalty0 1348--1360, 2001.
\newblock ISSN 0162-1459.
\newblock URL \url{http://www.jstor.org/stable/3085904}.

\bibitem[Fan and Chou(2016)]{Fan2016}
Miaolin Fan and Chun-An Chou.
\newblock Exploring stability-based voxel selection methods in {MVPA} using
  cognitive neuroimaging data: a comprehensive study.
\newblock \emph{Brain Informatics}, 3\penalty0 (3):\penalty0 193--203, April
  2016.
\newblock \doi{10.1007/s40708-016-0048-0}.

\bibitem[Febles et~al.(2022)Febles, Ortega, Sosa, and Sahli]{Febles2022}
Elsa~Santos Febles, Marlis~Ontivero Ortega, Michell~Vald{\'{e}}s Sosa, and
  Hichem Sahli.
\newblock Machine learning techniques for the diagnosis of schizophrenia based
  on event-related potentials.
\newblock \emph{Frontiers in Neuroinformatics}, 16, July 2022.
\newblock \doi{10.3389/fninf.2022.893788}.

\bibitem[Feng et~al.(2017)Feng, Sun, and Wang]{Feng2017}
Dexiang Feng, Min Sun, and Xueyong Wang.
\newblock A family of conjugate gradient methods for large-scale nonlinear
  equations.
\newblock \emph{Journal of Inequalities and Applications}, 2017\penalty0 (1),
  September 2017.
\newblock ISSN 1029-242X.
\newblock \doi{10.1186/s13660-017-1510-0}.

\bibitem[Fischl(2012)]{Fischl2012}
Bruce Fischl.
\newblock Freesurfer.
\newblock \emph{NeuroImage}, 62\penalty0 (2):\penalty0 774--781, August 2012.
\newblock ISSN 1053-8119.
\newblock \doi{10.1016/j.neuroimage.2012.01.021}.

\bibitem[Fletcher(1964)]{Fletcher1964}
R.~Fletcher.
\newblock Function minimization by conjugate gradients.
\newblock \emph{The Computer Journal}, 7\penalty0 (2):\penalty0 149--154,
  February 1964.
\newblock ISSN 1460-2067.
\newblock \doi{10.1093/comjnl/7.2.149}.

\bibitem[Franke et~al.(2010)Franke, Ziegler, Kl{\"{o}}ppel, and
  Gaser]{Franke2010}
Katja Franke, Gabriel Ziegler, Stefan Kl{\"{o}}ppel, and Christian Gaser.
\newblock Estimating the age of healthy subjects from t1-weighted {MRI} scans
  using kernel methods: Exploring the influence of various parameters.
\newblock \emph{{NeuroImage}}, 50\penalty0 (3):\penalty0 883--892, April 2010.
\newblock \doi{10.1016/j.neuroimage.2010.01.005}.

\bibitem[Friedman et~al.(2007)Friedman, Hastie, H{\"{o}}fling, and
  Tibshirani]{Friedman2007}
Jerome Friedman, Trevor Hastie, Holger H{\"{o}}fling, and Robert Tibshirani.
\newblock Pathwise coordinate optimization.
\newblock \emph{Annals of Applied Statistics 2007, Vol. 1, No. 2, 302-332},
  August 2007.
\newblock \doi{10.1214/07-AOAS131}.

\bibitem[Gao et~al.(2015)Gao, Ver~Steeg, and Galstyan]{Gao2014}
Shuyang Gao, Greg Ver~Steeg, and Aram Galstyan.
\newblock {Efficient Estimation of Mutual Information for Strongly Dependent
  Variables}.
\newblock In Guy Lebanon and S.~V.~N. Vishwanathan, editors, \emph{Proceedings
  of the Eighteenth International Conference on Artificial Intelligence and
  Statistics}, volume~38 of \emph{Proceedings of Machine Learning Research},
  pages 277--286, San Diego, California, USA, 09--12 May 2015. PMLR.
\newblock URL \url{https://proceedings.mlr.press/v38/gao15.html}.

\bibitem[Garcia and Marder(2017)]{Garcia2017}
Tanya~P. Garcia and Karen Marder.
\newblock Statistical approaches to longitudinal data analysis in
  neurodegenerative diseases: Huntington{\textquoteright}s disease as a model.
\newblock \emph{Current Neurology and Neuroscience Reports}, 17\penalty0 (2),
  February 2017.
\newblock ISSN 1534-6293.
\newblock \doi{10.1007/s11910-017-0723-4}.

\bibitem[Ghadimi and Lan(2013)]{Ghadimi2013}
Saeed Ghadimi and Guanghui Lan.
\newblock Stochastic first- and zeroth-order methods for nonconvex stochastic
  programming.
\newblock \emph{{SIAM} Journal on Optimization}, 23\penalty0 (4):\penalty0
  2341--2368, January 2013.
\newblock \doi{10.1137/120880811}.

\bibitem[Ghadimi and Lan(2015)]{Ghadimi2015}
Saeed Ghadimi and Guanghui Lan.
\newblock Accelerated gradient methods for nonconvex nonlinear and stochastic
  programming.
\newblock \emph{Mathematical Programming}, 156\penalty0 (1-2):\penalty0 59--99,
  February 2015.
\newblock \doi{10.1007/s10107-015-0871-8}.

\bibitem[Ghosh and Thoresen(2016)]{Ghosh2016}
Abhik Ghosh and Magne Thoresen.
\newblock Non-concave penalization in linear mixed-effects models and
  regularized selection of fixed effects.
\newblock \emph{AStA Advances in Statistical Analysis (2018), Volume 102, Issue
  2, pp 179--210}, July 2016.
\newblock \doi{10.1007/s10182-017-0298-z}.

\bibitem[Gilbert and Nocedal(1992)]{Gilbert1992}
Jean~Charles Gilbert and Jorge Nocedal.
\newblock Global convergence properties of conjugate gradient methods for
  optimization.
\newblock \emph{SIAM Journal on Optimization}, 2\penalty0 (1):\penalty0 21--42,
  February 1992.
\newblock ISSN 1095-7189.
\newblock \doi{10.1137/0802003}.

\bibitem[G{\"{o}}del(1931)]{Goedel1931}
Kurt G{\"{o}}del.
\newblock {{\"{U}}}ber formal unentscheidbare s{\"{a}}tze der principia
  mathematica und verwandter systeme i.
\newblock \emph{Monatshefte f{\"{u}}r Mathematik und Physik}, 38--38\penalty0
  (1):\penalty0 173--198, December 1931.
\newblock ISSN 1436-5081.
\newblock \doi{10.1007/bf01700692}.

\bibitem[G{\'{o}}mez-Verdejo et~al.(2019)G{\'{o}}mez-Verdejo,
  Parrado-Hern{\'{a}}ndez, and Tohka]{GomezVerdejo2019}
Vanessa G{\'{o}}mez-Verdejo, Emilio Parrado-Hern{\'{a}}ndez, and Jussi Tohka.
\newblock Sign-consistency based variable importance for machine learning in
  brain imaging.
\newblock \emph{Neuroinformatics}, 17\penalty0 (4):\penalty0 593--609, March
  2019.
\newblock \doi{10.1007/s12021-019-9415-3}.

\bibitem[Hager and Zhang(2006)]{Hager2006}
William Hager and Hongchao Zhang.
\newblock A survey of nonlinear conjugate gradient method.
\newblock 2, January 2006.

\bibitem[Hager and Zhang(2005)]{Hager2005}
William~W. Hager and Hongchao Zhang.
\newblock A new conjugate gradient method with guaranteed descent and an
  efficient line search.
\newblock \emph{SIAM Journal on Optimization}, 16\penalty0 (1):\penalty0
  170--192, January 2005.
\newblock ISSN 1095-7189.
\newblock \doi{10.1137/030601880}.

\bibitem[Hao et~al.(2020)Hao, Bao, Guo, Yu, Zhang, Risacher, Saykin, Yao, and
  Shen]{Hao2020}
Xiaoke Hao, Yongjin Bao, Yingchun Guo, Ming Yu, Daoqiang Zhang, Shannon~L.
  Risacher, Andrew~J. Saykin, Xiaohui Yao, and Li~Shen.
\newblock Multi-modal neuroimaging feature selection with consistent metric
  constraint for diagnosis of alzheimer{\textquotesingle}s disease.
\newblock \emph{Medical Image Analysis}, 60:\penalty0 101625, February 2020.
\newblock \doi{10.1016/j.media.2019.101625}.

\bibitem[He et~al.(2018)He, Xu, and Kang]{He2018}
Kevin He, Han Xu, and Jian Kang.
\newblock A selective overview of feature screening methods with applications
  to neuroimaging data.
\newblock \emph{{WIREs} Computational Statistics}, 11\penalty0 (2), September
  2018.
\newblock \doi{10.1002/wics.1454}.

\bibitem[Helmke(1994)]{Helmke1994}
Uwe Helmke.
\newblock \emph{Optimization and Dynamical Systems}.
\newblock Springer London, London, 1994.
\newblock ISBN 9781447134671.
\newblock \doi{10.1007/978-1-4471-3467-1}.

\bibitem[Hestenes and Stiefel(1952)]{Hestenes1952}
M.R. Hestenes and E.~Stiefel.
\newblock Methods of conjugate gradients for solving linear systems.
\newblock \emph{Journal of Research of the National Bureau of Standards},
  49\penalty0 (6):\penalty0 409, December 1952.
\newblock ISSN 0091-0635.
\newblock \doi{10.6028/jres.049.044}.

\bibitem[Hinton and Salakhutdinov(2006)]{Hinton2006}
G.~E. Hinton and R.~R. Salakhutdinov.
\newblock Reducing the dimensionality of data with neural networks.
\newblock \emph{Science}, 313\penalty0 (5786):\penalty0 504--507, July 2006.
\newblock ISSN 1095-9203.
\newblock \doi{10.1126/science.1127647}.

\bibitem[Hoheisel et~al.(2020)Hoheisel, Laborde, and Oberman]{Hoheisel2020}
Tim Hoheisel, Maxime Laborde, and Adam~M. Oberman.
\newblock A regularization interpretation of the proximal point method for
  weakly convex functions.
\newblock \emph{Journal of Dynamics \& Games}, 2020.
\newblock URL \url{https://api.semanticscholar.org/CorpusID:202607166}.

\bibitem[Hubbard and West(1995)]{Hubbard1995}
John~H. Hubbard and Beverly~H. West.
\newblock \emph{Differential Equations: A Dynamical Systems Approach}.
\newblock Springer New York, 1995.
\newblock ISBN 9781461241928.
\newblock \doi{10.1007/978-1-4612-4192-8}.

\bibitem[Hunt et~al.(2014)Hunt, Weissfeld, Boudreau, Aizenstein, Newman,
  Simonsick, Domelen, Thomas, Yaffe, and Rosano]{Hunt2014}
Megan J.~Olson Hunt, Lisa Weissfeld, Robert~M. Boudreau, Howard Aizenstein,
  Anne~B. Newman, Eleanor~M. Simonsick, Dane R.~Van Domelen, Fridtjof Thomas,
  Kristine Yaffe, and Caterina Rosano.
\newblock A variant of sparse partial least squares for variable selection and
  data exploration.
\newblock \emph{Frontiers in Neuroinformatics}, 8, 2014.
\newblock \doi{10.3389/fninf.2014.00018}.

\bibitem[Ibrahim et~al.(2012)Ibrahim, Abdallah, and N'Guessan]{Ibrahim2012}
Sidi~Zakari Ibrahim, Mkhadri Abdallah, and Assi N'Guessan.
\newblock A mixture of local and quadratic approximation variable selection
  algorithm in nonconcave penalized regression.
\newblock \emph{ARIMA}, 15:\penalty0 18, January 2012.

\bibitem[Ince et~al.(2016)Ince, Giordano, Kayser, Rousselet, Gross, and
  Schyns]{Ince2016}
Robin~A.A. Ince, Bruno~L. Giordano, Christoph Kayser, Guillaume~A. Rousselet,
  Joachim Gross, and Philippe~G. Schyns.
\newblock A statistical framework for neuroimaging data analysis based on
  mutual information estimated via a gaussian copula.
\newblock \emph{Human Brain Mapping}, 38\penalty0 (3):\penalty0 1541--1573,
  November 2016.
\newblock \doi{10.1002/hbm.23471}.

\bibitem[Ivanoska et~al.(2021)Ivanoska, Trivodaliev, Kalajdziski, and
  Zanin]{Ivanoska2021}
Ilinka Ivanoska, Kire Trivodaliev, Slobodan Kalajdziski, and Massimiliano
  Zanin.
\newblock Statistical and machine learning link selection methods for brain
  functional networks: Review and comparison.
\newblock \emph{Brain Sciences}, 11\penalty0 (6):\penalty0 735, May 2021.
\newblock \doi{10.3390/brainsci11060735}.

\bibitem[Jiang et~al.(2020)Jiang, Lu, Chen, Yao, Li, Zhang, and Guo]{Jiang2020}
Huiting Jiang, Na~Lu, Kewei Chen, Li~Yao, Ke~Li, Jiacai Zhang, and Xiaojuan
  Guo.
\newblock Predicting brain age of healthy adults based on structural {MRI}
  parcellation using convolutional neural networks.
\newblock \emph{Frontiers in Neurology}, 10, January 2020.
\newblock \doi{10.3389/fneur.2019.01346}.

\bibitem[Jin et~al.(2017)Jin, Netrapalli, and Jordan]{Jin2017}
Chi Jin, Praneeth Netrapalli, and Michael~I. Jordan.
\newblock Accelerated gradient descent escapes saddle points faster than
  gradient descent.
\newblock November 2017.
\newblock \doi{10.48550/arxiv.1711.10456}.

\bibitem[Jonsson et~al.(2019)Jonsson, Bjornsdottir, Thorgeirsson, Ellingsen,
  Walters, Gudbjartsson, Stefansson, Stefansson, and Ulfarsson]{Jonsson2019}
B.~A. Jonsson, G.~Bjornsdottir, T.~E. Thorgeirsson, L.~M. Ellingsen, G.~Bragi
  Walters, D.~F. Gudbjartsson, H.~Stefansson, K.~Stefansson, and M.~O.
  Ulfarsson.
\newblock Brain age prediction using deep learning uncovers associated sequence
  variants.
\newblock \emph{Nature Communications}, 10\penalty0 (1), November 2019.
\newblock \doi{10.1038/s41467-019-13163-9}.

\bibitem[Kafka and Wilke(2019)]{Kafka2019}
Dominic Kafka and Daniel Wilke.
\newblock Gradient-only line searches: An alternative to probabilistic line
  searches.
\newblock March 2019.
\newblock \doi{10.48550/ARXIV.1903.09383}.

\bibitem[Kanzow and Lechner(2020)]{Kanzow2020}
Christian Kanzow and Theresa Lechner.
\newblock Globalized inexact proximal newton-type methods for nonconvex
  composite functions.
\newblock \emph{Computational Optimization and Applications}, 78\penalty0
  (2):\penalty0 377--410, November 2020.
\newblock ISSN 1573-2894.
\newblock \doi{10.1007/s10589-020-00243-6}.

\bibitem[Karamardian(1976)]{Karamardian1976}
S.~Karamardian.
\newblock Complementarity problems over cones with monotone and pseudomonotone
  maps.
\newblock \emph{Journal of Optimization Theory and Applications}, 18\penalty0
  (4):\penalty0 445--454, April 1976.
\newblock ISSN 1573-2878.
\newblock \doi{10.1007/bf00932654}.

\bibitem[Khan et~al.(2007)Khan, Bandyopadhyay, Ganguly, Saigal, Erickson,
  Protopopescu, and Ostrouchov]{Khan2007}
Shiraj Khan, Sharba Bandyopadhyay, Auroop~R. Ganguly, Sunil Saigal, David~J.
  Erickson, Vladimir Protopopescu, and George Ostrouchov.
\newblock Relative performance of mutual information estimation methods for
  quantifying the dependence among short and noisy data.
\newblock \emph{Physical Review E}, 76\penalty0 (2):\penalty0 026209, August
  2007.
\newblock \doi{10.1103/physreve.76.026209}.

\bibitem[Khera et~al.(2017)Khera, Chaffin, Aragam, Emdin, Klarin, Haas,
  Roselli, Natarajan, and Kathiresan]{Khera2017}
Amit~V. Khera, Mark Chaffin, Krishna~G. Aragam, Connor~A. Emdin, Derek Klarin,
  Mary~E. Haas, Carolina Roselli, Pradeep Natarajan, and Sekar Kathiresan.
\newblock Genome-wide polygenic score to identify a monogenic risk-equivalent
  for coronary disease.
\newblock November 2017.
\newblock \doi{10.1101/218388}.

\bibitem[Kim et~al.(2018)Kim, Lee, and Kwon]{Kim2018}
Dongshin Kim, Sangin Lee, and Sunghoon Kwon.
\newblock A unified algorithm for the non-convex penalized estimation: The
  ncpen package.
\newblock November 2018.
\newblock \doi{10.48550/arxiv.1811.05061}.

\bibitem[Kim et~al.(2008)Kim, Choi, and Oh]{Kim2008}
Yongdai Kim, Hosik Choi, and Hee-Seok Oh.
\newblock Smoothly clipped absolute deviation on high dimensions.
\newblock \emph{Journal of the American Statistical Association}, 103\penalty0
  (484):\penalty0 1665--1673, 2008.
\newblock ISSN 0162-1459.
\newblock URL \url{http://www.jstor.org/stable/27640214}.

\bibitem[Kraskov et~al.(2004)Kraskov, St{\"{o}}gbauer, and
  Grassberger]{Kraskov2004}
Alexander Kraskov, Harald St{\"{o}}gbauer, and Peter Grassberger.
\newblock Estimating mutual information.
\newblock \emph{Physical Review E}, 69\penalty0 (6):\penalty0 066138, June
  2004.
\newblock \doi{10.1103/physreve.69.066138}.

\bibitem[L.~F.~Kozachenko(1987)]{L.F.Kozachenko1987}
N.~N.~Leonenko L.~F.~Kozachenko.
\newblock Sample estimate of the entropy of a~random vector.
\newblock \emph{Probl. Peredachi Inf.}, 23\penalty0 (2):\penalty0 9--16, 1987.
\newblock URL \url{http://mathscinet.ams.org/mathscinet-getitem?mr=908626}.

\bibitem[Lan(2011)]{Lan2011}
Guanghui Lan.
\newblock An optimal method for stochastic composite optimization.
\newblock \emph{Mathematical Programming}, 133\penalty0 (1-2):\penalty0
  365--397, January 2011.
\newblock \doi{10.1007/s10107-010-0434-y}.

\bibitem[Lee et~al.(2014)Lee, Sun, and Saunders]{Lee2014a}
Jason~D. Lee, Yuekai Sun, and Michael~A. Saunders.
\newblock Proximal newton-type methods for minimizing composite functions.
\newblock \emph{{SIAM} Journal on Optimization}, 24\penalty0 (3):\penalty0
  1420--1443, January 2014.
\newblock \doi{10.1137/130921428}.

\bibitem[Lee and Breheny(2015)]{Lee2015}
Sangin Lee and Patrick Breheny.
\newblock Strong rules for nonconvex penalties and their implications for
  efficient algorithms in high-dimensional regression.
\newblock \emph{Journal of Computational and Graphical Statistics}, 24\penalty0
  (4):\penalty0 1074--1091, October 2015.
\newblock \doi{10.1080/10618600.2014.975231}.

\bibitem[Lee et~al.(2016)Lee, Kwon, and Kim]{Lee2016}
Sangin Lee, Sunghoon Kwon, and Yongdai Kim.
\newblock A modified local quadratic approximation algorithm for penalized
  optimization problems.
\newblock \emph{Comput. Stat. Data Anal.}, 94\penalty0 (C):\penalty0 275--286,
  February 2016.
\newblock ISSN 0167-9473.
\newblock \doi{10.1016/j.csda.2015.08.019}.

\bibitem[Li(2022)]{Li2022}
Qiang Li.
\newblock Functional connectivity inference from {fMRI} data using multivariate
  information measures.
\newblock \emph{Neural Networks}, 146:\penalty0 85--97, February 2022.
\newblock \doi{10.1016/j.neunet.2021.11.016}.

\bibitem[Li et~al.(2016)Li, Jiang, Haupt, Arora, Liu, Hong, and Zhao]{Li2016}
Xingguo Li, Haoming Jiang, Jarvis Haupt, Raman Arora, Han Liu, Mingyi Hong, and
  Tuo Zhao.
\newblock On fast convergence of proximal algorithms for sqrt-lasso
  optimization: Don't worry about its nonsmooth loss function, 2016.

\bibitem[Liang et~al.(2022)Liang, Lee, Arefin, Dong, Walczak, Shi, Knoll, Ge,
  Ying, and Zhang]{Liang2022}
Zifei Liang, Choong~H Lee, Tanzil~M Arefin, Zijun Dong, Piotr Walczak, Song-Hai
  Shi, Florian Knoll, Yulin Ge, Leslie Ying, and Jiangyang Zhang.
\newblock Virtual mouse brain histology from multi-contrast mri via deep
  learning.
\newblock \emph{eLife}, 11, January 2022.
\newblock ISSN 2050-084X.
\newblock \doi{10.7554/elife.72331}.

\bibitem[Liem et~al.(2017)Liem, Varoquaux, Kynast, Beyer, Masouleh, Huntenburg,
  Lampe, Rahim, Abraham, Craddock, Riedel-Heller, Luck, Loeffler, Schroeter,
  Witte, Villringer, and Margulies]{Liem2017}
Franziskus Liem, Ga{\"{e}}l Varoquaux, Jana Kynast, Frauke Beyer,
  Shahrzad~Kharabian Masouleh, Julia~M. Huntenburg, Leonie Lampe, Mehdi Rahim,
  Alexandre Abraham, R.~Cameron Craddock, Steffi Riedel-Heller, Tobias Luck,
  Markus Loeffler, Matthias~L. Schroeter, Anja~Veronica Witte, Arno Villringer,
  and Daniel~S. Margulies.
\newblock Predicting brain-age from multimodal imaging data captures cognitive
  impairment.
\newblock \emph{{NeuroImage}}, 148:\penalty0 179--188, March 2017.
\newblock \doi{10.1016/j.neuroimage.2016.11.005}.

\bibitem[Linn et~al.(2016)Linn, Gaonkar, Doshi, Davatzikos, and
  Shinohara]{Linn2016}
Kristin~A. Linn, Bilwaj Gaonkar, Jimit Doshi, Christos Davatzikos, and
  Russell~T. Shinohara.
\newblock Addressing confounding in predictive models with an application to
  neuroimaging.
\newblock \emph{The International Journal of Biostatistics}, 12\penalty0
  (1):\penalty0 31--44, May 2016.
\newblock \doi{10.1515/ijb-2015-0030}.

\bibitem[Lord et~al.(2018)Lord, Sun, and Bollt]{Lord2018}
Warren~M. Lord, Jie Sun, and Erik~M. Bollt.
\newblock Geometric k-nearest neighbor estimation of entropy and mutual
  information.
\newblock \emph{Chaos: An Interdisciplinary Journal of Nonlinear Science},
  28\penalty0 (3), March 2018.
\newblock \doi{10.1063/1.5011683}.

\bibitem[Lubich et~al.(2006)Lubich, Wanner, and Hairer]{Lubich2006}
Christian Lubich, Gerhard Wanner, and Ernst Hairer.
\newblock \emph{Geometric Numerical Integration}.
\newblock Number v.31 in Springer Series in Computational Mathematics Ser.
  Springer Berlin / Heidelberg, Berlin, Heidelberg, 2\textsuperscript{nd} ed.
  edition, 2006.
\newblock ISBN 9783540306665.
\newblock Description based on publisher supplied metadata and other sources.

\bibitem[M.~Tsukada(2005)]{Tsukada2005}
M.~Kato M.~Tsukada, H.~Suyari.
\newblock On the probability distribution maximizing generalized entropies.
\newblock In \emph{Proceedings of 2005 Symposium on Applied Functional Analysis
  - Information Sciences and Related Fields}, pages 99--111, 2005.

\bibitem[Magri et~al.(2009)Magri, Whittingstall, Singh, Logothetis, and
  Panzeri]{Magri2009}
Cesare Magri, Kevin Whittingstall, Vanessa Singh, Nikos~K Logothetis, and
  Stefano Panzeri.
\newblock A toolbox for the fast information analysis of multiple-site {LFP},
  {EEG} and spike train recordings.
\newblock \emph{{BMC} Neuroscience}, 10\penalty0 (1), July 2009.
\newblock \doi{10.1186/1471-2202-10-81}.

\bibitem[Mar{\v{c}}enko and Pastur(1967)]{Marcenko1967}
V~A Mar{\v{c}}enko and L~A Pastur.
\newblock Distribution of eigenvalues for some sets of random matrices.
\newblock \emph{Mathematics of the USSR-Sbornik}, 1\penalty0 (4):\penalty0
  457--483, April 1967.
\newblock ISSN 0025-5734.
\newblock \doi{10.1070/sm1967v001n04abeh001994}.

\bibitem[Martino et~al.(2008)Martino, Valente, Staeren, Ashburner, Goebel, and
  Formisano]{Martino2008}
Federico~De Martino, Giancarlo Valente, No{\"{e}}l Staeren, John Ashburner,
  Rainer Goebel, and Elia Formisano.
\newblock Combining multivariate voxel selection and support vector machines
  for mapping and classification of {fMRI} spatial patterns.
\newblock \emph{{NeuroImage}}, 43\penalty0 (1):\penalty0 44--58, October 2008.
\newblock \doi{10.1016/j.neuroimage.2008.06.037}.

\bibitem[Mazumder et~al.(2011)Mazumder, Friedman, and Hastie]{Mazumder2011}
Rahul Mazumder, Jerome~H. Friedman, and Trevor Hastie.
\newblock Sparsenet: Coordinate descent with nonconvex penalties.
\newblock \emph{Journal of the American Statistical Association}, 106\penalty0
  (495):\penalty0 1125--1138, 2011.
\newblock ISSN 0162-1459.
\newblock URL \url{http://www.jstor.org/stable/23427579}.

\bibitem[Meckes(2021)]{Meckes2021}
Elizabeth Meckes.
\newblock The eigenvalues of random matrices.
\newblock \emph{IMAGE, the Bulletin of the International Linear Algebra
  Society, no. 65, pp. 9-22, 2020}, January 2021.
\newblock \doi{10.48550/arxiv.2101.02928}.

\bibitem[Mehta(2004)]{Mehta2004}
M.~L. Mehta.
\newblock \emph{Random matrices}.
\newblock Number 142 in Pure and applied mathematics (Academic Press).
  Elsevier, Amsterdam, 3\textsuperscript{rd} ed. edition, 2004.
\newblock ISBN 9780080474113.
\newblock Includes bibliographical references and indexes.

\bibitem[Mohr et~al.(2006)Mohr, Puls, Wrase, Heinz, Hochreiter, and
  Obermayer]{Mohr2006}
J.~Mohr, I.~Puls, J.~Wrase, A.~Heinz, S.~Hochreiter, and K.~Obermayer.
\newblock P-{SVM} variable selection for discovering dependencies between
  genetic and brain imaging data.
\newblock In \emph{The 2006 {IEEE} International Joint Conference on Neural
  Network Proceedings}. {IEEE}, 2006.
\newblock \doi{10.1109/ijcnn.2006.247207}.

\bibitem[Moon et~al.(1995)Moon, Rajagopalan, and Lall]{Moon1995}
Young-Il Moon, Balaji Rajagopalan, and Upmanu Lall.
\newblock Estimation of mutual information using kernel density estimators.
\newblock \emph{Physical Review E}, 52\penalty0 (3):\penalty0 2318--2321,
  September 1995.
\newblock \doi{10.1103/physreve.52.2318}.

\bibitem[Morduchovi{\v{c}}(2018)]{Morduchovic2018}
Boris~S. Morduchovi{\v{c}}.
\newblock \emph{Variational analysis and applications}.
\newblock Springer monographs in mathematics. Springer, Cham, Switzerland,
  softcover re-print of the hardcover 1\textsuperscript{st} edition 2018
  edition, 2018.
\newblock ISBN 9783030065133.
\newblock Literaturverzeichnis: Seiten 533-578.

\bibitem[Mordukhovich(2006{\natexlab{a}})]{Mordukhovich2006}
B.~Sh Mordukhovich.
\newblock \emph{Variational Analysis and Generalized Differentiation I: Basic
  theory}.
\newblock Number~1 in Variational analysis and generalized differentiation.
  Springer, Berlin ;, 2006{\natexlab{a}}.
\newblock ISBN 9783540312475.
\newblock Includes bibliographical references and indexes.

\bibitem[Mordukhovich(2006{\natexlab{b}})]{Mordukhovich2006a}
B.~Sh Mordukhovich.
\newblock \emph{Variational Analysis and Generalized Differentiation II:
  Applications}.
\newblock Number~2 in Variational analysis and differentiation. Springer, New
  York, 2006{\natexlab{b}}.
\newblock ISBN 9783540312468.
\newblock Includes bibliographical references and index.

\bibitem[Moreau(1965)]{Moreau1965}
J.J. Moreau.
\newblock Proximit{\'{e}} et dualit{\'{e}} dans un espace hilbertien.
\newblock \emph{Bulletin de la Soci{\'{e}}t{\'{e}} Math{\'{e}}matique de
  France}, 93:\penalty0 273--299, 1965.
\newblock URL \url{http://eudml.org/doc/87067}.

\bibitem[Nemirovsky et~al.(2023)Nemirovsky, Popiel, Rudas, Caius, Naci, Schiff,
  Owen, and Soddu]{Nemirovsky2023}
Idan~E. Nemirovsky, Nicholas J.~M. Popiel, Jorge Rudas, Matthew Caius, Lorina
  Naci, Nicholas~D. Schiff, Adrian~M. Owen, and Andrea Soddu.
\newblock An implementation of integrated information theory in resting-state
  {fMRI}.
\newblock \emph{Communications Biology}, 6\penalty0 (1), July 2023.
\newblock \doi{10.1038/s42003-023-05063-y}.

\bibitem[Nesterov(1983)]{Nesterov1983}
Y.~E. Nesterov.
\newblock A method for solving the convex programming problem with convergence
  rate $o(1/k^2)$.
\newblock \emph{Dokl. Akad. Nauk SSSR}, 269:\penalty0 543--547, 1983.
\newblock URL \url{https://ci.nii.ac.jp/naid/10029946121/en/}.

\bibitem[Nesterov(2004{\natexlab{a}})]{Nesterov2004a}
Yu. Nesterov.
\newblock Smooth minimization of non-smooth functions.
\newblock \emph{Mathematical Programming}, 103\penalty0 (1):\penalty0 127--152,
  December 2004{\natexlab{a}}.
\newblock \doi{10.1007/s10107-004-0552-5}.

\bibitem[Nesterov(2012)]{Nesterov2012}
Yu. Nesterov.
\newblock Efficiency of coordinate descent methods on huge-scale optimization
  problems.
\newblock \emph{{SIAM} Journal on Optimization}, 22\penalty0 (2):\penalty0
  341--362, January 2012.
\newblock \doi{10.1137/100802001}.

\bibitem[Nesterov(2004{\natexlab{b}})]{Nesterov2004}
Yurii Nesterov.
\newblock \emph{Introductory Lectures on Convex Optimization}.
\newblock Springer {US}, 2004{\natexlab{b}}.
\newblock \doi{10.1007/978-1-4419-8853-9}.

\bibitem[Neumaier et~al.(2024)Neumaier, Kimiaei, and Azmi]{Neumaier2024}
Arnold Neumaier, Morteza Kimiaei, and Behzad Azmi.
\newblock Globally linearly convergent nonlinear conjugate gradients without
  wolfe line search.
\newblock \emph{Numerical Algorithms}, February 2024.
\newblock ISSN 1572-9265.
\newblock \doi{10.1007/s11075-024-01764-5}.

\bibitem[Nikolova(2000)]{Nikolova2000}
Mila Nikolova.
\newblock Local strong homogeneity of a regularized estimator.
\newblock \emph{{SIAM} Journal on Applied Mathematics}, 61\penalty0
  (2):\penalty0 633--658, January 2000.
\newblock \doi{10.1137/s0036139997327794}.

\bibitem[Nocedal et~al.(2000)Nocedal, Wright, and Wright]{Nocedal2000}
J.~Nocedal, S.~Wright, and Stephen~J. Wright.
\newblock \emph{Numerical Optimization}.
\newblock Springer Series in Operations Research and Financial Engineering.
  Springer New York, New York, NY, second edition edition, 2000.
\newblock ISBN 9780387987934.
\newblock URL \url{https://books.google.ca/books?id=epc5fX0lqRIC}.

\bibitem[Odland(2018)]{Odland2018}
Tommy Odland.
\newblock tommyod/kdepy: Kernel density estimation in python, 2018.

\bibitem[P\'{a}l et~al.(2010)P\'{a}l, P\'{o}czos, and Szepesv\'{a}ri]{Pal2010}
D\'{a}vid P\'{a}l, Barnab\'{a}s P\'{o}czos, and Csaba Szepesv\'{a}ri.
\newblock Estimation of r\'{e}nyi entropy and mutual information based on
  generalized nearest-neighbor graphs.
\newblock In \emph{Proceedings of the 23\textsuperscript{rd} International
  Conference on Neural Information Processing Systems - Volume 2}, NIPS'10,
  pages 1849--1857, Red Hook, NY, USA, 2010. Curran Associates Inc.

\bibitem[Pan and Shang(2018)]{Pan2018}
Juming Pan and Junfeng Shang.
\newblock A simultaneous variable selection methodology for linear mixed
  models.
\newblock \emph{Journal of Statistical Computation and Simulation}, 88\penalty0
  (17):\penalty0 3323--3337, 2018.
\newblock \doi{10.1080/00949655.2018.1515948}.

\bibitem[Paquette et~al.(2020)Paquette, van Merri{\"{e}}nboer, Paquette, and
  Pedregosa]{Paquette2020}
Courtney Paquette, Bart van Merri{\"{e}}nboer, Elliot Paquette, and Fabian
  Pedregosa.
\newblock Halting time is predictable for large models: A universality property
  and average-case analysis.
\newblock June 2020.
\newblock \doi{10.48550/arxiv.2006.04299}.

\bibitem[Parnell et~al.(2020)Parnell, D{\"{u}}nner, Atasu, Sifalakis, and
  Pozidis]{Parnell2020}
Thomas Parnell, Celestine D{\"{u}}nner, Kubilay Atasu, Manolis Sifalakis, and
  Haralampos Pozidis.
\newblock Tera-scale coordinate descent on {GPUs}.
\newblock \emph{Future Generation Computer Systems}, 108:\penalty0 1173--1191,
  July 2020.
\newblock \doi{10.1016/j.future.2018.04.072}.

\bibitem[Pashayan et~al.(2015)Pashayan, Duffy, Neal, Hamdy, Donovan, Martin,
  Harrington, Benlloch, Olama, Shah, Kote-Jarai, Easton, Eeles, and
  Pharoah]{Pashayan2015}
Nora Pashayan, Stephen~W. Duffy, David~E. Neal, Freddie~C. Hamdy, Jenny~L.
  Donovan, Richard~M. Martin, Patricia Harrington, Sara Benlloch, Ali Amin~Al
  Olama, Mitul Shah, Zsofia Kote-Jarai, Douglas~F. Easton, Rosalind Eeles, and
  Paul~D. Pharoah.
\newblock Implications of polygenic risk-stratified screening for prostate
  cancer on overdiagnosis.
\newblock \emph{Genetics in Medicine}, 17\penalty0 (10):\penalty0 789--795,
  January 2015.
\newblock \doi{10.1038/gim.2014.192}.

\bibitem[Pe{\~{n}}a et~al.(1999)Pe{\~{n}}a, Rubio, and Serna]{Pena1999}
Ignacio Pe{\~{n}}a, Gonzalo Rubio, and Gregorio Serna.
\newblock Why do we smile? on the determinants of the implied volatility
  function.
\newblock \emph{Journal of Banking \& Finance}, 23\penalty0 (8):\penalty0
  1151--1179, August 1999.
\newblock ISSN 0378-4266.
\newblock \doi{10.1016/s0378-4266(98)00134-4}.

\bibitem[Pereda et~al.(2018)Pereda, Garc{\'{\i}}a-Torres, Meli{\'{a}}n-Batista,
  Ma{\~{n}}as, M{\'{e}}ndez, and Gonz{\'{a}}lez]{Pereda2018}
Ernesto Pereda, Miguel Garc{\'{\i}}a-Torres, Bel{\'{e}}n Meli{\'{a}}n-Batista,
  Soledad Ma{\~{n}}as, Leopoldo M{\'{e}}ndez, and Juli{\'{a}}n~J.
  Gonz{\'{a}}lez.
\newblock The blessing of dimensionality: Feature selection outperforms
  functional connectivity-based feature transformation to classify {ADHD}
  subjects from {EEG} patterns of phase synchronisation.
\newblock \emph{{PLOS} {ONE}}, 13\penalty0 (8):\penalty0 e0201660, August 2018.
\newblock \doi{10.1371/journal.pone.0201660}.

\bibitem[Polak and Ribiere(1969)]{Polak1969}
Eric Polak and G~Ribiere.
\newblock Note sur la convergence de m{\'{e}}thodes de directions
  conjugu{\'{e}}es.
\newblock \emph{ESAIM: Mathematical Modelling and Numerical Analysis -
  Mod{\'{e}}lisation Math{\'{e}}matique et Analyse Num{\'{e}}rique}, 3\penalty0
  (R1):\penalty0 35--43, 1969.
\newblock URL \url{http://eudml.org/doc/193115}.

\bibitem[Polyak(1964)]{Polyak1964}
B.T. Polyak.
\newblock Some methods of speeding up the convergence of iteration methods.
\newblock \emph{{USSR} Computational Mathematics and Mathematical Physics},
  4\penalty0 (5):\penalty0 1--17, January 1964.
\newblock \doi{10.1016/0041-5553(64)90137-5}.

\bibitem[Qian(1999)]{Qian1999}
Ning Qian.
\newblock On the momentum term in gradient descent learning algorithms.
\newblock \emph{Neural Networks}, 12\penalty0 (1):\penalty0 145--151, January
  1999.
\newblock \doi{10.1016/s0893-6080(98)00116-6}.

\bibitem[Quarteroni et~al.(2007)Quarteroni, Sacco, and Saleri]{Quarteroni2007}
Alfio Quarteroni, Riccardo Sacco, and Fausto Saleri.
\newblock \emph{Numerical Mathematics}.
\newblock Springer New York, 2007.
\newblock ISBN 9780387227504.
\newblock \doi{10.1007/b98885}.

\bibitem[Rabhi and Bouezmarni(2019)]{Rabhi2019}
Yassir Rabhi and Taoufik Bouezmarni.
\newblock Nonparametric inference for copulas and measures of dependence under
  length-biased sampling and informative censoring.
\newblock \emph{Journal of the American Statistical Association}, 115\penalty0
  (531):\penalty0 1268--1278, June 2019.
\newblock ISSN 1537-274X.
\newblock \doi{10.1080/01621459.2019.1611586}.

\bibitem[Racine(2022)]{Racine2022}
Jeffrey Racine.
\newblock A primer on regression splines, 2022.
\newblock URL
  \url{https://cran.r-project.org/web/packages/crs/vignettes/spline_primer.pdf}.
\newblock PDF document.

\bibitem[R{\'{e}}nyi(1959)]{Renyi1959}
A.~R{\'{e}}nyi.
\newblock On measures of dependence.
\newblock \emph{Acta Mathematica Academiae Scientiarum Hungaricae}, 10\penalty0
  (3--4):\penalty0 441--451, September 1959.
\newblock ISSN 1588-2632.
\newblock \doi{10.1007/bf02024507}.

\bibitem[Reshef et~al.(2011)Reshef, Reshef, Finucane, Grossman, McVean,
  Turnbaugh, Lander, Mitzenmacher, and Sabeti]{Reshef2011}
David~N. Reshef, Yakir~A. Reshef, Hilary~K. Finucane, Sharon~R. Grossman,
  Gilean McVean, Peter~J. Turnbaugh, Eric~S. Lander, Michael Mitzenmacher, and
  Pardis~C. Sabeti.
\newblock Detecting novel associations in large data sets.
\newblock \emph{Science}, 334\penalty0 (6062):\penalty0 1518--1524, December
  2011.
\newblock ISSN 1095-9203.
\newblock \doi{10.1126/science.1205438}.

\bibitem[Riahi and Qattan(2018)]{Riahi2018}
Mohamed~Kamel Riahi and Issam~Al Qattan.
\newblock Linearly convergent nonlinear conjugate gradient methods for a
  parameter identification problems.
\newblock June 2018.
\newblock \doi{10.48550/ARXIV.1806.10197}.

\bibitem[Riutort-Mayol et~al.(2022)Riutort-Mayol, B{\"{u}}rkner, Andersen,
  Solin, and Vehtari]{RiutortMayol2022}
Gabriel Riutort-Mayol, Paul-Christian B{\"{u}}rkner, Michael~R. Andersen, Arno
  Solin, and Aki Vehtari.
\newblock Practical hilbert space approximate bayesian gaussian processes for
  probabilistic programming.
\newblock \emph{Statistics and Computing}, 33\penalty0 (1), December 2022.
\newblock ISSN 1573-1375.
\newblock \doi{10.1007/s11222-022-10167-2}.

\bibitem[Rockafellar and Wets(2010)]{Rockafellar2010}
Ralph~Tyrrell Rockafellar and Roger J.-B. Wets.
\newblock \emph{Variational analysis}.
\newblock Number 317 in Die @Grundlehren der mathematischen Wissenschaften in
  Einzeldarstellungen. Springer, Heidelberg, corr. 3. printing. [softcover
  version of original hardcover edition 1998] edition, 2010.
\newblock ISBN 3642083048.

\bibitem[Ross(2014)]{Ross2014}
Brian~C. Ross.
\newblock Mutual information between discrete and continuous data sets.
\newblock \emph{{PLoS} {ONE}}, 9\penalty0 (2):\penalty0 e87357, February 2014.
\newblock \doi{10.1371/journal.pone.0087357}.

\bibitem[Ross(2019)]{Ross2019}
I.~M. Ross.
\newblock An optimal control theory for accelerated optimization.
\newblock February 2019.
\newblock \doi{10.48550/ARXIV.1902.09004}.

\bibitem[Roy(2021)]{Roy2021}
Arkaprava Roy.
\newblock Nonparametric group variable selection with multivariate response for
  connectome-based modeling of cognitive scores.
\newblock October 2021.
\newblock \doi{10.48550/ARXIV.2110.05641}.

\bibitem[Runcie and Crawford(2019)]{Runcie2019}
Daniel~E. Runcie and Lorin Crawford.
\newblock Fast and flexible linear mixed models for genome-wide genetics.
\newblock \emph{{PLOS} Genetics}, 15\penalty0 (2):\penalty0 e1007978, February
  2019.
\newblock ISSN 1553-7404.
\newblock \doi{10.1371/journal.pgen.1007978}.

\bibitem[Saad(2003)]{Saad2003}
Yousef Saad.
\newblock \emph{Iterative Methods for Sparse Linear Systems}.
\newblock Society for Industrial and Applied Mathematics, January 2003.
\newblock ISBN 9780898718003.
\newblock \doi{10.1137/1.9780898718003}.

\bibitem[Schelldorfer(2011)]{Schelldorfer2011}
J{\"{u}}rg Schelldorfer.
\newblock \emph{High-Dimensional Gaussian and Generalized Linear Mixed Models}.
\newblock PhD thesis, ETH Zurich, Z{\"{u}}rich, 2011.

\bibitem[Schl{\"o}gl et~al.(2002)Schl{\"o}gl, Neuper, and
  Pfurtscheller]{Schloegl2002}
A.~Schl{\"o}gl, C.~Neuper, and G.~Pfurtscheller.
\newblock Estimating the mutual information of an eeg-based brain-computer
  interface.
\newblock \emph{Biomedical Engineering}, 47\penalty0 (1-2):\penalty0 3--8,
  2002.

\bibitem[Shannon(1948)]{Shannon1948}
C.~E. Shannon.
\newblock A mathematical theory of communication.
\newblock \emph{Bell System Technical Journal}, 27\penalty0 (3):\penalty0
  379--423, July 1948.
\newblock \doi{10.1002/j.1538-7305.1948.tb01338.x}.

\bibitem[Sherif et~al.(2014)Sherif, Rioux, Rousseau, Kassis, Beck, Adalat, Das,
  Glatard, and Evans]{Sherif2014}
Tarek Sherif, Pierre Rioux, Marc-Etienne Rousseau, Nicolas Kassis, Natacha
  Beck, Reza Adalat, Samir Das, Tristan Glatard, and Alan~C. Evans.
\newblock Cbrain: a web-based, distributed computing platform for collaborative
  neuroimaging research.
\newblock \emph{Frontiers in Neuroinformatics}, 8, May 2014.
\newblock ISSN 1662-5196.
\newblock \doi{10.3389/fninf.2014.00054}.

\bibitem[Shi et~al.(2018)Shi, Du, Jordan, and Su]{Shi2018a}
Bin Shi, Simon~S. Du, Michael~I. Jordan, and Weijie~J. Su.
\newblock Understanding the acceleration phenomenon via high-resolution
  differential equations.
\newblock October 2018.
\newblock \doi{10.48550/arxiv.1810.08907}.

\bibitem[Shi and Shen(2005)]{Shi2005}
Zhen-Jun Shi and Jie Shen.
\newblock Convergence of descent method without line search.
\newblock \emph{Applied Mathematics and Computation}, 167\penalty0
  (1):\penalty0 94--107, August 2005.
\newblock ISSN 0096-3003.
\newblock \doi{10.1016/j.amc.2004.06.097}.

\bibitem[Silverman(1982)]{Silverman1982}
B.~W. Silverman.
\newblock Algorithm {AS} 176: Kernel density estimation using the fast fourier
  transform.
\newblock \emph{Applied Statistics}, 31\penalty0 (1):\penalty0 93, 1982.
\newblock \doi{10.2307/2347084}.

\bibitem[Simon et~al.(2013)Simon, Friedman, Hastie, and Tibshirani]{Simon2013}
Noah Simon, Jerome Friedman, Trevor Hastie, and Robert Tibshirani.
\newblock A sparse-group lasso.
\newblock \emph{Journal of Computational and Graphical Statistics}, 22\penalty0
  (2):\penalty0 231--245, April 2013.
\newblock \doi{10.1080/10618600.2012.681250}.

\bibitem[Snyman(1985)]{Snyman1985}
J~A Snyman.
\newblock Unconstrained minimization by combining the dynamic and conjugate
  gradient methods.
\newblock \emph{Quaestiones Mathematicae}, 8\penalty0 (1):\penalty0 33--42,
  January 1985.
\newblock ISSN 1727-933X.
\newblock \doi{10.1080/16073606.1985.9631898}.

\bibitem[Snyman(2004)]{Snyman2004}
J.~A. Snyman.
\newblock A gradient-only line search method for the conjugate gradient method
  applied to constrained optimization problems with severe noise in the
  objective function.
\newblock \emph{International Journal for Numerical Methods in Engineering},
  62\penalty0 (1):\penalty0 72--82, 2004.
\newblock ISSN 1097-0207.
\newblock \doi{10.1002/nme.1189}.

\bibitem[Sofer et~al.(2014)Sofer, Dicker, and Lin]{Sofer2014}
Tamar Sofer, Lee Dicker, and Xihong Lin.
\newblock Variable selection for high dimensional multivariate outcomes.
\newblock \emph{Statistica Sinica}, 2014.
\newblock \doi{10.5705/ss.2013.019}.

\bibitem[Sontag(1998)]{Sontag1998}
Eduardo~D. Sontag.
\newblock \emph{Mathematical Control Theory}.
\newblock Springer eBook Collection. Springer New York, New York, NY, second
  edition edition, 1998.
\newblock ISBN 9781461205777.
\newblock \doi{10.1007/978-1-4612-0577-7}.

\bibitem[Spall(2012)]{Spall2012}
James~C. Spall.
\newblock Cyclic seesaw process for optimization and identification.
\newblock 154\penalty0 (1):\penalty0 187--208, March 2012.
\newblock \doi{10.1007/s10957-012-0001-1}.

\bibitem[Speed(2011)]{Speed2011}
Terry Speed.
\newblock A correlation for the 21\textsuperscript{st} century.
\newblock \emph{Science}, 334\penalty0 (6062):\penalty0 1502--1503, December
  2011.
\newblock ISSN 1095-9203.
\newblock \doi{10.1126/science.1215894}.

\bibitem[Sprungk(2017)]{Sprungk2017}
Bj{\"{o}}rn Sprungk.
\newblock \emph{Numerical methods for Bayesian inference in Hilbert spaces}.
\newblock Universit{\"{a}}tsverlag Chemnitz, Chemnitz, 2017.

\bibitem[Stein and Shakarchi(2003)]{Stein2003}
Elias~M. Stein and Rami Shakarchi.
\newblock \emph{Fourier Analysis: An Introduction}, volume~1.
\newblock Princeton University Press, Princeton, 15. druck edition, 2003.
\newblock ISBN 9780691113845.
\newblock Hier auch sp{\"{a}}ter erschienene, unver{\"{a}}nderte Nachdrucke.

\bibitem[Steuer et~al.(2002)Steuer, Kurths, Daub, Weise, and
  Selbig]{Steuer2002}
R.~Steuer, J.~Kurths, C.~O. Daub, J.~Weise, and J.~Selbig.
\newblock The mutual information: Detecting and evaluating dependencies between
  variables.
\newblock \emph{Bioinformatics}, 18\penalty0 (suppl{\_}2):\penalty0 S231--S240,
  October 2002.
\newblock \doi{10.1093/bioinformatics/18.suppl_2.s231}.

\bibitem[Su et~al.(2014)Su, Boyd, and Cand\`{e}s]{Su2014}
Weijie Su, Stephen Boyd, and Emmanuel~J. Cand\`{e}s.
\newblock A differential equation for modeling nesterov's accelerated gradient
  method: Theory and insights.
\newblock In \emph{Proceedings of the 27\textsuperscript{th} International
  Conference on Neural Information Processing Systems - Volume 2}, NIPS'14,
  pages 2510--2518, Cambridge, MA, USA, 2014. MIT Press.

\bibitem[Sun and Zhang(2001)]{Sun2001}
Jie Sun and Jiapu Zhang.
\newblock Global convergence of conjugate gradient methods without line search.
\newblock \emph{Annals of Operations Research}, 103\penalty0 (1/4):\penalty0
  161--173, 2001.
\newblock ISSN 0254-5330.
\newblock \doi{10.1023/a:1012903105391}.

\bibitem[Suresh et~al.(2022)Suresh, Newton, Everett, Lin, and
  Duerstock]{Suresh2022}
Shruthi Suresh, David~T. Newton, Thomas~H. Everett, Guang Lin, and Bradley~S.
  Duerstock.
\newblock Feature selection techniques for a machine learning model to detect
  autonomic dysreflexia.
\newblock \emph{Frontiers in Neuroinformatics}, 16, August 2022.
\newblock \doi{10.3389/fninf.2022.901428}.

\bibitem[Tibshirani(1996)]{Tibshirani1996}
Robert Tibshirani.
\newblock Regression shrinkage and selection via the lasso.
\newblock \emph{Journal of the Royal Statistical Society. Series B
  (Methodological)}, 58\penalty0 (1):\penalty0 267--288, 1996.
\newblock ISSN 0035-9246.
\newblock URL \url{http://www.jstor.org/stable/2346178}.

\bibitem[Tibshirani et~al.(2011)Tibshirani, Bien, Friedman, Hastie, Simon,
  Taylor, and Tibshirani]{Tibshirani2011}
Robert Tibshirani, Jacob Bien, Jerome Friedman, Trevor Hastie, Noah Simon,
  Jonathan Taylor, and Ryan~J. Tibshirani.
\newblock Strong rules for discarding predictors in lasso-type problems.
\newblock \emph{Journal of the Royal Statistical Society: Series B (Statistical
  Methodology)}, 74\penalty0 (2):\penalty0 245--266, November 2011.
\newblock \doi{10.1111/j.1467-9868.2011.01004.x}.

\bibitem[Trefethen and Bau(2022)]{Trefethen2022}
Lloyd~N. Trefethen and David Bau.
\newblock \emph{Numerical Linear Algebra, Twenty-fifth Anniversary Edition}.
\newblock Society for Industrial and Applied Mathematics, January 2022.
\newblock ISBN 9781611977165.
\newblock \doi{10.1137/1.9781611977165}.

\bibitem[Tsai et~al.(1999)Tsai, Fisher, Wible, Wells, Kim, and
  Willsky]{Tsai1999}
Andy Tsai, John~W. Fisher, Cindy Wible, William~M. Wells, Junmo Kim, and
  Alan~S. Willsky.
\newblock Analysis of functional {MRI} data using mutual information.
\newblock In \emph{Medical Image Computing and Computer-Assisted Intervention
  {\textendash} {MICCAI}'99}, pages 473--480. Springer Berlin Heidelberg, 1999.
\newblock \doi{10.1007/10704282_51}.

\bibitem[Tsallis(1988)]{Tsallis1988}
Constantino Tsallis.
\newblock Possible generalization of boltzmann-gibbs statistics.
\newblock \emph{Journal of Statistical Physics}, 52\penalty0 (1-2):\penalty0
  479--487, July 1988.
\newblock \doi{10.1007/bf01016429}.

\bibitem[van~den Berg et~al.(2023)van~den Berg, Gameiro, Lessard, James, and
  Mischaikow]{Berg2023}
J.B. van~den Berg, M.~Gameiro, J.-P. Lessard, J.D.~Mireles James, and
  K.~Mischaikow.
\newblock Ordinary differential equations: A constructive approach.
\newblock Lecture notes available from Florida Atlantic University, 2023.
\newblock URL
  \url{https://cosweb1.fau.edu/~jmirelesjames/ODE_course/lectureNotes_version3.pdf}.

\bibitem[Victor(2002)]{Victor2002}
Jonathan~D. Victor.
\newblock Binless strategies for estimation of information from neural data.
\newblock \emph{Physical Review E}, 66\penalty0 (5):\penalty0 051903, November
  2002.
\newblock \doi{10.1103/physreve.66.051903}.

\bibitem[Vignat and Plastino(2005)]{Vignat2005}
C.~Vignat and A.~Plastino.
\newblock The p-sphere and the geometric substratum of power-law probability
  distributions.
\newblock \emph{Physics Letters A}, 343\penalty0 (6):\penalty0 411--416, August
  2005.
\newblock \doi{10.1016/j.physleta.2005.05.027}.

\bibitem[Vignat and Plastino(2007)]{Vignat2007b}
C.~Vignat and A.~Plastino.
\newblock Scale invariance and related properties of q-gaussian systems.
\newblock \emph{Physics Letters A}, 365\penalty0 (5-6):\penalty0 370--375, June
  2007.
\newblock \doi{10.1016/j.physleta.2007.02.003}.

\bibitem[Vignat and Plastino(2009)]{Vignat2009}
C.~Vignat and A.~Plastino.
\newblock Why is the detection of $q-$gaussian behavior such a common
  occurrence?
\newblock \emph{Physica A: Statistical Mechanics and its Applications},
  388\penalty0 (5):\penalty0 601--608, March 2009.
\newblock ISSN 0378-4371.
\newblock \doi{10.1016/j.physa.2008.11.001}.

\bibitem[Vignat et~al.(2004)Vignat, {Hero III}, and Costa]{Vignat2004}
C~Vignat, A.O {Hero III}, and J.A Costa.
\newblock About closedness by convolution of the tsallis maximizers.
\newblock \emph{Physica A: Statistical Mechanics and its Applications},
  340\penalty0 (1-3):\penalty0 147--152, September 2004.
\newblock \doi{10.1016/j.physa.2004.04.001}.

\bibitem[Visscher et~al.(2017)Visscher, Wray, Zhang, Sklar, McCarthy, Brown,
  and Yang]{Visscher2017}
Peter~M. Visscher, Naomi~R. Wray, Qian Zhang, Pamela Sklar, Mark~I. McCarthy,
  Matthew~A. Brown, and Jian Yang.
\newblock 10 years of gwas discovery: Biology, function, and translation.
\newblock \emph{The American Journal of Human Genetics}, 101\penalty0
  (1):\penalty0 5--22, July 2017.
\newblock ISSN 0002-9297.
\newblock \doi{10.1016/j.ajhg.2017.06.005}.

\bibitem[Wang(2006)]{Wang2006}
Cheng-jing Wang.
\newblock Some remarks on conjugate gradient methods without line search.
\newblock \emph{Applied Mathematics and Computation}, 181\penalty0
  (1):\penalty0 370--379, October 2006.
\newblock ISSN 0096-3003.
\newblock \doi{10.1016/j.amc.2006.01.040}.

\bibitem[Wang et~al.(2018)Wang, Aragam, and Xing]{Wang2018a}
Haohan Wang, Bryon Aragam, and Eric~P. Xing.
\newblock Variable selection in heterogeneous datasets: A truncated-rank sparse
  linear mixed model with applications to genome-wide association studies.
\newblock \emph{Methods}, 145:\penalty0 2--9, 2018.
\newblock ISSN 1046-2023.
\newblock \doi{10.1016/j.ymeth.2018.04.021}.
\newblock URL
  \url{http://www.sciencedirect.com/science/article/pii/S1046202317304917}.
\newblock Data mining methods for analyzing biological data in terms of
  phenotypes.

\bibitem[Wang et~al.(2013)Wang, Kim, and Li]{Wang2013}
Lan Wang, Yongdai Kim, and Runze Li.
\newblock Calibrating nonconvex penalized regression in ultra-high dimension.
\newblock \emph{Annals of Statistics 2013, Vol. 41, No. 5, 2505-2536}, November
  2013.
\newblock \doi{10.1214/13-AOS1159}.

\bibitem[Watt(2020)]{Watt2020}
Jeremy Watt.
\newblock \emph{Machine learning refined}.
\newblock Cambridge University Press, New York, second edition edition, 2020.
\newblock ISBN 1108690939.
\newblock Includes bibliographical references and index.

\bibitem[Wigner(1955)]{Wigner1955}
Eugene~P. Wigner.
\newblock Characteristic vectors of bordered matrices with infinite dimensions.
\newblock \emph{The Annals of Mathematics}, 62\penalty0 (3):\penalty0 548,
  November 1955.
\newblock ISSN 0003-486X.
\newblock \doi{10.2307/1970079}.

\bibitem[Wigner(1958)]{Wigner1958}
Eugene~P. Wigner.
\newblock On the distribution of the roots of certain symmetric matrices.
\newblock \emph{The Annals of Mathematics}, 67\penalty0 (2):\penalty0 325,
  March 1958.
\newblock ISSN 0003-486X.
\newblock \doi{10.2307/1970008}.

\bibitem[Wu(2011)]{Wu2011a}
Qing-jun Wu.
\newblock A nonlinear conjugate gradient method without line search and its
  global convergence.
\newblock In \emph{2011 International Conference on Computational and
  Information Sciences}. IEEE, October 2011.
\newblock \doi{10.1109/iccis.2011.45}.

\bibitem[Xiong and Shang(2019)]{Xiong2019}
Jingwei Xiong and Junfeng Shang.
\newblock A penalized approach to mixed model selection via cross-validation.
\newblock \emph{Communications in Statistics - Theory and Methods}, 0\penalty0
  (0):\penalty0 1--27, 2019.
\newblock \doi{10.1080/03610926.2019.1669806}.

\bibitem[Yang et~al.(2017)Yang, Zeng, Goddard, Wray, and Visscher]{Yang2017}
Jian Yang, Jian Zeng, Michael~E Goddard, Naomi~R Wray, and Peter~M Visscher.
\newblock Concepts, estimation and interpretation of {SNP}-based heritability.
\newblock \emph{Nature Genetics}, 49\penalty0 (9):\penalty0 1304--1310,
  September 2017.
\newblock \doi{10.1038/ng.3941}.

\bibitem[Yang et~al.(2024)Yang, Asgharian, and Bhatnagar]{Yang2024}
Kai Yang, Masoud Asgharian, and Sahir Bhatnagar.
\newblock Accelerated gradient methods for sparse statistical learning with
  nonconvex penalties.
\newblock \emph{Statistics and Computing}, 34\penalty0 (1), January 2024.
\newblock ISSN 1573-1375.
\newblock \doi{10.1007/s11222-023-10371-8}.

\bibitem[Yang and Zou(2014)]{Yang2014}
Yi~Yang and Hui Zou.
\newblock A fast unified algorithm for solving group-lasso penalize learning
  problems.
\newblock \emph{Statistics and Computing}, 25\penalty0 (6):\penalty0
  1129--1141, August 2014.
\newblock \doi{10.1007/s11222-014-9498-5}.

\bibitem[Yu et~al.(2015)Yu, Won, Lee, Lim, and Yoon]{Yu2015}
Donghyeon Yu, Joong-Ho Won, Taehoon Lee, Johan Lim, and Sungroh Yoon.
\newblock High-dimensional fused lasso regression using
  majorization{\textendash}minimization and parallel processing.
\newblock \emph{Journal of Computational and Graphical Statistics}, 24\penalty0
  (1):\penalty0 121--153, January 2015.
\newblock \doi{10.1080/10618600.2013.878662}.

\bibitem[Yu and Peng(2017)]{Yu2017}
Yongchao Yu and Jigen Peng.
\newblock The moreau envelope based efficient first-order methods for sparse
  recovery.
\newblock \emph{Journal of Computational and Applied Mathematics},
  322:\penalty0 109--128, October 2017.
\newblock ISSN 0377-0427.
\newblock \doi{10.1016/j.cam.2017.03.014}.

\bibitem[Zhang(2010)]{Zhang2010a}
Cun-Hui Zhang.
\newblock Nearly unbiased variable selection under minimax concave penalty.
\newblock \emph{Annals of Statistics 2010, Vol. 38, No. 2, 894-942}, February
  2010.
\newblock \doi{10.1214/09-AOS729}.

\bibitem[Zhang et~al.(2010)Zhang, Ersoz, Lai, Todhunter, Tiwari, Gore,
  Bradbury, Yu, Arnett, Ordovas, and Buckler]{Zhang2010}
Zhiwu Zhang, Elhan Ersoz, Chao-Qiang Lai, Rory~J Todhunter, Hemant~K Tiwari,
  Michael~A Gore, Peter~J Bradbury, Jianming Yu, Donna~K Arnett, Jose~M
  Ordovas, and Edward~S Buckler.
\newblock Mixed linear model approach adapted for genome-wide association
  studies.
\newblock \emph{Nature Genetics}, 42\penalty0 (4):\penalty0 355--360, March
  2010.
\newblock \doi{10.1038/ng.546}.

\bibitem[Zhou(2009)]{Zhou2009}
Guangming Zhou.
\newblock A descent algorithm without line search for unconstrained
  optimization.
\newblock \emph{Applied Mathematics and Computation}, 215\penalty0
  (7):\penalty0 2528--2533, December 2009.
\newblock ISSN 0096-3003.
\newblock \doi{10.1016/j.amc.2009.08.058}.

\bibitem[Zou and Hastie(2005)]{Zou2005}
Hui Zou and Trevor Hastie.
\newblock Regularization and variable selection via the elastic net.
\newblock \emph{Journal of the Royal Statistical Society. Series B (Statistical
  Methodology)}, 67\penalty0 (2):\penalty0 301--320, March 2005.
\newblock ISSN 13697412, 14679868.
\newblock \doi{10.1111/j.1467-9868.2005.00503.x}.
\newblock URL \url{http://www.jstor.org/stable/3647580}.

\bibitem[Zou and Li(2008)]{Zou2008}
Hui Zou and Runze Li.
\newblock One-step sparse estimates in nonconcave penalized likelihood models.
\newblock \emph{The Annals of Statistics}, 36\penalty0 (4):\penalty0
  1509--1533, 2008.
\newblock ISSN 0090-5364.
\newblock \doi{10.48550/arxiv.0808.1012}.
\newblock URL \url{http://www.jstor.org/stable/25464679}.

\end{thebibliography}

\end{document}